\documentclass[11pt,reqno]{amsart}
\usepackage[a4paper,margin=27mm]{geometry}
\usepackage{amsmath,amssymb,mathtools,mathrsfs}
\usepackage{enumitem,booktabs,array}
\usepackage{needspace}
\IfFileExists{lmodern.sty}{\usepackage{lmodern}\usepackage{microtype}}{}
\usepackage[colorlinks=true,linkcolor=blue,citecolor=blue,urlcolor=blue]{hyperref}

\numberwithin{equation}{section}

\theoremstyle{plain}
\newtheorem{theorem}{Theorem}[section]
\newtheorem{lemma}[theorem]{Lemma}
\newtheorem{proposition}[theorem]{Proposition}
\newtheorem{corollary}[theorem]{Corollary}

\newtheorem{mainthm}{Theorem}

\newtheorem{maincor}[mainthm]{Corollary}
\theoremstyle{definition}
\newtheorem{definition}[theorem]{Definition}
\newtheorem{example}[theorem]{Example}

\theoremstyle{remark}
\newtheorem{remark}[theorem]{Remark}

\newcommand{\C}{\mathbb C}\newcommand{\Z}{\mathbb Z}
\newcommand{\R}{\mathbb R}\newcommand{\Q}{\mathbb Q}
\newcommand{\A}{\mathbb A}\newcommand{\PP}{\mathbb P}
\newcommand{\Gm}{\mathbb G_m}
\newcommand{\HH}{\mathbb H}
\newcommand{\GL}{\mathrm{GL}}\newcommand{\SL}{\mathrm{SL}}
\newcommand{\PGL}{\mathrm{PGL}}
\newcommand{\WG}{\mathsf W}
\newcommand{\der}{\mathrm{der}}\newcommand{\ab}{\mathrm{ab}}

\newcommand{\IC}{\operatorname{IC}}
\newcommand{\IH}{\operatorname{IH}}
\newcommand{\Dol}{\mathrm{Dol}}
\newcommand{\B}{\mathrm B}
\newcommand{\MD}{M_{\Dol}}
\newcommand{\MB}{M_{\B}}
\newcommand{\sslash}{\mathbin{/\mkern-6mu/}}
\newcommand{\ptau}{{}^{p}\!\tau}
\newcommand{\pH}{{}^{p}\!\mathcal H}
\newcommand{\Cc}{C^{\circ}}
\newcommand{\fg}{\mathfrak g}\newcommand{\ft}{\mathfrak t}
\newcommand{\Nil}{\mathcal N}
\newcommand{\Bfl}{\mathcal B}
\DeclareMathOperator{\Spec}{Spec}
\DeclareMathOperator{\Hom}{Hom}\DeclareMathOperator{\End}{End}
\DeclareMathOperator{\tr}{tr}\DeclareMathOperator{\Gr}{Gr}
\DeclareMathOperator{\id}{id}\DeclareMathOperator{\Ad}{Ad}
\DeclareMathOperator{\ad}{ad}
\DeclareMathOperator{\Lie}{Lie}\DeclareMathOperator{\rank}{rank}
\DeclareMathOperator{\Pic}{Pic}\DeclareMathOperator{\Supp}{Supp}
\DeclareMathOperator{\im}{im}\DeclareMathOperator{\codim}{codim}
\DeclareMathOperator{\fib}{fib}\DeclareMathOperator{\Cone}{Cone}
\DeclareMathOperator{\Flag}{Flag}\DeclareMathOperator{\sgn}{sgn}
\DeclareMathOperator{\ch}{ch}\DeclareMathOperator{\res}{res}
\DeclareMathOperator{\pr}{pr}\DeclareMathOperator{\Ext}{Ext}

\hypersetup{pdftitle={The P=W conjecture for reductive groups},pdfauthor={Tao Su}}

\begin{document}

\title[$P=W$ for reductive groups]{The $P=W$ conjecture for reductive groups}
\author{Tao Su}
\email{sutao@bimsa.cn}
\address{Beijing Institute of Mathematical Sciences and Applications, Beijing, China}
\date{}

\begin{abstract}
Let $G$ be a connected complex reductive group and let $C$ be a compact
Riemann surface with a nonempty finite set $D$ of punctures. Tame
nonabelian Hodge theory gives a homeomorphism $\Psi$ from a moduli space
of parabolic $G$-Higgs bundles on $(C,D)$ to a $G$-character variety of
$C\setminus D$ with prescribed local monodromy. The $P=W$ conjecture of
de~Cataldo, Hausel and Migliorini predicts that $\Psi$ identifies the
perverse filtration of the Hitchin map with the weight filtration of the
character variety, after doubling the index.

For every $G$ and every genus, we prove this identity in rational
intersection cohomology for nonempty character varieties whose local
monodromies are constrained to prescribed conjugacy-class closures,
assuming that the prescribed semisimple parts form an intrinsically
generic tuple. Repeated eigenvalues are allowed. When the moduli spaces are
rationally smooth, for instance when every prescribed class is
semisimple, or when only the semisimple part of the local monodromy is
prescribed, or for any mixture of these two conditions, the identity
holds in ordinary rational cohomology.

The proof constructs a degree-two operator $E$ (cup product with the
Goldman symplectic class) and a degree $-2$ operator $F$ (obtained by
adding a flag at a moving point of the curve, winding its eigenvalues,
and removing the flag again). We show that the leading actions of $E$ satisfy relative hard Lefschetz for
the perverse filtration and curious hard Lefschetz for the weight
filtration, that $F$ lowers the perverse index by two and the weight
index by four, and that
$[[F,E],E]=-2E$. An elementary lemma on $\mathfrak{sl}_2$-triples then
forces both filtrations to coincide with the eigenspace filtration of
$[E,F]$.
\end{abstract}

\maketitle
{\small\tableofcontents}

\section{Introduction}\label{sec:intro}

\subsection{Two algebraic structures on one topological space}
\label{subsec:intro-setting}

Throughout, $C$ is a smooth connected projective complex curve of genus
$g$, $D=\{p_1,\ldots,p_s\}\subset C$ is a set of $s\ge1$ points, and
$\Cc=C\setminus D$. Let $G$ be a connected complex reductive group.

\emph{The Betti side.} The fundamental group $\pi_1(\Cc)$ is generated by
$a_1,b_1,\ldots,a_g,b_g,c_1,\ldots,c_s$, where $c_j$ is a loop around
$p_j$, subject to the single relation $\prod_i[a_i,b_i]\prod_jc_j=1$.
Fix closed conjugation-invariant subsets $\overline{\mathcal C}_1,
\ldots,\overline{\mathcal C}_s\subset G$ (closures of conjugacy classes,
see \S\ref{subsec:intro-results}) and a central element
$\zeta\in Z(G)\cap G_{\der}$. The \emph{character variety} is the affine GIT quotient
\begin{equation}\label{eq:intro-MB}
 \MB=\Bigl\{(A_i,B_i,M_j)\in G^{2g}\times\textstyle\prod_j
   \overline{\mathcal C}_j:\ \prod_{i=1}^g[A_i,B_i]\prod_{j=1}^sM_j
   =\zeta\Bigr\}\sslash G,
\end{equation}
where $G$ acts by simultaneous conjugation and $[A,B]=ABA^{-1}B^{-1}$.
The points of $\MB$ are represented by the closed $G$-orbits of these
monodromy tuples. For $\zeta=1$, Richardson's criterion identifies these with
conjugacy classes of completely reducible representations of
$\pi_1(\Cc)$ \cite[Theorem~16.4]{Ric88}. Here complete reducibility means that every parabolic containing the
image has a Levi subgroup containing it. A direct proof over
$\C$ appears in Lemma~\ref{lem:general-closed-direct}.
The
quotient is affine, and Deligne's theory equips its cohomology with a
weight filtration $W_\bullet$ \cite{Del71,Del74}.

\emph{Central twists.} Because $D\ne\varnothing$, a twist $\zeta$ can
be absorbed into one puncture. Fix $p_0\in D$. Replacing $M_{p_0}$ by
$\zeta^{-1}M_{p_0}$ identifies~\eqref{eq:intro-MB} with the untwisted
character variety ($\zeta=1$) whose boundary condition at $p_0$ is the
translate $\zeta^{-1}\overline{\mathcal C}_{p_0}$, again a
conjugacy-class closure. For $\zeta\ne1$, the Dolbeault partner and the
map $\Psi$ below are, by definition, those of this untwisted problem,
composed with the inverse identification on the Betti side. Translation
by a central element changes neither the adjoint monodromy nor the
genericity condition of Definition~\ref{def:intro-generic}, and it does
not require $\zeta$ to be the exponential of an integral cocharacter.

\emph{The Dolbeault side.} A \emph{parabolic $G$-Higgs bundle} on $(C,D)$
is a principal $G$-bundle $E$ with a reduction to a parabolic subgroup
at each $p_j$ (a ``parabolic structure'' with real weights) together
with a \emph{Higgs field}
$\varphi\in H^0(C,\ad(E)\otimes K_C(D))$ whose residues are compatible
with the parabolic structure; see \S\ref{subsec:moduli} for the precise
(parahoric) definition. The moduli space $\MD$ of polystable objects is
a quasi-projective variety. Evaluating a basis of homogeneous invariant
polynomials on $\varphi$ defines the \emph{Hitchin map}
\[
 h:\MD\longrightarrow\mathcal A,
\]
which is proper with affine target (Hitchin \cite{Hit87a,Hit87b} and
Nitsure \cite[Theorem~6.1]{Nit91} for vector bundles on closed curves;
Proposition~\ref{prop:pg-master} and Theorem~\ref{thm:fine-master} for
the spaces used here). The
decomposition theorem applied to $h$ yields the \emph{perverse
filtration} $P_\bullet$ on the (intersection) cohomology of $\MD$
(Definition~\ref{def:perverse-filtration}).

\emph{Nonabelian Hodge theory.} By the work of Hitchin, Donaldson,
Corlette and Simpson in the compact case \cite{Hit87a,Don87,Cor88,Sim92},
Simpson for $\GL_n$ on punctured curves \cite{Sim90}, Biquard
\cite{Biq97} and Biquard--Garc\'{\i}a-Prada--Mundet i Riera for general
$G$ \cite{BGPM20}, there is a homeomorphism
\begin{equation}\label{eq:intro-Psi}
 \Psi:\MD\longrightarrow\MB,
\end{equation}
provided the parabolic weights and residues of the Higgs field are
matched with the eigenvalues of the local monodromies by the explicit
formulas recalled in \S\ref{subsec:nah}. The map $\Psi$ is real analytic
on smooth loci but not algebraic. Consequently $\Psi^*$ identifies the
rational cohomology of the two spaces, but the two algebraic structures
produce two \emph{a priori} unrelated filtrations on it.

\subsection{The $P=W$ conjecture}\label{subsec:intro-pw}

De~Cataldo, Hausel and Migliorini \cite{dCHM12} proposed that the two
filtrations agree after doubling the perverse index:
\begin{equation}\label{eq:PW}
 P_k H^j(\MD,\Q)=\Psi^*W_{2k}H^j(\MB,\Q)=\Psi^*W_{2k+1}H^j(\MB,\Q)
 \qquad(j,k\in\Z).
\end{equation}
This is an equality of subspaces, under the specific map $\Psi$, and not
merely an equality of the dimensions of the associated graded pieces.
For singular moduli spaces one replaces $H^j$ by rational intersection
cohomology $\IH^j$ (recalled in \S\ref{subsec:sheaves}); this
formulation is recorded in \cite[Question~4.1.7]{dCMa20} and
studied in \cite[Conjecture~1.2]{FM22}.

The conjecture originated in the computation of mixed Hodge polynomials
of character varieties by Hausel and Rodriguez-Villegas \cite{HRV08},
which predicted a ``curious'' symmetry of the weight filtration, and in
the topology of Hitchin systems. The following results are relevant
here.
\begin{itemize}[leftmargin=2em]
\item For $G=\GL_n$ on closed curves, with $n$ and the degree coprime,
\eqref{eq:PW} was proved for $n=2$ by de~Cataldo--Hausel--Migliorini
\cite{dCHM12}, for genus two by de~Cataldo--Maulik--Shen \cite{dCMS22},
and in general by Maulik--Shen \cite{MS24} and
Hausel--Mellit--Minets--Schiffmann \cite{HMMS}; a third proof was given
by Maulik--Shen--Yin \cite{MSY25}.
\item Hausel--Mellit--Minets--Schiffmann also treat moduli spaces of
stable full-flag parabolic $\GL_n$-Higgs bundles for generic stability
conditions \cite[\S\S8.3--8.6, Corollary~8.18]{HMMS}.
\item For $\SL_n$ with $n$ prime and degree coprime to $n$,
de~Cataldo--Maulik--Shen reduced the conjecture on closed curves to
$\GL_n$ \cite[Theorem~1.5 and Remark~1.6]{dCMS22SL} and proved it
in genus two \cite[Theorem~0.1]{dCMS22SL}. Combined with the general
$\GL_n$ theorem, this gives every genus $g\ge2$
\cite[Theorem~0.3]{MS24}.
Felisetti--Mauri proved the intersection-cohomology comparison for the
untwisted $\GL_n$ and $\SL_n$ character varieties of genus one in
arbitrary rank, and of genus two in rank two \cite[Main theorem, p.~857]{FM22}.
\item On the Betti side, the curious hard Lefschetz property conjectured
by Hausel--Letellier--Rodriguez-Villegas \cite{HLRV11} for generic
semisimple $\GL_n$-character varieties of punctured curves was proved by Mellit
\cite{Mel25}. Shende computed the weights of tautological classes
\cite{She17}.
\end{itemize}
Punctured curves and general reductive groups introduce two features
absent from the coprime $\GL_n$ case. First, when eigenvalues of the
local monodromy coincide, the boundary geometry is governed by nilpotent
orbit closures and Springer resolutions, and the moduli spaces are in
general singular; compare Letellier \cite{Let15} and Ballandras
\cite{Bal23} for $\GL_n$. Second, a change of the global form of $G$
(for instance from $\SL_n$ to $\PGL_n$) introduces finite central
lifting data and possibly disconnected centralizers.

\subsection{Main results}\label{subsec:intro-results}

\emph{Boundary conditions.} Every $M\in G$ has a Jordan decomposition
$M=SV$ with $S$ semisimple, $V$ unipotent and $SV=VS$. We prescribe, at
each puncture $p$, the closure $\overline{\mathcal C}_p$ of a conjugacy
class $\mathcal C_p$; the semisimple parts of elements of
$\overline{\mathcal C}_p$ are all conjugate to a fixed semisimple
$S_p\in G$. Three cases are singled out:
\begin{itemize}[leftmargin=2em]
\item \emph{semisimple}: $\mathcal C_p$ is the class of $S_p$ (then
$\overline{\mathcal C}_p=\mathcal C_p$);
\item \emph{full spectral}: $\overline{\mathcal C}_p$ is the set of all
$M$ with semisimple part conjugate to $S_p$ (this is the closure of the
class of $S_pu$, for $u$ regular unipotent in $Z_G(S_p)^0$);
\item \emph{closed Jordan}: the general case.
\end{itemize}
In terms of $\log V$, a closed Jordan condition requires $\log V$ to lie
in a prescribed nilpotent orbit closure in $\Lie Z_G(S_p)^0$, together
with all its translates by the possibly disconnected group $Z_G(S_p)$.

\emph{Genericity.} The central element $\zeta$ and the classes $S_p$
must satisfy the obvious necessary condition: the image of
$\zeta^{-1}\prod_pS_p$ in the abelianization $G_\ab=G/G_\der$ is $1$.
For a parabolic subgroup $P\subset G$, write $L_P=P/R_u(P)$ for its
Levi quotient and $L_P^\ab=L_P/[L_P,L_P]$ for the abelianization of
$L_P$; these notions are recalled in \S\ref{subsec:groups}.

\begin{definition}[Intrinsic genericity]\label{def:intro-generic}
The data $(\zeta,(S_p)_{p\in D})$ are \emph{intrinsically generic} if for
every proper parabolic subgroup $P\subsetneq G$ and every choice of
elements $s_p\in L_P$ with $s_p$ conjugate in $G$ to $S_p$ (for all
$p\in D$), the image of $\zeta^{-1}\prod_ps_p$ in $L_P^\ab$ is not $1$.
\end{definition}

Genericity is a condition on the tuple of semisimple parts; no single
puncture is required to have distinct eigenvalues. Genericity guarantees that
every representation in~\eqref{eq:intro-MB} is irreducible, i.e.\ is not
contained in a proper parabolic subgroup (Lemma~\ref{lem:generic}).

\begin{example}[$\GL_n$]\label{ex:intro-GLn}
Let $G=\GL_n$ and $\zeta=1$. Let $c_{p1},\ldots,c_{pk_p}$ be the distinct
eigenvalues of $S_p$, with multiplicities $m_{p1},\ldots,m_{pk_p}$. The
compatibility condition is $\prod_{p,a}c_{pa}^{m_{pa}}=1$. Proper
parabolic subgroups are stabilizers of partial flags, and it suffices to
test maximal ones, i.e.\ stabilizers of a subspace of dimension $r$,
$0<r<n$. The data are intrinsically generic if and only if
\begin{equation}\label{eq:intro-GLn-generic}
 \prod_{p,a}c_{pa}^{\,e_{pa}}\ne1\quad\text{whenever}\quad
 0\le e_{pa}\le m_{pa},\quad \sum_ae_{pa}=r\ \ (p\in D),\quad 0<r<n.
\end{equation}
This is the genericity condition of Hausel--Letellier--Rodriguez-Villegas
\cite[Definition~2.1.1]{HLRV11}.
\end{example}

\emph{Matching parameters.} We match the Dolbeault weights and residues
with the Betti data by the dictionary of \S\ref{subsec:nah}, retaining all
compatible topological components. We call these parameters
\emph{compatible}. With these conventions, the main theorem is as
follows.

\Needspace{18\baselineskip}
\begin{mainthm}[Generic $P=W$]\label{thm:A}
Let $G$ be a connected complex reductive group, $C$ a smooth projective
curve of arbitrary genus, and $D\subset C$ a nonempty finite set. Let
$(\zeta,(\overline{\mathcal C}_p)_{p\in D})$ be closed Jordan boundary
conditions whose semisimple data are intrinsically generic, and assume
$\MB\ne\varnothing$. Let $\MD$ be the moduli space of polystable
parabolic $G$-Higgs bundles with compatible parameters and
$\Psi:\MD\to\MB$ the tame nonabelian Hodge homeomorphism. Then for all
$j,k$,
\begin{equation}\label{eq:thmA}
 P_k\IH^j(\MD,\Q)=\Psi^*W_{2k}\IH^j(\MB,\Q)=\Psi^*W_{2k+1}\IH^j(\MB,\Q),
\end{equation}
where the perverse filtration is normalized by the defect of the Hitchin
map (Definition~\ref{def:perverse-filtration}).
\end{mainthm}

Theorem~\ref{thm:A} is proved in \S\ref{sec:proof}, by
Theorem~\ref{thm:boundary-descent} and the central descent of
\S\ref{subsec:central-quotient}. We record the ordinary-cohomology and
equivariant consequences separately.

\begin{maincor}[Ordinary $P=W$]\label{cor:ordinary}
Under the hypotheses of Theorem~\ref{thm:A}, suppose that $\MB$
(equivalently $\MD$) is rationally smooth. Then \eqref{eq:PW} holds in
ordinary rational cohomology, with the same normalization of the
perverse filtration. In particular, this conclusion holds whenever the
boundary condition at each puncture is either semisimple or full
spectral, allowing arbitrary mixtures of these two conditions.
\end{maincor}

\begin{maincor}[Finite symmetry summands]\label{cor:equivariant}
Under the hypotheses of Theorem~\ref{thm:A}, let $\Gamma$ be a finite
group acting algebraically on the boundary problem and on $\MD$ and
$\MB$, compatibly with $\Psi$, and over the Hitchin base on $\MD$.
For every idempotent $e\in\Q[\Gamma]$, \eqref{eq:thmA} holds on the
images of $e$, with the filtrations induced from the full intersection
cohomology groups. Under the additional hypothesis of
Corollary~\ref{cor:ordinary}, the same conclusion holds in ordinary
cohomology. For an abelian group, individual character summands are
taken over a splitting field, or summed over Galois orbits over $\Q$,
as in \S\ref{subsec:finite}.
\end{maincor}

\begin{remark}[Central twists]\label{rem:central-twists}
For a finite central subgroup $Z$, the full twisting group $H^1(\Cc,Z)$
acts on the complete union of lifted boundary sectors used in descent;
on a fixed boundary tuple only its stabilizer subgroup acts.
Corollary~\ref{cor:equivariant} applies to these actions, as
explained in \S\ref{subsec:central-quotient}.
\end{remark}

Corollary~\ref{cor:ordinary} is proved in
\S\ref{subsec:rational-smooth}; Corollary~\ref{cor:equivariant} follows
from the preservation of both filtrations by the finite projectors of
\S\ref{subsec:finite}. Several remarks on the scope of these results
are in order.

\begin{remark}\label{rem:intro-scope}
(1) The group $G$ is arbitrary: neither simply connected, adjoint, nor
of type $A$. The moduli spaces include all topological types of
$G$-bundles compatible with the boundary data, and no component of
$\MB$ or $\MD$ is discarded.

(2) After choosing the growth weights and logarithmic extensions,
the residue theorem~\eqref{eq:fuchs} constrains the degrees detected by
characters of $G$. These equations do not detect every torsion
topological class of a semisimple group. All compatible topological
classes and central lifting sectors are retained; the character degrees
cannot be varied independently of the chosen parameters.

(3) Nonemptiness and genericity are used; nongeneric problems, and
problems with strictly semistable points, are outside the scope of this
paper.

(4) Corollary~\ref{cor:ordinary} includes singular full spectral
spaces: they are rationally smooth because the nilpotent cone is
\cite[\S2.3 and its proof in \S\S2.6--2.9]{BM83} (Proposition~\ref{prop:springer}). A smaller nilpotent
orbit closure need not be rationally smooth, which is why the general
statement uses $\IH$.
\end{remark}

\begin{maincor}[General linear and special linear groups]\label{cor:B}
Every nonempty generic semisimple $\GL_n$-character variety of $\Cc$,
and every nonempty generic semisimple fixed-determinant
$\SL_n$-character variety of $\Cc$, satisfies \eqref{eq:PW} in ordinary
cohomology. The genus and $n$ are arbitrary, repeated eigenvalues at a
puncture are allowed, and the identity holds on every character summand
of the compatible finite group actions of
Corollary~\ref{cor:equivariant}. Full spectral and mixed
semisimple/full spectral conditions satisfy the same identity in
ordinary cohomology, and general closed Jordan conditions satisfy
\eqref{eq:thmA}.
\end{maincor}

Corollary~\ref{cor:B} is proved at the end of \S\ref{sec:proof}.
The full-flag parabolic cases overlap with the setup of
\cite[\S8, Corollary~8.18]{HMMS}, with the stability and residue
conditions specified there. For instance, a semisimple
$\GL_n$ problem with multiplicities $m_{pa}$ has dimension
\begin{equation}\label{eq:intro-dim}
 \dim\MB=(2g-2)n^2+2+\sum_{p\in D}\Bigl(n^2-\sum_am_{pa}^2\Bigr),
\end{equation}
and the fixed-determinant $\SL_n$ problem has dimension
$(2g-2)(n^2-1)+\sum_p(n^2-\sum_am_{pa}^2)$; compare
\cite[Theorem~2.1.5]{HLRV11}. For regular semisimple boundary classes,
the same dimension count follows from
Lemma~\ref{lem:stable-zero-section}.

\subsection{Strategy of the proof}\label{subsec:intro-strategy}

\emph{An $\mathfrak{sl}_2$ criterion.} Let $V$ be a finite-dimensional
vector space with two finite increasing filtrations $P_\bullet$ and
$L_\bullet$, and let $E\in\End V$ satisfy $E(P_k)\subset P_{k+2}$ and
$E(L_k)\subset L_{k+2}$. Suppose that $E$ induces hard Lefschetz
isomorphisms $\Gr_{N-r}\xrightarrow{\sim}\Gr_{N+r}$ for both
filtrations, with the same center $N$. These hypotheses alone do not
force $P=L$. Section~\ref{sec:linear} shows that they do if in addition
there is a single $F\in\End V$ with
\begin{equation}\label{eq:intro-F}
 F(P_k)\subset P_{k-2},\qquad F(L_k)\subset L_{k-2},\qquad
 [[F,E],E]=-2E.
\end{equation}
Indeed, then $(E,[E,F],F)$ is an $\mathfrak{sl}_2$-triple and both
filtrations equal $\bigoplus_{a\le k-N}\ker([E,F]-a)$
(Theorem~\ref{thm:common-lowering}). We apply this with
$P_\bullet$ the perverse filtration and $L_k=W_{2k}$. This strategy is
the one of \cite{HMMS}, where both filtrations are shown to coincide with
the filtration of an $\mathfrak{sl}_2$-triple built from tautological
classes and Hecke operators; our operators are constructed differently,
in a way that works uniformly for all $G$.

\emph{The raising operator $E$} (\S\ref{sec:raising}) is cup product
with the class of the Goldman symplectic form on $\MB$, normalized by an
invariant form $\kappa$ on $\fg$. For regular semisimple boundary data we
work on a \emph{unipotent stabilization} $f:\widetilde\MB\to\MB$
(Definition~\ref{def:stabilization}): over the rigidified stacks it is a
vector bundle with fiber the unipotent radical $U_B$ of a Borel subgroup.
Following the Bruhat-cell constructions of Mellit \cite{Mel25} for
$\GL_n$ and of de~Amorin \cite{dA25} in general, $\widetilde\MB$ has a
finite closed filtration with pieces $(\C^*)^t\times(\C^e/F)$, $F$
finite. Every statement we need about these coordinates is proved in
\S\ref{subsec:stabilized-bruhat}; we use neither the quotient
statements of \cite{dA25} nor a product description of strata of $\MB$
itself. For very generic $\GL_n$ data,
\cite[Theorem~2.10 and Appendix~C]{Su26} gives stable torus--affine
descriptions: each stratum $Z$ satisfies
$Z\times\A^{\dim U_B}\cong(\C^*)^a\times\A^b$.
The affine factor is not canceled. On each stabilized piece the
pulled-back Goldman class is the cup-product pairing on the parabolic
cohomology of a local system of Weyl reflections, hence nondegenerate.
This gives curious hard Lefschetz on $\widetilde\MB$, and the
compactly supported trace of $f$ transfers it to $\MB$. Thus the
weight-$4$ part of $E$ is Lefschetz in every root system. On the
Dolbeault side,
the Goldman class differs from the first Chern class of a determinant
line bundle, which is relatively ample for the Hitchin map, by a
combination of first Chern classes of line bundles attached to the
punctures; we show that these boundary classes raise the perverse
filtration by at most one. Hence the perverse-degree-$2$ part of $E$
satisfies relative hard Lefschetz.

\emph{The lowering operator $F$} (\S\ref{sec:lowering}) comes from a
moving auxiliary puncture $q\in\Cc$. A Gysin map $B$ adds a scalar flag
at $q$; a trace map $A$ removes it; and $AB=(-1)^{b}|\WG|$, where
$b=\dim G/B$ and $\WG$ is the Weyl group. Winding the eigenvalues at $q$
around a cocharacter loop $\gamma$ gives a unipotent monodromy with
normalized logarithm $N_\gamma$. After integrating the moving point over
$C$, the operator
\[
 F^{\rm raw}=\tr_C\circ A_C\circ N_\gamma^2\circ B_C\circ\mathrm{unit}_C
\]
lowers weights by four because $N_\gamma$ lowers them by two. A
computation of characteristic classes shows $[[F^{\rm raw},E],E]
=c\cdot E$ for an explicit nonzero constant $c$
(Lemma~\ref{lem:whole-quadratic-return}).

\emph{The perverse bound for $F$} (\S\ref{sec:perverse}) is the heart of
the argument. We realize the same operator by a correspondence over the
Hitchin base, by degenerating the moving-point construction along the
Hodge parameter $\lambda\to0$ (from flat connections to Higgs fields). An
annular estimate proves that the degenerated correspondence is proper,
and that its central fiber lies over the diagonal of the Hitchin base.
A degree-zero morphism of the Hitchin direct image over $\mathcal
A\times C$, followed by the trace over $C$, lowers the perverse index by
two.

\emph{Descent} (\S\ref{sec:proof}). The previous steps are carried out on
a \emph{fine reference} problem: regular semisimple boundary classes,
full flags, and a simply connected derived group. The Springer
decomposition, applied in families constructed in Section~\ref{sec:comparison-families},
extracts the intersection cohomology of every closed Jordan condition with explicit
shifts. A finite central quotient (\S\ref{subsec:central-lifts})
returns to the given group $G$.

\begin{example}[The rank-one case]\label{ex:intro-rank-one}
Let $G=\C^*$ and fix compatible puncture monodromies. Then
$\MB\cong(\C^*)^{2g}$, with coordinates the monodromies $x_i,y_i$ along
$a_i,b_i$, and $H^*(\MB)$ is the exterior algebra on the classes
$\xi_i=[d\log x_i]$, $\eta_i=[d\log y_i]$, all of weight two (and Hodge
type $(1,1)$). Hence $H^j(\MB)$ is pure of weight $2j$. The partner
$\MD$ is a translate of $T^*\Pic^d(C)$, and $h:\MD\to H^0(C,K_C)\cong
\C^g$ is the projection. With $K=Rh_*\Q[g]$, the summand
$H^j(\Pic^d(C))\otimes\Q_{\C^g}[g][-j]$ of $K$ has perverse index $j$, so
$P_kH^j=H^j$ for $k\ge j$ and $0$ for $k<j$. Thus \eqref{eq:PW} holds.
The $\mathfrak{sl}_2$ criterion is visible here: $E=\sum_i\xi_i\eta_i$
(the Goldman form) and $F=\sum_i\iota_{\eta_i^\vee}\iota_{\xi_i^\vee}$
(contraction with the dual bivector) satisfy $[E,F]=j-g$ on $H^j$, and
both filtrations are the filtrations by cohomological degree.
\end{example}

\subsection{Organization}

Section~\ref{sec:background} recalls the definitions and standard tools,
and Section~\ref{sec:linear} proves the $\mathfrak{sl}_2$ criterion.
Section~\ref{sec:comparison-families} constructs the residue families and
identifies the comparison maps. Section~\ref{sec:raising} treats the
raising operator $E$; Section~\ref{sec:lowering} constructs the lowering
operator $F$, the compact
coefficient and proves the weight estimate. The characteristic-class
and commutator calculations are in Section~\ref{sec:cup-comparisons}.
Section~\ref{sec:perverse} constructs the correspondences and proves
the perverse bound for $F$. Section~\ref{sec:proof} combines these arguments,
including closed Jordan conditions, central descent and fixed
determinants.

The logical chain is:
\begin{center}
\begin{tabular}{l}
Theorem~\ref{thm:common-leading} (the operator $E$)
$+$ Proposition~\ref{prop:betti-weight} and
Lemma~\ref{lem:whole-quadratic-return} ($F$ on $H^*(\MB)$)\\
\quad$+$ Theorem~\ref{thm:rc-perverse} (perverse degree of $F$)
$+$ Theorem~\ref{thm:common-lowering} ($\mathfrak{sl}_2$ criterion)\\
\quad$\Longrightarrow$ Theorem~\ref{thm:fine-comparison} (fine
references)\\
\quad$\Longrightarrow$ Theorem~\ref{thm:boundary-descent} (closed Jordan
conditions) $\Longrightarrow$ Theorem~\ref{thm:A}.
\end{tabular}
\end{center}
The last two steps are Springer extraction (\S\ref{subsec:extraction})
and central descent (\S\ref{subsec:central-quotient}).

Two appendices contain supporting proofs. Appendix~\ref{app:analytic}
contains the analytic inputs to harmonic transport: the harmonic inverse
at a nonzero real Hodge parameter, uniform prolongation frames and
compactness. Appendix~\ref{app:calculations} contains the group-theoretic,
Springer-fiber and Bruhat calculations. The family, coefficient and
correspondence constructions are in the main text.

\subsection*{Conventions}
Varieties are over $\C$; cohomology has rational coefficients unless
stated otherwise; ``dimension'' means complex dimension. Shifts are
cohomological: $H^i(K[a])=H^{i+a}(K)$. The Tate structure $\Q(-1)$ has
Hodge type $(1,1)$ and weight $2$. We write $\WG$ for the Weyl group and
reserve $W_\bullet$ for weight filtrations. The stacks of intrinsically generic moduli problems have finite
stabilizers after rigidification by $Z(G)$; their rational statements
descend from finite quotient charts by averaging
(Lemma~\ref{lem:finite-quotient}). Auxiliary classifying and equivariant
Springer stacks may have positive-dimensional affine stabilizers.
Their coherent operations are formed in the stack formalism recalled in
\S\ref{subsec:coefficient-enhancement}, before pullback to the generic moduli charts.

\section{Background}\label{sec:background}

This section fixes notation and states, in the form used later, the
standard results on which the proof relies. Readers familiar with a
subsection may skip it; every later use refers back to a numbered
statement here.

\subsection{Reductive groups}\label{subsec:groups}

Our reference is \cite[Chapters~6--8]{Spr98}. For $\GL_n$ the reader
may keep in mind: $T$ the diagonal matrices, $B$ the upper triangular
matrices, $U_B$ the upper unitriangular matrices.

A connected linear algebraic group is \emph{reductive} if it has
no nontrivial connected normal unipotent subgroup. Here an element is
\emph{unipotent} if all its eigenvalues in a faithful matrix
representation are one. A \emph{torus} is an algebraic group isomorphic
to $(\C^*)^r$. Fix a
maximal torus $T$ of $G$ and write
\[
 X^*(T)=\Hom(T,\Gm),\qquad X_*(T)=\Hom(\Gm,T),\qquad
 \langle\ ,\ \rangle:X^*(T)\times X_*(T)\to\Z
\]
for the character lattice, the cocharacter lattice and their perfect
pairing. The \emph{roots} $\Phi\subset X^*(T)$ are the nonzero weights of
$T$ on $\fg=\Lie G$; to each root $\alpha$ corresponds a coroot
$\alpha^\vee\in X_*(T)$, and $(X^*(T),\Phi,X_*(T),\Phi^\vee)$ is the
\emph{root datum} of $G$. We always use the lattices of $G$ itself, not
their rational spans: a coroot may be divisible in $X_*(T)$, and a
surjection of tori may have disconnected kernel. The \emph{Weyl group}
is $\WG=N_G(T)/T$.

A \emph{Borel subgroup} $B\supset T$ is a maximal connected solvable
subgroup; it determines positive roots $\Phi^+$ and simple roots. Put
\[
 b=|\Phi^+|=\dim G/B=\dim U_B,
\]
where $U_B=R_u(B)$ is the unipotent radical. A \emph{parabolic subgroup}
is a closed subgroup $P$ containing a Borel subgroup; then $G/P$ is
projective, and $P=L\ltimes R_u(P)$ where $R_u(P)$ is the unipotent
radical and $L$ is a \emph{Levi subgroup}. We write $L_P=P/R_u(P)$ and
$L_P^\ab=L_P/[L_P,L_P]$, a torus. For $\GL_n$, parabolic subgroups are
the stabilizers of partial flags and $L_P^\ab$ records the determinants of
the graded pieces.

A \emph{pinning}, relative to $T\subset B$, chooses compatible
parametrizations of the simple root subgroups. The associated root
homomorphisms $\SL_2\to G$ send the diagonal torus through the coroots
$\alpha_i^\vee$. With
$n_i=\left(\begin{smallmatrix}0&-1\\1&0\end{smallmatrix}\right)$ and
$x_i(z)=\left(\begin{smallmatrix}1&z\\0&1\end{smallmatrix}\right)$ in the
image of the $i$th copy, one has $n_i^2=\alpha_i^\vee(-1)$ and, for $e\ne0$,
\begin{equation}\label{eq:rank-one}
 n_ix_i(e)n_i=x_i(-e^{-1})\,\alpha_i^\vee(e^{-1})\,n_i\,x_i(-e^{-1}).
\end{equation}
Representatives of Weyl group elements therefore carry constant torus
factors, which we keep track of. The \emph{Bruhat decomposition}
$G=\bigsqcup_{w\in\WG}Bn_wB$ is a finite partition into locally closed
pieces, and for a simple reflection $s$
\begin{equation}\label{eq:bruhat-rule}
 (BsB)(BwB)=BswB\ \ \text{if}\ \ \ell(sw)>\ell(w),\qquad
 (BsB)(BwB)\subset BswB\cup BwB\ \ \text{otherwise}
\end{equation}
\cite[Lemma~8.3.7]{Spr98}, where $\ell$ is the length function. Recording,
along a word in simple reflections, which alternative occurs at each
letter gives the \emph{Bruhat walks} of \S\ref{subsec:torus-coordinates};
these are closely related to Deodhar's distinguished subexpressions
\cite{Deo85}.

For $S\in G$, its \emph{centralizer} is
$Z_G(S)=\{g\in G:gS=Sg\}$. If $S$ is semisimple, this group is
reductive and may be disconnected; $Z_G(S)^0$ contains a maximal torus.
If $G_\der$ is simply connected, Steinberg's theorem \cite{Ste75}
gives connectedness; Lemma~\ref{lem:semisimple-centralizer-direct}
proves the complex case used here. In particular, a regular semisimple
$S\in T$ then has $Z_G(S)=T$. For other global forms one must keep the
full centralizer, including its component group.

Fix a faithful representation $\rho:G\hookrightarrow\SL(V)$ (append the
inverse determinant to a faithful representation in $\GL(V')$) and let
$\kappa(x,y)=\tr(\rho(x)\rho(y))$ be the trace form of $\rho$ on $\fg$.
The form $\kappa$ is invariant, nondegenerate, and rational on
$X_*(T)\otimes\Q$, and it is
orthogonal for the decomposition of $\fg$ into the center and the simple
factors of $\fg_\der$.

\subsection{Nilpotent cones and Springer theory}\label{subsec:springer}

Let $L$ be a connected reductive group with Lie algebra $\mathfrak l$,
Borel subgroup $B_L$ and Weyl group $\WG_L$. The \emph{nilpotent cone}
$\Nil_L\subset\mathfrak l$ is the zero fiber of the adjoint quotient
$\mathfrak l\to\mathfrak l\sslash L$. The cone $\Nil_L$ is a finite union of
\emph{nilpotent orbits} (adjoint $L$-orbits). The \emph{Springer
resolution} and the \emph{Springer fiber} of $N\in\Nil_L$ are
\[
 \mu:\widetilde\Nil_L=L\times^{B_L}\Lie R_u(B_L)\longrightarrow\Nil_L,
 \qquad
 \Bfl_N=\mu^{-1}(N)=\{\,B'\in L/B_L:\ N\in\Lie B'\,\}.
\]
Thus $\mu$ replaces a nilpotent element by that element together with a
Borel subalgebra containing it; for $\GL_n$, by that element together
with an invariant full flag. The \emph{Grothendieck--Springer
resolution} $\widetilde{\mathfrak l}=\{(x,B'):x\in\Lie B'\}\to\mathfrak l$
is a $\WG_L$-Galois cover over the regular semisimple locus; its
multiplicative analog is $\widetilde L=\{(g,B'):g\in B'\}\to L$.

\begin{proposition}[Springer theory]\label{prop:springer}\leavevmode
\begin{enumerate}[label=(\roman*)]
\item\label{it:kostant} $\Nil_L$ is an irreducible normal complete
intersection \cite{Kos63}; an alternative proof over $\C$ is given in
Lemma~\ref{lem:nilpotent-cone} below. Consequently, by Zariski's main theorem,
every Springer fiber is connected.
\item\label{it:semismall} The map $\mu$ is proper and birational. Every
nilpotent orbit is a \emph{relevant stratum}, meaning that
$2\dim\Bfl_N=\codim(L\cdot N\subset\Nil_L)$
\cite[\S2.1]{BM83}, \cite[Proposition~6.7(ii)]{Spr76}.
We give a direct dimension argument in \S\ref{app:springer}.
In particular $\mu$ is semismall. The decomposition theorem gives
\[
 R\mu_*\Q[\dim\Nil_L]=\bigoplus_{(\mathcal O,\mathcal L)}
 V_{\mathcal O,\mathcal L}\otimes\IC(\overline{\mathcal O},\mathcal L),
\]
a sum over nilpotent orbits $\mathcal O$ and irreducible
$L$-equivariant local systems $\mathcal L$ on $\mathcal O$, with
finite-dimensional multiplicity spaces. The multiplicity of the constant
local system on $\mathcal O\ni N$ is
$V_{\mathcal O}=H^{2b_{\mathcal O}}(\Bfl_N,\Q)^{A_L(N)}$, where
$b_{\mathcal O}=\dim\Bfl_N$ and $A_L(N)=\pi_0Z_L(N)$; it is nonzero,
since $A_L(N)$ permutes the top-dimensional irreducible components of
$\Bfl_N$ and the sum of their classes is invariant
\cite[\S\S1.2, 1.7 and 2.2]{BM83}.
\item\label{it:springer-action} Restriction from the
Grothendieck--Springer resolution defines an action of $\WG_L$ on
$R\mu_*\Q$ \cite[\S\S2.2, 2.6--2.7]{BM83}. The spaces
$V_{\mathcal O,\mathcal L}$ are zero or irreducible $\WG_L$-modules, and distinct
occurring complex pairs give inequivalent representations
\cite[Theorem~1.13 and its proof in \S3.11]{Spr78}. We use the
restriction convention: the zero orbit corresponds to sign and the
regular orbit to the trivial representation. Springer's original
normalization assigns sign to the regular orbit
\cite[\S7.4]{Spr76}; irreducibility is unchanged by a sign twist.
The Fourier-transform convention also differs from restriction by
tensoring with sign \cite[Theorem~1.1]{AHJR14}.
\item\label{it:rational-smooth} The $\WG_L$-invariant summand of
$R\mu_*\Q$ is $\Q_{\Nil_L}$; hence $\Nil_L$ is rationally smooth
\cite[\S2.3 and its proof in \S\S2.6--2.9]{BM83}.
\item\label{it:DLP} For every $N$, $H^{\rm odd}(\Bfl_N,\Z)=0$ and
$H^{\rm even}(\Bfl_N,\Z)$ is torsion free \cite{DLP88}; moreover
$H^1(\Bfl_N,\Q)=0$ and $H^2(\Bfl_N,\Q)$ is a sum of copies of $\Q(-1)$
(Proposition~\ref{prop:springer-low-degree}).
\end{enumerate}
\end{proposition}

In type $A$ the Springer fiber dimension of a nilpotent of Jordan type
$\lambda$ is $\sum_j\binom{\lambda'_j}{2}$, where $\lambda'$ is the
conjugate partition. For a $2\times2$ block, the regular orbit has
$b=0$ and the zero orbit has $b=1$.

\subsection{Mixed Hodge structures}\label{subsec:mhs}

A \emph{(rational) mixed Hodge structure} is a finite-dimensional
$\Q$-vector space $V$ with an increasing \emph{weight filtration}
$W_\bullet$ and a decreasing \emph{Hodge filtration} $F^\bullet$ on
$V_\C$, such that each $\Gr^W_mV$ is pure of weight $m$: its Hodge
components have types $(p,q)$ with $p+q=m$. Deligne constructed
functorial mixed Hodge structures on $H^j(X)$ and $H^j_c(X)$ for every
complex variety $X$ \cite{Del71,Del74}. Morphisms of mixed Hodge
structures $f:V\to V'$ are \emph{strict} for both filtrations
\cite[Th\'eor\`eme~2.3.5]{Del71}:
\[
 \begin{aligned}
 f(W_mV)&=f(V)\cap W_mV',\\
 f_\C(F^pV_\C)&=f_\C(V_\C)\cap F^pV'_\C.
 \end{aligned}
\]
Thus kernels, images and idempotent summands carry the induced
filtrations, and passing to either associated graded is exact.

The \emph{Tate structure} $\Q(-1)$ is one-dimensional of type $(1,1)$;
tensoring with $\Q(-b)$ raises weights by $2b$. A mixed Hodge structure
is \emph{Hodge--Tate} if every $\Gr^W_m$ has only types $(p,p)$; then
odd weight-graded pieces vanish. For a smooth $X$ of dimension $d$,
Poincar\'e duality $H^j(X)\cong H^{2d-j}_c(X)^\vee(-d)$ is an isomorphism
of mixed Hodge structures, so $W_\ell H^j(X)$ is the annihilator of
$W_{2d-\ell-1}H_c^{2d-j}(X)$.

For an endomorphism $T$ of a filtered space $(V,W_\bullet)$ we say that
$T$ has \emph{weight degree at most $a$} if $T(W_k)\subset W_{k+a}$ for
all $k$; then $T$ induces $\Gr^W_k\to\Gr^W_{k+a}$, its \emph{leading
part}.

\subsection{Constructible complexes and the perverse filtration}
\label{subsec:sheaves}

For a variety $X$ let $D^b_c(X)$ be the bounded derived category of
sheaves of $\Q$-vector spaces with constructible cohomology. We use the
six operations $f^*,Rf_*,Rf_!,f^!,\otimes,R\mathcal{H}om$ and their
standard adjunctions and base-change morphisms \cite{KS90}; the
base-change maps are used as isomorphisms under the properness or local
acyclicity hypotheses stated at each application. When $f$ is proper, $Rf_!=Rf_*$, and \emph{proper base change} identifies the
stalk of $Rf_*K$ at $y$ with $R\Gamma(f^{-1}(y),K)$.

\emph{Perverse sheaves.} The middle perversity $t$-structure on
$D^b_c(X)$ has heart the abelian category of perverse sheaves: $K$ is
perverse if $\dim\Supp\mathcal H^i(K)\le-i$ and the same holds for its
Verdier dual \cite[\S2.1]{BBD82}. We write $\ptau_{\le k}$ for the
truncation functors and $\pH^k$ for perverse cohomology. On a smooth
$X$ of dimension $m$, $\Q_X[m]$ is perverse.

\emph{Intersection cohomology.} For $X$ of pure dimension $d$ with dense
smooth open $j:X^\circ\hookrightarrow X$, the \emph{intersection complex}
$\IC_X=j_{!*}\Q_{X^\circ}[d]$ is the unique perverse extension of
$\Q_{X^\circ}[d]$ with no nonzero subobject or quotient supported on
$X\setminus X^\circ$ \cite[Corollary~1.4.25]{BBD82}, and
\[
 \IH^j(X)=\HH^{j-d}(X,\IC_X)
\]
is the intersection cohomology of Goresky--MacPherson \cite{GM83}.
Similarly $\IC(\overline{\mathcal O},\mathcal L)$ denotes the
intermediate extension of a local system $\mathcal L$. We say that $X$
is \emph{rationally smooth} if $\IC_X\cong\Q_X[d]$; then
$\IH^j(X)=H^j(X)$, compatibly with mixed Hodge structures. Finite
quotients of smooth varieties are rationally smooth
(Lemma~\ref{lem:finite-quotient}).

\emph{Mixed Hodge modules.} Saito's mixed Hodge modules enrich suitable constructible complexes
with weight and Hodge filtrations, and lift the six operations to this
setting \cite[Introduction, Theorem~0.1 and \S4]{Sai90}.
Rational realization is exact and faithful on the abelian category of
mixed Hodge modules, and conservative on its bounded derived category:
a morphism whose realization is an isomorphism is an isomorphism.
Derived realization need not be faithful, so equality of derived
morphisms or their homotopies is established in the enhanced coefficient
diagrams, rather than inferred from equality of their rational images. Every complex of
geometric origin used below (constant sheaves, intersection complexes,
their direct images) is the realization of a mixed Hodge module complex,
and $\HH^*(X,\IC_X)$ receives a mixed Hodge structure in this way.
Saito's theory extends to algebraic stacks \cite{Tub25}.

A \emph{pure} Hodge module has a single weight. Polarizable pure Hodge
modules decompose by \emph{strict support}: a summand supported on an
irreducible closed subset $Z$ has no nonzero subobject or quotient
supported on a proper closed subset of $Z$. Such a summand is the
intermediate extension of a polarized variation on a dense smooth open
of $Z$. Maps between different strict supports vanish, and a map on a
common strict support is determined on that dense open
\cite[\S5.1.3, Lemma~5.1.10 and Corollary~5.1.11]{Sai88}.

\begin{theorem}[Decomposition theorem and relative hard Lefschetz]
\label{thm:decomposition}
Let $h:X\to\mathcal A$ be a projective map of varieties, with $X$ of
pure dimension $d$. Each $\pH^i(Rh_*\IC_X)$ underlies a polarizable pure
Hodge module of weight $d+i$ and decomposes by strict support.
Moreover, $Rh_*\IC_X$ splits as the direct sum of these objects
shifted by $[-i]$.
If $\eta$ is the first Chern class of an $h$-ample line bundle, then for
$r$ as in Definition~\ref{def:perverse-filtration}
and every $s\ge0$, cup product with $\eta^s$ induces isomorphisms
\[
 \eta^s:\pH^{r-s}(Rh_*\IC_X[-r])\xrightarrow{\ \sim\ }
 \pH^{r+s}(Rh_*\IC_X[-r])(s).
\]
\end{theorem}

The general decomposition theorem is due to
Beilinson--Bernstein--Deligne--Gabber
\cite[Theorems~6.2.5 and~6.2.10]{BBD82}. We use the projective
Hodge-theoretic form of Saito
\cite[Theorem~5.3.1 and Remark~5.3.12]{Sai88}; see also
\cite[Theorems~1.6.1 and~3.1.16]{dCM09}.
The analytic curve input is Zucker's Hodge theorem for polarized
coefficients, including its filtered comparison
\cite[Theorem~7.12, Corollary~7.13 and Proposition~9.1]{Zuc79}.
The local estimates use Schmid's one-variable orbit and limiting Hodge
theorems \cite[Theorems~4.9, 5.13 and~6.16]{Sch73}.
Relative hard Lefschetz gives the splitting: split off the two extreme
perverse cohomology objects using the extreme Lefschetz isomorphism,
correct their inclusions to make the retractions orthogonal, and
induct on the cohomological amplitude. This argument also works in the
Hodge category with the Tate twists retained. Thus every perverse
truncation is a direct summand. The chosen splitting is not canonical;
the truncation morphism is.

\begin{definition}[Perverse filtration]\label{def:perverse-filtration}
Let $h:X\to\mathcal A$ be proper, $X$ of pure dimension $d$. The
\emph{defect} of $h$ is $r=\dim(X\times_{\mathcal A}X)-d$. Put
$K_h=Rh_*\IC_X[-r]$ and define
\[
 P_k\IH^j(X)=\im\Bigl[\HH^{j-d+r}(\mathcal A,\ptau_{\le k}K_h)
 \longrightarrow\HH^{j-d+r}(\mathcal A,K_h)=\IH^j(X)\Bigr].
\]
If $X$ is smooth this reads $K_h=Rh_*\Q_X[d-r]$ and
$P_kH^j(X)=\im[\HH^{j-d+r}(\ptau_{\le k}K_h)\to H^j(X)]$.
\end{definition}

The definition is intrinsic; no splitting is chosen.  If $X$ is rationally smooth, the normalized
filtration is concentrated in $[0,2r]$: on a target stratum $S$ with
fiber dimension $e$ one has $\dim S+2e\le d+r$, and proper base change
bounds the stalk cohomology of $Rh_*\Q_X[d]$ above by $2e-d$.
The perverse support condition gives the upper bound $r$ before
normalization, and duality gives the lower bound $-r$.
For every pure-dimensional $X$ and projective $h$,
Theorem~\ref{thm:decomposition} gives isomorphisms
$\eta^s:\Gr^P_{r-s}\IH^j(X)\xrightarrow{\sim}
 \Gr^P_{r+s}\IH^{j+2s}(X)$.
For the Hitchin maps in this paper $\dim X=2N$ and $r=N$
(Theorem~\ref{thm:fine-master}), so the center is $N$.

An endomorphism $T$ of $\IH^*(X)$ has \emph{perverse degree at most $a$}
if $T(P_k)\subset P_{k+a}$. A stronger, sheaf-level notion is used in
\S\ref{subsec:boundary-classes}: a morphism $\tau:K_h\to K_h[m]$ has
\emph{strong perverse degree at most $a$} if it maps
$\ptau_{\le k}K_h$ into $(\ptau_{\le k+a}K_h)[m]$ compatibly with the
truncation maps; this implies perverse degree at most $a$ on
hypercohomology.

\subsection{Gysin maps, traces, correspondences and gluing}
\label{subsec:correspondences}

Let $i:Z\hookrightarrow X$ be a closed embedding of smooth varieties of
codimension $c$. The complex orientation gives
$i^!\Q_X\cong\Q_Z[-2c](-c)$, hence the \emph{Gysin map}
$i_*:H^j(Z)(-c)\to H^{j+2c}(X)$, with \emph{self-intersection formula}
$i^*i_*x=c_c(N_{Z/X})\smile x$. For a smooth proper map $\pi:E\to X$ of
relative dimension $b$, integration over the fibers is a \emph{trace}
$R\pi_*\Q_E[2b](b)\to\Q_X$. We use Deligne's period convention
$\Q(a)=(2\pi i)^a\Q$. Under the de Rham comparison, the trace
$H_c^{2d}(X,\Q(d))\to\Q$ for a smooth $d$-fold is
$(2\pi i)^{-d}\int_X$. Trace and Gysin maps retain the displayed Tate
twists; their filtered de Rham and topological realizations agree
\cite[Propositions~2.5.3 and~2.5.5 and Remark~2.5.7]{Sai88}.
Pullback followed by Gysin map adds an incidence condition; restriction
followed by trace removes it. The same
statements hold on smooth Deligne--Mumford stacks with finite
stabilizers, by averaging on finite quotient charts.

A \emph{cohomological correspondence} from $K\in D^b_c(X)$ to
$L\in D^b_c(Y)$ consists of a space $Z$, maps $c_1:Z\to X$ and $c_2:Z\to
Y$ with $c_2$ proper, and a morphism $c_1^*K\to c_2^!L$. By adjunction
it induces $\HH^*(X,K)\to\HH^*(Y,L)$. Knowing the action of a
correspondence on fibers does not determine it; the coefficient morphism
is part of the data.

\emph{Gluing over a curve.} For a map of complexes $u:\mathsf D\to\mathsf
R$, its \emph{homotopy fiber} is $\fib(u)=\Cone(u)[-1]$, realized as
$\mathsf D\oplus\mathsf R[-1]$ with differential $(x,y)\mapsto(dx,u(x)-dy)$.
If a curve is covered by an open surface $O$ and small disks around
$D$, Mayer--Vietoris identifies its cochains with
\[
 \fib\Bigl(C^*(O)\oplus C^*(\text{disks})\xrightarrow{\ \res-\,\mathrm{att}\
 }C^*(\text{boundary circles})\Bigr).
\]
To define an endomorphism of this complex one must give maps on $O$, on
the disks and on the circles, \emph{and} homotopies comparing the two
composites on each circle. The following example shows that the
homotopies carry information.

\begin{example}[What the pieces can miss]\label{ex:collar}
Let $u\in H^2(C,\Q)$ with $\int_Cu=1$. The map $u:\Q_C\to\Q_C[2]$
restricts to zero on $O$ and on each disk (both have $H^2=0$), but its
trace over $C$ is $1$. In the homotopy-fiber model the whole of $u$ is
carried by a homotopy on the boundary circles (a ``collar term'').
\end{example}

\subsection{Finite group actions}\label{subsec:finite}

\begin{lemma}\label{lem:finite-quotient}
Let a finite group $\Gamma$ act algebraically on a variety $X$, with
quotient $q:X\to Y=X/\Gamma$.
\begin{enumerate}[label=(\roman*)]
\item $H^*(Y)=H^*(X)^\Gamma$ and $H^*_c(Y)=H^*_c(X)^\Gamma$ as mixed
Hodge structures. If $\Gamma$ acts linearly on $\A^e$ then
$H_c^*(\A^e/\Gamma)=\Q(-e)[-2e]$.
\item If $X$ is smooth, $Y$ is rationally smooth.
\item If $X$ is of pure dimension $d$, then
$(q_*\IC_X)^\Gamma\cong\IC_Y$; if $X$ is rationally smooth, the same
holds for $\Q_X[d]$.
\end{enumerate}
\end{lemma}

\begin{proof}
Since $q$ is finite and proper, proper base change computes its
derived stalks on finite fibers; hence $R^iq_*=0$ for $i>0$ on sheaves
and $q_*$ is exact, and the averaging operator
$|\Gamma|^{-1}\sum_\gamma\gamma$ is an exact projector onto
$(q_*\Q_X)^\Gamma=\Q_Y$; this gives (i), and the formula for $\A^e$
follows because algebraic actions preserve the complex orientation. For
(iii), a finite map is small, so $q_*\IC_X$ is the intermediate extension
of its restriction to a dense open set over which $q$ is a covering;
there the invariant part is the constant sheaf on each component of $Y$,
and intermediate extension commutes with the exact projector. Statement
(ii) is (iii) for $\IC_X=\Q_X[d]$.
\end{proof}

When $\Gamma$ is abelian, the isotypic summand for a character
$\chi:\Gamma\to E^*$ ($E$ a cyclotomic field) is cut out by
$e_\chi=|\Gamma|^{-1}\sum_\gamma\chi(\gamma)^{-1}\gamma$; over $\Q$ one
sums $e_\chi$ over a Galois orbit. These projectors are algebraic, hence
preserve $W_\bullet$ by strictness. A finite direct image is perverse
$t$-exact, so if $\Gamma$ acts over the Hitchin base, the projectors
commute with perverse truncation and preserve $P_\bullet$. An ample
class can be replaced by the sum of its $\Gamma$-translates, so relative
hard Lefschetz holds on every isotypic summand. When an action permutes
connected components, the projectors are applied to the union of all
components; choosing one component would discard summands.

\subsection{Families}\label{subsec:families}

A \emph{local system} is a locally constant sheaf of finite-dimensional
vector spaces; parallel transport along paths identifies its stalks up
to homotopy. To compare filtrations in a family $f:\mathcal X\to S$ one
needs (a) that $R^jf_*\Q$ is a local system whose stalks are the
cohomology of the fibers (\emph{fiber base change}), and (b) that the
filtration steps are sub-local systems. When $f$ is proper, (a) is proper
base change plus local constancy. Our families of moduli spaces are not
proper over their parameter spaces, and smoothness alone does not give
(a); we prove (a) separately (Lemma~\ref{lem:basechange}).

For (b) on the weight side, the cohomology of a smooth mixed Hodge
module complex has a flat weight filtration. Indeed, the regular
connection underlying each weight subobject is a subconnection, and
its polarized pure quotients are variations of Hodge structure
\cite[Lemmas~5.1.10 and~5.2.12]{Sai88}.
The stronger equivalence with admissible variations is given by
\cite[Theorem~3.27 and the following algebraic remark]{Sai90}.
Only flatness and polarized pure quotients are needed below; the scalar
comparison will take place at an ordinary point of the extended
coefficient diagram (Lemma~\ref{lem:scalar-coefficient-extension}).

For (b) on the perverse side we use local acyclicity. Let
$\phi:Y\to S$ be a morphism to a smooth curve, choose $0\in S$, and
write $i:Y_0=\phi^{-1}(0)\hookrightarrow Y$. The \emph{nearby-cycle
complex} $\psi_\phi K$ records the cohomology of nearby fibers in small
neighborhoods of points of $Y_0$. The \emph{vanishing-cycle complex}
$\varphi_\phi K$ fits into the canonical specialization triangle
\[
 i^*K\longrightarrow\psi_\phi K\longrightarrow\varphi_\phi K
 \longrightarrow i^*K[1]
\]
\cite[\S\S3.4.10--3.4.12]{Sai88}. These are our unshifted conventions.
For the Hodge-module coefficients used below, $\psi_\phi[-1]$ and
$\varphi_\phi[-1]$ are exact on the perverse heart
\cite[(2.2.6) and Proposition~2.5]{Sai90}.
The complex $K$ is \emph{universally locally acyclic (ULA)} over $S$ if,
after every base change, its vanishing cycles along every curve vanish.
A smooth morphism has ULA constant sheaf, and proper direct image over
$S$ preserves ULA. For a
proper map over a smooth connected curve with semisimple coefficient,
lisse total direct image and vanishing $\varphi$, the perverse filtration
on fibers is locally constant \cite[Theorem~3.2.1(iii)]{dCMa20}.

Finally we need a topological recognition principle. A
\emph{rational cohomology manifold} of dimension $n$ is a locally compact
Hausdorff space $X$ of finite cohomological dimension whose dualizing
complex is $\omega_X\cong\mathfrak o_X[n]$ for a rank-one local system
$\mathfrak o_X$ (the orientation sheaf); equivalently its local groups
$H^i(X,X\setminus x;\Q)$ are $\Q$ for $i=n$ and $0$ otherwise.
The recognition principle needed here is proved below.

\begin{lemma}[Invariance of domain]\label{lem:domain-invariance}
Let $X,Y$ be locally contractible rational cohomology manifolds of the
same dimension $n$, without boundary. Then every continuous injective map
$f:X\to Y$ is open. In particular a continuous bijection $X\to Y$ is a
homeomorphism.
\end{lemma}

\begin{proof}
Let $x\in X$, and let $K$ be a compact neighborhood of $x$. Then $f|_K$
is a homeomorphism onto its image. Choose an open neighborhood $W$ of
$f(x)$ disjoint from the compact set $f(K\setminus\operatorname{int}K)$
and with $\mathfrak o_Y|_W$ trivial. Then $A=f(K)\cap W$ is closed in $W$
and homeomorphic to the open subset $f^{-1}(W)\cap\operatorname{int}K$ of
$X$, so $A$ is a rational cohomology $n$-manifold. For the closed
inclusion $i:A\hookrightarrow W$,
$i^!\Q_W[n]\cong i^!\omega_W\cong\omega_A\cong\mathfrak o_A[n]$; thus
$i^!\Q_W$ is a rank-one local system in degree $0$, and its stalk at
$f(x)$, which is the stalk of the sheaf $\underline\Gamma_A\Q_W$ of
sections supported in $A$, is nonzero. A nonzero germ of
$\underline\Gamma_A\Q_W$ at $f(x)$ is a locally constant function near
$f(x)$, nonzero on every neighborhood and supported in $A$; on a
connected neighborhood it is a nonzero constant, so that neighborhood
lies in $A\subset f(X)$. Applying this to $f$ restricted to arbitrary
open subsets shows that $f$ is open.
\end{proof}

The no-boundary hypothesis is the reason why real parameters below are
taken in open intervals.

\subsection{Moduli spaces}\label{subsec:moduli}

\emph{Betti spaces.} For an affine variety $Y$ with $G$-action,
$Y\sslash G=\Spec\C[Y]^G$ is the affine GIT quotient defined by invariant
regular functions. The space~\eqref{eq:intro-MB} is the coarse space of
the quotient stack $[\widetilde{\MB}/G]$ of the representation variety
$\widetilde\MB$; the stack remembers automorphisms. Since $Z(G)$ acts
trivially, we \emph{rigidify}, i.e.\ divide the automorphism groups by
$Z(G)$. When the boundary data are generic, the rigidified stack has
finite stabilizers:

\begin{lemma}\label{lem:generic}
Assume the data of Definition~\ref{def:intro-generic} are intrinsically
generic. Then every point of $\widetilde\MB$ is \emph{irreducible}: its
image is contained in no proper parabolic subgroup. Consequently its
$G$-orbit is closed and its stabilizer is finite modulo $Z(G)$.
\end{lemma}

\begin{proof}
Suppose all $A_i,B_i,M_p$ lie in a proper parabolic $P$, and choose a
Levi factor $L_P$. Under the Levi projection, the semisimple part of
$M_p$ has image $s_p\in L_P$ which is $G$-conjugate to $S_p$.
The unipotent part maps to $1$ in the torus $L_P^\ab$, as do all
commutators. Projecting the surface relation to $L_P^\ab$ therefore
gives $\prod_p\pi(s_p)=\pi(\zeta)$, contrary to intrinsic genericity.
The tuple is irreducible. Lemma~\ref{lem:irreducible-closed-direct}, the
irreducible case of Richardson's closed-orbit criterion, now gives both
closedness of its orbit and finiteness of its stabilizer modulo $Z(G)$.
\end{proof}

\emph{Parabolic and parahoric Higgs bundles.} A principal
$G$-bundle is an algebraic family of $G$-torsors, locally a product with
$G$ in the \'etale topology; a torsor is a space with a free transitive
$G$-action. For $\GL_n$, principal bundles are frame bundles of vector
bundles. For a vector bundle, a
\emph{parabolic structure} at $p$ is a weighted flag in the fiber $E_p$,
with weights in $[0,1)$; the \emph{parabolic degree} is
$\deg E+\sum_{p,a}m_{pa}\vartheta_{pa}$, where $\vartheta_{pa}$ are the
weights with multiplicities $m_{pa}$, and a parabolic Higgs bundle is
\emph{stable} if every proper Higgs-invariant subbundle has smaller
parabolic slope. For general $G$, subbundles are replaced by reductions to
parabolic subgroups compatible with the Higgs field, and slopes by
degrees of the associated character line bundles \cite{BGPM20,KSZ24}.

Rational weights admit two equivalent descriptions. On a \emph{root
curve} $\mathcal C$ the point $p$ with weight denominator $e$ is replaced
by the chart $[\Spec\C[w]/\mu_e]$ with $z=w^e$; a bundle on this chart is
a $\mu_e$-equivariant bundle on the $w$-disk, and the action on the fiber
at $w=0$ is its \emph{inertia homomorphism}
\cite{Cad07,BV12}, \cite[Proposition~6.1]{BMW12}. For a rational point
$\theta\in X_*(T)\otimes\Q$, the connected subgroup generated by the
torus and root lattices is the usual \emph{Bruhat--Tits parahoric}
\cite{BS15,KSZ24}. We use the full root-chart automorphism group of
Lemma~\ref{lem:central-full-local}, which can have additional components
when $G$ is not simply connected. These components are retained; the
root-curve description is not identified with only the connected
parahoric. The Lie algebra lattice of this full group and that of the
pro-unipotent radical of the connected parahoric subgroup are
\begin{equation}\label{eq:lattices}
 \mathcal P=\ft[[z]]\oplus\bigoplus_\alpha
   z^{\lceil\alpha(q)\rceil}\fg_\alpha[[z]],\qquad
 \mathcal S=z\ft[[z]]\oplus\bigoplus_\alpha
   z^{\lfloor\alpha(q)\rfloor+1}\fg_\alpha[[z]],
\end{equation}
where $q=-\theta$ is the \emph{growth weight}. A logarithmic Higgs field
is a section of $\mathcal P\otimes K_C(D)$; it is \emph{strongly
parabolic} if it lies in $\mathcal S\otimes K_C(D)$. A nondegenerate
invariant form identifies $\mathcal P^\vee=\mathcal S(D)$ and
$\mathcal S^\vee=\mathcal P(D)$, because $-\lceil x\rceil=\lfloor-x\rfloor$.

\emph{$\lambda$-connections.} For $\lambda\in\C$, a
\emph{$\lambda$-connection} on a vector bundle is a $\C$-linear map
$\nabla_\lambda:\mathcal E\to\mathcal E\otimes K_C(D)$ with
$\nabla_\lambda(fs)=\lambda\,df\otimes s+f\nabla_\lambda s$; for
$\lambda=0$ this is a Higgs field and for $\lambda\ne0$,
$\lambda^{-1}\nabla_\lambda$ is a logarithmic connection
\cite{Sim97}. Under a gauge transformation $g$, the connection form
changes by $A\mapsto\Ad(g^{-1})A+\lambda g^{-1}dg$. The moduli spaces of
semistable principal $\lambda$-connections form a family over $\A^1_\lambda$
(the \emph{Hodge moduli space}) with a $\Gm$-action scaling $\lambda$ and
the Higgs field \cite{FHZ26}.

\emph{Residues and the residue theorem.} At a puncture with growth weight
$q$, the regularizing frame $e_{\rm reg}=ez^{q}$ changes the connection form
$A$ to $\Ad(z^{-q})A+\lambda q\,dz/z$. We call the residue in this frame
the \emph{affine residue} $R_p$. For a character $\chi\in X^*(G)$, let
$\deg_\chi$ be the degree of the associated line bundle and put
$\delta(\chi)=\deg_\chi-\sum_pd\chi(q_p)$ (the parabolic degree). The
residue theorem applied to the associated line connection gives the
\emph{Fuchs relation}
\begin{equation}\label{eq:fuchs}
 \sum_{p\in D}d\chi(R_p)=-\lambda\,\delta(\chi)\qquad(\chi\in X^*(G)).
\end{equation}
Relation~\eqref{eq:fuchs} is the only obstruction to deforming the
central part of an object
(\S\ref{subsec:A-restrictions}).

\subsection{The tame nonabelian Hodge correspondence}\label{subsec:nah}

For $\GL_n$, consider a semisimple block of a local monodromy with
eigenvalue $c\ne0$. Choose a real number $\vartheta$ (a lift of the
phase) and write
\begin{equation}\label{eq:nah-dictionary}
 c=e^{-2\pi i\vartheta}e^{-4\pi i\beta},\qquad
 \beta=\frac{i}{4\pi}\log|c|\in i\R .
\end{equation}
The corresponding parabolic Higgs bundle has weight $\vartheta$ (growth
weight $a=-\vartheta$) and Higgs residue eigenvalue $\beta$; changing
$\vartheta$ by an integer changes the extension across the puncture but
not the local system. These are the parameters $\alpha=-\vartheta$,
$s=\beta$ of \cite[Proposition~6.3, Table~1 and Proposition~7.7]{BGPM20}
(arXiv v4 numbering), in which the filtered local system has weight
$0$. In Mochizuki's Kashiwara--Malgrange--Sabbah--Simpson (KMS)
notation \cite[\S2.1.6 and \S\S7.2.1--7.2.2, Definition~7.2]{Moc07}, a real growth
weight $a$ and a semisimple Higgs residue $\beta$ produce, at the real
Hodge parameter $\tau$, the growth weight and residue of the
$\tau$-connection
\begin{equation}\label{eq:kms}
 p_\tau=a,\qquad e_\tau=(1+\tau^2)\beta-\tau a;
\end{equation}
indeed $p(\lambda)=a+2\operatorname{Re}(\lambda\bar\beta)$ and
$e(\lambda)=\beta-\lambda a-\lambda^2\bar\beta$, and for $\beta\in i\R$
and $\lambda=\tau\in\R$ these reduce to~\eqref{eq:kms}. At $\tau=1$ the
residue $2\beta-a$ satisfies $\exp(-2\pi i(2\beta-a))=c$. For general
$G$ the same formulas are applied to cocharacters, with centralizer
orbits in place of eigenvalue blocks.

\begin{theorem}[Tame nonabelian Hodge correspondence]\label{thm:NAH}
For compatible parameters there is a homeomorphism
$\Psi:\MD\to\MB$ which sends a polystable parabolic $G$-Higgs bundle to
the monodromy of the flat connection $\bar\partial_h+\partial_h+\varphi
+\varphi^{\dagger_h}$ of its adapted harmonic metric $h$; it matches
stable objects with irreducible representations
\cite{Sim90} ($\GL_n$), \cite[Theorems~8.1 and 11.4, \S9]{Biq97},
\cite[Table~1 and Theorem~7.8]{BGPM20}.
\end{theorem}

The map $\Psi$ is not algebraic, and does not respect mixed Hodge
structures. In the proof we need $\Psi$ in families and on strata of
nilpotent orbit types; this is Theorem~\ref{thm:harmonic-transport}.

\emph{The Hitchin map.} Let $p_1,\ldots,p_\ell$ be homogeneous generators
of $\C[\fg]^G$, of degrees $d_1,\ldots,d_\ell$. Evaluation on $\varphi$
gives
\[
 h:\MD\longrightarrow\mathcal A_{K_C(D)}=\bigoplus_iH^0(C,K_C(D)^{d_i}),
\]
and fixing the residues restricts $h$ to an affine subspace
$\mathcal A\subset\mathcal A_{K_C(D)}$. For a vector bundle these are the
coefficients of the characteristic polynomial of $\varphi$. The image
dimension and the defect of $h$ for the spaces used here are computed in
Theorem~\ref{thm:fine-master}.

\emph{Terminology for families.} A \emph{fine} boundary datum at $p$
consists of a full flag (a Borel reduction) at $p$ together with an
ordered Cartan value: an element of $T$ (Betti side) or of $\ft$ (residue
side). A \emph{weak} fine datum allows, in addition, an arbitrary
unipotent (resp.\ nilpotent) part preserving the flag. A \emph{residue
master} is the family of fine weak moduli spaces obtained by letting
the ordered Cartan values vary with fixed stability weights, subject
to~\eqref{eq:fuchs}. A \emph{fine reference} is a member of such a
family with regular semisimple boundary values, used to anchor the
comparison. A fine reference is \emph{separating} if no proper parabolic
subgroup satisfies a character relation among the boundary values of the
reference: the
multiplicative relation of Definition~\ref{def:intro-generic} for the
monodromies along every parameter path used, and its additive analog
\eqref{eq:additive-generic} for the Higgs residues. The families are
constructed in Section~\ref{sec:comparison-families}.

\section{A common lowering operator determines the filtrations}
\label{sec:linear}

Two filtrations that both satisfy hard Lefschetz for the same operator
$E$, with the same center, need not coincide: on
$V=\Q e_0\oplus\Q e_1\oplus\Q e_2$ with $E(e_0)=e_2$, $E(e_1)=E(e_2)=0$,
both $L_0=\langle e_0\rangle\subset L_1=\langle e_0,e_1\rangle\subset
L_2=V$ and $L'_0=\langle e_0+e_1\rangle\subset L'_1=\langle e_0,e_1
\rangle\subset L'_2=V$ are Lefschetz with center $1$. This section shows
that the ambiguity disappears once a single lowering operator lowers
both filtrations and satisfies one commutator identity. The argument
uses no splitting of the decomposition theorem; it is a filtered version
of \cite[Proposition~8.3 and Lemma~8.4]{HMMS}; the proof below
is independent of the general complete-reducibility theorem.

We write $[A,B]=AB-BA$. An \emph{$\mathfrak{sl}_2$-triple} $(e,h,f)$
satisfies $[e,f]=h$, $[h,e]=2e$, $[h,f]=-2f$. A \emph{Lefschetz string} is a basis $v_0,\ldots,v_m$ on which
$ev_i=v_{i+1}$, $fv_i=i(m-i+1)v_{i-1}$ and $hv_i=(2i-m)v_i$, with
$v_{-1}=v_{m+1}=0$. These formulas directly satisfy the three relations.

Throughout this section $V$ is a finite-dimensional vector space over a
field of characteristic zero, with a finite increasing filtration
$L_\bullet$. An endomorphism $T$ has \emph{filtration degree at most $a$}
if $T(L_k)\subset L_{k+a}$ for all $k$. Let $E$ have filtration degree at
most $2$, and write $e:\Gr^L_k\to\Gr^L_{k+2}$ for its leading part.

\begin{definition}\label{def:lefschetz-module}
The pair $(L_\bullet,E)$ is a \emph{filtered Lefschetz module with center
$N$} if $e^r:\Gr^L_{N-r}V\to\Gr^L_{N+r}V$ is an isomorphism for every
$r\ge0$. If $V=\bigoplus_jV^j$ is graded and $E:V^j\to V^{j+2}$, we
mean $e^r:\Gr^L_{N-r}V^j\xrightarrow{\sim}
\Gr^L_{N+r}V^{j+2r}$ in each cohomological degree $j$.
\end{definition}

Then $\Gr^LV$ is a representation of $\mathfrak{sl}_2$ in which $e$ acts
as the raising operator and $h$ acts on $\Gr^L_kV$ by $k-N$: on
$\Gr^L_{N-m}$ the kernel of $e^{m+1}$ is a complement to
$e\,\Gr^L_{N-m-2}$: the restriction of $e^{m+1}$ to this
last subspace is an isomorphism onto $\Gr^L_{N+m+2}$, by the Lefschetz
isomorphism $e^{m+2}$ on $\Gr^L_{N-m-2}$. Iterating gives the strings
above and the indicated operators $f,h$.

\begin{lemma}[Negative-degree injectivity]\label{lem:negative-injectivity}
Let $(L_\bullet,E)$ be a filtered Lefschetz module. An endomorphism $T$
of filtration degree at most $-2$ with $[E,[E,T]]=0$ is zero. An
endomorphism of filtration degree at most $-1$ commuting with $E$ is
zero.
\end{lemma}

\begin{proof}
Scalar extension reduces the linear identities to $\C$: the finitely
many coefficients lie in a finitely generated characteristic-zero field
which embeds in $\C$. Give the primitive strings a positive Hermitian
form with $e^*=f$, and take their orthogonal direct sum. On the full
endomorphism space, including maps between different strings, the
Hilbert--Schmidt form then makes $u=\ad(e)$ and $v=\ad(f)$ adjoints.
For a homogeneous endomorphism $x$ of weight $a$, $[u,v]=\ad(h)$ gives
\[
                  \|ux\|^2-\|vx\|^2=-a\|x\|^2.
\]
Thus $u$ is injective for $a<0$, and $u^2$ is injective for $a<-2$.
If $a=-2$ and $u^2x=0$, put $y=ux$. Then $y$ has weight zero and
$uy=0$, so the same identity gives $vy=0$. Hence
$\|y\|^2=\langle x,vy\rangle=0$, and negative-weight injectivity
gives $x=0$.

If $T\ne0$, choose the least filtration degree $a$ for which
$T(L_k)\subset L_{k+a}$ for every $k$. The leading symbol of $T$ on
$\Gr^LV$ is then nonzero of weight $a$. The assumed commutator puts that symbol in
the injective kernel of $u^2$ for $a\le-2$, or of $u$ for $a<0$,
respectively. This is a contradiction.
\end{proof}

\begin{proposition}[Canonical normalization]\label{prop:canonical-lowering}
Let $(L_\bullet,E)$ be a filtered Lefschetz module with center $N$, and
let $F$ have filtration degree at most $-2$ with
\begin{equation}\label{eq:double-commutator}
 [[F,E],E]=-2E.
\end{equation}
Then $(E,H,F)$ with $H=[E,F]$ is an $\mathfrak{sl}_2$-triple, $H$ is
diagonalizable with integer eigenvalues, and
\begin{equation}\label{eq:cartan-splitting}
 L_kV=\bigoplus_{a\le k-N}\ker(H-a).
\end{equation}
Such an $F$ is unique, and every endomorphism preserving $L_\bullet$ and
commuting with $E$ commutes with $F$ and $H$.
\end{proposition}

\begin{proof}
Let $f$ be the degree $-2$ part of $F$ on $\Gr^LV$ and $f_0$ the lowering
operator of the $\mathfrak{sl}_2$-structure on $\Gr^LV$. Taking leading
parts in~\eqref{eq:double-commutator} gives $[[f,e],e]=-2e=[[f_0,e],e]$,
so $f=f_0$ by Lemma~\ref{lem:negative-injectivity} applied on $\Gr^LV$.
Hence $H$ preserves $L_\bullet$ and acts on $\Gr^L_k$ by $k-N$.
Equation~\eqref{eq:double-commutator} says $[H,E]=2E$. The element
$R=[H,F]+2F$ has filtration degree at most $-2$, and the Jacobi identity
gives $[E,R]=[[E,H],F]+[H,[E,F]]+2[E,F]=-2H+0+2H=0$; so $R=0$ by the
lemma, and $(E,H,F)$ is an $\mathfrak{sl}_2$-triple.

The product $\prod_k(H-(k-N))$ over the occupied graded degrees $k$
annihilates $V$ (apply the factors in decreasing order of $k$; each kills
the top remaining graded piece). This polynomial has distinct roots, so $H$ is
diagonalizable and its spectral projectors, being polynomials in $H$,
preserve $L_\bullet$. If $v\ne0$ is an $H$-eigenvector with eigenvalue
$a$ and $k$ is minimal with $v\in L_k$, the image of $v$ in $\Gr^L_k$ is
nonzero and $H$ acts there by $k-N$, so $k=a+N$. This
proves~\eqref{eq:cartan-splitting}.

If $F'$ is another such operator, $F-F'$ has filtration degree at most
$-2$ and is killed by $(\ad E)^2$, hence vanishes. If $\Pi$ preserves
$L_\bullet$ and commutes with $E$, then $[F,\Pi]$ has filtration degree
at most $-2$ and $(\ad E)^2[F,\Pi]=[[[F,E],E],\Pi]=-2[E,\Pi]=0$, so
$[F,\Pi]=0$ and then $[H,\Pi]=0$.
\end{proof}

\begin{theorem}[Common-lowering criterion]\label{thm:common-lowering}
Let $P_\bullet$ and $L_\bullet$ be two finite filtrations of $V$ such
that $(P_\bullet,E)$ and $(L_\bullet,E)$ are filtered Lefschetz modules
with the same center $N$. If a single endomorphism $F$ lowers both
filtrations by two and satisfies~\eqref{eq:double-commutator}, then
$P_\bullet=L_\bullet$. The same holds on the image of every idempotent
that preserves both filtrations and commutes with $E$.
\end{theorem}

\begin{proof}
By Proposition~\ref{prop:canonical-lowering} both filtrations equal
$\bigoplus_{a\le k-N}\ker([E,F]-a)$. An idempotent $\Pi$ as in the
statement commutes with $[E,F]$ by the last assertion of the
proposition, so it preserves the eigenspace decomposition.
\end{proof}

In the application $L_k=W_{2k}$, odd weight-graded pieces vanish, and the
theorem gives $P_k=W_{2k}=W_{2k+1}$. All operators act on one vector
space, identified through the given map $\Psi$; an abstract isomorphism
of associated graded spaces would not suffice.

\begin{example}[Symplectic tori]\label{ex:torus-sl2}
Let $V=\Lambda^\bullet(\Q^{2g})^\vee$ be the cohomology of
$(\C^*)^{2g}$, with symplectic basis $\xi_i,\eta_i$ of degree one, and
$L_k=\bigoplus_{j\le k}\Lambda^j$. Let $E$ be multiplication by
$\omega=\sum_i\xi_i\eta_i$ and $F=\sum_i\iota_{\eta_i^*}\iota_{\xi_i^*}$
contraction with the dual bivector. Then $[E,F]=j-g$ on $\Lambda^j$; this
is the classical Lefschetz decomposition of an exterior algebra, and
$(L_\bullet,E)$ is a filtered Lefschetz module with center $g$. For the
torus of Example~\ref{ex:intro-rank-one}, both $P_\bullet$ and
$W_{2\bullet}$ equal $L_\bullet$.
\end{example}


\section{Families and comparison maps}\label{sec:comparison-families}\label{app:families}

Pointwise nonabelian Hodge theory does not identify the two filtrations;
we need identifications induced by maps of families. This section
establishes the following three properties of families of moduli spaces,
which are used throughout Sections~\ref{sec:raising}--\ref{sec:proof}.
They are theorems proved below, not additional hypotheses.

\begin{enumerate}[label=(F\arabic*)]
\item\label{F:betti} (Theorem~\ref{thm:weak-betti}.) Over the
intrinsically generic locus of the ordered Cartan parameter space of the
weak Betti family, ordinary and compactly supported cohomology form
local systems with canonical fiber base change, and the weight
filtration consists of flat sub-local systems. All compatible central
lifting sectors are included.
\item\label{F:dol} (Theorem~\ref{thm:fine-master}.) The rational fine
Dolbeault residue master (\S\ref{subsec:nah}), after rigidification, is smooth with finite stabilizers,
has a proper relative Hitchin map, and its fiber restrictions identify
the perverse filtrations. The constant sheaf of the master is ULA over the residue
base, and the Hitchin defect is half the fiber dimension.
\item\label{F:harm} (Theorem~\ref{thm:harmonic-transport}.) Continuous
families of harmonic metrics identify the relevant paths of Betti and
Dolbeault fibers with the nonabelian Hodge map, compatibly with finite
projectors and, along paths with fixed graded centralizer, with the
Springer multiplicities giving the intersection complexes of prescribed
nilpotent orbit closures.
\end{enumerate}

The section also constructs the residue masters and compactification
pairs used in~\ref{C:master}. We first establish smoothness over the residues
(\S\ref{subsec:A-restrictions}) and construct the masters and their
compactification pairs (\S\ref{subsec:B-masters}). Fiber inclusions then
give ordinary base change (Lemma~\ref{lem:basechange}), proper Hitchin
direct images give perverse transport (Theorem~\ref{thm:fine-master}),
and Riemann--Hilbert in families gives flat weights
(Theorem~\ref{thm:weak-betti}). The fixed-parameter harmonic inverse
(\S\ref{app:harmonic}) and harmonic compactness identify the resulting
transports with $\Psi$, also on intersection cohomology
(Theorem~\ref{thm:harmonic-transport}).

Moduli spaces are coarse spaces of centrally rigidified stacks; local
inertia homomorphisms, topological components and central lifting sectors
are kept, and actions permuting sectors act on their union.

\emph{The fine weak Hodge master.} A rational parahoric type is
\emph{fine} if the inertia centralizer at every puncture is $T$. Choose
stability weights outside all numerical walls, subject to the central
degree equations. The \emph{fine weak Hodge master} is the family of
parahoric $\lambda$-connections with full boundary flags and variable
ordered Cartan residues, the nilpotent part preserving the flag being
unrestricted. Residues are affine-normalized as in~\S\ref{subsec:moduli}:
on a ramified disk, the equivariant residue divided by the ramification
index, in the regularizing frame $ez^q$. The parameter space $\Lambda$ is
the affine space cut out by~\eqref{eq:fuchs}. In the masters used below
all character parabolic degrees vanish, $\delta(\chi)=0$, while the
ordinary degrees $\deg_\chi$ keep their prescribed values. We write $a$
for the real growth weight of the harmonic metric and $q$ for a nearby
rational weight used to build an algebraic family; changing $q$ is part of
the construction and does not change the endpoint problem.

\begin{lemma}[Rational charts for real weights]\label{lem:rational-chart}
For intrinsically generic full spectral data, the real weights $a$ admit
arbitrarily close rational weights $q$ with the same parahoric lattices,
the same full (extended affine) stabilizers, the same central degree
equations, and the same stable locus.
\end{lemma}

\begin{proof}
Work in a bounded affine chart. Root exponents in~\eqref{eq:lattices} are
constant on the facet of $a$. The full stabilizer is described by the
finitely many equations $wa-a=\mu$ ($w\in\WG$, $\mu\in X_*(T)$) that can
occur in the chart; keep those satisfied at $a$ and avoid acquiring new
ones. With the central degree equations these define a rational affine
face containing $a$, with the same local gauge group and component action.
For each conjugacy type of maximal parabolic there are finitely many
character tests and relative positions, with degrees in a fixed discrete
lattice, so only finitely many wall equations $d+\ell(q)=0$ meet the chart
($\ell$ includes the sign of the growth-weight contribution). For a
reduction compatible with the Higgs residue $\beta$, the residue theorem
gives $\ell(\beta)=0$, and intrinsic genericity excludes $d+\ell(a)=0$
(exponentiating both would give a forbidden character product). Preserving
these finitely many signs is an open condition in the face, and rational
points are dense there.
\end{proof}

\subsection{Smoothness over the residues}
\label{subsec:A-restrictions}

At fixed $\lambda$ and fixed ordered Cartan residues the deformation
complex of the fine master is
\begin{equation}\label{eq:fine-complex}
 \bigl[\mathcal P\xrightarrow{d_{\nabla_\lambda}}\mathcal S\otimes
 K_C(D)\bigr],
\end{equation}
with $\mathcal P,\mathcal S$ as in~\eqref{eq:lattices}. A \emph{framing}
of the torsor of the reductive quotient $L$ of the parahoric at each
puncture removes the gauge transformations that are nontrivial in $L$;
prescribing the whole residue matrix, the complex becomes
\begin{equation}\label{eq:framed-complex}
 \bigl[\mathcal S\xrightarrow{d_{\nabla_\lambda}}\mathcal S\otimes
 K_C(D)\bigr].
\end{equation}
The residue evaluation on framed charts is therefore smooth (see below),
which allows smooth pullback of nilpotent cones and orbit closures in
$\Lie L$.

\emph{Obstructions at a Higgs point.}
At $\lambda=0$ the differential is $\mathcal O_C$-linear.
Using $\mathcal P^\vee=\mathcal S(D)$ and
$\mathcal S^\vee=\mathcal P(D)$, curve Serre duality identifies the dual
of $\mathbb H^2$ of either fixed-residue complex above with the kernel
of $[\Phi,-]$ on $H^0(\mathcal P)$. These are \emph{unframed}
infinitesimal automorphisms; a boundary frame does not remove the
central kernel of this dual complex. For the variable whole-residue
framed complex
$[\mathcal S\to\mathcal P\otimes K_C(D)]$ the dual kernel instead lies
in $H^0(\mathcal S)$. The distinction follows directly by dualizing the
two terms, and is important for the residue evaluation.

On the stable locus the component of the infinitesimal automorphism
algebra in $\fg_{\der}$ vanishes. For the zero-residue Higgs points
needed below, this also follows from the compact, finite-equivariant
degree-zero correspondence: after passage to the ramified cover the
Higgs field has no pole. Stability excludes a proper equivariant
parabolic reduction of degree zero, so the full equivariant reductive
monodromy group lies in no proper parabolic. Its centralizer is
reductive; a noncentral torus in that centralizer would put the
monodromy in a proper Levi subgroup. Its Lie algebra is therefore
$\mathfrak z(\fg)$. Finite residual automorphisms are allowed.

The central obstruction is precisely the line-connection residue
equation~\eqref{eq:fuchs}. For a line bundle of degree $d$ the
logarithmic Atiyah cocycle is $\lambda\,d\log g_{ij}$, and the residue
sequence sends a residue tuple $(r_p)$ to $\sum_pr_p$ in
$H^1(K_C)=\C$. Thus the obstruction is
$\lambda d+\sum_pr_p$. On its zero hyperplane, solutions form locally
a torsor under $H^0(K_C)$ over a Picard component. This is smooth also
at zero residues. Apply the calculation in a character basis of the
central torus. A chosen finite central lifting sector adds no
infinitesimal obstruction, since finite \'etale data and their cocycle
identities lift uniquely over nilpotent thickenings. The same
calculation on the ramified cover descends by exact finite invariants.
Consequently the master is smooth over $\Lambda$ at zero-residue Higgs
points.

\emph{Propagation to the full master.}
The existence and zero-scaling-limit construction in
Proposition~\ref{prop:pg-master} uses boundedness, canonical
Harder--Narasimhan strata and semistable reduction, independently of
smoothness over the residue space. Every point has, after a finite
parameter extension if necessary, a scaling limit with $\lambda=0$
and all ordered residues zero. The chosen stability chamber makes
this limit stable. The smooth locus of the finite-presentation
residue morphism is open and scaling invariant. Its complement is
closed and invariant; a point of that complement would have a
zero-scaling limit in the complement, contrary to the preceding
calculation. This proves smoothness over the entire compatible
affine residue space, without dualizing an
$\mathcal O_C$-nonlinear connection differential.

\emph{Whole-residue evaluation on the stable coarse-type charts.}
For an arbitrary retained parahoric type, work on the chart with a frame
in every reductive quotient.  At fixed whole residues and fixed $\lambda$
its complex is \eqref{eq:framed-complex}.  On a curve its obstruction
space is the cokernel of the map on $H^1$ induced by
$d_{\nabla_\lambda}:\mathcal S\to\mathcal S\otimes K_C(D)$.
Under the ordinary Serre pairings on the two bundles, the transpose
of this map is, up to sign, the map on $H^0$ induced by
$d_{\nabla_\lambda}:\mathcal P\to\mathcal P\otimes K_C(D)$.
Indeed, resolve both complexes by their Dolbeault double complexes.
Pair Dolbeault representatives by $\kappa$ and integrate:
invariance cancels the bracket terms, and the two differential terms
sum, up to the differential-form sign, to $\lambda\partial$ of the
pairing.  This term integrates to zero by Stokes.
The pairing of a local section of $\mathcal S$ with one of
$\mathcal P$ lies in $\mathcal O_C(-D)$, so the logarithmic terms in
the Serre pairings have no poles.  The calculation also descends on
taking finite deck invariants.
Thus the dual obstruction space consists of the unframed horizontal
infinitesimal automorphisms.  For a stable Higgs point its derived
Lie algebra part vanishes.  For $\lambda\ne0$ on the generic paths
used here it vanishes because the monodromy lies in no proper
parabolic.  The remaining central equations are precisely
\eqref{eq:fuchs}, by the character calculation above.  More explicitly,
the splitting $\fg=\mathfrak z(\fg)\oplus[\fg,\fg]$ is respected
by $d_{\nabla_\lambda}$.  For a square-zero extension the derived
obstruction vanishes and the character maps detect the entire
central obstruction, since their differentials span the dual center.
On the compatible central target the character-line deformations lift
by the logarithmic Atiyah calculation.  The discrete central lifting
data lift uniquely across nilpotent thickenings.
Thus the framed evaluation, including variation of $\lambda$, is
smooth at these points over its central compatibility space.  This
argument uses stability at the point under consideration and applies
also when the stability weights lie on numerical walls.

\subsection{Residue masters, compactification and transport}
\label{subsec:B-masters}

For the constructions in this subsection, cocharacters are in $X_*(T)$ for the
given group $G$, the fine rational weights are those of
Lemma~\ref{lem:rational-chart} and \S\ref{subsec:A-transport} (outside the
walls, with inertia centralizer $T$), all character parabolic degrees
vanish, and the affine residues satisfy~\eqref{eq:fuchs}. The rational
root presentation gives a finite ramified logarithmic lift of each
parahoric chart, with its inertia action and allowed lattice changes; all
constructions are made on the lift and descended. Statements about
smoothness concern the rigidified stacks, which have finite stabilizers;
their coarse spaces may have finite quotient singularities.

A $\Gm$-variety is \emph{semiprojective} if its fixed locus is proper and
every orbit has a limit as the scaling parameter tends to $0$
\cite{HRV15}. The inputs are: existence, Harder--Narasimhan
stratification and semistable reduction for principal meromorphic
$\lambda$-connections \cite[Theorems~2.27 and~4.47,
Proposition~4.48]{FHZ26} (numbering of arXiv v2); properness of the
principal Hitchin map \cite[Corollary~6.21 and Remark~6.22]{AHLH23};
semistable reduction for $\Theta$-stratified stacks
\cite[Theorem~6.5]{AHLH23}; the local root-lattice and character-degree
calculations for logahoric Higgs fields
\cite[\S2.1 and Lemmas~4.10--4.11]{KSZ24}; and the relative quotient construction of
\cite[Appendix~B, Theorem~7.2.1]{HLRV11} and
\cite[Theorem~3.2.1 and Lemma~3.5.1]{dC21}. The full root-chart arrows are
those of Lemma~\ref{lem:central-full-local}; reductions and stability are
tested directly on the retained root-curve and equivariant objects.
The ordinary semistable moduli spaces identify objects up to
S-equivalence, that is, isomorphism of their polystable graded objects.
For logarithmic connections the Lie algebroid in the vector-bundle
construction is $T_{C/S}(-D)$, dual to $\Omega^1_{C/S}(D)$
(cf.\ \cite[\S2.4]{FHZ26}). For semistable reduction we use the
characteristic-zero slice argument of \cite[Remark~6.4 and
Lemma~6.10]{AHLH23} on the total stack over $\C$. This also works over
the moving parameter variety: maps from $[\A^1_T/\Gm]$ to that variety
factor through $T$, and a modified family over a discrete valuation
ring has the same base map by separatedness.
We retain the full root-chart arrows in the equivariant stack below. For $\lambda$-connections the
finite-cover construction uses the affine residue shift
$A\mapsto\Ad(g^{-1})A+\lambda g^{-1}dg$, with the local inertia
isomorphisms and their cocycle identities; the local formulas hold in any
genus.

For an arbitrary global form of $G$, the boundedness argument uses
representations of highest weight $m_j\omega_j$, where $\omega_j$ is a
fundamental weight in the rational character space, extended by zero
on the central directions, and $m_j>0$ is chosen
so that $m_j\omega_j\in X^*(T)$. These representations are representations
of $G$ itself. In the proof of \cite[Lemma~4.10(1)]{HLFH23}, the bound
on the highest-weight subbundle then bounds
$m_j\langle d_{E_P},\omega_j\rangle$; dividing by $m_j$ gives the
required bound on $\langle d_{E_P},\omega_j\rangle$, where
$d_{E_P}(\chi)=\deg E_P(\chi)$. The finite set of
multipliers only changes the uniform constants. This is the integral
multiple convention also used in the proof of Lemma~A.9 of that source, and it
preserves the subsequent boundedness and cone-optimization arguments.

\begin{proposition}[Geometry of the fine weak master]\label{prop:pg-master}
The nonempty fine weak master has a quasi-projective coarse space with
zero scaling limits and proper fixed locus, and its ordered residue map to
the affine space~\eqref{eq:fuchs} is smooth and surjective. The same holds
relatively over $\Cc$ when a disjoint auxiliary section $q\in\Cc$ with a
sufficiently small fine weight is added; then the fixed locus is proper
over $\Cc$.
\end{proposition}

\begin{proof}
\emph{The cover, including a moving section.}
Let the desired local ramification orders be $e_i$ and put
$n=\operatorname{lcm}_i e_i$. For $n=1$ the identity cover suffices.
Otherwise write $m_i=n/e_i$, with $p_i$ a fixed or moving section,
and let $\Delta=\sum_i m_i p_i$ on the product curve over the
parameter chart. Take a line bundle $L$ on $C$ so positive that
$L^n(-\Delta)$ is very ample on each fiber, has vanishing $H^1$,
and has positive degree. On a chosen fiber choose a section of
this bundle with simple zeros disjoint from the marked points.
Such sections form a nonempty open: evaluation at each marking is
surjective, and vanishing of a first jet at an unspecified point
imposes two conditions with one free point. The fiberwise first-jet
evaluation is surjective by increasing $\deg L$ if necessary.
Cohomology and base change extends this choice to a section on a
Zariski neighborhood of the parameter. Shrink so that its residual
zero divisor remains reduced, \'etale over the base, and disjoint
from the markings. Its product with the canonical section of
$\mathcal O(\Delta)$ is a section $s$ of $L^n$, whose divisor is
$\Delta+E$ with $E$ the residual divisor.

In the cyclic algebra given by $y^n=s$ use the finite algebra
\[
 \bigoplus_{j=0}^{n-1}L^{-j}
 \left(\left\lfloor\frac{j(\Delta+E)}{n}\right\rfloor\right).
\]
Here multiplication is the usual product of powers of $y$, with
$y^n=s$. This is the normalized cyclic covering algebra. Indeed,
at a branch of multiplicity $m=n/e$, after extracting a root of a
unit \'etale locally and using a relative parameter $x$, the
equation is $y^n=x^m$. It normalizes to $m$ disjoint charts
$x=v^e$ up to constant units. Equivalently, in the $j$th
eigenspace a coefficient of $y^j$ may have pole of order at
most $\lfloor j/e\rfloor$ in $x$, which gives exactly the displayed
algebra. Away from the branch divisor the cover is \'etale.
Thus the covering curve is smooth and projective over the
parameter, and carries its action of $\mu_n$. Its fibers are
connected since the polynomial is Eisenstein at any simple point of
$E$. At each marking the inertia order is exactly $e_i$, with the
specified root coordinate action. The additional points of $E$
have trivial inertia on the moduli objects and no pole. In a
faithful associated bundle a local frame at such a point averages
to an invariant full frame, by triviality of the inertia on its
fiber. Invariant regular coefficients and defining tensors
descend in this frame, and an invariant regular one-form has the
form $v^{n-1}b(v^n)\,dv$, regular downstairs. Consequently the
extra branch points introduce no new bundle, connection, or arrow
downstairs. On the fixed or moving root chart the quotient by the
cyclic deck action gives the required finite-cover presentation.

\emph{Ordinary and equivariant semistability.}
On the smooth projective Galois cover use the ordinary meromorphic
$\lambda$-connection stack in fixed degree.
In characteristic zero its semistable substack has a quasi-projective
good moduli space, and the canonical destabilizing reduction of an
unstable object is unique
\cite[Theorem~2.27, Proposition~4.46 and Theorem~4.47]{FHZ26}.
Here a \emph{good moduli morphism} $\pi$ has
$\pi_*\mathcal O=\mathcal O$ and exact pushforward on
quasi-coherent sheaves \cite[Definition~4.1]{Alp13}.
Every deck isomorphism preserves that reduction, and uniqueness gives
its cocycle. For an equivariant reduction its upstairs character degrees
are the cover degree times the root-curve degrees: a root line of
order $e$ has $e$th power the point line, while coarse line bundles
pull back with the ordinary degree multiplier. The signs of the
destabilizing tests therefore agree. It follows that equivariant
semistability is exactly ordinary semistability upstairs.

\emph{Finite forgetting.} Apply Lemma~\ref{lem:finite-deck-forgetting}
to the bounded ordinary semistable meromorphic $\lambda$-connection
stack on the cover and also to its regular Higgs counterpart.
Choose a faithful $G$-module $V_0$ and a deck-linearized relatively
ample line on the cover (pull back an ample line from the coarse
curve). For $n$ uniformly large on this finite-type stack,
$E(V_0)(n)$ is fiberwise generated, with vanishing $H^1$, and its
space of relative sections is locally free of fixed rank on each
degree component. Frames of that space form an $H$-torsor, for
$H=\mathrm{GL}(W)$ of this rank. An automorphism fixing the frame
is trivial by generation and faithfulness. Thus the frame space is
an algebraic space, and the deck action on it commutes with $H$.
Around each deck-fixed point of the ordinary good moduli space take
the intersection of the translates of an affine open neighborhood.
It is affine since the moduli space is separated. The frame space
over its saturated preimage is affine: the frame torsor is affine,
the good moduli morphism is cohomologically affine, and a
quasi-compact quasi-separated cohomologically affine algebraic
space over an affine scheme is affine \cite[Proposition~3.3]{Alp13}.
Call this frame scheme $R$.
In these coordinates the stack of coherent deck structures is the
disjoint union of $[R_\sigma/Z_H(\sigma(\Gamma))]$.
Indeed the action on the framed section space gives a representation
of $\Gamma$, which in families is locally a fixed representation
by its isotypic multiplicity bundles. There are finitely many
possibilities in fixed rank. On choosing a compatible frame the
remaining change of frame is its centralizer. Lemma~\ref{lem:finite-deck-forgetting}
therefore makes forgetting deck data finite on good moduli spaces.
The assertion is local on the target, including the open complement
of its deck-fixed locus, where the source is empty. Local inertia
types are locally constant since homomorphisms from finite cyclic
groups into $G$ are, up to conjugacy, finitely many torus torsion
types in a bounded order. Their conditions are open and closed in
families (in characteristic zero each such conjugacy class is a
connected component, up to components of the type space).
Retained degree components are likewise open and closed. These
cuts preserve finiteness.

We will also need finiteness of the inertia morphism after
rigidification. Frame the section spaces separately for the irreducible
summands of the faithful representation. Each deck centralizer $J$ is
then a product of general linear groups, and $Z(G)$ is a fixed subgroup
of its scalar center. The group $J'=J/Z(G)$ acts on the retained
affine chart with finite stabilizers by fine stability. To check
properness of this action, take two integral points over a discrete
valuation ring and a generic change of frame between them. After a
finite extension lift the change to $J$ and write it as
$a\nu(\pi)b$ by elementary divisors, with $a,b$ integral units and
$\pi$ a uniformizer. Embed the affine chart equivariantly as a closed
subscheme of a $J'$-module. In a weight basis for $\nu$, integrality
of the second point forces the negative-weight coordinates of the
first, after change by $b$, to vanish in the special fiber. Its
$\nu$-limit there belongs to the chart and is fixed by $\nu$.
Finite stabilizers force the projection of $\nu$ to $J'$ to be
trivial. The change of frame therefore extends over the valuation
ring, since integrality in an affine group can be tested after the
extension. Thus the action is proper, and its inertia is both proper
and affine, hence finite.

\emph{Quasi-projectivity and zero scaling limits.}
The finitely many fixed-action quotients above give the equivariant
good moduli spaces on the invariant affine charts. The universal
property glues them on overlaps. Finite forgetting to the ordinary
quasi-projective moduli space gives quasi-projectivity; central
rigidification leaves the coarse space unchanged. Apply ordinary
semistable reduction \cite[Proposition~4.48]{FHZ26} to the naive
scaling DVR family $(E,t\nabla_\lambda)$ over $t\lambda$.
The resulting ordinary semistable extension gives a coarse zero limit,
after a finite DVR extension if needed. Since forgetting is finite,
valuative properness extends the equivariant generic coarse orbit as
well, with its retained open-and-closed inertia and degree data.
Embed the coarse space in a projective closure.
The unique projective limit of the original orbit agrees after
ramification with the retained semistable point, so it lies in the
open coarse locus already over the original DVR. Fine stability gives
a stable representative.

\emph{The fixed locus.} At a scaling-fixed point $\lambda$ vanishes
and all ordered residues vanish. On the cover the logarithmic
residue belongs to the Lie algebra of the inertia centralizer.
For fine inertia this is the Cartan, so the whole residue vanishes
and the lifted Higgs field is regular. All positive-weight Hitchin
invariants vanish as well. The equivariant regular nilpotent cone
is the closed zero locus of the pullback of the ordinary Hitchin
coefficients in the equivariant regular moduli space, with the
retained inertia and degree cuts. By finite forgetting it is proper
over the base: the ordinary degree-zero principal Hitchin morphism
on the compact cover is projective by \cite[Set-up~2.0.1 and
Proposition~2.2.2]{dC21}. To check the semiharmonic hypothesis,
decompose a faithful representation into irreducibles. The connected
center acts by scalars on each, so principal semistability in
characteristic zero gives semistability of each associated Higgs bundle
\cite[Proposition~2.32]{FHZ26}. Each has degree zero because its
determinant is a character. Their direct sum is therefore semistable
of degree zero, and higher Chern classes vanish on the curve.
Conversely the semiharmonic condition gives semistability in the
adjoint representation and hence principal semistability
\cite[Corollary~2.33]{FHZ26}. These are fiberwise conditions on
families in the prescribed degree sector. The fine parahoric flags, including an auxiliary flag,
are encoded by the fixed fine inertia actions on the root curve.
The scaling-fixed locus is closed in this proper space. Additional
trivial-inertia points have no pole and impose the regularity condition.

\emph{Smoothness and surjectivity.} Smoothness over the residues is
\S\ref{subsec:A-restrictions}. For fixed markings the image of a nonempty master is
open, invariant under scaling, and contains the origin by the
zero-limit property. It is therefore the entire residue parameter space.

\emph{Relative version.} Choose the auxiliary weight in the real span of
the coroots, so that $d\chi(q_{\rm aux})=0$ for all characters; the
auxiliary flag then changes neither $\deg_\chi$ nor $\delta(\chi)$, and a
small weight crosses no wall. For $u\in\Cc$, use the relative cyclic cover constructed above,
with the auxiliary section disjoint from $D$. The cover is a smooth
projective family on a Zariski neighborhood of $u$. The section-frame
construction is relative and makes forgetting deck data finite there.
Relative projectivity in \cite[Proposition~2.2.2]{dC21} gives properness
of the ordinary zero Hitchin fiber, hence of the retained fixed locus.
These local statements give the required relative properness.
The points of $\Cc$ with nonempty master form an open set by smoothness.
This set is also the image of the fixed locus, since each scaling orbit
has a fixed zero limit over the same point of $\Cc$. Relative properness
makes that image closed. Thus a nonempty moving master meets every
fiber of the connected curve $\Cc$, and the same contraction argument
gives surjectivity on each fiber.

To obtain one relatively ample line on the entire moving-point family,
choose a fixed marking $p_0\in D$. Write $p_1,\ldots,p_a$ for the other
markings, including the moving one, and $N$ for the least common multiple
of the required orders. For $N>1$ choose $N$th roots $L_i$ of the
degree-zero lines $\mathcal O(p_i-p_0)$ on the family of curves.
These choices exist after a finite \'etale Galois cover of $\Cc$:
rigidify along $p_0$ (the restrictions there are constant) and pull back
multiplication by $N$ on the Jacobian, with its rigidified Poincar\'e line.
All $p_i$ stay disjoint from $p_0$ on this base. Adjoin the $N$th roots
of the tautological rational sections of $\mathcal O(p_i-p_0)$ in
the corresponding $L_i$ and normalize. The covering group is
$(\mu_N)^a$: valuations at $p_1,\ldots,p_a$ make the Kummer classes
independent modulo $N$. Away from these sections and $p_0$ the cover is
\'etale. At $p_i$ a single factor has equation $v^N=x$ up to a unit.
At $p_0$ all factors have pole of order one, so their ratios are roots
of units and the remaining factor has the same local chart after taking
its inverse. After extracting roots of units \'etale locally, the relative covering
algebra on each branch chart is a product of the algebras
$\mathcal O[v]/(v^N-x)$. These smooth relative charts and their
finite \'etale complements pull back under arbitrary base change. The cover is finite over the
projective relative curve and has inertia order $N$ at every marking.
Identify each cyclic inertia group with $\mu_N$ through its action on
the local parameter $v$. At $p_0$ the diagonal deck element
$(\zeta,\ldots,\zeta)$ multiplies the inverse-root parameter by
$\zeta^{-1}$, so its tangent-character label is $\zeta^{-1}$.
For desired order $e\mid N$ prescribe the original inertia type
along $z\mapsto z^{N/e}$ in these tangent-character coordinates,
with trivial action of its kernel. Averaging a faithful frame over this
kernel descends bundles and their tensors in $w=v^{N/e}$. An invariant
logarithmic form $b(w)\,dv/v$ becomes
$(e/N)b(w)\,dw/w$, giving exactly the smaller root chart for
$\lambda$-connections as well. These are the same retained inertia cuts
as before. For $N=1$ take the identity cover.

Thus on a single finite \'etale Galois base change the relative master
is finite over the ordinary quasi-projective space by
Lemma~\ref{lem:finite-deck-forgetting} applied to section frames, with its
proof applying to meromorphic $\lambda$-connections as well as regular
Higgs bundles. The original root-curve family and its coarse space
descend by their moduli characterization. Tensor the Galois translates
of a relatively ample line upstairs. The resulting relatively ample line
has the permutation descent datum. Its descended line is relatively
ample, as can be checked after this faithfully flat base change.
\end{proof}

The compactification below is relative to all residue coordinates; the
coordinate $z$ is a compactification coordinate, independent of
$\lambda$.

\begin{lemma}[Compactification pair]\label{lem:pg-pair}
Let $f:\mathcal Z\to B\times\A^r$ be a smooth semiprojective rigidified
master as above ($B$ a point or $\Cc$, scaling trivial on $B$), with a
linearized relatively ample line. There is a smooth proper orbifold pair
$(\overline{\mathcal Z},\mathcal E)\to B\times\A^r$, with both
$\overline{\mathcal Z}$ and the divisor $\mathcal E$ smooth over the base,
whose complement is $\mathcal Z$, including the fiber $\lambda=0$.
\end{lemma}

\begin{proof}
Let $d_i>0$ be the weights of the coordinates of $f$, $S=B\times\A^r_w$,
and
\[
 \mathcal Z'=\{(x,z,w):f_i(x)=z^{d_i}w_i\},\qquad t\cdot(x,z,w)=(tx,tz,w).
\]
This is the base change of $f$ along $(b,z,w)\mapsto(b,(z^{d_i}w_i)_i)$,
hence smooth over $B\times\A^1_z\times\A^r_w$; both $\mathcal Z'$ and its
divisor $\{z=0\}=f^{-1}(B\times0)\times\A^r_w$ are smooth over $S$.
On coarse spaces the zero limits are $(x_0,0,w)$, and the fixed locus is
the coarse fixed locus of the master times $\A^r_w$, proper over $S$.
The necessary linearization exists on every normal quasi-projective
coarse chart. Indeed, for a normal integral variety $X$ the map
$\mathrm{Cl}(X)\to\mathrm{Cl}(X\times\Gm)$ is an isomorphism:
localizing divisors at the generic point of $X$ removes the horizontal
part because $\C(X)[t,t^{-1}]$ is a PID, and its units are $c t^n$
with $c\in\C(X)^*$. Thus only vertical divisors and the original
vertical principal relations remain. The invertible sheaf
$a^*L\otimes p^*L^{-1}$ on $X\times\Gm$, for action $a$ and projection
$p$, comes from a divisorial sheaf on $X$; by restriction at the
identity this sheaf is trivial. Normalize the resulting isomorphism
at the identity. Its cocycle defect is a unit on $X\times\Gm^2$,
of the form $u(x)t^i s^j$, and equals one on both identity sections,
so is identically one. This gives a linearization of an ample line
on each component. The compactifying quotient is defined by the
action and the no-limit-at-infinity open, hence the local quotients
agree under base change to overlapping moving-point charts.

First take the tame coarse
space $Z'_c$ of $\mathcal Z'$. It is normal, and a power of the
linearized ample line descends. On the open $V_c\subset Z'_c$ of
points without a limit at infinity, the relative quotient theorem
gives a separated proper geometric quotient $Q=V_c/\Gm\to S$
\cite[\S11, Theorems~11.1--11.2]{Sim97},
\cite[Theorem~3.2.1]{dC21}. Here is the valuative argument for
properness. Locally over $S$ the linearization embeds $Z'_c$ as an
invariant open in its projective closure, with coordinate weights
$m_i$. Its fixed locus is closed in the closure by relative
properness. For a generic orbit representative with coordinates $x_i$
over a discrete valuation field, scaling by a scalar of valuation $s$
retains precisely the coordinates minimizing $v(x_i)+s m_i$ in the
special fiber. As $s$ decreases, the minimizing weights give a chain
of fixed points; each breakpoint gives one orbit joining adjacent
points. The first point lies in $Z'_c$ and the last outside, since
they specialize the two generic endpoints. An orbit ending inside
$Z'_c$ lies there by openness and invariance, so its zero limit is
also inside. Exactly one breakpoint therefore runs from inside to
outside. Its rational value is realized after finite ramification,
and any two rescalings giving a limit in $V_c$ differ by a unit.
This gives existence and uniqueness of the limiting orbit. Let
$\mathcal V$ be the inverse image of $V_c$ in $\mathcal Z'$.

The action on the variety $V_c$ has finite stabilizers. On an invariant
affine chart, invert the homogeneous coordinates $a_i$ nonzero at the
chosen point. If $d$ is the gcd of their degrees, a Laurent monomial
$u=\prod a_i^{e_i}$ with $\sum e_i\deg a_i=d$ gives the slice
$\Gm\times^{\mu_d}\{u=1\}$. Thus $[V_c/\Gm]\to Q$ is proper on finite
cyclic quotient charts. The tame coarse morphism $\mathcal V\to V_c$ is
proper \cite[\S3, Theorem~3.2 and Corollary~3.3(a)]{AOV08}; by base
change, so is $[\mathcal V/\Gm]\to[V_c/\Gm]$. Their composite is proper
and retains the original finite stabilizers, which need not be cyclic.
Smoothness of $[\mathcal V/\Gm]$ and its boundary over $S$ descends from
the smooth atlases $\mathcal V$ and $\mathcal V\cap\{z=0\}$.
Every point with $z\ne0$ lies in $\mathcal V$, and
$x\mapsto[(x,1,f(x))]$ identifies this open with $\mathcal Z$; its
complement is the smooth divisor
$[\mathcal V\cap\{z=0\}/\Gm]$, independent of $w$.
\end{proof}

\begin{proposition}[Transport stationary at the Higgs endpoint]
\label{prop:pg-transport}
For any finite family of piecewise smooth paths
$(\lambda,q,s)\mapsto(\lambda,q,a_\lambda+\lambda\rho(s),\ldots)$ in the
parameter space there are transports of the master and of its
compactification pair that are the identity at $\lambda=0$ and have
displacement $O(\lambda)$ on compact charts. They realize the cochain
continuation maps with their composition homotopies, and their graphs are
proper over both endpoints.
\end{proposition}

\begin{proof}
Choose a horizontal splitting of the tangent sequence of the pair, tangent
to $\mathcal E$, average it on finite charts and patch. The lift of the
path velocity is linear in $\lambda\rho'(s)$, hence equals $\lambda
V_{\lambda,s}$ with $V$ bounded on compact charts; properness makes the
flow complete, and integrating gives the $O(\lambda)$ bound and the identity
at $\lambda=0$. The flows are homeomorphisms of the complement, so both
projections of their graphs are proper. A flow trivialization identifies
the evaluation maps of the family with continuation; convex interpolation
of splittings over finite simplices gives the composition faces,
compatibly with finite symmetries.
\end{proof}

\subsection{Fiber restriction and base change}\label{subsec:A-basechange}

\begin{lemma}[Fiber restriction and base change]\label{lem:basechange}
Let $f:\mathcal H\to\Lambda$ be a nonempty fine master, or its base change
to a linear coincidence subspace, with coarse space $H$, and
$r=\dim\Lambda$. Then the inclusions of fibers induce isomorphisms
\[
 H^*(H,\Q)\xrightarrow{\sim}H^*(H_a,\Q),\qquad
 R\Gamma(H,\Q)_\Lambda\xrightarrow{\sim}Rf_*\Q_H;
\]
in particular the stalk-to-fiber maps of $Rf_*\Q_H$ are isomorphisms.
The complexes $Rf_!\Q_H$ likewise have locally constant cohomology
and the canonical stalks $R\Gamma_c(H_a,\Q)$. All these assertions
are compatible with automorphisms and finite projectors.
\end{lemma}

\begin{proof}
By smoothness the image of $f$ is open; it contains the origin (zero
limits, Proposition~\ref{prop:pg-master}) and is invariant under the
positive-weight scaling, hence equals $\Lambda$; the same holds after
restriction to a nonempty linear coincidence subspace.

By Lemma~\ref{lem:pg-pair} there is a smooth proper orbifold pair
$(\overline{\mathcal H},\mathcal E)\to\Lambda$ with complement
$\mathcal H$. Horizontal lifts tangent to $\mathcal E$, averaged over
finite charts and patched by a partition of unity, have complete flows
over compact paths; they trivialize the pair, hence $H$, over the
contractible $\Lambda$, and the resulting isomorphism is induced by the
fiber inclusion. To identify the derived direct image, let
$V=R\Gamma(H,\Q)$ with evaluation $V_\Lambda\to Rf_*\Q_H$. In a finite
quotient chart, normal purity gives $j_a^!\Q_H=\Q_{H_a}[-2r](-r)$, and
exceptional base change identifies the costalk of evaluation at $a$ with
$V\to R\Gamma(H_a,\Q)$ shifted by $[-2r]$, an isomorphism by the above.
All costalks of the cone vanish, so evaluation is an isomorphism; this
step is needed because $f$ is not proper. The trivializations of the pair
give the compact-support local systems and their fiber base change.
There is no unshifted comparison between the compactly supported
cohomology of the total space and that of a fiber: globally the product
over $\Lambda\cong\A^r$ gives
$R\Gamma_c(H,\Q)\cong R\Gamma_c(H_a,\Q)[-2r](-r)$.
\end{proof}

\subsection{Hitchin fiber dimensions and perverse transport}
\label{subsec:A-perverse}

\begin{theorem}[Geometry of the fine master]\label{thm:fine-master}
The rigidified fine master is smooth over its residue space, and its weak
Higgs part has a proper joint Hitchin map. For a nonempty fixed-residue
fiber of dimension $2N$, every nonempty Hitchin fiber is pure of
dimension $N$ and the Hitchin defect is $N$. The coarse constant sheaf is
ULA over the residue space, and the fiber restrictions identify the
normalized perverse filtrations, also on finite projector summands and
linear coincidence restrictions.
\end{theorem}

\begin{proof}
\emph{Idea.} Properness follows from the $\Gm$-action contracting to a
proper nilpotent cone. Poisson-commutativity gives a lower bound on
fiber dimensions. The symplectic weight count bounds the core above,
and proper scaling specialization gives the upper bound on all fibers.
Perverse transport follows from local acyclicity over the residue base,
combined with the fiber base change of Lemma~\ref{lem:basechange}.

\emph{Properness.} Apply the properness criterion
\cite[Definition~B.1 and Proposition~B.5]{FHZ26} at $\lambda=0$: the
source has zero limits, the coefficient-and-residue target has positive
weights, and the fiber over the origin is the fixed-type nilpotent cone,
which is proper (Proposition~\ref{prop:pg-master}). The joint map $(h,b)$
factors through $\mathcal A_{\rm mer}\times_{\prod_p\ft/\WG}\prod_p\ft$,
which is finite over $\mathcal A_{\rm mer}$; so $h$ is proper.
The coarse Higgs master is a closed slice of the quasi-projective
coarse master of that proposition. Hence the joint Hitchin map is
projective: an immersion of its source into a projective space,
together with the proper map, is a proper immersion into the product
with its target and therefore a closed immersion. The same argument
works over a moving-base chart and after a linear coincidence
restriction. The Springer maps used later are likewise projective,
being closed flag incidences over their targets.

\emph{Dimensions.} The fixed-residue stable Higgs stack is a symplectic
reduction of the cotangent stack of framed parahoric bundles: the duality
$\mathcal P^\vee=\mathcal S(D)$ identifies cotangent directions with weak
Higgs fields, and the frame-changing torus has the ordered residues as
moment map. The Hitchin functions Poisson-commute: in a Dolbeault chart the
Hamiltonian vector field of $\eta\,p(\Phi)$ is $(\eta\,dp(\Phi),0)$
because $[dp(\Phi),\Phi]=0$, and the Serre-dual coefficient lies in
$\mathcal P(-D)\subset\mathcal S$, so the field is allowed and descends.
Hence $h$ has image of dimension at most $N$, and every fiber has
dimension at least $N$. At zero residues the zero fiber is the
\emph{semiprojective core}: the points whose scaling orbit also has a
limit at infinity \cite{HRV15}. At an orbifold fixed point lift the scaling action after an isogeny
$t\mapsto t^m$. The symplectic form upstairs has weight $m$, so its
integer tangent weights pair as $u,m-u$. The tangent dimension with
weights $u\le0$, including the fixed directions, is at most $N$.
The attracting stratum for the limit at infinity has this dimension, so
all strata of the core have dimension at most $N$. Upper semicontinuity
for the proper joint map, along scaling orbits to the origin, therefore
bounds all fibers above by $N$. Since the image has dimension at most
$N$, the fiber-dimension inequality bounds every irreducible component
of every fiber below by $N$. Thus the fibers are pure of dimension $N$,
the image has dimension $N$, and
$\dim(H_b\times_{\mathcal A}H_b)-\dim H_b=N$.

\emph{Perverse transport.} On a smooth chart over the residue base $B$ the
constant sheaf is ULA; finite direct image and invariants preserve this.
Put $r=\dim B$, $\mathbf H=(h,b):Y\to\mathcal A\times B$ and
$K=R\mathbf H_*\Q_Y[N+r]$. Proper direct image preserves ULA, and
successive coordinate cuts with the perverse vanishing-cycle criterion
\cite[Proposition~3.1.4 and Corollary~3.1.6]{dCMa20} give
\begin{equation}\label{eq:perverse-restrict}
 i_a^*(\ptau_{\le k}K)[-r]\cong\ptau_{\le k}(i_a^*K[-r])
 =\ptau_{\le k}Rh_{a*}\Q_{Y_a}[N].
\end{equation}
The remaining push $\mathcal A\times B\to B$ is not proper. By the
decomposition theorem $\ptau_{\le k}K$ is a direct summand of $K$, so
$R\pr_{B*}\ptau_{\le k}K$ is a summand of $Rb_*\Q_Y[N+r]$, which has fiber
base change and lisse cohomology by Lemma~\ref{lem:basechange}; naturality
restricts the base-change isomorphism to the summand. Taking images of the
canonical truncation morphism and using~\eqref{eq:perverse-restrict}
identifies the fibers with $P_kH^j(Y_a)$. The splitting is used only for
this verification; the truncation morphism and its image are canonical.
\end{proof}

\subsection{Riemann--Hilbert in families and flat weights}
\label{subsec:A-betti}

Let $S_B$ be the space of ordered Cartan monodromies satisfying the
central product equations, and let $S_B^{\rm gen}$ be its intrinsically
generic locus. The weak Betti family over this locus has a compatible full
flag at each puncture and variable unipotent part; we call the extra flag
an \emph{enrichment}. Boalch's classification identifies enriched
monodromy with logarithmic parahoric connections, including automorphisms
\cite[Theorems~A and~D]{Boa11}. This comparison is local on $S_B^{\rm gen}$;
it is not one Riemann--Hilbert identification over the whole affine
residue space for a fixed $q$.

The numerical weight of Boalch's enriched local system must be
specified. On a logarithm chart write the Cartan monodromy as
$\exp(2\pi i r)$, so the operator residue at $\lambda=1$ is $e=-r$ and
$R=q-r$. Boalch's apartment parameters are $\theta=-q$ and
$\phi=\operatorname{Re}(r-q)=-\operatorname{Re}R$. When $\phi$ lies in the
chosen open Borel chamber, its parabolic is that Borel; the numerical
weight is recovered from $(q,r)$ and adds no choice to the flag; this
weight need not be zero. The zero-filtration harmonic inverse of
\S\ref{app:harmonic} is not used in this comparison.
We need only the following local dependence on parameters.

\begin{lemma}[Riemann--Hilbert in families]\label{lem:enriched-rh}
For a rational fine parahoric type, the enriched local Riemann--Hilbert
correspondence is holomorphic in the connection and its ordered residues,
on a chart where the real part of the affine-normalized residue of the
horizontal equation stays in one open Weyl chamber. Integral root
resonances of the raw residue are allowed.
\end{lemma}

\begin{proof}
Write the horizontal equation as $dy=Ay$. For growth weight $q$, the gauge
$z^{-q}$ gives $A=(b(s)+R(z,s))dz/z$ with $R=O(|z|^\epsilon)$, uniformly on
a small chart (nonconstant affine-root terms have positive, discrete
degree, and the reductive quotient is toral). Put $t=-\log z$. On the
opposite big cell of $G/B$ the flag equation is
$u'=D(s)u+F(t,u,s)$ with $\|F(t,0,s)\|+\|\partial_uF\|\le
Ce^{-\epsilon\operatorname{Re}t}$, where the eigenvalues of $D(s)$ are
positive roots of $b(s)$ with real parts $\ge\delta>0$. The decaying
solution is the fixed point of
$u(t,s)=-\int_0^\infty e^{-vD(s)}F(t+v,u(t+v,s),s)\,dv$, a contraction in
the norm $\sup e^{\eta\operatorname{Re}t}\|u\|$ ($0<\eta<\epsilon$) after
moving the starting circle toward the puncture; iteration gives
holomorphic dependence, and parallel transport extends the flag to a base
point. A turn around the puncture conjugates the coefficients by
$\exp(-2\pi iq)$, which preserves the limiting Borel, so by uniqueness the
flag is monodromy invariant and covariant under parahoric frame changes.
In the normal form $z^\tau z^{\sigma+N}$ of
\cite[Theorem~D and its proof in \S3]{Boa11}, the fine condition puts the
nilpotent part in the unipotent radical of this Borel, whose flag has
exactly this asymptotic property, and the automorphism comparison there
identifies stabilizers. No inverse of a resonant root operator is used.
\end{proof}

\begin{theorem}[The weak Betti family]\label{thm:weak-betti}
On the intrinsically generic locus of the ordered full-flag weak Betti
family, ordinary and compactly supported rational cohomology are local
systems with fiber base change, and the weight filtrations are flat
sub-local systems. This includes collisions at several punctures, all
finite character summands, and unions permuted by component groups.
\end{theorem}

\begin{proof}
\emph{Idea.} Locally on the generic parameter space, Riemann--Hilbert
identifies the Betti family with a slice of a fine Hodge master, whose
fiber restrictions are isomorphisms by Lemma~\ref{lem:basechange}. Local
constancy is a topological statement, so it transfers along this
homeomorphism; flatness of the weight filtration then comes from the
algebraic structure of the Betti family alone.

Fix $s\in S_B^{\rm gen}$ and a sufficiently small analytic logarithm
chart $U$ around $s$, with ordered logarithms $r_p(u)$. Choose fine
rational weights $q_p$ outside the numerical walls and with full torus
inertia, satisfying $\sum_pd\chi(q_p-r_p)=0$ for central characters.
This is possible because $\sum_pd\chi(r_p)$ is an integer constant on
$U$, and a rational central correction realizes these finitely many
values. Retain all compatible finite topological sectors. Coroot shifts
of the logarithms change neither monodromy nor central degree nor the
$G$-topological sector; a large multiple of $2\rho^\vee$ puts
$\phi_p=\operatorname{Re}(r_p-q_p)$ in the chamber of the ordered Borel.
Shrink $U$ so that all these strict inequalities persist. Define
\[
 s_q:U\longrightarrow\Lambda_q,\qquad
 u\longmapsto\bigl(1,(q_p-r_p(u))_p\bigr),
 \qquad H_{q,U}=H_q^{\rm an}\times_{\Lambda_q}s_q(U).
\]
Here $H_q\to\Lambda_q$ is the fine Hodge master for this choice of $q$;
only its pullback to this chart is compared with Betti objects.
Lemma~\ref{lem:enriched-rh} and holomorphic dependence of parallel
transport give a holomorphic map $H_{q,U}\to X_B^{\rm an}|_U$ over $U$.
Boalch's classification, with its automorphism comparison, glues the
local parahoric extensions to the flat bundle on $\Cc$ and makes this
map bijective on each fiber. Genericity makes every representation irreducible, so no
object is lost to the stability condition. Both total spaces are rational
cohomology manifolds of the same dimension (residue smoothness on one
side; smooth enriched peripheral spaces and vanishing infinitesimal
centralizer modulo the center on the other, compare
\cite[Theorem~B]{Boa11}), so Lemma~\ref{lem:domain-invariance} makes the
bijection a homeomorphism.

The proper smooth compactification pair from
Lemma~\ref{lem:basechange} pulls back to $s_q(U)$. The local
trivializations of the pulled-back pair, or the same evaluation and costalk argument, identify
$Rf_{q,U*}\Q$ with a locally constant complex and give the canonical
stalk-to-fiber isomorphisms. Under the homeomorphism over $U$ these are
exactly the corresponding maps for $X_B|_U\to U$. Proper pair transport
also gives the compact-support statement. Such charts cover
$S_B^{\rm gen}$, and the conclusions glue: being a local system and being
an isomorphism for the canonical base-change map are local properties.
Changing $q$ or a logarithm lift on an overlap requires no identification
of their numerical filtrations. Both constructions prove the same
intrinsic assertions for the Betti morphism. This proves local constancy
over the entire generic Betti parameter locus, not a global
trivialization; its local systems may have monodromy. No assertion at a
nongeneric parameter follows from this argument.
Finally the Betti family is algebraic: applying the
mixed Hodge module operations to its constant complex, the cohomology
modules have lisse realization over the smooth parameter space.
Their weight submodules are subconnections, giving the flat weight
filtration described in \S\ref{subsec:families}. The base-change morphism in mixed
Hodge modules is an isomorphism because its rational realization is, so
the fiber mixed Hodge structure is Deligne's. (This does not assert that
Riemann--Hilbert preserves mixed Hodge structures.)
\end{proof}

\subsection{Refinements and transport of intersection complexes}
\label{subsec:A-transport}

For a harmonic Higgs bundle $(\bar\partial,\Phi,h)$ and real $\tau$, the
\emph{preferred evaluation} is the $\tau$-connection
$\bar\partial_\tau=\bar\partial+\tau\Phi^\dagger$,
$\nabla_\tau=\tau\partial_h+\Phi$.

A residue refinement is made in the reductive quotient $L$ of the coarse
parahoric: the preimage of a Borel of $L$ has flag quotient $L/B_L$, and
once the ordered semisimple labels are fixed the residual flags are full
flags of their centralizers. If angular weights tie but imaginary labels
differ, the spectral ordering supplies the partial flag canonically.

The refinement must be the full flag incidence over \emph{every} stable
coarse object. After Lemma~\ref{lem:rational-chart}, choose a rational
polarization for the refined problem: the relevant walls $d+\ell(q)=0$
are proper hyperplanes on the face preserving the joint blocks (a wall
containing the face would give $d+\ell(a)=0$ and $\ell(\beta)=0$, hence a
forbidden character product), and a refinement smaller than the minimal
numerical gap preserves coarse stability for every residual flag. With
coarse and fine weights $q_{\rm c},q_{\rm f}$, the incidence equations use
$R_{\rm f}=R_{\rm c}+\lambda(q_{\rm f}-q_{\rm c})$.
For a root $\alpha$ of the coarse reductive quotient, write
$m=\alpha(q_{\rm c})\in\Z$ and $e=R_{\rm c}-\lambda q_{\rm c}$.
The root belongs to the graded centralizer precisely when
\[
 \alpha(e)+\lambda m=0,
 \qquad\text{equivalently}\qquad
 \alpha(R_{\rm f})=\lambda\alpha(q_{\rm f}-q_{\rm c}).
\]
These are linear equations in the fine master. Restrict to the open
part of their intersection where no further root enters the centralizer.
On this part the local Weyl group is fixed. By
\eqref{eq:framed-complex} the framed residue evaluation is smooth, so
forgetting the residual flags is a smooth pullback of a product of
Springer maps; its decomposition and restriction maps are the Springer
ones. For an orbit $\mathcal O\ni N$ in $\mathfrak l$ with
$b=\dim\Bfl_N$, the multiplicity of the constant local system is
$V_{\mathcal O}=H^{2b}(\Bfl_N,\Q)^{A_L(N)}$
(Proposition~\ref{prop:springer}\ref{it:semismall}), and changing a local
transport by $Z_L(N)$ acts trivially on it. Forgetting flags preserves
invariant polynomials, so all maps lie over the common Hitchin base.

\begin{theorem}[Harmonic transport of boundary conditions]
\label{thm:harmonic-transport}
For a generic full spectral principal problem with real growth weights
$a$ and imaginary semisimple labels $\beta$, preferred evaluation
identifies the real Hodge family over an open interval containing
$[0,1]$ with the product of the Higgs space and the interval; at $\tau=1$
it is the tame nonabelian Hodge map $\Psi$. This identification preserves
the nilpotent orbit of each connected centralizer (and the union of
translates for a disconnected one). On every constant-coefficient closed
Jordan multiplicity, the ratio of the master restrictions is the map on
intersection cohomology induced by $\Psi$. These identifications commute
with finite projectors, Springer extraction, and the perverse transport
of Theorem~\ref{thm:fine-master}.
\end{theorem}

\begin{proof}
\emph{Idea.} For each real $\tau$, the harmonic metric of a Higgs
bundle defines a $\tau$-connection (preferred evaluation). This map is
continuous in $\tau$ (Proposition~\ref{prop:preferred-continuity}) and,
for $\tau\ne0$, invertible by Theorem~\ref{thm:curve-inverse}.
Invariance of domain upgrades the resulting continuous bijection to a
homeomorphism of families. The semisimple labels fix the centralizer. On its Lie algebra
the pointwise comparison of nilpotent weight filtrations fixes the orbit
dimension. Continuity then keeps the residue in the same orbit.

\emph{Fiberwise bijectivity.} Let $q$ be the rational polarization from
Lemma~\ref{lem:rational-chart}, $\mathscr X$ the full spectral coarse
family and $\mathscr Y\to\mathscr X$ the fine Springer incidence. By
\eqref{eq:kms} the raw residue of the preferred operator is $e_\tau$, and
its affine residue in the $q$-master is
\begin{equation}\label{eq:actual-path}
 R_\tau=e_\tau+\tau q=(1+\tau^2)\beta+\tau(q-a).
\end{equation}
At $\tau=0$ the stable loci for $a$ and $q$ agree. For $\tau\ne0$ the flat
exponents are $-a+(\tau+\tau^{-1})\beta$; for a proper parabolic
character, a nonzero imaginary residue sum stays nonzero and a zero sum
keeps the nonintegral angular obstruction, so the representation is
irreducible and stable for $q$.

For the inverse, first work with the \emph{coarse} parahoric object in
$\mathscr X$ and endow it with the real growth label $a$ of the harmonic
metric.
Lemma~\ref{lem:rational-chart} identifies its principal parahoric
lattices and component actions with those in the $q$-chart. This is a
change of the numerical growth labels on the same principal object;
in associated representations the filtrations are reindexed by $a$,
not left at $q$. Irreducibility makes stability independent of this
change when $\tau\ne0$, and the central degree equations are preserved.
Indeed the $q$-labeled object has filtered-local-system weight
\[
 q+\operatorname{Re}(e_\tau/\tau)=q-a,
\]
which need not vanish. It is the object relabeled by $a$ for which
$a+\operatorname{Re}(e_\tau/\tau)=0$ for every weight of every
representation of $G$. Thus its associated filtered local systems have
trivial filtration weights. The underlying local systems are semisimple. Indeed, if the Zariski
closure $H$ of the monodromy had a nontrivial unipotent radical $U$,
the Borel--Tits canonical parabolic would satisfy
$H\subset N_G(U)\subset P(U)\ne G$ \cite[Proposition~3.1(ii)]{BT71}, contradicting
irreducibility. Thus $H^0$ is reductive; in characteristic zero the
finite component group preserves complete reducibility, so every
associated representation is semisimple.
Each irreducible summand has parabolic degree zero. By
Theorem~\ref{thm:curve-inverse} these filtered $\tau$-flat objects carry
adapted $\tau$-harmonic metrics, giving an exact tensor functor
$\Xi_{\tau,0}$ from them to tame harmonic bundles. Composing with the
associated-bundle functor $\operatorname{Rep}(G)\to\{\text{filtered
$\tau$-flat bundles}\}$ transports the defining tensors and their
relations, so the frames carrying the fixed tensor tuple to its transported
sections reconstruct the principal torsor on the punctured curve.

To obtain a harmonic reduction of this torsor, fix a tensor frame at a
basepoint. Write $H$ for the stabilizer of the tuple and $h$ for the
faithful vector metric. The defining tensors are harmonic morphisms:
they are parallel for the Chern connection $D_h$ and are annihilated
by the Higgs field $\Phi$ and its adjoint. Thus Chern holonomy belongs
to $K_J=H\cap U(h)$. Use the metric downstairs on the punctured
curve and pull it back on a covering chart, so the deck arrows are
isometries as well. The group $J=H\cap H^{\dagger_h}$ is closed under
Hermitian adjoint, so $\operatorname{Lie}J$ is the complexification
of $\operatorname{Lie}K_J$ (take the skew-Hermitian and Hermitian
parts). After parallel transport to the basepoint, each coefficient
of $\Phi$ and its adjoint belongs to $\operatorname{Lie}J$.
Choose a maximal compact $K_H$ of $H$ containing $K_J$. By conjugacy
of maximal compacts, $K_H$ preserves an $H$-translate $h'$ of the
chosen model metric. Extend this metric and the $K_H$-reduction by
$D_h$-transport. Both adjoints agree on $\operatorname{Lie}J$, since
both act by minus the identity on $\operatorname{Lie}K_J$. Thus $D_h$
is still the Chern connection, the adjoint of $\Phi$ is unchanged,
and the principal reduction gives the same preferred operators.
The positive endomorphism $h^{-1}h'$ is $D_h$-parallel. Its constant
eigenvalues give equivalent growth norms on this punctured curve.

Lemma~\ref{lem:real-growth-tensor} proves tensor compatibility of
real-weight prolongations in every representation of $G$, and therefore
extends the reduction at the real growth weight $a$ of the harmonic
metric. Only then do we express the principal object in the rational
chart $q$; rational tensor reconstruction is
\cite[Proposition~6.1]{BMW12}. Finite inertia actions are morphisms in
the same tensor category. For a recovered Higgs label $(a_0,\beta_0)$ and its
$\tau$-connection label $(p,e)$ the general KMS formulas give
$e+\tau p=(1+\tau^2)\beta_0$. Here $p=a$ and
$p+\operatorname{Re}(e/\tau)=0$, so $\operatorname{Re}\beta_0=0$.
Hence $a_0=p=a$ and $\beta_0=(e+\tau a)/(1+\tau^2)=\beta$.
Hence preferred evaluation is bijective on
fibers, and it is continuous by Proposition~\ref{prop:preferred-continuity}.

A fine refinement may change the coarse facet and is not subjected to
this same-object relabeling. If its weight is $q_{\rm f}$, its path is
$R_{{\rm f},\tau}=e_\tau+\tau q_{\rm f}$, and forgetting flags maps it to
$R_\tau=e_\tau+\tau q$ on $\mathscr X$. The harmonic inverse just proved
is applied after this forgetful map. The proper Springer incidence and
its multiplicity complexes are transported by the fixed-centralizer
and intersection-complex argument below; no zero-filtration assertion
is made for the $q_{\rm f}$-labeled fine object.

\emph{Fixed centralizers.} For a root $\alpha$ of the parahoric reductive
quotient the rational chart preserves the affine-root equation, so
$\alpha(q)=\alpha(a)$ and $\alpha(R_\tau)=(1+\tau^2)\alpha(\beta)$; since
$1+\tau^2\ne0$, the connected graded centralizer is independent of
$\tau$. On a frame chart the full spectral condition is a smooth pullback
of one fixed nilpotent cone, which is rationally smooth
(Proposition~\ref{prop:springer}\ref{it:rational-smooth}); finite quotients
and smooth factors preserve this. So source and target, over an open
interval containing $[0,1]$, are locally contractible rational cohomology
manifolds of the same dimension without boundary, and
Lemma~\ref{lem:domain-invariance} makes the continuous bijection a
homeomorphism. At $\tau=1$ it is $\Psi$ by construction.

\emph{Orbit types.} Apply the pointwise comparison of
Lemma~\ref{lem:pointwise-kms} to the adjoint harmonic bundle.
The Lie algebra $\mathfrak l$ of the common centralizer is the joint
graded summand of growth and semisimple residue both zero. Besides the
Cartan subalgebra, this summand contains precisely the roots $\alpha$
with $m=\alpha(a)\in\Z$ and $\alpha(\beta)=0$. Indeed, multiplying a
root frame by $z^m$ shifts its growth from $m$ to zero and adds $\tau m$
to its residue, leaving $(1+\tau^2)\alpha(\beta)$.

On this summand the nilpotent endomorphism is $\ad_{\mathfrak l}N_\tau$.
Write $d_j$ for the dimension of grade $j$ of its nilpotent weight
filtration centered at zero. A Jordan block of length $s$ contributes
one dimension in degrees $1-s,3-s,\ldots,s-1$, so
$\dim\ker\ad_{\mathfrak l}N_\tau=d_0+d_1$. The lemma gives
constancy of both numbers. Thus the orbit dimension is fixed. In a
local frame the residues vary continuously, and each connected-group
nilpotent orbit is open and closed in the union of the finitely many
orbits of its dimension, since its boundary has smaller dimension.
A path in such a frame therefore stays in one orbit. This applies
also to a marked connected orbit when continuous frames have been chosen.

For the full centralizer orbit there is also a pointwise argument,
which covers the two very even orbits in type $D$. Put
$s=\exp(2\pi i(a-2\beta))$ and $L=Z_G(s)$. In every representation
$V$, an eigenvalue $c$ of $s$ determines the growth modulo $\Z$ by its
phase and the imaginary Higgs label by its absolute value. In a frame
of this joint graded summand the nilpotent acts by $d\rho(N_\tau)$:
shifting a torus-weight frame to the chosen growth representative adds
only a scalar to its diagonal residue. Choose Jacobson--Morozov
cocharacters $h_0,h_\tau$ lowering their respective nilpotents by two
in the sense that $\operatorname{Ad}(h_i(x))N_i=x^{-2}N_i$,
and conjugate them into $T$ inside $L^\circ$. Their weight multiplicities
on each $V_c$ are exactly the $d_j$ for that block. Hence
$\chi_V(s h_\tau(x))=\chi_V(s h_0(x))$ for all $V$ and $x\in\C^*$.
By highest-weight triangularity these characters span
$\C[X^*(T)]^{\WG}$, whose functions separate the finite Weyl orbits on $T$.
Finitely many closed loci, one for each Weyl element, therefore cover
$\C^*$; irreducibility gives one Weyl element valid for all $x$. At
$x=1$ this element fixes $s$, so its representative in $L$ identifies
the two cocharacters. For their common value $h$, the $Z_{L^\circ}(h)$-orbit
of either nilpotent is open in the weight space $\mathfrak l_{-2}(h)$:
its tangent space is $[\mathfrak l_0(h),N]=\mathfrak l_{-2}(h)$ by
$\mathfrak{sl}_2$ representation theory. The two open orbits in this
vector space coincide. This gives the full $L$-orbit, hence the
permitted union of connected orbits when $L$ is disconnected.

\emph{Intersection complexes.} For such a condition the homeomorphism
identifies the family with a product over the interval and carries the
constant top-stratum line to itself; the intersection complex, being the
intermediate extension of that line, is the product intersection complex.
The smooth-pullback Springer construction identifies the multiplicity with
this complex with shift $2b$ and twist $(-b)$, and its restriction to
fibers follows from the smooth-pullback incidence square and proper base
change. So the ratio of the fine-master restrictions on each such summand
is the $\IH$ map of $\Psi$.

\emph{Defects of smaller strata.} The dense stratum of a closure has
dimension $2(N-b)$ and Poisson-commuting Hitchin functions, so its Hitchin
fibers have dimension $\ge N-b$. Over a stratum with Springer fiber
dimension $b'$, the inverse image of a coarse Hitchin fiber lies in a
fine Hitchin fiber of dimension $N$; hence that part of the coarse fiber
has dimension at most $N-b'$. Here $b'=b$ on the dense stratum and
$b'\ge b$ on every smaller stratum. Thus every fiber of the closure has
dimension at most $N-b$, and its image has dimension $N-b$ by the dense
stratum calculation. Applying the fiber-dimension inequality on each
irreducible component over its image gives the reverse bound for every
irreducible fiber component. The fibers are therefore pure of dimension
$N-b$, and the defect is $N-b$. Removing the Springer degree $2b$, the
weight $2b$ and the perverse index $b$ gives the normalized filtrations of
the smaller problem, and the perverse restriction argument of
Theorem~\ref{thm:fine-master} applies to every proper projector summand.
\end{proof}


\section{The raising operator}\label{sec:raising}

The criterion of \S\ref{sec:linear} needs one operator $E$ whose leading
parts satisfy relative hard Lefschetz for $P_\bullet$ and curious hard
Lefschetz for $W_{2\bullet}$, with the same center. We take $E$ to be cup
product with the class of the Goldman symplectic form. The Betti
argument computes the pulled-back class on torus factors of pieces of
$\widetilde\MB$, descends the cohomological statement by trace, and then
extends the result when eigenvalues collide. On the Dolbeault side the
same class differs from a relatively ample class by boundary classes,
and the decisive point is that the boundary classes raise the perverse
filtration by at most one.

Subsections~\ref{subsec:torus-coordinates}--\ref{subsec:collisions}
treat the Betti side and prove curious hard Lefschetz
(Theorem~\ref{thm:regular-chl} and \S\ref{subsec:collisions});
\S\S\ref{subsec:ample}--\ref{subsec:boundary-classes} treat the
Dolbeault side; Theorem~\ref{thm:common-leading} combines the two.

\subsection{Set-up}\label{subsec:raising-setup}

We use the faithful representation and invariant trace form fixed in
\S\ref{subsec:groups}. Goldman defined, for closed surfaces, a symplectic
form on the smooth locus of the character variety from $\kappa$ and the
cup product in group cohomology \cite[\S\S1.4, 1.7--1.8]{Gol84}. For
punctured surfaces with fixed boundary classes, we use the reduced form
of the quasi-Hamiltonian fusion product of handle doubles and conjugacy
classes. The compact-group construction and its surface pairing are
\cite[Theorems~5.1, 6.1 and~9.3]{AMM98}; the algebraic complex reductive
realization is \cite[Theorem~2.2 and Corollaries~3.6--3.7]{Saf17}. Concretely,
at an irreducible representation $\rho$ the tangent space is the
parabolic cohomology $H^1_{\rm par}(\Cc,\fg_{\Ad\rho})$ (classes whose
restrictions to the boundary circles are trivial, see
\S\ref{subsec:torus-coordinates}), and the form is
\[
 \omega(x,y)=\bigl\langle\kappa(\tilde x\smile y),[\Sigma,\partial\Sigma]
 \bigr\rangle,
\]
where $\Sigma$ is the compact surface obtained from $C$ by removing small
open disks around the punctures and $\tilde x\in H^1(\Sigma,\partial
\Sigma;\fg_{\Ad\rho})$ is any lift of $x$. We call the de~Rham class of $\omega$, with one
fixed period normalization, the \emph{Goldman class}, and write $E$ for
cup product with it.

By Lemma~\ref{lem:generic}, under intrinsic genericity all representations
are irreducible with stabilizers finite modulo $Z(G)$. This holds on the
whole fiber, not only on a stable open subset; the rank computations
below use this.

The proof uses the three family properties \ref{F:betti}--\ref{F:harm}
established in Section~\ref{sec:comparison-families}:
property~\ref{F:betti} transports curious hard Lefschetz,
property~\ref{F:dol} transports the perverse bound on boundary classes,
and property~\ref{F:harm} identifies the two cup actions.

\subsection{Torus coordinates and parabolic cohomology}
\label{subsec:torus-coordinates}

We first assume that $G_\der$ is simply connected and that every boundary
class is regular semisimple. Then every $S_p$ has centralizer exactly
$T$ (\S\ref{subsec:groups}), which is the ``nice'' condition of
de~Amorin \cite{dA25}. The general case is reduced to this one by the
central cover of \S\ref{subsec:central-lifts} and by the collision
argument of \S\ref{subsec:collisions}. Put $r=\dim T$ and $z=\dim Z(G)$.

\emph{Unipotent stabilization.} Absorb the twist as in
\S\ref{subsec:intro-setting}, order the punctures so that the last one,
$p_s$, has regular boundary value $C_s\in T$, and put $S=T/Z(G)$. Let
\[
 R=\Bigl\{(A_i,B_i,M_1,\ldots,M_{s-1})\in G^{2g}\times\textstyle
 \prod_{j<s}\overline{\mathcal C}_j:\ \prod_i[A_i,B_i]\,M_1\cdots M_{s-1}
 \,C_s=1\Bigr\},
\]
the relation space with the last monodromy fixed to $C_s$; the torus $T$
acts on $R$ by simultaneous conjugation. Since $Z_G(C_s)=T$, every
$G$-orbit of tuples meets $R$ in one $T$-orbit, so $[R/S]$ is the
rigidified quotient stack of~\eqref{eq:intro-MB} and $\MB=R\sslash S$.

\begin{definition}\label{def:stabilization}
The \emph{unipotent stabilization} of $\MB$ is
$\widetilde\MB=(R\times U_B)\sslash S$, where $S$ acts on $U_B$ by
conjugation, together with the map $f:\widetilde\MB\to\MB$ induced by
the projection $R\times U_B\to R$.
\end{definition}

Equivalently, the map $(r,u)\mapsto u\,r\,u^{-1}$ identifies
$R\times U_B$ with the space of tuples $(A_i,B_i,M_1,\ldots,M_{s-1})$ whose surface product
lies in $U_BC_s^{-1}$: the map $n\mapsto C_s^{-1}nC_sn^{-1}$ is an
automorphism of $U_B$ because $C_s$ is regular, so every element of
$U_BC_s^{-1}$ is uniquely of the form $uC_s^{-1}u^{-1}$. This is the form
used in \S\ref{subsec:stabilized-bruhat}; compare
\cite[Lemma~11.1]{dA25}.

The exponential map $\Lie U_B\to U_B$ is an $S$-equivariant isomorphism
of varieties, so the stack map
\[
 f_{\rm st}:[R\times U_B/S]\longrightarrow[R/S]
\]
is a representable vector bundle of rank $b$, and its oriented trace is
an isomorphism $Rf_{{\rm st},!}\Q=\Q[-2b](-b)$. The coarse-space maps of
both stacks have rational direct image $\Q$, by averaging over the
finite stabilizers (Lemma~\ref{lem:finite-quotient}). Composing these
proper direct images gives, for the coarse map $f$,
\begin{equation}\label{eq:stabilization-trace}
 Rf_!\Q_{\widetilde\MB}=\Q_{\MB}[-2b](-b),
\end{equation}
an identity of mixed Hodge complexes compatible with cup product by
classes pulled back from $\MB$. The coarse map $f$ itself need not be
an affine-space bundle: over a point with finite stabilizer $F$ its fiber
is $U_B/F$.

\emph{Bruhat walks.} For a word $(i_L,\ldots,i_1)$ in the simple
reflections put $b_i(e)=n_ix_i(e)$. The \emph{Bruhat walk} of
$(e_1,\ldots,e_L)$ is the sequence $(p_0,\ldots,p_L)$ in $\WG$ with
$p_0=1$ and $b_{i_m}(e_m)\cdots b_{i_1}(e_1)\in Bn_{p_m}B$. By
\eqref{eq:bruhat-rule}, if $s_{i_m}p_{m-1}>p_{m-1}$ then
$p_m=s_{i_m}p_{m-1}$ (an \emph{up} step); otherwise
$p_m\in\{s_{i_m}p_{m-1},p_{m-1}\}$ (a \emph{down} or \emph{stay} step);
compare \cite[\S\S5.4, 7.3]{Mel25} and \cite[Lemmas~8.1--8.2]{dA25}. For
the closed walks used here $p_0=p_L=1$, so there are $a$ up steps, $a$
down steps and $k$ stay steps, with $L=2a+k$. There is a coordinate
$e_m$ at each up and stay step, nonzero at stays, and the set of points
with a given walk is $(\C^*)^k\times\C^a$. Write $S_\bullet$ and
$D_\bullet$ for the stay and down positions and
$w_{>m}=s_{i_L}\cdots s_{i_{m+1}}$. In coordinates rescaled by the
preceding Cartan factor, \eqref{eq:rank-one} gives
\begin{equation}\label{eq:word-data}
 \mu(e)=\kappa_\bullet\prod_{m\in S_\bullet}
   w_{>m}\bigl(\alpha_{i_m}^\vee(e_m^{-1})\bigr),\qquad
 t\cdot e_m=\alpha_{i_m}(p_{m-1}(t))\,e_m,
\end{equation}
where $\mu$ is the Cartan part of the product and
$\kappa_\bullet=\prod_{m\in D_\bullet}w_{>m}(\alpha_{i_m}^\vee(-1))$
records the Tits-lift corrections. In the raw coordinates the torus
acts through the word prefix rather than the Bruhat prefix; the
triangular monomial change to the coordinates of~\eqref{eq:word-data} is
given below. Combining the word with the handle coordinates gives a
\emph{torus equation}
\begin{equation}\label{eq:torus-equation}
 \Phi:T^{2g}\times(\C^*)^k\longrightarrow T_\der,
\end{equation}
independent of the remaining $A_{\rm ext}+a$ affine coordinates, where
$A_{\rm ext}=b(2g+s-1)-L/2$. \S\ref{subsec:stabilized-bruhat}
proves the required $T$-equivariant handle, puncture and walk charts
directly from Bruhat factorizations; the construction follows
\cite[\S\S3--8]{dA25}, which extends \cite[\S\S5.4, 7.3--7.4]{Mel25}.
The remaining affine coordinates are $T$-homogeneous, so their zero
section preserves every Cartan level. Only these constructions are taken
from \cite{dA25}: the rank of $\Phi$ is proved under intrinsic
genericity in Lemma~\ref{lem:stable-zero-section}, and every finite
kernel is computed from the integral character map.

\emph{Raw and rescaled stay coordinates.} Number stays in their walk
order. Put $\beta_j=p_{m_j-1}^{-1}\alpha_{i_{m_j}}$,
$R_j=s_{\beta_j}$ and $E_j=R_j\cdots R_1$. If $a_j$ is the raw stay
coordinate and $h_j$ the accumulated Cartan factor, then, up to fixed
Cartan multipliers,
\[
 h_j=\beta_j^\vee(a_j^{-1})h_{j-1},\qquad
 t\cdot a_j=\beta_j(E_{j-1}(t))a_j.
\]
Set $a_j=e_j\beta_j(h_{j-1})$. This is a triangular monomial
automorphism of the stay coordinates, and in the new coordinates
\begin{equation}\label{eq:stay-rescaling}
 t\cdot e_j=\beta_j(t)e_j,\qquad
 h_j=\beta_j^\vee(e_j^{-1})R_j(h_{j-1}),\qquad
 v_j=R_jv_{j-1}-q_j,
\end{equation}
where $v_j=d\log h_j$ and $q_j=\beta_j^\vee d\log e_j$.
The word is considered after collecting the handle Cartan factors into
a separate factor $D$; hence $h_0=1$. With
$W=R_k\cdots R_1$ and $R_{<j}=R_1\cdots R_{j-1}$,
\[
 d\log\mu=v_k=W\sum_jR_{<j}q_j.
\]
The Weyl identity
$w_{>m_j}\alpha_{i_{m_j}}^\vee=-W R_{<j}\beta_j^\vee$
recovers \eqref{eq:word-data}. Thus the moment map $\mu$ and the gauge
action in~\eqref{eq:word-data} are written in the same coordinates.

\emph{Parabolic cohomology.} Let $\Sigma_{g,k}$ be a compact oriented
genus-$g$ surface with $k$ boundary circles and $\mathbb V$ a local
system on it. The \emph{parabolic cohomology}
$H^1_{\rm par}(\Sigma_{g,k},\mathbb V)$ is the image of
$H^1(\Sigma_{g,k},\partial\Sigma_{g,k};\mathbb V)\to H^1(\Sigma_{g,k};
\mathbb V)$: classes whose restriction to each boundary circle is
trivial. For the adjoint local system, $H^1_{\rm par}$ is the tangent
space to deformations with fixed boundary conjugacy classes.

\begin{lemma}\label{lem:parabolic-lattice}
Let $V=X_*(T)\otimes\Q$. The rational complex of the linearized gauge
action $D_0$ and the linearized torus equation $D_1$ on a cell is
\begin{equation}\label{eq:parabolic-complex}
 V\xrightarrow{D_0}V^{2g}\oplus\bigoplus_{j=1}^k\im(1-R_j)
 \xrightarrow{D_1}V,
\end{equation}
the cochain complex computing $H^1_{\rm par}(\Sigma_{g,k},\mathbb V_W)$
for the local system $\mathbb V_W$ with handle monodromies
$A_i,B_i\in\WG$ and boundary monodromies the reflections
$R_j=s_{\beta_j}$, $\beta_j=p_{m_j-1}^{-1}\alpha_{i_{m_j}}$ at the stay
steps $m_j$; these satisfy $\prod_i[A_i,B_i]R_1\cdots R_k=1$.
\end{lemma}

\begin{proof}
Put $w_m=s_{i_m}\cdots s_{i_1}$ and $E_m=p_m^{-1}w_m$. This is unchanged
at up and down steps and becomes $R_jE_{m-1}$ at a stay; since
$p_0=p_L=1$ we get $R_k\cdots R_1=\prod_i[A_i,B_i]$, the stated relation.
By~\eqref{eq:word-data} the derivative at a stay is
$W R_{<j}\beta_j^\vee\,d\log e_j$, with $W=\prod_i[A_i,B_i]$ and
$R_{<j}=R_1\cdots R_{j-1}$, and the corresponding gauge derivative is
$\beta_j^\vee\beta_j(v)=(1-R_j)v$.

Write $C_i=[A_i,B_i]$ and $C_{<i}=\prod_{\ell<i}C_\ell$. The handle torus
factor is
\[
 A_i(t_i)\,A_iB_i(u_it_i^{-1})\,A_iB_iA_i^{-1}(u_i^{-1});
\]
with $a_i=A_i\,d\log t_i$, $b_i=B_i\,d\log u_i$ its derivative is
$(1-A_iB_iA_i^{-1})a_i+(A_i-C_i)b_i$.
The isolated twisted handle has gauge derivative
$((1-A_i)C_i^{-1}v,(1-B_i)C_i^{-1}v)$; convolution with its word,
whose Weyl product is $C_i$, replaces $v$ by $C_iv$.
Moving each handle torus factor left past the preceding word factors
conjugates its derivative by $C_{<i}$
(Lemma~\ref{lem:H-charts}). Hence
\begin{gather*}
 D_0(v)=((1-A_i)v,(1-B_i)v,(1-R_j)v),\\
 D_1(a,b,q)=\sum_iC_{<i}\bigl((1-A_iB_iA_i^{-1})a_i+(A_i-C_i)b_i\bigr)
   +W\sum_jR_{<j}q_j .
\end{gather*}
These are the coboundary and the Fox-calculus relation maps of the
standard presentation of $\pi_1(\Sigma_{g,k})$, with boundary cocycles
restricted to $\im(1-R_j)$. In particular $D_1D_0=0$ and the middle
cohomology is $H^1_{\rm par}$.
\end{proof}

The invariant form $\kappa$ gives a perfect alternating pairing on
$H^1_{\rm par}$: Poincar\'e--Lefschetz duality pairs
$H^1(\Sigma,\partial\Sigma;\mathbb V)$ perfectly with
$H^1(\Sigma;\mathbb V)$, the transpose of the natural map $J$ between
them is $-J$, so $\ker J$ is the annihilator of $\im J$. With
$h^0=\dim H^0(\Sigma,\mathbb V_W)$ (the monodromy invariants), the
Euler characteristic and the boundary sequence give
\begin{equation}\label{eq:parabolic-dimension}
 \dim H^1_{\rm par}(\Sigma_{g,k},\mathbb V_W)
   =(2g-2)r+\sum_j\rank(1-R_j)+2h^0,
\end{equation}
and both maps in~\eqref{eq:parabolic-complex} have rank $r-h^0$.

To pass from equivariant charts to their quotients, we keep track of
the finite kernel of the gauge action. The following elementary quotient calculation
also explains why an isogeny need not be split.

\begin{lemma}[The quotient of an equivariant stabilized cell]
\label{lem:H-quotient}
Let $S$ and $K$ be complex algebraic tori, let $V$ be a linear
$S$-representation, and let $\alpha:S\to K$ have finite kernel $F$.
Let $S$ act on $V\times K$ through its representation on $V$ and
translation by $\alpha$ on $K$. With $Q=K/\alpha(S)$ there are
noncanonical isomorphisms
\[
 [(V\times K)/S]\simeq[V/F]\times Q,
 \qquad (V\times K)/S\simeq(V/F)\times Q
\]
of the quotient stack and of its coarse space, respectively.
The finite group $F$ is diagonalizable and acts diagonally after a
choice of basis of $V$. The same description applies separately to
the components of a torsor under a possibly disconnected torus kernel.
\end{lemma}
\begin{proof}
Writing $L=\alpha(S)=S/F$, the exact sequence
$1\to L\to K\to Q\to1$ splits by the character lattices of tori.
Hence $K\simeq L\times Q$. Induction from $F$ to $S$ identifies
$[(V\times(S/F))/S]$ with $[V/F]$, proving the assertions.
We do not split the isogeny $S\to S/F$. A finite subgroup of a complex
torus is diagonalizable, so its linear representations decompose into
characters. At $(0,k)$ the stabilizer is $F$; at $(v,k)$ it is $F_v$.
\end{proof}

\begin{lemma}[Stratification of the stabilization]
\label{lem:stable-zero-section}
The unipotent stabilization $\widetilde\MB$ has a finite filtration by
closed subvarieties whose successive differences are finite disjoint
unions of pieces
\[
 (\C^*)^t\times(\C^e/F),\qquad
 t=k+(2g-2)r+2z,\quad e=A_{\rm ext}+a,\quad
 t+2e=\dim\MB+2b,
\]
one family of pieces for each closed Bruhat walk, with $a$ up steps and
$k$ stay steps. Here $F$ is a finite diagonalizable group acting
linearly on $\C^e$: the kernel of the integral gauge map on a component
of the corresponding fiber of $\Phi$. (This group is invisible if
integral lattices are replaced by their real spans.)
\end{lemma}

\begin{proof}
\S\ref{subsec:stabilized-bruhat} gives $T$-equivariant product
charts $\C^e\times\Phi^{-1}(c)$ (Corollary~\ref{cor:H-charts}). The zero
section of the affine factor consists of points of $R\times U_B$, i.e.\
of representations. By Lemma~\ref{lem:generic} these representations
have stabilizer finite modulo $Z(G)$. Consequently the torus-coordinate
stabilizer has dimension $z$, so $h^0=z$ in
\eqref{eq:parabolic-complex}. Both differentials have rank $r-z$.
Thus $\Phi$ is surjective onto $T_{\der}$ on a nonempty cell, with
possibly disconnected kernel, and its fibers have dimension
$2gr+k-(r-z)$.

On a connected component $K$ of this fiber, translate to the identity
component and write $\lambda:T/Z(G)\to K$ for the gauge translation.
The kernel $F$ of $\lambda$ is finite. With $Q=K/\lambda(T/Z(G))$,
Lemma~\ref{lem:H-quotient} gives
\[
 [\C^e\times K/(T/Z(G))]\simeq[\C^e/F]\times Q.
\]
Taking coarse spaces gives the asserted product.
Subtracting $r-z$ gives the formula for $t$; substituting $L=2a+k$
gives $t+2e=\dim\MB+2b$. The closed filtration follows from the outer
Bruhat order and the prefix maps, as proved in
Lemma~\ref{lem:H-filtration}. No partition of $\MB$ with the same
product shape is needed.
\end{proof}

\subsection{Identifying the Goldman form}\label{subsec:goldman-form}

Only a cohomology class on each stabilized piece is needed for
Lefschetz. In particular, Proposition~\ref{prop:goldman-pairing} does not
require every affine-containing term of a differential form to vanish. The fusion
formula is \cite[Theorem~6.1]{AMM98}; \S\ref{subsec:stabilized-bruhat}
checks its compatibility with each Bruhat coordinate change.
The following identity explains why such a local replacement can be
checked inside the complete surface product.

\begin{lemma}[Bar potential under local replacements]\label{lem:H-block-bars}
Put $B(f,g)=\kappa(f^{-1}df\wedge dg\,g^{-1})$ and
\[
 \Omega(f_1,\ldots,f_n)=\sum_{j>1}B(f_1\cdots f_{j-1},f_j).
\]
Replacing a block of factors by another block with the same product
changes the full $\Omega$ by the difference of the two internal bar sums.
This difference is independent of the position of the block.
In particular, simultaneously expanding a factor and its inverse, even
when they are separated by other factors, makes canceling corrections.
\end{lemma}
\begin{proof}
The cocycle identity
$B(f,g)+B(fg,h)=B(g,h)+B(f,gh)$ gives
\[
 \Omega(\ldots,a,b,\ldots)
 =\Omega(\ldots,ab,\ldots)+B(a,b).
\]
Iteration proves the block assertion. The inverse identity gives
$\Omega(f_n^{-1},\ldots,f_1^{-1})=-\Omega(f_1,\ldots,f_n)$,
which proves the final assertion without moving either block past any
intervening factor.
\end{proof}

\begin{proposition}[Goldman class on a stabilized piece]
\label{prop:goldman-pairing}
On every piece $Q\times(\C^e/F)$ of
Lemma~\ref{lem:stable-zero-section}, the restriction of the pulled-back
Goldman class is the pullback of a constant logarithmic class on $Q$. On the rational
cocharacter space of $Q$, this class is the parabolic cup pairing on
$\ker D_1/\operatorname{im}D_0$, with the fixed Goldman normalization.
In particular the torus form is nondegenerate, in every root system.
\end{proposition}

\begin{proof}
The affine zero section induces an isomorphism on ordinary rational
cohomology, including after a finite linear quotient. We calculate there,
using the bar potential $\Omega$ of Lemma~\ref{lem:H-block-bars}. The
fusion form is $\Omega/2$. Lemmas~\ref{lem:H-handle-bars} and~\ref{lem:H-puncture-bars}
therefore identify the original surface form with the form in the
stabilized coordinates, retaining the three handle Cartan
corrections and the fixed-conjugacy-class form.

For the remaining word, Lemma~\ref{lem:H-word-bars} shows that up and
down steps add no variable Cartan bar, while a stay adds
$\kappa(d\log\beta_j^\vee(a_j^{-1})\wedge v_{j-1})$.
Here $a_j$ is the raw coordinate. At the closed endpoint its Bruhat
primitive is a $U_B$--$B$ bar and vanishes. Hence
\[
 \Omega_{\rm word}
 =\sum_j\kappa(-\beta_j^\vee d\log a_j\wedge v_{j-1})
 =\sum_j\kappa(v_{j-1}\wedge q_j),
\]
where \eqref{eq:stay-rescaling} and self-adjointness of $R_j$ give the
second equality. The separate terminal bar is
$\kappa(d\log D\wedge d\log\mu)=0$ on $D\mu=\text{constant}$.

Choose normalizer lifts and consider handles
$n_{A_i}=\dot A_it_i$, $n_{B_i}=\dot B_iu_i$ and reflection factors
$n_j=\beta_j^\vee(e_j)\dot R_j$. The right logarithmic derivative of a
product is $u_{ab}=u_a+\rho(a)u_b$, and that of an inverse is
$u_{a^{-1}}=-\rho(a)^{-1}u_a$. Thus the surface relation has derivative
$D_1$, simultaneous conjugation has derivative $D_0$, and its bar sum
is exactly the preceding handle-and-word sum. Constant lift factors
change the level only.

Half the polygon bar sum is the antisymmetrized relative cup product.
For a reflection boundary $R$, the relative-lift correction is
$\tfrac12\kappa(\xi,(R^{-1}-R)\eta)=0$ because $R^2=1$.
The fixed final Cartan value has zero boundary variation and can be
capped. The reduced form is therefore the pairing on
$\ker D_1/\operatorname{im}D_0$. Poincar\'e--Lefschetz duality proves
nondegeneracy. No change of the normalization of $E$ is made in this
calculation.
\end{proof}

\subsection{Curious hard Lefschetz at regular boundary values}
\label{subsec:regular-chl}

On $(\C^*)^t$ with a nondegenerate constant logarithmic two-form
$\omega$, $H^j$ is pure of weight $2j$ and type $(j,j)$, and $\omega$ has
weight $4$; Example~\ref{ex:torus-sl2} shows that $\omega$ is Lefschetz
with center $t/2$ for the filtration $W_{2\bullet}$. The common middle
weight $t+2e=\dim\MB+2b$ of all pieces allows us to glue.

\begin{theorem}\label{thm:regular-chl}
Let $\MB$ be a nonempty intrinsically generic character variety with
regular semisimple boundary classes, and $d=\dim\MB$. Then $H^*(\MB)$ is
Hodge--Tate, the Goldman class has weight degree at most $4$, and its
leading part $\ell$ satisfies \emph{curious hard Lefschetz}: for $a\ge0$
\[
 \ell^a:\Gr^W_{d-2a}H^j(\MB,\C)\xrightarrow{\ \sim\ }
 \Gr^W_{d+2a}H^{j+2a}(\MB,\C),
\]
and likewise for $H_c^*$ and on the summands cut out by commuting finite
projectors.
\end{theorem}

\begin{proof}
All pieces in this argument are pieces of $\widetilde\MB$.
On a piece $(\C^*)^t\times(\C^e/F)$ the torus factor satisfies
exterior-algebra Lefschetz with compact middle weight $t$, and
$H^*_c(\C^e/F)=\Q(-e)[-2e]$ by Lemma~\ref{lem:finite-quotient}; so each
piece is Hodge--Tate with middle weight $t+2e=d+2b$. Cup product with the
pulled-back Goldman class commutes with the localization sequences of the closed
filtration of Lemma~\ref{lem:stable-zero-section}, and strictness of
weights and the five lemma glue the Lefschetz isomorphisms, as in
\cite[Proposition~2.19]{Mel25}. The cup-compatible trace \eqref{eq:stabilization-trace} identifies
$H_c^{j+2b}(\widetilde\MB)$ with $H_c^j(\MB)(-b)$, changing the middle
weight from $d+2b$ to $d$. Poincar\'e duality on
the rationally smooth $\MB$ gives the statement for $H^*$.

For general $G$ use the cover $G^\sharp=Z(G)^0\times G_\der^{\rm sc}\to
G$, keep every compatible lift of the boundary values and of $\zeta$,
and apply the above on each central level. The finite group of
\S\ref{subsec:central-lifts} preserves the Goldman form (though not
each cell), and exact invariants descend the Lefschetz isomorphisms.
\end{proof}

\subsection{Collisions of eigenvalues and Springer multiplicities}
\label{subsec:collisions}

When eigenvalues coincide we keep a full flag at each puncture and allow
an arbitrary unipotent part preserving it. This produces a smooth family
over the space of ordered boundary values, containing the regular
values. Two facts are needed: the weight-$4$ part of the Goldman class is
flat in this family, and its cup product commutes with the Springer
projectors.

\emph{The weak family as a quasi-Hamiltonian reduction.} The space
$G\times_{U_B}B$, where $u\in U_B$ acts by $(g,x)\mapsto(gu^{-1},uxu^{-1})$,
remembers a decorated flag in $G/U_B$ and the Cartan part of the
boundary element. The space $G\times_{U_B}B$ is the (non-affinized)
multiplicative implosion,
with a quasi-Hamiltonian $G\times T$-structure. The axioms and fusion
formulas of \cite[Definition~2.2]{AMM98} have the algebraic complex
reductive realization of \cite[Theorem~2.2, Corollaries~3.6--3.7 and
\S\S3.5--3.6]{Saf17}. Fuse the handle doubles with one copy per
puncture and reduce by $G$ at the surface relation. A stabilizer fixes
the decorated $G/U_B$ frame, so lies in a conjugate of $U_B$ as well as
in the centralizer of the underlying irreducible representation. The
latter is finite modulo the center, and its intersection with $U_B$
is trivial in characteristic zero. The reduction $P$ is therefore
smooth; its remaining torus-valued quasi-Hamiltonian structure has a
closed nondegenerate two-form. The residual torus $K=T^s/Z(G)$ acts locally
freely, and $P/K$ is the weak full-flag family, with unrestricted
unipotent boundary parts. Local logarithms of the Cartan values give an
ordinary moment map $P\to\mathfrak k^*$.

\emph{Variation of the class.} Choose a basis of characters of $K$,
connection forms $\theta_i$ on the $K$-torsor and logarithmic moment
coordinates $\mu_i$. The form $\omega_P+\sum_id\mu_i\wedge\theta_i$ is
basic and its differential is $-\sum_id\mu_i\wedge F_i$, $F_i$ the
curvatures. Hence, in the Gauss--Manin connection of the cohomology
local system (property~\ref{F:betti}),
\begin{equation}\label{eq:DH}
 \nabla[\Omega_t]=-\sum_id\mu_i\otimes b_{i,t},\qquad
 b_{i,t}=c_1(\mathcal L_i),
\end{equation}
where $\mathcal L_i$ is the line bundle associated with the $i$th
character; this is the Duistermaat--Heckman variation formula
\cite[Theorem~1.1 and \S2, (2.3)--(2.10)]{DH82} in the present setting.
The cited global theorem assumes a compact torus and a proper momentum
map. Here the same local curvature computation uses the cohomological
transport of property~\ref{F:betti}. Since $c_1(\mathcal L_i)$ has
weight two, the weight-$4$ part of the Goldman class is \emph{horizontal}.
This weight-$4$ part is rational after period normalization. At a regular anchor,
$Rf_{{\rm st},*}\Q=\Q$ for the stack vector bundle, and the proper
coarse maps have rational direct image $\Q$. Their composition gives
$\Q_{\MB}\xrightarrow{\sim}Rf_*\Q_{\widetilde\MB}$ in mixed Hodge
complexes; in particular $f^*$ is an isomorphism on ordinary cohomology
as mixed Hodge structures. Each stabilized piece
$Q\times(\C^e/F)$ has the ordinary cohomology of its affine zero section
$Q$. The logarithmic formula there is rational after the fixed period
normalization: $\kappa$ is rational on $X_*(T)\otimes\Q$, and the
level, gauge quotient and reflections are defined over their rational
cocharacter spaces.

Choose a dense open $O\subset\widetilde\MB$ whose components lie in
maximal-dimensional stabilized pieces, by removing the finitely many
stratum boundaries and closures of lower-dimensional strata. The
preceding rational classes therefore restrict rationally to $O$.
The space $\widetilde\MB$ is rationally smooth, being the coarse space
of the smooth finite-inertia vector-bundle stack. Put
$m=\dim\widetilde\MB$ and $Z=\widetilde\MB\setminus O$.
Duality identifies
$H_Z^2(\widetilde\MB,\Q)=H^{\rm BM}_{2m-2}(Z,\Q)(-m)$;
this group is spanned by the fundamental classes of the
codimension-one components of $Z$ and is pure of weight two.
The localization sequence and strictness thus make restriction to $O$
injective on $\Gr^W_4H^2$. This proves rationality downstairs through
$f^*$, and it propagates along the rational weight local systems.

\emph{Springer modifications.} Let $\pi:Y\to X$ be a connected
\emph{Springer modification}: the map from a weak full-flag family to the
corresponding full spectral family, which over each boundary point is a
product of Springer fibers on a frame chart with ordered labels
and connected centralizer (Proposition~\ref{prop:springer}). Its coarse
fiber is the quotient of this product by the finite effective stabilizer.
For example, in $\PGL_2$ the class of $\operatorname{diag}(-1,1)$
has centralizer $N_G(T)$ and two invariant lines with the same Cartan
value. The Springer fiber of its connected centralizer $T$ is a point.
We make the calculation on the cover first and descend the full union
of sectors in \S\ref{subsec:central-quotient}. By
Proposition~\ref{prop:springer}\ref{it:kostant} and
Proposition~\ref{prop:springer-low-degree}, $R^0\pi_*\Q=\Q$,
$R^1\pi_*\Q=0$, and fiber $H^2$ is pure Tate. The Leray sequence then
makes $H^2(X)\to H^2(Y)$ injective with cokernel of weight at most two
(evaluate a global section of $R^2\pi_*\Q$ at one point of each stratum
to embed the cokernel into sums of fiber $H^2$). Hence the weight-$4$
part of the class on $Y$ comes from $X$, and its cup product commutes
with every endomorphism of $R\pi_*\Q$, in particular with the Springer
projectors.

The weak family meets the regular locus in every generic parameter
component, and fiber base change preserves nonemptiness. Horizontality
transports Theorem~\ref{thm:regular-chl} to all collision values. With
the Springer decomposition
\[
 R\pi_*\Q_Y[d]=\bigoplus_\lambda V_\lambda\otimes
 \IC_{Z_\lambda}(\mathcal A_\lambda)(-b_\lambda),\qquad
 2b_\lambda=d-\dim Z_\lambda,
\]
the summand appears in $H^{j+2b}(Y)$ as $\IH^j(Z,\mathcal A)(-b)$, and
removing the weight shift $2b$ moves the center from $d$ to $\dim Z$.
This gives curious hard Lefschetz on every occurring multiplicity. For a
disconnected centralizer we argue on the connected-centralizer cover with
all component-group translates included and then take invariants.

\subsection{A relatively ample reference class}\label{subsec:ample}

The perverse side needs a relatively ample class for the Hitchin map. At
a reference point with finite-order monodromy, both the Goldman class and
an ample class come from the same universal Chern class on a compact
covering curve. We first construct this reference and then measure the
change of the Goldman class along the residues.

The torus factor is separated by the finite central quotient
(Proposition~\ref{prop:central-Dolbeault-quotient}) and treated by
Example~\ref{ex:intro-rank-one}; in this subsection the group is
semisimple. By the absorbed-twist convention in the Introduction,
the reference is an untwisted local system. Consider a
rational fine strongly parabolic reference (zero ordered Cartan residue)
whose peripheral monodromy is regular semisimple of finite order; its
existence in each required degree is part of the family construction
(Section~\ref{sec:comparison-families}). Choose a finite cover $\widehat C\to C$ of
degree $M$, branched over $D$, on which the whole peripheral monodromy
becomes trivial. Both the flat bundle and its harmonic Higgs bundle then
extend holomorphically over $\widehat C$.

\begin{proposition}\label{prop:pure-ample-reference}
There is a rational class $\alpha_B\in H^2(\MB)$, pure Tate of weight
$4$, equal to a fixed nonzero multiple of the Goldman class, whose
pullback $\alpha_D=\Psi^*\alpha_B$ equals $M^{-1}c_1(\mathcal L_\rho)$
for a line bundle $\mathcal L_\rho$ some positive power of which descends
to a Hitchin-relatively ample line bundle on the coarse space.
\end{proposition}

\begin{proof}
On the stack put $\mathcal L_\rho=(\det Rq_*\mathcal U_\rho)^{-1}$, where
$\mathcal U_\rho$ is the universal bundle associated with $\rho$ on the
compact cover and $q$ the projection. Decompose $V$ by central characters
$\nu$; the associated bundles are flat, hence of degree zero, and a
central element acts on the determinant line by
$\prod_\nu\nu(z)^{\dim V_\nu(1-\widehat g)}=\det\rho(z)^{1-\widehat g}=1$.
A further power kills the finite effective stabilizers. By
Grothendieck--Riemann--Roch,
$c_1(\mathcal L_\rho)=-\int_{\widehat C}\ch_2(\mathcal U_\rho)$.

By the Chern--Weil computation of Biswas--Schumacher
\cite[Proposition~6.1, equation~(27)]{BS06} (arXiv v1 numbering), this
class is represented by the positive $L^2$ (Petersson--Weil) pairing on
harmonic infinitesimal deformations, up to an exact correction. The
principal harmonic deformation complex injects isometrically into the
associated one, since $d\rho$ has an equivariant trace-orthogonal
splitting; only polystability of the associated bundle is used. Harmonic
families are constructed in deck-equivariant gauge charts on
$\widehat C$; the kernel of the linearized equation consists of parallel
adjoint sections commuting with the Higgs field, i.e.\ infinitesimal
automorphisms of the stable principal object, which are central. Hence
every nonzero principal tangent vector has positive norm, top
self-intersections on positive-dimensional subvarieties of Hitchin fibers
are positive, and the Nakai--Moishezon criterion gives relative
ampleness, also through finite orbifold charts.

On the Betti side put $\alpha_B=-M^{-1}\int_{\widehat C}\ch_2(\mathcal
U_{\rho,B})$; pullback to the cover is algebraic on character stacks.
Shende's theorem \cite[Theorem, p.~2]{She17} shows that slanting a class of
$H^{2k}(BG)$ against a surface cycle yields a class of Hodge type
$(k,k)$; with $k=2$ this makes $\alpha_B$ pure Tate of weight $4$. Its
tangent pairing is the normalized Goldman cup pairing. The
harmonic metric identifies the two universal smooth bundles, so
$\alpha_D=M^{-1}c_1(\mathcal L_\rho)$ without using~\eqref{eq:PW}.
Pullback along the central gerbe is injective in degree two and finite
rigidification preserves rational cohomology.
\end{proof}

\subsection{The perverse degree of boundary classes}
\label{subsec:boundary-classes}

By~\eqref{eq:DH}, moving the residues changes the Goldman class by
first Chern classes of boundary character lines. We show that their cup
products raise the perverse filtration by at most one. First we prove
vanishing on every Hitchin fiber for a dense set of residues, then use
local acyclicity to reach all residues.

Fix a rational fine Dolbeault residue master. A character tuple
$\chi=(\chi_p)$ of $T^s/Z(G)$ defines a line bundle by applying $\chi_p$
to the Cartan torsor at $p$; let $b_\chi$ be its rational first Chern
class. For a maximal proper parabolic, the rational characters trivial on
$Z(G)$ form a ray of dominant characters; choose $\xi$ on each of the
finitely many rays up to conjugacy. Consider the dense set of residues
$(\beta_p)$ that are regular and satisfy
\begin{equation}\label{eq:additive-generic}
 \sum_p(w_p\xi)(\beta_p)\ne0
\end{equation}
for all such $\xi$ and all Weyl elements $w_p$; its complement is a finite
union of hyperplanes.

\begin{lemma}\label{lem:boundary-fibre-zero}
For residues satisfying~\eqref{eq:additive-generic}, every $b_\chi$
restricts to zero on every Hitchin fiber, including singular fibers.
\end{lemma}

\begin{proof}
Fix a Hitchin value $a$. The \emph{cameral cover} over the open set where
the Higgs field is regular semisimple parametrizes ordered Cartan
representatives; it is nonempty because the residues are regular. Let
$\WG_a\subset\WG$ be the monodromy group of a component. A rational
$\WG_a$-invariant cocharacter outside the center would define an
invariant parabolic reduction over the function field, which extends
over $C$ by properness of $G/P$; applying the dominant character to the
Higgs field gives a logarithmic differential with residue sum zero,
contradicting~\eqref{eq:additive-generic}. Hence
$X_*(T)_\Q^{\WG_a}=X_*(Z(G)^0)_\Q$.

On the fiber stack, the cameral locus carries the diagonalizing
$T$-torsor of the universal Higgs bundle. Evaluating a character line at
two points of a connected cameral component gives the same Chern class,
so $c_{w\chi}=c_\chi$ for $w\in\WG_a$, and the evaluation factors through
coinvariants; by the previous paragraph, averaging kills the derived
part. At a puncture use a root disk $u^r=z$ with $\varphi=\psi(u)du/u$,
$\psi(0)=r\beta_p$ regular; the disk can be chosen uniformly on the
fiber, since the discriminant depends only on $a$. The ordered Cartan
branch and its torsor extend over the disk; evaluation at $0$ is the
boundary line, at a nonzero point the cameral evaluation, and their
Chern classes agree. Thus derived boundary characters vanish on the
fiber. For the central part write $\chi_p=\chi_p^\der+\lambda_p$ with
$\lambda_p\in X^*(G)_\Q$; diagonal central cancellation gives
$\sum_p\lambda_p=0$, and all central lines evaluated at one common point
cancel. Pullback along the central gerbe is injective in rational degree
two, and averaging descends to the coarse fiber.
\end{proof}

\begin{lemma}\label{lem:boundary-degree-one}
Every boundary character class $b_\chi$ has strong perverse degree at
most one (Definition~\ref{def:perverse-filtration}) on the whole fine
residue master.
\end{lemma}

\begin{proof}
Choose the decomposition in the derived category of mixed Hodge modules,
using projective relative hard Lefschetz. Write
$K_h\cong\bigoplus_i M_i[-i]$, with $M_i$ pure of weight $w+i$.
The algebraic boundary line defines
$K_h\to K_h(1)[2]$. Its $(i,j)$ component belongs to
$\Ext^{2+i-j}(M_i,M_j(1))$, so it vanishes for $j>i+2$.
For $j=i+2$ it is a morphism of pure Hodge modules of equal weight.
Let $I$ be its image. Polarizable pure Hodge modules of fixed weight
are semisimple \cite[Corollary~5.2.13]{Sai88}. On each strict
support the polarization splits a Hodge subvariation on a smooth dense
open by its orthogonal complement. These projections extend across
the boundary with the Hodge filtration by the filtered
minimal-extension formulas \cite[Proposition~3.2.2, formulas
(3.2.2.1)--(3.2.2.2)]{Sai88}. Splitting the surjection
onto $I$ and the inclusion of $I$ into $M_{i+2}(1)$ gives a factorization
\[
 \mathrm{id}_{I[-i]}\colon I[-i]\longrightarrow K_h
 \xrightarrow{\ \smile\, b_\chi\ }K_h(1)[2]\longrightarrow I[-i].
\]
At every point over a residue satisfying~\eqref{eq:additive-generic}
with regular labels, proper base change identifies the middle stalk map
on cohomology with cup product on the Hitchin fiber.
Lemma~\ref{lem:boundary-fibre-zero} makes this map zero, so every
cohomology group of the stalk of $I$ vanishes there.

As a shifted direct summand of the joint direct image, $I$ is ULA over
the residue space. Through any remaining residue take a small disk
whose punctured part avoids the finite union of excluded hyperplanes.
After base change to this disk the nearby cycles of $I$ are zero.
Local acyclicity identifies the special restriction with these nearby
cycles, so $I$ vanishes there too.
Thus the canonical top maps $M_i\to M_{i+2}(1)$ vanish. To see that
the resulting truncation maps are compatible, put
$D_k=(\ptau_{\ge k+2}K_h)(1)[2]$. Its lowest possible perverse degree
is $k$. A map from $\ptau_{\le k}K_h$ to $D_k$ is determined by its
map in perverse degree $k$, which here is the vanishing top map.
Also $\Hom(\ptau_{\le k}K_h,D_k[-1])=0$. The truncation triangle
therefore gives a unique lift to $(\ptau_{\le k+1}K_h)(1)[2]$.
The lifts commute with the truncation maps by the same orthogonality.
This is strong perverse degree at most one. Iterated ULA coordinate restriction carries
the same bound to every residue fiber, including collision fibers.
\end{proof}

The determinant line of Proposition~\ref{prop:pure-ample-reference}
extends over the fixed-weight master using the same universal bundle on
the cover: each central-character determinant has degree zero by its
prescribed parabolic character degree. Scaling the Higgs field linearizes
this line and contracts the Hitchin base to $0$; since the fiber over $0$
is ample by the proposition and ampleness is open, a power of the line is
relatively ample on the whole master.

The weight-$4$ part of the Goldman class already commutes with Springer
projectors. The lowering argument needs this for the whole class,
including its weight-two corrections.

\begin{lemma}\label{lem:whole-goldman-descent}
Let $\pi:Y\to X$ be a connected generic Springer modification. The whole
complex Goldman class on $Y$ descends uniquely to $H^2(X,\C)$, and its cup
product commutes with all endomorphisms of $R\pi_*\Q_Y$.
\end{lemma}

\begin{proof}
Fix a coarse point of $X$ and a representative object. The fine flags of that object
form a product $F$ of Springer fibers with finite effective automorphism
group $H$, and the coarse fiber is $F/H$; keeping the object fixed gives
a morphism from $F$ to the smooth fine moduli stack. By
Proposition~\ref{prop:springer-low-degree}, $H^0(F)=\Q$, $H^1(F)=0$ and
$H^2(F)$ is pure Tate of weight two. Hence $R^0\pi_*\Q=\Q$ and
$R^1\pi_*\Q=0$. Choose a proper surjection $F_0\to F$ from a smooth
projective variety, for example the disjoint union of resolutions of
the irreducible components. By \cite[Proposition~8.2.5]{Del74}, pullback
embeds $\Gr^W_2H^2(F)=H^2(F)$ into $H^2(F_0)$. The pulled-back Goldman
form is a holomorphic two-form on $F_0$, so its class has type $(2,0)$;
the class also comes from the Tate group $H^2(F)$, so has type $(1,1)$.
This class is therefore zero on $F_0$, and injectivity gives zero on $F$. Finite-quotient
invariance gives the same conclusion on $F/H$.
The Leray sequence
$0\to H^2(X,\C)\to H^2(Y,\C)\to H^0(X,R^2\pi_*\C)$ then gives existence
and uniqueness of the descended class. A class pulled back from $X$
commutes with every endomorphism of $R\pi_*\Q_Y$ by the projection
formula. For disconnected centralizers include all central-cover
components and take invariants.
\end{proof}

For a general reductive group take the product of the semisimple
operators with the torus operator of Example~\ref{ex:intro-rank-one}:
the trace form is orthogonal for $\fg=\mathfrak z\oplus\fg_\der$,
product filtrations have additive centers, and the finite central
quotient preserves both operators.

\begin{theorem}\label{thm:common-leading}
For the problems of Theorem~\ref{thm:A}, the normalized Goldman cup
product $E$ has perverse degree at most $2$, with leading part satisfying
relative hard Lefschetz, and weight degree at most $4$, with leading part
satisfying curious hard Lefschetz; both have center $N=\frac12\dim\MB$.
The leading parts are rational. They commute with projectors from
finite symmetries preserving the Goldman class and, after the comparison of
Theorem~\ref{thm:harmonic-transport}, commute with the Springer
projectors, with the shifts of \S\ref{subsec:collisions}.
\end{theorem}

\begin{proof}
On each fine master, Proposition~\ref{prop:pure-ample-reference} at the
master's own rational unitary reference provides the pure class and the
ample line; the fine weights are chosen off the finitely many unit-phase
character walls so that this reference is generic. A regular datum is
connected to such a reference by a separating imaginary residue, a short
angular segment, and a contraction of the imaginary residue at rational
separated angular weights; the testing-character margins and root gaps
keep the stable problems and the centralizers fixed along this path, and
property~\ref{F:harm} supplies the comparison squares.

At the reference, the boundary character lines on the two sides agree
under $\Psi$, being associated with the same universal smooth bundle;
they extend as flat sections along both paths. Integrating~\eqref{eq:DH}
in a logarithm chart gives, on both sides and with the same
coefficients,
\[
 T_t\alpha_{\rm ref}=q[\Omega_t]+\sum_ia_i(t)\,b_{i,t},
\]
with $T_t$ parallel transport and $q\ne0$ the period normalization. On
the Dolbeault side $T_t\alpha_{\rm ref}$ is the class of the relatively
ample line, and each $b_{i,t}$ has strong perverse degree at most one
(Lemma~\ref{lem:boundary-degree-one}); so the perverse-degree-$2$ part of
$E$ is the relatively ample Lefschetz operator, and
Theorem~\ref{thm:decomposition} applies. On the Betti side each
$b_{i,t}$ has weight two, so the weight-$4$ part of $E$ is the horizontal
curious Lefschetz operator of \S\ref{subsec:collisions}. No ampleness
across different stability chambers is used.

Lemma~\ref{lem:whole-goldman-descent} makes $E$ commute with all Springer
projectors, and property~\ref{F:harm} transports this to Dolbeault
cohomology. On a multiplicity with Springer fiber dimension $b$ the
summand is $\IH^j(Z)(-b)\subset H^{j+2b}(Y)$, with perverse shift $b$ and
weight shift $2b$; removing them gives the center of the smaller orbit
closure. Finally take the union of compatible central lifts and finite
invariants.
\end{proof}


\section{The lowering operator on Betti cohomology}\label{sec:lowering}

We now construct the lowering counterpart of $E$. The idea is to add a
flag at a moving point $q\in\Cc$, wind the eigenvalues at $q$ around a
loop, and then remove the flag and integrate over $q$. Winding produces a
unipotent monodromy whose logarithm lowers weights; adding and removing
the flag are algebraic and preserve weights. To integrate over the
compact curve $C$ (and not only over $\Cc$) the construction must be
attached at each puncture; these attachments also produce the boundary
terms in the commutator with $E$.

We work at a regular separating reference and first for semisimple $G$;
the central torus and finite central descent are treated in
\S\ref{sec:proof}. Fix $T\subset B\subset G$ and $b=\dim G/B$. Three
different motions occur:
\begin{itemize}[leftmargin=2em]
\item the \emph{spatial braid} moves $q$ around a puncture, with the
auxiliary eigenvalue $t\in T$ fixed;
\item the \emph{auxiliary winding} moves $t$ around a cocharacter loop
$\gamma_\nu(a)=t_*\exp(-2\pi ia\nu)$, $a\in[0,1]$, $\nu\in X_*(T)$;
\item the \emph{scalar approach} moves $t$ to $1$, where the auxiliary
point can be filled in and the original problem reappears.
\end{itemize}
Only the auxiliary winding defines the logarithm $N_\gamma$.

\emph{Terminology.} The \emph{physical points} are the punctures in $D$;
$q\in\Cc$ is the \emph{auxiliary} point. A \emph{cap} is a small disk
around a physical point, containing $q$ when $q$ approaches that point;
its boundary annulus is a \emph{collar}. A \emph{carrier} is the space
$Z$ supporting a cohomological correspondence
(\S\ref{subsec:correspondences}), together with its two maps. A
\emph{prism} is a parametrized homotopy in a coefficient diagram,
together with its boundary faces. The \emph{Cartan open} $\Omega\subset
T$ is the set of auxiliary values $t$ for which every fused value
$S_pw^{-1}(t)$ ($p\in D$, $w\in\WG$) is regular and the resulting
boundary data are intrinsically generic.

\subsection{Adding and removing a scalar flag}\label{subsec:scalar-flag}

At $q\in\Cc$ introduce an extra puncture carrying a full flag and an
ordered Cartan value $t\in T$ of the monodromy, and let $Y_t(q)$ be the
resulting weak character variety (\S\ref{subsec:nah}). At $t=1$ the locus
where the auxiliary monodromy is itself $1$ is
\[
 i:E_q\hookrightarrow Y_1(q),\qquad
 \pi:E_q=\Flag(\mathcal E|_q)\longrightarrow\MB.
\]
Here the evaluation $G$-torsor $\mathcal E|_q$ is defined on the
unrigidified moduli stack. The center acts trivially on $G/B$, so the
associated flag bundle $E_q$ descends to the rigidified stack; the
$G$-torsor itself need not descend. In the present semisimple case the
rigidification $r$ is a finite central gerbe. Locally its fibers are
$BZ(G)$, and finite-group averaging gives $Rr_*\Q=\Q$, hence pullback
identifies the rational cochain objects used here. We form universal torsors and their
characteristic operations before rigidification and descend the rational
operations through this equivalence. Pullback along $r$ commutes with
the Gysin maps and proper flag traces used below, since their squares
are base-change squares. A character line descends after taking a power
trivial on $Z(G)$; dividing its first Chern class by that power defines
the rational class. This convention also applies to the capped bundles
of \S\ref{subsec:surface-class}.

The map $i$ is a regular embedding of codimension $b$ with normal bundle
$T_\pi^*$, the relative cotangent bundle of the flag bundle. Define
\begin{equation}\label{eq:scalar-legs}
 B=i_*\pi^*,\qquad A=\pi_*i^*,\qquad AB=c\cdot\id,\qquad
 c=(-1)^b|\WG|.
\end{equation}
The value of $c$ is the self-intersection formula
$\int_{G/B}c_b(T^*(G/B))=(-1)^b\chi(G/B)=(-1)^b|\WG|$. For $G=\SL_2$
the cotangent line of $\PP^1$ has degree $-2$, so $AB=-2$.

\subsection{The triangular source}\label{subsec:triangular}

To attach the moving construction near a puncture we need a source that
also exists for $t\ne1$. Let $p$ be a puncture with regular separating
monodromy $S\in T$, $Z_G(S)=T$. Let $\nu:\mathcal R\to G$ be the framed
surface word omitting the factors at $p$ and $q$. Locally
\begin{gather*}
 Y_t=\{(r,U,F):\nu(r)=U^{-1}S^{-1},\ (U,F)\in\widetilde G_t\}/T,\\
 \widetilde G_t=\{(u,gB):g^{-1}ug\in B,\ \overline{g^{-1}ug}=t\in B/U_B\}.
\end{gather*}
The fixed-flag locus $L(t)=\{F=B\}$ relates this problem to the main
problem in which the monodromies at $p$ and $q$ are fused into $St$.

\begin{lemma}[The triangular source]\label{lem:triangular-source}
On the open set where $St$ is regular and the boundary data are
intrinsically generic,
\[
 L(t)\cong\bigl[\nu^{-1}((St)^{-1})\times U_B\,/\,T\bigr].
\]
The projection $p_t:L(t)\to\MB(St)$ to the fused main problem is an
$\A^b$-bundle, and $L(t)\hookrightarrow Y_t$ is a closed regular embedding
of codimension $b$; hence $G(t)=(i_{L(t)})_*p_t^*$ is defined on the
whole main cohomology. At $t=1$,
\begin{equation}\label{eq:triangle-scalar}
 A\circ G(1)=\id.
\end{equation}
The spatial braid around $p$ acts trivially on $Rp_{t*}\Q$, compatibly
with its powers and with the Gysin map.
\end{lemma}

\begin{proof}
For regular $c\in T$ the map $n\mapsto c^{-1}ncn^{-1}$ on $U_B$ is solved
in increasing root height: its component in the $\alpha$-coordinate has
linear coefficient a nonzero multiple of $\alpha(c)-1$, and successive
solution gives an algebraic inverse, also in families. So every
$SU\in StU_B$ is uniquely $n(St)n^{-1}$, and replacing $r$ by
$n^{-1}rn$ gives the product presentation.

The cokernel of the differential of the main word map is dual to the
infinitesimal stabilizer. Here is the complex-algebraic verification;
\cite[Proposition~4.1(3)]{AMM98} is the compact quasi-Hamiltonian version.
Use an invariant nondegenerate bilinear form and right logarithmic
differentials. For a conjugacy-class factor $M$, the annihilator of its
differential is $\fg^M$. For a commutator $[A,B]$, differentiation in
$A$ and $B$ shows that an annihilator $X$ is fixed by $ABA^{-1}$ and
satisfies $\Ad_{A^{-1}}X=\Ad_{[A,B]^{-1}}X$; these two equations are
equivalent to $X\in\fg^A\cap\fg^B$. In a product word, transport $X$
by the inverse of each preceding product. Each preceding factor fixes
the transported $X$, so these conditions give exactly the common
infinitesimal stabilizer of the word variables. At the fused value
$(St)^{-1}$ this is also the stabilizer of the fused representation,
which is zero by genericity for semisimple $G$. Thus the word map is
smooth along this incidence. Pulling back the fixed point $B\in G/B$
in the Grothendieck chart gives a regular embedding with complex
orientation. At $t=1$, the scalar
incidence and the fixed-flag incidence are the zero section and one
cotangent fiber in the Springer chart; they meet transversally in a
section of $E_q\to\MB$, which gives~\eqref{eq:triangle-scalar} with
positive orientation by refined Gysin base change.

The braid acts on the product presentation by
$[r',n]\mapsto[(St)^{-1}r'(St),n]=[r',(St)n(St)^{-1}]$. The second form
fixes the main source; the gauge $(St)^m$ compares the two presentations
of the $m$-th power, and $(St)^{m+n}=(St)^m(St)^n$ supplies the
composition comparisons. The affine-bundle unit is an equivalence of
direct images on this action diagram, and its inverse, unit and counit
give the trivialization.
\end{proof}

We use the Springer action of $\WG$ obtained by specializing the deck
action of the Grothendieck--Springer resolution; in this convention the
zero orbit carries the sign representation
(Proposition~\ref{prop:springer}\ref{it:springer-action}). On the proper
Springer direct image,
\begin{equation}\label{eq:scalar-sign}
 BA=c\,e_-,\qquad e_-=|\WG|^{-1}\sum_{w\in\WG}\sgn(w)\,w .
\end{equation}
Indeed the zero-orbit multiplicity is the one-dimensional top cohomology
of $G/B$, and the scalar Gysin and trace maps are its inclusion and
projection, with composite $c$. We form these maps in the universal proper
Springer diagram and pull back by the framed word map; proper base change
and the smoothness just proved identify them with $A$ and $B$.

\subsection{The puncture packet}\label{subsec:packet}

The ordered auxiliary eigenvalues force us to keep one main source for
each $w\in\WG$. Let $H=R\Gamma(\MB,\Q)[-2b](-b)$ be the shifted main
complex. On the regular Cartan open, the deck map $d_w:Y_{w^{-1}t}\to Y_t$
gives $\Phi_w(t)=d_w\,G(w^{-1}t)$, whose source has fused monodromy
$Sw^{-1}(t)$.

\begin{proposition}[The puncture packet]\label{prop:finite-packet}
At the scalar limit put $P=\Q[\WG]\otimes H$ and
\[
 \Phi([w]\otimes x)=wG(x),\qquad
 b_P(x)=(-1)^b\sum_w\sgn(w)[w]\otimes x,\qquad
 a_P([w]\otimes x)=\sgn(w)x .
\]
Then
\begin{equation}\label{eq:packet-endpoints}
 \Phi\,b_P=B,\qquad A\,\Phi=a_P,\qquad a_Pb_P=c,
\end{equation}
and the two reductions of $A\Phi b_P$ to $c$ are joined by a specified
homotopy. On the $w$-summand, the auxiliary winding is the winding of the
fused monodromy $Sw^{-1}(t)$.
\end{proposition}

\begin{proof}
The maps $\Phi_w(t)$ are algebraic over one Cartan open, and proper
Grothendieck--Springer specialization gives $wG$ at the scalar limit.
By~\eqref{eq:scalar-sign} and~\eqref{eq:triangle-scalar},
$\Phi b_P=ce_-G=BAG=B$ and $A\Phi=a_P$; summing signs gives $a_Pb_P=c$.
The braid fixes the auxiliary group element and acts on the main word
variables independently of the flag, so deck relabelings, Springer maps
and Gysin orientations are equivariant on a common action diagram. Both
scalar reductions use the same universal Springer relation and the same
transverse square, and their common composite is
$ABAG\to cAG\to c$; \S\ref{subsec:coefficient-enhancement}
explains why this is a homotopy before passing to cohomology. Continuation
of the source family $\MB(Sw^{-1}(t))$ identifies its loop with the
auxiliary loop.
\end{proof}

A \emph{coefficient diagram} below means a finite diagram of complexes,
chain maps, and specified homotopies between composites (and
homotopies between those); this is what is needed to glue over a curve
(Example~\ref{ex:collar}).

\subsection{The compact coefficient and the weight estimate}
\label{subsec:compact-coefficient}

Choose disjoint disks around the punctures. On $\Cc$ the moving
construction gives a coefficient $\mathcal V_t$ (the cohomology of
$Y_t(q)$ as $q$ varies). At $p$ use the packet
$P_{p,t}=\bigoplus_{w\in\WG}R\Gamma(\MB(S_pw^{-1}(t)),\Q)[-2b](-b)$.
The braid-equivariant maps $\Phi_w$ give an attachment
$a_{p,t}:P_{p,t}\to R\Gamma(S^1_p,\mathcal V_t)$: insert the constant
circle cochain on the trivially acting main source, then apply the
equivariant map. Define the \emph{compact coefficient}
\begin{equation}\label{eq:compact-coefficient}
 \mathfrak K_t=\fib\Bigl(R\Gamma(\Cc,\mathcal V_t)\oplus\bigoplus_p
 P_{p,t}\xrightarrow{\ \res-a\ }\bigoplus_pR\Gamma(S^1_p,\mathcal V_t)
 \Bigr).
\end{equation}
Thus the attachments and their homotopies enter the differential of
$\mathfrak K_t$.

For the weight computation, the curve is represented by a finite cellular
complex in which every cell is assigned $\Q(0)$; we write
$C^*(C,\Q(0))$. The rational realization of $C^*(C,\Q(0))$ computes
$H^*(C)$, but the assigned Hodge structures are part of this coefficient construction (all
of weight zero), not the mixed Hodge structure of the curve;
\S\ref{subsec:coefficient-enhancement} constructs it together with
the attachments.

The weight estimate uses the integral monodromy lattices of the
geometric vertex variations: a pure Tate quotient has a definite
polarization, so its integral monodromy group is finite. The finite
coefficient cones inherit their monodromy eigenvalues from these
vertices. A common power then acts trivially on every pure weight
quotient, so its logarithm lowers weights by at least two. The integral
lattices are essential; Hodge--Tate type alone would not imply this
finiteness. The work consists in building the compact
coefficient~\eqref{eq:compact-coefficient}, with all its attachments,
inside mixed Hodge modules.

\subsection{The coefficient diagram and its weights}
\label{subsec:coefficient-enhancement}

We use the $\infty$-categorical enhancement of mixed Hodge modules and its
six operations on stacks \cite[\S2.1, Theorems~4.1.1 and~4.1.8,
and Propositions~4.2.5--4.2.6]{Tub25}.
For an algebraic function $f$ to $\A^1$, coherent proper and smooth
exchange for unipotent nearby cycles $\Upsilon_f$ is
\cite[Theorem~4.4.3]{Tub25}. Write
$L_n=(u\mapsto u^n)_*\Q$ on $\Gm$, with its cyclic permutation action.
Total nearby cycles are the filtered colimit of
$\Upsilon_f(M\otimes f_\eta^*L_n)$, where $f_\eta$ is the restriction
over $\Gm$. For $n\mid m$ the transition pulls functions on the $n$
sheets back to the $m$ sheets. Each monodromy eigenvalue is paired with
its inverse character once $n$ is sufficiently divisible. The deck
action records the finite part. For example the sign local system
on $\Gm$ has nonzero total fiber and zero untwisted unipotent fiber
\cite[(5.7.2)--(5.7.6)]{SaiMix}. This yields the cyclic-cover construction
of \cite[Definition~4.4.8]{Tub25}: under $t=u^n$ the logarithm
variation pulls back by $\log t=n\log u$, and finite proper base
change at the special fiber and the projection formula on the
punctured fiber give the displayed tensor factor
\cite[\S1.1 and Proposition~1.3]{SaiExt}. At each stage the projection
formula for a proper output map and pullback of $L_n$ on a smooth
source square give the unipotent exchange maps compatibly with these
transitions. Taking their colimit gives total exchange. Under $t=u^n$,
the monodromy logarithm on the $u$-cover is $n$ times the original
logarithm, so the comparison uses division by $n$. The scheme
comparisons with Saito are triangulated, and shifted total nearby
cycles preserve bounded constructibility and are perverse exact
\cite[Propositions~4.4.6 and~4.4.10 and Theorem~4.4.11]{Tub25}.
An enhancement keeps homotopies between composites (and higher ones),
whereas the derived category keeps only their classes. The moduli
charts used in the finite coefficient diagram are
of finite type and have finite effective inertia after central
rigidification. Their Gysin embeddings, affine-bundle maps and Springer
maps are representable. Universal evaluation and equivariant Springer
stacks may have positive-dimensional affine stabilizers; their relations
are formed before pullback to these moduli charts. We do not take their
possibly unbounded equivariant global cohomology as a coefficient of the
compact complex. Rational realization is conservative on the bounded
constructible objects used, and the affine-bundle unit is an equivalence
in the enhanced category. Bounded constructibility is checked on the finite-type
charts, then descends through finite inertia with rational coefficients;
these bounds also apply to all the finite cellular cones below.
The enhancement supplies coherent functors for the geometric diagram
constructed here. Conservativity detects equivalences, while
faithfulness on the perverse heart checks the universal Springer
relations. We do not lift arbitrary rational homotopies through rational
realization. The characteristic-class homotopies in \S\ref{subsec:physical-lines} and \S\ref{subsec:surface-class}
are constructed separately in rational cochains.

On a finite quotient chart $[U/\Gamma]$, rational invariants are the image
of the exact projector $|\Gamma|^{-1}\sum_{\gamma\in\Gamma}\gamma$;
they introduce no higher group cohomology. A finite collection of
representable operations and cellular cones therefore has a common
cohomological bound. For an affine bundle $p$, local triviality and
conservativity applied to the adjunction unit give $K\simeq Rp_*p^*K$.
Nearby-cycle exchange is used as an isomorphism for proper output maps
and smooth pullbacks; for a general source map we retain the natural
comparison arrow.

\emph{1. Words give the spatial action.} Remove the physical disks, choose
an outer base point with paths to tangential base points at the
punctures. The punctured surface retracts onto a graph with free
generators $a_1,\ldots,a_r$, $r=2g+|D|-1$; with the auxiliary point
removed we use the fundamental \emph{groupoid} on these base points.
Point-pushing isotopies act on it, and a physical circle acts by its word
$w_p$. A framed representation is given by finitely many group elements,
the relation and boundary decorations; substituting path words uses
products and inverses in $G$, so the action and its composition identities
are algebraic, and its nerve is an algebraic diagram. The outer base object is a section of the moving-puncture surface
bundle. A pointed graph is a $K(\pi,1)$, and the space of its pointed
self-equivalences has discrete homotopy type: maps into a $1$-type have
no higher homotopy, and a natural automorphism of a connected groupoid
functor that is the identity at the base object is the identity along
every path. The word action therefore models the moving-point family
with its higher comparisons. Gauges, including central gauges, are
retained as evaluated $G$-words before rigidification.

\emph{2. The physical-circle attachment.} In the two-holed collision disk
put $C_{\rm tot}=MU$; the full twist conjugates $M,U$ by $C_{\rm tot}$ and
sends the tangential path $\gamma_p$ to $C_{\rm tot}\gamma_p\delta_p^{-1}$.
For the $w$-source write $C(t)=Sw^{-1}(t)$, $SU=nC(t)n^{-1}$ and
$r'=n^{-1}rn$; the two presentations
$\beta^m(r',n)=(C(t)^{-m}r'C(t)^m,n)$ and
$\beta_0^m(r',n)=(r',C(t)^mnC(t)^{-m})$ are compared by the gauge $C(t)^m$,
with composition faces from $C(t)^{m+n}=C(t)^mC(t)^n$. The second fixes
the main source, and the affine-bundle unit gives the trivial source
action used by the Gysin attachment.

\emph{3. Auxiliary evaluation and the Euler class.} Keep the whole local
system on the auxiliary boundary with its Borel reduction (locally
$[B/B]$). Changing the auxiliary path by $\delta_q^k$ fixes $(U,F)$ but
changes evaluation by the gauge $U^k$ and a character line by $t^k$. On
the graph the auxiliary tangent circle bundle has a section; relative to
positive frames on the caps its transition windings $n_p$ satisfy
$\sum_pn_p=\deg T_C=2-2g$ (clutching for the tangent line), so the total
evaluation face is $U^{\sum n_p}$ with character multiplier
$t^{\sum n_p}$. At $U=1$ these gauges are identities, but on the weak
fiber $t=1$ the element $U$ may be nontrivial unipotent; $U^k$ is a
comparison of boundary evaluation, not a global gauge of the outside
representation: $U$ need not centralize its other holonomies. Changing
the graph section changes the local rotation data by a coboundary;
it leaves their total Euler class unchanged.

\emph{4. A finite cellular formula.} With the left-action convention, the
open and circle cochain models are
$C_{\Cc}(V)=[V\xrightarrow{(\mathsf T(a_j)-1)_j}V^{\oplus r}]$ and
$C_w(V)=[V\xrightarrow{\mathsf T(w)-1}V]$. For a word $w=l_1\cdots l_k$
define the Fox derivatives $\mathsf J_{w,j}$ (add
$\mathsf T(l_1\cdots l_{i-1})$ when $l_i=a_j$, subtract
$\mathsf T(l_1\cdots l_{i-1})\mathsf T(a_j^{-1})$ when $l_i=a_j^{-1}$);
then $\mathsf T(w)-1=\sum_j\mathsf J_{w,j}(\mathsf T(a_j)-1)$,
$\mathsf J_{uv}=\mathsf J_u+\mathsf T(u)\mathsf J_v$ and
$\mathsf J_{u^{-1}}=-\mathsf T(u^{-1})\mathsf J_u$. Restriction to $w_p$
is the identity in degree $0$ and $(v_j)\mapsto\sum_j\mathsf J_{w_p,j}v_j$
in degree $1$. The physical gauge and affine unit give $P_p\to C_{w_p}(V)$
with a nullhomotopy of $(\mathsf T(w_p)-1)\Phi$, and the compact object is
$\fib\bigl(C_{\Cc}(V)\oplus\bigoplus_pP_p\to\bigoplus_pC_{w_p}(V)\bigr)$.

This nullhomotopy is obtained before rational realization. The gauges
$C(t)^m$ compare the two source actions, with composition comparisons
given by their multiplication identity. The second source action is
strictly over the fused main problem, so the affine-bundle unit
identifies its direct image coherently with the trivial action on
$P_p$. Naturality of the complex-oriented Gysin map makes $\Phi$
equivariant for this trivial action and the circle action on the weak
coefficient. Applying circle cochains gives
\[
 P_p\longrightarrow P_p^{h\mathbb Z}
 \longrightarrow V^{h\mathbb Z}
 =\fib\bigl(\mathsf T(w_p)-1:V\longrightarrow V\bigr),
\]
where the first map is the constant circle cochain. In the two-cell
circle model this is precisely $\Phi$ together with the specified
nullhomotopy. The same geometric power faces give its composition
comparisons. Thus no arbitrary rational homotopy is lifted through
realization.

\emph{5. The common scalar homotopy.} The universal Springer objects
$\mathsf K=R\mu_*\Q_{T^*(G/B)}[2b]$ and $\mathsf Z_0=i_{0*}\Q(-b)$ lie in
the heart of equivariant perverse Hodge modules. The semismall
decomposition and the deck-specialization convention are those of
\cite[\S\S1.6--1.9, 2.2, 2.6]{BM83}; the zero-orbit summand is the
one-dimensional top flag cohomology with sign action. The scalar Gysin
and trace identify its inclusion and projection, and self-intersection
fixes their normalization, giving $AB=c$, $BA=ce_-$ and
$Aw=\sgn(w)A$. These maps lift to Hodge modules by the geometric
operations; faithfulness of rational realization on the heart then
checks the relations there. Mapping
spaces in a heart are discrete ($\Hom(\mathsf P,\mathsf Q[-k])=0$ for
$k>0$), so these relations have unique composition faces. The triangular
source need not be perverse; its relation $AG=1$ comes from the transverse
Gysin square. Thus $\Phi b_P$ reduces by $ce_-G\to BAG\to B$, and after
applying $A$ both reductions are the pasting $ABAG\to cAG\to c$. At $U=1$, filling the auxiliary meridian makes the main spatial
coefficient constant, including the physical tangential objects. These
maps and their common face give $A_CB_C=c$ on the compact fiber.

\emph{6. Cell weights and the cross class.} The constant cellular model
has terms $\Q^{1+|D|}$, $\Q^{r+|D|}$ and $\Q^{|D|}$ in degrees $0$, $1$
and $2$, each assigned $\Q(0)$; its cohomology has dimensions $1,2g,1$, all of weight
zero, with orientation $C^*(C,\Q(0))\to\Q(0)[-2]$. For a physical
character, the forward multiplier $g(t)^{-1}=\chi(w^{-1}t)^{-1}$
defines the Kummer extension of $\Q(0)$ by $\Q(1)$. On a logarithm
cover it is given by $W_{-2}=\Q e_1$, $W_0=\Q e_0\oplus\Q e_1$ and
$F^0=\C\bigl(e_0-(\log g/2\pi i)e_1\bigr)$; with $F^{-1}$ the whole bundle and $F^1=0$. A change
$\log g\mapsto\log g+2\pi ik$ accompanies $e_0\mapsto e_0+ke_1$.
This variation also has an algebraic construction. Take the homotopy
fiber of the two restriction maps from the cohomology of $\Gm$ over
the parameter base to the sections $1$ and $g$, keeping the sections
separately labeled even where they coincide. Its first relative
cohomology, twisted by $(1)$, is an extension of $\Q(0)$ by $\Q(1)$.
Integrating $dx/(2\pi i x)$ along a path from $g$ to $1$ gives the
displayed sign. Thus the coefficient belongs to the geometric Hodge
diagram by the algebraic operations of \cite[\S4]{Sai90}.

At a normal-crossing boundary write $\log g=\sum_i m_i\log z_i+h$ with
$h$ holomorphic. In the Deligne frame
$\widetilde e_0=e_0-(\sum_i m_i\log z_i/2\pi i)e_1$, $e_1$, the
extended $F^0$ is generated by $\widetilde e_0-(h/2\pi i)e_1$.
It projects to $1$ in $\Gr^W_0$, so the Hodge and weight graded bundles
are locally free on every boundary stratum. Each residue lowers $W$ by
two and acts trivially on its graded pieces. Hence $W$ is the relative
monodromy filtration for every sum of the residues, including successive
relative filtrations and the multiplicity-weighted sums on curve
pullbacks. These computations verify admissibility for this variation.
Tensored with the circle cell, it is the mixed Chern class of the clutch,
so the parameter--circle cross class survives in the weight construction.

\emph{7. Continuation and logarithms.} Spatial word automorphisms are
morphisms of mixed Hodge modules over the Cartan open; auxiliary winding
is parallel transport of the smooth Hodge coefficient, preserving flat
$W$ but possibly not the Hodge filtration of an endpoint.

\begin{lemma}[A finite logarithm on the diagram]\label{lem:bounded-diagram-log}
Let $\mathsf E$ be an object of a diagram category $\operatorname{Fun}(I,
D^b(\Q))$ with cohomology in degrees $[d_0,d_1]$. If $H^i(f)=0$ for all
$i$, then $f^{d_1-d_0+1}=0$. In particular, if $(T^m-1)^s$ kills all
cohomology, then $(T^m-1)^{s(d_1-d_0+1)}=0$, and
$m^{-1}\log(T^m)$ is a finite polynomial on the whole diagram.
\end{lemma}

\begin{proof}
The composite of $f$ with the top cohomology projection
$\mathsf E\to H^{d_1}(\mathsf E)[-d_1]$ is zero by $t$-structure
orthogonality, so $f$ factors through $\tau_{\le d_1-1}\mathsf E$; repeat
$d_1-d_0+1$ times. Finite exponential and logarithm identities give
$\log\Ad(T^m)=\ad\log(T^m)$.
\end{proof}

This is nilpotence in the homotopy category of the whole diagram,
not strict nilpotence of an arbitrary chain representative. Choose
one common polynomial bound for its finitely many object types; the
result retains the original winding and all its attaching faces.
The weight estimate for this logarithm is the separate variation
argument in Proposition~\ref{prop:betti-weight}.

\begin{lemma}[Extension through the scalar parameter]
\label{lem:scalar-coefficient-extension}
The finite compact-coefficient diagram extends through the auxiliary
Cartan value $t=1$, with its word actions, ordered deck operations,
physical-circle attachments and their composition faces. Its cohomology
is lisse there, with the fiber base change of property~\ref{F:betti}.
Its scalar fiber is the compact diagram with the Springer operations
and normalized scalar Gysin maps constructed above.
\end{lemma}

\begin{proof}
The point is to extend the deck operations on the proper coefficient
before forming the compact fiber. Write $\chi:G\to T/\WG$,
$d=\dim G$, and
\[
 \widetilde G=\{(U,F):U\text{ preserves }F\},\qquad
 X=(G\times_{T/\WG}T)_{\rm red}.
\]
Remembering the ordered Cartan value gives a projective map
$\widetilde\mu:\widetilde G\to X$. Every compatible $(U,t)$ is reached:
choose a Borel containing the commuting Jordan parts of $U$, then
choose the required Weyl label. Thus $\widetilde\mu$ is surjective,
and $X$ is irreducible because $\widetilde G$ is irreducible.
The map is birational over the dense
regular semisimple open where the Weyl action on $T$ is free.
The map $\mu:\widetilde G\to G$ is small
\cite[Corollary~1.7 and Proposition~4.2]{Sho88}.
Here is the dimension argument. Put $b=\dim(G/B)$ and let $S_n$ be a
conjugation-invariant stratum with $\dim\mu^{-1}(U)=n>0$.
Then $\dim S_n+n=b+\dim(S_n\cap B)$. In
$B\times_G\widetilde G$, the Bruhat cell indexed by $w$ has dimension
$\ell(w)$ and its group fiber has dimension $\dim B-\ell(w)$.
Each incidence stratum is irreducible of dimension $\dim B$ and
contains regular semisimple elements. Removing those elements gives
$\dim(S_n\cap B)+n<\dim B$, hence $\dim S_n+2n<d$.
Since $X\to G$ is finite, the same strict fiber-dimension inequality
holds for $\widetilde\mu$. Smallness and proper duality now give
\[
 \mathsf K_X:=R\widetilde\mu_*\Q_{\widetilde G}[d]=\IC_X
\]
in the Hodge heart: the strict support and cosupport inequalities
characterize the intermediate extension.

The Weyl group acts on $X$ by changing the ordered Cartan value and
fixing $U$. Intermediate extension extends the regular deck
identifications uniquely to $\mathsf K_X$. The same uniqueness gives
the conjugation structure and its interchange with the Weyl action:
on $G\times X$ use the conjugation and projection pullbacks with the
common smooth shift $[d]$, and on $G^r\times X$ use $[rd]$.
These are identities between heart objects before word-map pullback.
Pull back along the actual framed word map and use proper base change
for $\widetilde\mu$, then push to the Cartan base. This defines the
semilinear ordered deck operations through $1$ without requiring a
deck automorphism of a singular flag resolution. Their scalar
restriction is the prescribed Springer action, including the sign
convention of \S\ref{subsec:lie-springer}.

Choose a Weyl-invariant neighborhood of $1$ in which every fused
$C_w(t)=S_pw^{-1}(t)$ remains regular and the fused data remain
generic. The weak boundary tuple is intrinsically generic at $1$:
appending the semisimple value $1$ adds no character-product
obstruction. The triangular source maps and their Gysin arrows are
defined there. Their physical-circle comparison is the same
algebraic gauge $C_w(t)^m$; the identities for its powers supply the
same composition faces at $1$. The word transformations are evaluated
conjugations of $U$, so their interchange with the extended deck
operations is the pulled-back conjugation face above. Fox restriction
maps are finite sums of these word operations. Thus the attachment
nullhomotopies and their pastings extend before the graph, circle and
compact fibers are formed. The boundary evaluation $U^k$ is retained
on the unipotent weak fiber. The character multipliers have value $1$
at the scalar point and a nonsingular local logarithm, so the Kummer
coefficients extend too. These are the original geometric Gysin and
gauge maps; their shifts, orientations and trace normalizations are
unchanged.

Property~\ref{F:betti} supplies lisse cohomology and scalar fiber base
change for the weak and fused vertices. Finite cohomology exact
sequences then give lissity of the entire compact fiber. Consequently
its ordinary scalar fiber identifies naturally with nearby fibers,
compatibly with all the extended arrows. At $1$ the attachments are
$wG(1)$; the universal Springer relations and the transverse square
$AG=1$ give exactly the common scalar face constructed above.
\end{proof}

\subsection{Weight decrease and the returned operator}\label{subsec:weight-return}

\begin{proposition}[Weight decrease]\label{prop:betti-weight}
Let $T_\gamma$ be the auxiliary winding acting on~\eqref{eq:compact-coefficient}.
A power of $T_\gamma$ is unipotent, and the normalized logarithm
$N_\gamma$ of that power lowers weights by two on the cohomology of the
whole diagram. The scalar maps extend to maps of compact coefficients
\[
 B_C:C^*(C,\Q(0))\otimes H\to\mathfrak K_{\lim},\qquad
 A_C:\mathfrak K_{\lim}\to C^*(C,\Q(0))\otimes H,\qquad A_CB_C=c .
\]
Consequently
\begin{equation}\label{eq:raw-return}
 F_\nu^{\rm raw}=\tr_C\circ A_C\circ N_\gamma^2\circ B_C\circ
 \mathrm{unit}_C
\end{equation}
has cohomological degree $-2$ and weight degree at most $-4$ on
$H^*(\MB)$.
\end{proposition}

\begin{proof}
Choose a finite graph retract of $\Cc$ and keep the restrictions to the
puncture circles. Word maps, Gysin maps, affine-bundle units and ordered
deck maps are algebraic over a common Cartan open, and
\S\ref{subsec:coefficient-enhancement} constructs the braid action
and attaching homotopies before applying the mixed Hodge module
operations; so the finite complex has mixed Hodge module coefficients.
By property~\ref{F:betti}, the main and weak cohomologies are local
systems on the parameter region, hence so is the cohomology of the total
complex. Its weight filtration is flat and its pure quotients are
polarized variations, as in \S\ref{subsec:families}. The vertex
coefficients are Hodge--Tate
(Theorem~\ref{thm:regular-chl}) and the cellular terms have weight zero,
so the total cohomology is Hodge--Tate.

The geometric vertex variations have invariant integral lattices.
To see this without assuming integral base change on the whole parameter
region, choose a common dense algebraic open $S^\circ$ of the smooth
rational-lisse region $S$ where integral constructibility and fiber
comparison hold for the finitely many vertices. Their cohomology modulo
torsion is an integral local system there, inside the given rational
local system. The map $\pi_1(S^\circ)\to\pi_1(S)$ is surjective, since
loops can be perturbed off the removed complex subset of positive
codimension. Hence the rational monodromy of $\pi_1(S)$ preserves the
same lattice; the kernel acts trivially because it acts trivially on
its rationalization. A rational invariant summand inherits a full
invariant lattice by intersection, and a flat weight quotient inherits
the corresponding quotient lattice modulo torsion. Each pure Tate weight quotient has
a definite polarization, so its integral monodromy group is finite.
Consequently the eigenvalues of the chosen winding on the finitely many
vertex cohomologies are roots of unity. The graph, circle and compact
cones are finite monodromy-equivariant homotopy-fiber constructions.
Their cohomology exact sequences show, successively, that eigenvalues
on the total cohomology occur among those on the vertices.
Choose a common power having only eigenvalue
$1$ on total cohomology. On each pure Tate quotient this power is a
unipotent isometry of a definite form, hence is the identity. Its
logarithm strictly lowers the flat weight filtration and, since odd
weight quotients vanish, lowers weights by at least two. The parameter
region contains the whole winding loop.

At the scalar approach, the main sources become $H$ and the weak
coefficient has its scalar restriction isomorphism; proper Springer
specialization and refined Gysin base change identify the endpoint
squares of~\eqref{eq:packet-endpoints}, whose common homotopy glues the
maps on the compact complex and gives $A_CB_C=c$.
Lemma~\ref{lem:scalar-coefficient-extension} realizes this comparison
as the ordinary fiber of the extended compact diagram. The weight
subconnections remain flat through that point, so the transported
$N_\gamma$ lowers the same $W$ by two. A sufficiently small local loop
around this scalar point has trivial monodromy; $N_\gamma$ remains the
transported logarithm of the original auxiliary winding. The curve unit and the
orientation trace $C^*(C,\Q(0))\to\Q(0)[-2]$ have weight degree zero, the
two logarithms lower weights by four, and the scalar shifts and twists
cancel. If a logarithm is needed on the derived diagram itself,
Lemma~\ref{lem:bounded-diagram-log} supplies a finite polynomial. The
logarithm is taken rationally: parameter continuation preserves flat $W$
but not the Hodge filtration of an endpoint.
\end{proof}

The weight estimate gives the raw returned operator. The next section
fixes its normalization by comparing the surface and physical boundary
classes on the same coefficient diagram.

\section{Characteristic classes and the commutator}
\label{sec:cup-comparisons}\label{subsec:commutator}\label{app:cup}

The compact coefficient and its weights were constructed in
\S\ref{subsec:coefficient-enhancement}. Here we compare bundles and
characteristic classes on that diagram. We first explain the return
calculation and fix its normalization, then prove the physical-line and
surface comparisons used in the calculation. These comparisons are made
in rational cochains: a homotopy between cup operations need not preserve
mixed Hodge structures at an endpoint. Example~\ref{ex:collar} explains
why the overlap and composition homotopies must be retained.
Here $\WG$ denotes the Weyl group.

\subsection{The quadratic return and its normalization}\label{subsec:quadratic-return}

The weight estimate does not determine the normalization of the
lowering operator; the commutator with $E$ does. Put
$q_\kappa(x)=\kappa(x,x)/2$. Let $e_0$ be the surface characteristic
class of the capped principal bundle (the slant product of $q_\kappa$ of
the universal bundle with the fundamental class of $C$) and
$E_0=m_{e_0}$ its cup product. To relate this particular extension
to the Goldman class, first use the finite-order reference of
Proposition~\ref{prop:pure-ample-reference}. Let
$\exp(2\pi i a_p)$ be the physical monodromy in its ordered Cartan,
and let $\xi_p$ be the vector of boundary character classes, so that
$\chi(\xi_p)=c_1(\mathscr L_{\chi,p})$. At a point of the compact
cover with $z=w^m$ one has $ma_p\in X_*(T)$. Relative to the upstairs
flat frame, the pullback cap frame is $w^{-ma_p}$. Thus on a Cartan
character line the pullback cap adds the divisor
$m\langle\chi,a_p\rangle[\widehat p]$. The point $\widehat p$ is
fixed as the moduli point varies, so its normal line is trivial
and the point divisor has square zero. Compute $q_\kappa=\ch_2(\mathcal U_\rho)$
in the common Cartan reduction. For a weight $\mu$ of $\rho$ the
change in the degree-four term is
$m\langle\mu,a_p\rangle\,c_1(\mathcal U_{\rho,\mu})[\widehat p]$,
where the weight bundle is restricted to the point divisor. Summing
over the distinct weights, with their ranks accounted for by $c_1$,
gives $m\kappa(a_p,\xi_p)[\widehat p]$. The exterior bundles agree, and
the ramification indices over $p$ sum to the cover degree $M$.
With the sign of $\alpha_B$ in Proposition~\ref{prop:pure-ample-reference}
this gives
\[
 e_0=-\alpha_B+\sum_{p\in D}\kappa(a_p,\xi_p).
\]
Along a chosen ordered logarithm path the topological cap and its
character lines extend over the family by clutching with the chosen
logarithms. Their classes are therefore flat. The variation formula
\eqref{eq:DH} changes the Goldman class by boundary classes along
this path. Consequently at the regular reference
\begin{equation}\label{eq:whole-raiser}
 E=a_\kappa E_0+\sum_it_i\,\delta_i,
\end{equation}
where $a_\kappa\ne0$ is a normalization constant and the $\delta_i$ are
cup products with physical boundary character classes
(\S\ref{subsec:boundary-classes}).

\emph{The collar mechanism.} Fix a puncture $p$ and an integral physical
character $\chi$. On the $w$-entry of its packet, comparison with the
line from the packet source introduces a relative line
$\mathscr P_{g_w^{-1}}$ with forward braid multiplier
$g_w(t)^{-1}$, where $g_w(t)=\chi(w^{-1}t)$
(Proposition~\ref{prop:equivariant-picard}). Along the auxiliary loop
$\gamma_\nu(a)=t_*\exp(-2\pi ia\nu)$, this multiplier is a constant
times $\exp(2\pi ik_wa)$, with $k_w=\langle\chi,w^{-1}\nu\rangle$.
Thus, if $\alpha$ is the positive auxiliary-circle class and $\theta_p$
the class of the spatial circle oriented by the forward braid,
\[
 c_1(\mathscr P_{g_w^{-1}})=k_w\,\alpha\smile\theta_p.
\]
This is a nonzero cross class when $k_w\ne0$
(Proposition~\ref{prop:mixed-obstruction}).

On the auxiliary universal cover, a logarithm $\ell_w$ of $g_w^{-1}$
gives the collar homotopy comparing the two cup actions. The change
$\ell_w(a+1)-\ell_w(a)=2\pi ik_w$ makes that homotopy change, up to
homotopy, by
$k_w\theta_p\smile a_{p,w}$, where $a_{p,w}$ is the circle attachment
of this packet entry (Proposition~\ref{prop:collar-defect}). The cup
operator on the homotopy fiber~\eqref{eq:compact-coefficient} contains
this homotopy in its collar component. On the scalar compact coefficient,
write $\widetilde\delta$ for the resulting lift of
$\delta=m_{c_1(\mathscr L_{\chi,p})}$. The full deck difference of $\widetilde\delta$ is
therefore
\begin{gather*}
 T_\gamma\widetilde\delta T_\gamma^{-1}-\widetilde\delta
 =J_\chi:=\sum_w k_w\,j_p(\theta_p\smile a_{p,w})\pr_{p,w},\\
 [N_\gamma,\widetilde\delta]=J_\chi,\quad [N_\gamma,J_\chi]=0.
\end{gather*}
Here $j_p$ is the collar connecting map and $\pr_{p,w}$ the disk
projection. After a common unipotent power $M$, the first deck
difference is $MJ_\chi$ and the second vanishes; division by $M$ gives
the displayed identity for the normalized logarithm $N_\gamma$.
These are identities on the whole compact
coefficient, including its attachments and homotopies. The compatible
scalar endpoints and composition faces are proved in
Propositions~\ref{prop:scalar-line-endpoints}
and~\ref{prop:compact-cup-lift}. This collar defect produces the
physical boundary terms in the return calculation below. Integral
multiples clear finite descent obstructions; division then gives
rational physical character classes.

\begin{lemma}[The quadratic return]\label{lem:whole-quadratic-return}
Suppose $G_\der$ is simple of rank $r$, and let $\nu\in X_*(T)$ be a
nonzero allowed cocharacter. Put
\[
 c_\nu=4(-1)^b|\WG|\,\kappa(\nu,\nu)/r.
\]
Then $[[F_\nu^{\rm raw},E_0],E_0]=c_\nu E_0$, and the normalized operator
$F_B=-2F_\nu^{\rm raw}/(a_\kappa c_\nu)$ satisfies
\begin{equation}\label{eq:whole-commutator}
 [[F_B,E],E]=-2E .
\end{equation}
If $\fg_\der$ has several simple factors, write $\kappa=\sum_i\kappa_i$
and $E_0=\sum_iE_{0,i}$, choose $\nu_i$ supported in factor $i$ of rank
$r_i$, and put $c_i=4(-1)^b|\WG|\kappa_i(\nu_i,\nu_i)/r_i$ and
$F_B=-\sum_i2F_{\nu_i}^{\rm raw}/(a_\kappa c_i)$; then
\eqref{eq:whole-commutator} holds. (In all cases $b$ and $|\WG|$ are those
of the full flag variety of $G$.)
\end{lemma}

\begin{proof}
\emph{Idea.} The auxiliary winding changes the capped bundle by a
lattice modification along the moving point, so it changes the surface
class $e_0$ by a quadratic polynomial in the winding number, whose linear
coefficient $y$ is an auxiliary flag class. Hence the double commutator
of $N_\gamma^2$ with $E_0$ is multiplication by $2y^2$. Removing the flag
turns $y^2$ into $c_\nu$ times the quadratic class evaluated at the moving
point (flag Euler formula and a Weyl-group average), and integrating over
the moving point gives back $e_0$. The boundary classes produce collar
terms which, after normalization, recombine with the surface term
into~\eqref{eq:whole-commutator}.

We use identities on the complete compact coefficient diagram, including
its disk and collar components, proved in
\S\S\ref{subsec:physical-lines}--\ref{subsec:scalar-collar-comparison} below.

\emph{Surface term.} First suppose the root system is simple. For the capped bundle, the integer deck
transformation of the auxiliary loop acts by
\[
 T^ke_0=e_0+ky+\tfrac12k^2z,\qquad y=-\kappa(\nu,\xi_q),\qquad
 z=\kappa(\nu,\nu)\,c_1(T_C)
\]
(Proposition~\ref{prop:compact-surface-comparison}), where $\xi_q$ is the
vector of first Chern classes of the auxiliary Cartan character lines.
Semilinearity and the finite logarithm give $[N_\gamma,E_0]=m_y$,
$[N_\gamma,m_y]=m_z$ and $[m_y,E_0]=0$, hence
$[N_\gamma^2,E_0]=2m_yN_\gamma+m_z$ and
\[
 [[N_\gamma^2,E_0],E_0]=2m_y^2 .
\]
The scalar maps are $E_0$-linear. Flag self-intersection and the Weyl
quadratic identity give
\[
 2A_Cm_y^2B_C=c_\nu\,m_{q_\kappa(\mathcal E|_{x=q})},\qquad
 \sum_{w\in\WG}\kappa(w\nu,\xi)^2=\frac{2|\WG|\,\kappa(\nu,\nu)}{r}\,
 q_\kappa(\xi);
\]
the second identity holds because the reflection representation of a
simple Weyl group is irreducible, so
$\sum_w\kappa(w\nu,\cdot)^{\otimes2}$ is a multiple of $\kappa$, and the
multiple is found by taking the $\kappa$-trace. The first combines this
with the flag Euler formula
$\pi_*(c_b(T^*_\pi)f(\xi))=(-1)^b\sum_wf(w\xi)$; the surface comparison proves it on
vertices, edges and composition triangles. Integrating over the curve
(the slant product of $q_\kappa(\mathcal E|_{x=q})$ over $q\in C$ is
$e_0$) gives $[[F_\nu^{\rm raw},E_0],E_0]=c_\nu E_0$.

\emph{Boundary terms.} For a physical character $\chi$ at $p$,
Proposition~\ref{prop:compact-cup-lift} and the sign calculation~\eqref{eq:returned-boundary} give
\[
 [F_\nu^{\rm raw},m_{c_1(\mathscr L_\chi)}]
 =2(-1)^b\sum_{w\in\WG}\langle\chi,w^{-1}\nu\rangle\,L_{p,w^{-1}\nu},
\]
where $L_{p,\mu}$ is the logarithm of the physical monodromy obtained by
the lattice modification $e'=ez^{-\mu}$; it sends $e_0$ to
$\kappa(\mu,\xi_p)$ and kills every physical character class. The Weyl
average and the normalization give
\[
 [F_B,\delta_i]=-L_i,\qquad [L_i,a_\kappa E_0]=\delta_i,\qquad
 [L_i,\delta_j]=0 .
\]
Since $[E_0,\delta_i]=0$, the Jacobi identity gives
$[[F_B,a_\kappa E_0],\delta_i]=[a_\kappa E_0,L_i]=-\delta_i$. Expanding
$[[F_B,E],E]$ for~\eqref{eq:whole-raiser}: the surface term is
$a_\kappa^2\cdot(-2/(a_\kappa c_\nu))\cdot c_\nu E_0=-2a_\kappa E_0$, each
of the two mixed terms is $-\sum_it_i\delta_i$, and the terms quadratic
in the $\delta_i$ vanish. This is~\eqref{eq:whole-commutator}.

For several simple factors, let $\kappa_i$ be the restriction of
$\kappa$ to factor $i$, extended by zero on the other factors.
Winding by $\nu_i$ differentiates only the $i$th characteristic class.
The full flag trace averages over all of $\WG$, so
Proposition~\ref{prop:compact-surface-comparison} gives
$[[F_{\nu_i}^{\rm raw},E_0],E_0]=c_i E_{0,i}$.
The boundary Weyl sum likewise selects the $i$th component of a
character. Summing the normalized returns therefore gives $-2E$. Integral
multiples place the cocharacters in the allowed lattice; $F^{\rm raw}$ and
$c_\nu$ both scale quadratically, so $F_B$ is unchanged.
\end{proof}

\begin{example}\label{ex:SL2-constants}
For $G=\SL_2$, $\rho$ the standard representation and $\nu=\alpha^\vee
=\operatorname{diag}(1,-1)$: $b=1$, $|\WG|=2$, $r=1$, $\kappa(\nu,\nu)=2$,
so $c=-2$ and $c_\nu=-16$.
\end{example}


\subsection{Physical logarithms and character lines}
\label{subsec:physical-lines}

The boundary part of the Goldman class comes from physical character
lines. The commutator of these classes with the compact return is controlled by a line
bundle on the product of the auxiliary parameter circle and the spatial
braid circle.

At a regular separating puncture, the Hecke translation $t_{p,\mu}$ uses
the lattice frame $e'=ez^{-\mu}$ and shifts the affine residue by
$-\lambda\mu$; its finite-jet normalization involves the factors
$\alpha(v_p)+\lambda c_{\alpha,j}$ with $c_{\alpha,j}\in\mathbb R$
and $\Im\alpha(v_p)\ne0$, uniformly invertible on the real interval. Let
$L_{p,\mu}$ be its logarithm after a common unipotent power, divided by
that power. The negative-exponential auxiliary winding closes by
$j_{-\mu}$ with return $q_\mu=f^{-1}j_{-\mu}$ and central translation
$t_{p,-\mu}$, so its logarithm on the packet source is
$\widehat L_{p,\mu}=-L_{p,\mu}$. Translations commute and are additive in
$\mu$. At a fixed physical point the graph has zero self-normal class, and
\begin{equation}\label{eq:physical-log-action}
 L_{p,\mu}(e_0)=\kappa(\mu,\xi_p),\qquad L_{p,\mu}(\xi_{p'})=0
\end{equation}
for every physical character class $\xi_{p'}$.

\emph{The line on the braid diagram.} Over the Cartan open $\Omega$, fix
an integral physical character $\chi$ (integral multiples clear finite
descent obstructions; divide afterward). For one Weyl entry put
$C(t)=Sw^{-1}(t)$, $R_C=\nu^{-1}(C(t)^{-1})$, $Z=R_C\times N_w$,
$X_C=[R_C/T]$, $L=[Z/T]$. Here $N_w$ is the copy of $U_B$ in the
$w$-entry, and the packet arrow is still $d_wG(w^{-1}t)$.
The maps $p:L\to X_C$ and $i:L\to Y$ are the affine bundle
and the regular embedding. Our positive character line is
\begin{equation}\label{eq:line-convention}
 \mathscr L_\chi=[Y^{\rm fr}\times\C/T],\qquad
 a\cdot(x,z)=(axa^{-1},\chi(a)z).
\end{equation}
On the slice $M=S$, the point push acts on the tangential path by
$\gamma_p\mapsto C\gamma_p\delta_p^{-1}$, so the frame of the Hom line
between outer and physical evaluations changes by $\chi(C)/\chi(S)$, and a
vector in the new frame has the inverse multiplier.

\begin{proposition}[Equivariant line comparison]\label{prop:equivariant-picard}
On the triangular braid diagram over $\Omega$,
\[
 \mathscr L_{\chi,\beta}|_L\cong p^*\mathscr L_{\chi,X}\otimes
 q^*\mathscr P_{g^{-1}},\qquad g(t)=\chi(C(t))/\chi(S)=\chi(w^{-1}t),
\]
where $q:L\to\Omega\times B\Z$ records the parameter and the braid action,
the main line has trivial braid action, and $\mathscr P_{g^{-1}}$ has
forward multiplier $g(t)^{-m}$ for the $m$th braid power.
\end{proposition}

\begin{proof}
The triangular parametrization identifies $i^*\mathscr L_\chi$ with
$p^*\mathscr L_{\chi,X}$ before the braid is introduced. The gauge $C^m$
identifies the two action diagrams $\beta^m,\beta_0^m$; relative to the
trivial main-line comparison of $\beta_0$, the path-frame ratio is $g^m$
and the forward vector multiplier is $g^{-m}$. The identities
$C^{m+n}=C^mC^n$ and $g^{-(m+n)}=g^{-m}g^{-n}$ give the composition faces.
\end{proof}

\begin{proposition}[The parameter--circle cross class]
\label{prop:mixed-obstruction}
For the auxiliary loop $\gamma_\nu(a)=t_*\exp(-2\pi ia\nu)$ and the
physical circle oriented by the forward braid clutch, with positive
degree-one classes $\alpha,\theta$,
\[
 c_1(\mathscr P_{g^{-1}})=k\,\alpha\smile\theta,\qquad
 k=\langle\chi,w^{-1}\nu\rangle .
\]
If $k\ne0$ this class is nonzero on the base torus and on the braid space
over the loop.
\end{proposition}

\begin{proof}
Along the loop, $g^{-1}$ is a constant times $\exp(2\pi ika)$, so the
clutching line has Chern number $k$; with the unitary clutch $f=g^{-1}$
the connection $bf^{-1}df$ on the cut circle has curvature
$db\wedge f^{-1}df$ and $c_1=(d\log f/2\pi i)\wedge db$ integrates to $k$
in the auxiliary--physical order. A constant factor has zero winding. The
$\beta_0$ mapping torus is an affine bundle over $X_C\times S^1$, and the
unit of the nonempty finite component system injects base degree-one
classes; K\"unneth with the circle gives nonvanishing.
\end{proof}

\begin{proposition}[The collar cup homotopy]\label{prop:collar-defect}
Let $a_w$ be the circle attachment of one packet entry. The line
comparison gives
$m_{c_1(\mathscr L_\chi)}G_w=G_wm_{c_1(\mathscr L_{\chi,X})}
+c_1(\mathscr P_{g^{-1}})G_w$. A logarithm $\ell$ of $g^{-1}$ on the
auxiliary universal cover gives a degree-one collar homotopy $h_w$ with
\begin{equation}\label{eq:collar-homotopy}
 dh_w+h_wd=m_{c_1(\mathscr L_\chi)}a_w-a_wm_{c_1(\mathscr L_{\chi,X})},
\end{equation}
and for the positive auxiliary deck transformation $\tau$,
$\tau^*h_w-h_w=k\,\theta\smile a_w$ with $k=\langle\chi,w^{-1}\nu\rangle$,
up to homotopy; these commute with the affine unit, the Gysin map and the
deck relabelings.
\end{proposition}

\begin{proof}
Additivity of $c_1$ and the Thom projection formula give the first
identity. On the cut circle the section $s(a,b)=\exp(-b\ell(a))$
satisfies the clutching relation and trivializes the relative line;
transgressing gives $h_w$. Since $\ell(a+1)-\ell(a)=2\pi ik$, the deck
ratio is $\exp(-2\pi ikb)$, i.e.\ the primitive $\ell\,db/2\pi i$ changes by
$k\,db$, which gives the deck difference.
\end{proof}

\begin{proposition}[Scalar endpoints]\label{prop:scalar-line-endpoints}
At auxiliary scalar incidence the physical character line is pulled back
from the main problem on the full moving-point groupoid, its collar
comparison is the scalar specialization of
Proposition~\ref{prop:equivariant-picard}, and both scalar maps and their
common homotopy are modules for its cup product.
\end{proposition}

\begin{proof}
Quotienting the fundamental groupoid of $\Cc\setminus\{q\}$ by the
auxiliary meridians gives that of $\Cc$ (van Kampen), naturally in $q$, so
a representation with auxiliary monodromy $1$ extends uniquely across $q$
with its physical flags. In triangular coordinates the scalar braid is
$(r,F)\mapsto(S^{-1}rS,F)$ and the relative multiplier is $g(1)^{-1}=1$.
So $B=i_*\pi^*$ and $A=\pi_*i^*$ satisfy the projection formula for this
line; at the transverse intersection $n=1$ the physical frame comparison is
the identity, and both endpoint reductions are $ABAG\to cAG\to c$ as
module maps.
\end{proof}

\begin{proposition}[The compact cup operator]\label{prop:compact-cup-lift}
For $\delta=m_{c_1(\mathscr L_{\chi,p})}$ on the main coefficient there is
a degree-two operator $\widetilde\delta$ on the scalar compact coefficient
with $\widetilde\delta B_C=B_C(1\otimes\delta)$,
$A_C\widetilde\delta=(1\otimes\delta)A_C$, and
\begin{equation}\label{eq:compact-cup-commutator}
 [N,\widetilde\delta]=J_\chi:=\sum_wk_w\,j_p(\theta_p\smile a_{p,w})
 \pr_{p,w},\qquad k_w=\langle\chi,w^{-1}\nu\rangle,\qquad[N,J_\chi]=0,
\end{equation}
where $\pr_{p,w}$ is disk projection, $j_p$ the collar connecting map and
$\theta_p$ the forward-braid circle class.
\end{proposition}

\begin{proof}
Write the compact cone as $\mathsf K=(\mathsf O\oplus\mathsf P)\oplus
\mathsf T[-1]$ with $d(o,p,v)=(do,dp,\res o-ap-dv)$. The cup map is
diagonal on $\mathsf O\oplus\mathsf P$ with third component
$\delta_{\mathsf T}v+\sum_wh_wp_{p,w}$; \eqref{eq:collar-homotopy} is its
chain-map condition. At another puncture $p'$ the frame at $p$ is
unchanged data and no relative clutch occurs. Choosing the branch of
$\log\chi(w^{-1}t)^{-1}$ vanishing at $t=1$ along the scalar arc makes the
trivializations agree with the filling comparisons, which gives the
endpoint identities. Open and disk cup maps descend to the parameter
family, so the deck difference is the collar term $J_\chi$, which commutes
with packet monodromy. In this deck convention
$T_\gamma\widetilde\delta T_\gamma^{-1}-\widetilde\delta=J_\chi$; after a
common unipotent power $M$ the difference is $MJ_\chi$ and the next
adjoint difference vanishes, so the finite logarithm
gives~\eqref{eq:compact-cup-commutator}.
\end{proof}

\emph{Signs.} Let $h$ be the counterclockwise point push. The braid is
precomposition $\beta(\rho)=\rho\circ h_*$, inverse to geometric transport,
so $\theta_\beta=-\theta_{\rm ccw}$. For the restriction-minus-attachment
cone, $\tr_C(j_p\theta_{\rm ccw})=1$ (use the representative
$d\rho_O\wedge\theta_{\rm ccw}$ with $\rho_O$ increasing outward), hence
$\tr_C(j_p\theta_\beta)=-1$. At the cochain level, inversion is
$C_T\to C_{T^{-1}}$, $(a,b)\mapsto(a,-T^{-1}b)$ for
$d(a,b)=(da,(T-1)a-db)$, so it reverses the constant circle generator.
Forward auxiliary transport satisfies $P_{a+1}=P_aT_\gamma$, so deck pullback is $\Ad(T_\gamma)$. Commuting
$\delta$ through the scalar endpoints, the middle term is
$[N^2,\widetilde\delta]=2J_\chi N$; on the collar the trace of
$A_Cj_p\theta_\beta$ is $-a_P$, disk projection gives
$N^P([w]x)=[w](-L_{p,w^{-1}\nu}x)$, and composing with $B_C$ gives $b_P$.
The two minus signs cancel and the packet signs give $(-1)^b$:
\begin{equation}\label{eq:returned-boundary}
 [F_\nu^{\rm raw},\delta]=2(-1)^b\sum_{w\in\WG}\langle\chi,w^{-1}\nu
 \rangle\,L_{p,w^{-1}\nu}.
\end{equation}
For a simple root system of rank $r$ write $\chi(\xi_p)=a_\kappa
\kappa(\mu,\xi_p)$; the Weyl sum in~\eqref{eq:returned-boundary} is
$|\WG|a_\kappa\kappa(\nu,\nu)L_{p,\mu}/r$, so with $c_\nu$ and $F_B$ as in
Lemma~\ref{lem:whole-quadratic-return},
\begin{equation}\label{eq:normalized-boundary}
 [F_B,m_{a_\kappa\kappa(\mu,\xi_p)}]=-L_{p,\mu}.
\end{equation}
(The adjoint $\alpha$-root line has Chern class $\alpha(\xi_p)$ in this
convention, and roots span the derived rational characters.) For several
simple factors use $c_i$ and $F_B$ as in the lemma; orthogonality keeps
only the $i$th component of $\mu$.

\begin{proposition}[Physical equivariance]\label{prop:physical-equivariance}
$[F_\nu^{\rm raw},L_{p,\mu}]=0$. Consequently, for
$y_i=a_\kappa\kappa(\mu_i,\xi_{p_i})$, $\delta_i=m_{y_i}$ and
$L_i=L_{p_i,\mu_i}$,
\[
 [F_B,\delta_i]=-L_i,\qquad[F_B,L_i]=0,\qquad[[F_B,\delta_i],\delta_j]=0,
\]
and for every polynomial $f$,
$[F_B,m_{f(y)}]=-\sum_im_{(\partial_if)(y)}L_i$.
\end{proposition}

\begin{proof}
Build the same diagram over the joint ordered Cartan base of physical
parameters $S$ and auxiliary parameter $t$; the scalar endpoint is $t=1$
for every $S$. At the separating anchor all relevant physical and fused
Levi character products have nonzero logarithmic modulus and physical
roots have modulus $\ne1$; multiplying by the compact phase tori and by a
small regular radial translate of the auxiliary torus preserves these
finitely many margins (and avoids the old-node coincidences
$w'A_p=A_p$, $w'\ne1$, $A_p=\log|S_p|$). Physical and auxiliary loops are
then coordinate loops of a product, so their continuations commute; after
a common power $N$ is a morphism of the physical local system, and so is
$\tr_CA_CN^2B_C\mathrm{unit}_C$. Taking the finite logarithm of the
physical continuation gives the first claim. With
$[L_i,\delta_j]=[L_i,L_j]=0$ and~\eqref{eq:normalized-boundary} the rest
follows by the Leibniz rule.
\end{proof}

\subsection{The surface class on the compact diagram}
\label{subsec:surface-class}

The surface computation has the same gluing issue, but now the moving
graph $\Delta=\{x=q\}\subset C\times C$ has a nontrivial normal line. We
construct one comparison of capped bundles on the whole curve; its deck
derivative gives the auxiliary flag class, its second derivative the
tangent Euler class, and a radial Thom homotopy identifies the derivative
with evaluation at the moving point, including the collar homotopy.

\begin{proposition}[The surface comparison]\label{prop:compact-surface-comparison}
Fix a simple factor $i$ of rank $r_i$ and an allowed nonzero $\nu$
supported on it. Put $c_\nu=4(-1)^b|\WG|\kappa(\nu,\nu)/r_i$ and write
$q_{\kappa_i}(v)=\kappa_i(v,v)/2$ and let $e_{0,i}$ be its surface
slant over $C$. On the compact coefficient, the
capped surface class $e=e_0$ satisfies
\[
 T^ke=e+ky+\tfrac12k^2z,\qquad y=-\kappa(\nu,\xi_q),\qquad
 z=\kappa(\nu,\nu)c_1(T_C),
\]
the scalar maps are $E_0$-linear, and with the maps and trace
of~\eqref{eq:raw-return},
\[
 2A_Cm_y^2B_C=c_\nu\,m_{q_{\kappa_i}(\mathcal E|_{x=q})},\qquad
 [[F_\nu^{\rm raw},m_{e_0}],m_{e_0}]=c_\nu\,m_{e_{0,i}},
\]
as identities of the rational coefficient diagram.
\end{proposition}

\begin{proof}
Keep the spatial coordinate $x$ separate from the moving point $q$. A
\emph{cap} is the extension of the flat bundle over a puncture given by a
chosen logarithm and lattice. All torsors and classifying maps in this
proof are formed on the unrigidified stacks; only their rational
characteristic operations are descended, as in
\S\ref{subsec:scalar-flag}. The collar, scalar and normal-deformation
diagrams are pulled back there before this descent, so the same
comparison faces descend together.

\emph{1. Characteristic operations.} Enlarge the coefficient diagram by the
sources of its correspondences, their normal deformations and scalar
faces, and pull this diagram to the auxiliary logarithm cover. Let
$\mathscr R$ be the homotopy limit there of its rational characteristic
cochain algebras, using the universal classifying maps and their
connection resolutions. Auxiliary winding acts externally on this
diagram by deck pullback, and hence semilinearly on its modules. Pullback, Gysin and trace are $\mathscr
R$-modules by the projection formula, so the compact fiber is an
$\mathscr R$-module and even characteristic operations commute on it. This does not
assert that the compact coefficient is itself a ring. The scalar
endpoints will be modules for the particular surface class constructed
below; they are not assumed linear for arbitrary auxiliary characters.

\emph{2. The capped bundle cylinder.} In a character line of holonomy
$\exp(2\pi ia)$ a flat frame and its cap frame satisfy
$e(b+1)=e(b)\exp(2\pi ia)$ and $s=e\exp(-2\pi iaf)$ with
$f=(2\pi i)^{-1}\log z$; for $e=e_{\rm hol}H(z)z^{-A}$, $a=-A$, this gives
the ordinary lattice $s=e_{\rm hol}H(z)$. On the $w$-disk use the main
bundle of $\MB(S_pw^{-1}t)$ on the whole curve, with physical logarithm
$a_p+w^{-1}\ell$. On the collar, contract the triangular coordinate by
$n_u=\eta_+(u)n\eta_+(u)^{-1}$ ($\eta_+$ positive on positive roots), a
$T$-equivariant contraction commuting with conjugation by $C$; pulling back
gives a bundle comparison on the whole curve, and connections with given
exterior restriction form an affine space. At $n=1$ the split and fused
caps lie in one Cartan torsor: on an annulus around $p$ and $q$ put
$f=z_q/z_p$ (winding zero) and match $e_{\rm split}=s_{\rm split}
z_p^az_q^\ell$ with $e_{\rm fused}=s_{\rm fused}z_p^{a+\ell}$, so
$s_{\rm split}=s_{\rm fused}\exp(-\ell\log f)$, extended by a collar cutoff. Choose the logarithm with value zero at
$q=p$; the permitted cutoffs and their interpolations form a convex
space. The comparison keeps the physical-cap braid frame ratio $t$
and the auxiliary-cap ratio $S_p$, with inverse forward vector
multipliers. At collision only the commuting endpoint $n=1$ is filled;
the retained cylinder compares this endpoint with the noncommuting
family. The auxiliary deck changes $\ell$ to $\ell-k\nu$, with relative
factor $f^{k\nu}$ and fused factor $z_p^{k\nu}$, product $z_q^{k\nu}$:
this is the new cap frame. For the associated $\chi$-line the factor is the
lattice $\mathcal O(-k\langle\chi,\nu\rangle\Delta)$: the split and
fused factors are
$\mathcal O(-k\langle\chi,\nu\rangle(\Delta-P))$ and
$\mathcal O(-k\langle\chi,\nu\rangle P)$, respectively, where
$P=\{p\}\times D_q$. This is the sign in the quadratic formula below.

Pull this same evaluation torsor and cylinder to the deformation to
the normal cone of the scalar embedding. On the central normal bundle
they are pulled back from the source torsor. The scalar and triangular
normal restrictions are the original transverse square at $n=1,\ell=0$,
where the comparison fills the identity auxiliary meridian. Thus the
Thom transformation carries the same $AG=1$ face used in the packet.
The map $d_w$ reorders the same eigenvalue/logarithm pairs and preserves
the Lie logarithm of the unflagged $U$. The cap comparison is therefore
pulled back from the flag-forgetting diagram with this logarithm data.
Along the unshifted principal scalar arc its logarithm has the canonical
nilpotent value. On caps a Borel connection suffices, because $\kappa(\ft,\mathfrak n)=\kappa(\mathfrak n,
\mathfrak n)=0$. Applying the quadratic class and integrating over $x$
gives $e\in H^2\mathscr R$, with edge and triangle components
\begin{equation}\label{eq:HK}
 H=\kappa(\xi,\Lambda)+\tfrac12\kappa(\Lambda,d\Lambda),\qquad
 K_{012}=\tfrac12\kappa(h_1,h_2),\qquad dK_{012}=H_{01}+H_{12}-H_{02},
\end{equation}
where $\xi$ is the Cartan characteristic cocycle, $\Lambda=a\theta_p+
\ell\theta_q$ the connection change on an edge and $h_1,h_2$ successive
changes. The differential includes the parameter differentials
$da$ and $d\ell$. These are form components of the rational classifying
comparison, whose path-power and deck identities also supply the
composition and tetrahedron faces; logarithmic forms alone are not
being declared rational cochains.

\emph{3. The shifted fused cap and the Euler class.} On $C_x\times D_q$
put $D'=\Delta-P$. Since the fused baseline has already changed to
$\xi-k\nu P$,
\begin{gather*}
 q_\kappa(\xi-k\nu\Delta)-q_\kappa(\xi-k\nu P)=-k\kappa(\nu,\xi)D'
 +k^2\kappa(\nu,\nu)\bigl(PD'+\tfrac12D'^2\bigr),\\
 PD'+\tfrac12D'^2=\tfrac12(\Delta^2-P^2).
\end{gather*}
In the relative Thom complex of the collision disk, integration over
$x$ gives $\pi_*(PD')=[p]$ and $\pi_*(D'^2)=c_1(T_C)
-2[p]$, so the point terms cancel. For global gluing one must retain the full Borel logarithm, not only
its Cartan quotient. The translated cap calculations take place over
$\Omega^{\mathrm{aux,rs}}=\Omega\cap T^{\mathrm{rs}}$, where
$T^{\mathrm{rs}}=\{t:\alpha(t)\ne1\text{ for every root }\alpha\}$.
On its ordered logarithm cover write $U=u t u^{-1}$ with $u\in U_B$
and set $A=\Ad_u\ell$. Thus $\exp(2\pi iA)=U$, and the translated
branches $\ell-k\nu$ are defined on this regular domain.
For the scalar approach use instead the unshifted principal logarithm
near $\ell=0$, with its canonical nilpotent value at unipotent $U$.
Only this principal cap and its collar, normal and scalar faces extend
to the weak scalar fiber; no logarithm with prescribed nonzero integral
Cartan part is asserted there. On these respective domains use
$s_i=e_i\exp(-A_i\log z_i)$. Choose lifts of the angular changes with
$\tau_{ij}=(2\pi i)^{-1}\log u_{ij}$ and
$\tau_{ij}+\tau_{jk}-\tau_{ik}=n_{ijk}$.
Let $g_{ij}$ be the corresponding path gauges. Concatenation of the
lifted paths and point~3 of \S\ref{subsec:coefficient-enhancement} give
$g_{ij}g_{jk}=U_i^{n_{ijk}}g_{ik}$ and
$g_{ij}A_jg_{ij}^{-1}=A_i$. Therefore
$G_{ij}=\exp(-2\pi i\tau_{ij}A_i)g_{ij}$ satisfies
$G_{ij}G_{jk}=G_{ik}$. An integral change of angular lift changes the
two factors by inverse meridian powers. These transitions define the
capped $G$-torsor with its specified gauge comparisons. The Cartan transition of the cap is
$a_{ij}\exp(-\ell\log u_{ij})$ and the projected cancellation is
$t^{-n_{ijk}}t^{n_{ijk}}=1$. On the principal scalar face $A_i$ is
nilpotent and its full exponential cancellation is still required:
$t=1$ does not imply $U=1$.
After $\ell\mapsto\ell-k\nu$ the transitions change by $u_{ij}^{k\nu}$,
which is $N_\Delta^{-k\nu}$; hence $\xi_q'=\xi_q-k\nu c_1(N_\Delta)$ and
$T^ky=y+k\kappa(\nu,\nu)c_1(T_C)$, i.e.\ $z=\kappa(\nu,\nu)c_1(T_C)$.

\emph{4. Thom localization.} The relative line $\mathcal O(\Delta)$ with
its section has first Chern transformation the Thom class $u_\Delta$,
unique up to no positive-degree ambiguity since
$H^{-j}(\Delta,\Q)=0$ for $j>0$. With $r$ the radial retraction onto the
graph and $K_\xi$ the radial prism, $\xi-r^*i_\Delta^*\xi=dK_\xi$, so
\[
 \pi_*(\xi u_\Delta)\simeq i_\Delta^*\xi\quad\text{with homotopy }
 \pi_*(K_\xi u_\Delta),
\]
the projection formula with homotopy; similarly
$u_\Delta^2\simeq i_{\Delta!}c_1(N_\Delta)$, and the same for $P$. These normal and scalar restrictions retain the cap
Cartan quotient: at a regular cap it comes from the auxiliary flag, on
a dummy disk from the physical cap, and at the principal scalar face
from the Borel quotient of the nilpotent cap. No diagonalization of
the nilpotent connection is used. The identity-entry transverse scalar--triangular
square carries their common quotient and section. The specialized Weyl
entries on the proper direct image are compared in
Lemma~\ref{lem:scalar-collar-return}.

For the transverse pair $(\Delta,P)$ use the divisor functions
$z_x-z_q,z_x$, their product complex Thom orientation and the same
connection/tubular resolution. The quotient has divisor $\Delta-P$; the cap frame uses the negative
associated-line lattice specified above. Pulling the singular radial prisms to
this two-divisor normal diagram gives their pastings, with the
original scalar and triangular normal faces. Applying the quadratic
class gives precisely the $PD'$ polarization and the triangle term;
no new collar homotopy is chosen after localization.
The graph modification tensors the
$\chi$-line by $\mathcal O(-k\langle\chi,\nu\rangle\Delta)$, so the
quadratic change is $-k\kappa(\nu,\xi)u_\Delta+\frac12k^2\kappa(\nu,\nu)
u_\Delta^2$, and integrating identifies its linear term with the diagonal
flag class together with its radial, collar, normal and scalar
homotopies. Applying this to both $\Delta$ and $P$ identifies the
derivative of the integrated cap homotopy with the difference of their
evaluated flag homotopies. This is the required comparison between
spatial integration and diagonal evaluation. On the $w$-disk the linear term is $-\kappa(w^{-1}\nu,\xi_p)$
and the fixed graph has zero normal class; so $z$ may have a nonzero
collar component although its disk vertices vanish. The identity is extended from the regular anchor to the scalar
coefficient using only the unshifted principal cap. Its
conjugation-compatible logarithm extends to the canonical nilpotent
value, and the ordinary lattice and flag classes extend with it.
Let $r:\mathfrak K_*\xrightarrow{\sim}\mathfrak K_{\lim}$ be the
restriction equivalence of the retained diagram. The extending principal
cap identifies its operations with $rE_*r^{-1}$,
$rm_{y,*}r^{-1}$ and $rm_{z,*}r^{-1}$. Put $T_{\lim}=rT_*r^{-1}$.
Conjugating the regular-anchor identity gives
\[
 T_{\lim}^kE_{\lim}T_{\lim}^{-k}
 =E_{\lim}+k m_{y,\lim}+\tfrac12k^2m_{z,\lim}.
\]
This is transport of a rational cochain identity through the weak-family
restrictions, proper Springer specialization and finite cones.
No translated resonant lattice is extended at that endpoint.

\emph{5. From the deck formula to commutators.} Put $N=M^{-1}\log T^M$
(Lemma~\ref{lem:bounded-diagram-log}) and $E_0=m_e$. Semilinearity gives
$T^kE_0T^{-k}=m_{T^ke}$; comparing $T^{Mk}e=e+Mky+\frac12M^2k^2z$ with
$\exp(Mk\ad N)E_0$ for all $k$ gives $[N,E_0]=m_y$, $[N,m_y]=m_z$, and
$[[N,E_0],E_0]=[m_y,E_0]=0$ by commutativity of even operations.

\emph{6. The scalar return.} Lemma~\ref{lem:scalar-collar-return}
identifies the two scalar degree-four operations on the compact main
coefficient, using the attaching homotopy of $Y=[N,E_0]$ constructed
in steps~2--5. The scalar legs are linear for the unflagged principal
$E_0$ by the scalar normal pullback and projection formula. Hence the
algebra in step~5 gives
\[
 [[F_\nu^{\rm raw},m_{e_0}],m_{e_0}]
 =\tr_C\,2A_CY^2B_C\,\mathrm{unit}_C=c_\nu m_{e_{0,i}},
\]
since integration of the evaluated $q_{\kappa_i}$ over $q$ gives $e_{0,i}$.
\end{proof}

\begin{example}[Why the scalar comparison uses correspondences]
\label{ex:scalar-deck-flag}
For $G=\SL_2$ put
$U(t,u)=\left(\begin{smallmatrix}t&u\\0&t^{-1}\end{smallmatrix}\right)$
and $d=t-t^{-1}$. The old ordered flag is $\C e_1$; the other eigenline
for $d\ne0$ is $\langle(-u,d)\rangle$. Its graph closure has an exceptional
$\PP^1$ over $(t,u)=(1,0)$. On that fiber the new tautological line is
$\mathcal O(-1)$, whereas the old line is constant. Individual flag
lines therefore cannot be identified through deck specialization.
The invariant quadratic flag formula nevertheless survives:
$h^2=x^2$ in $H_T^*(\PP^1)$, with $h$ the tautological class and $x$
the torus character. The proof below compares the specified
cohomological operations and their collar homotopies.
\end{example}

\subsection{The scalar collar comparison}\label{subsec:scalar-collar-comparison}

\begin{lemma}\label{lem:scalar-collar-return}
With the scalar legs of Proposition~\ref{prop:betti-weight} and the
transported linear coefficient $Y=[N,E_0]$ from steps~2--5 above,
\[
 2A_CY^2B_C=c_\nu m_{q_{\kappa_i}(\mathcal E|_{x=q})}
\]
on the compact scalar main coefficient.
\end{lemma}

\begin{proof}
We compare the actual attaching homotopies in rational cochains. At an
open scalar vertex, $i^*Y$ is the auxiliary flag character on
$\Flag(\mathcal E|_q)$. The self-intersection formula gives
$\pi_*(c_b(T_\pi^*)f(\xi))=(-1)^b\sum_w f(w\xi)$: on the universal
flag over $BT$ the ratio of the cotangent and tangent Euler classes
at each fixed point is $(-1)^b$. The covectors $\kappa(w^{-1}\nu,-)$ vanish on all other
simple factors. On the $i$th Cartan the reflection representation is
irreducible. Taking the $\kappa_i$-trace gives
\[
 \sum_{w\in\WG}\kappa(w^{-1}\nu,-)^{\otimes2}
 =\frac{|\WG|\kappa(\nu,\nu)}{r_i}\,\kappa_i.
\]
This proves the vertex equality universally over $BG$. For example, on
$\SL_2\times\SL_2$ with the standard trace on each factor, winding by
the first coroot returns $32x_1^2$ in $H^4(BT)$, where $x_1$ is the
first diagonal character. The second winding returns $32x_2^2$. On each packet disk the same
Weyl sum gives the equality over $BT$. Choose these vertex homotopies
universally and pull them back. We compare their restrictions on a
single collar.

\emph{The local correspondence.}
Use the separating fine reference with $Z_G(S)=T$ from
\S\ref{subsec:triangular}, before finite central descent. Write $P=SU$
and $C=C_w(t)=Sw^{-1}t$. Near the auxiliary value $U=1$, the fused
element $P$ stays in the regular semisimple neighborhood of $S$ with
its local ordering. The full centralizer condition makes this ordering
a local branch of $T^{\rm rs}\to T^{\rm rs}/\WG$. The physical
logarithm $a_p$ consequently extends to a conjugation-compatible
logarithm of $P$ there. Write $P^a=\exp(a\log P)$, $a\in\mathbb R$.
On the principal auxiliary chart define $U^a$ by the full logarithm
of step~3; for unipotent $U$ this uses its nilpotent logarithm.
On the triangle $P=n C n^{-1}$ one has
$P^a=nC^a n^{-1}$ with
$C^a=\exp(2\pi ia(a_p+w^{-1}\ell))$.

The universal weak chart $[\widetilde G_t/T]$ and the main chart $BT$
map to $\mathcal B=[G/G]$ by $P$ and $C$. The triangle $[U_B/T]$
maps to both over $\mathcal B$ by the specified conjugating gauge $n$.
Pullback along the inverse framed word
$[\mathcal R/G]\to\mathcal B$ gives the moduli correspondence, with
$r'=n^{-1}rn$. Give the weak chart the identity on $(U,F)$ with
base gauge $P^a$ in $\mathcal B$. Its action after word pullback is
$(r,U,F)\mapsto(P^{-a}rP^a,U,F)$. On the triangle its two
presentations are
\[
 (r',n)\mapsto(C^{-a}r'C^a,n),\qquad
 (r',n)\mapsto(r',C^anC^{-a}).
\]
Their comparison is the gauge $C^a$. Its power law, the corresponding
$T$-gauge on the main $BT$, and the affine-bundle unit give the trivial
action on the packet source. Restriction to $a\in\mathbb Z$ gives the
specified forward braid, gauges and composition faces of
\S\ref{subsec:coefficient-enhancement}. Here $\mathbb R$ is the additive topological group with its ordinary
topology. These continuous actions are formed on rational coefficient
diagrams, rather than on mixed Hodge modules. The outgoing scalar map restricts
to $U=1$, the incoming map starts there, and the proper deck
correspondences preserve actual $U$. Their scalar pasting can therefore
be computed on this germ before the final projection.

\emph{The angular comparison of the linear term.}
In steps~2--4 the discrete deck difference is
$kY+k^2m_z/2$. Subtracting its normal Euler term, the radial graph
prism evaluates the auxiliary character at $x=q$ and the fused
character at $x=p$. To keep its actual path, set $p=0$,
$q=q_0\exp(2\pi ib)$, with a fixed disk tangent for $x-q$.
At $n=1$ the split developing frame in each torus character is
$x^{a_p}(x-q)^{w^{-1}\ell}$. Removing the auxiliary lattice factor
at $x=q$ leaves $q^{a_p}$; removing the physical factor at $x=0$
leaves $(-q)^{w^{-1}\ell}$. Their angular factors relative to $b=0$
are $S^b$ and $(w^{-1}t)^b$, respectively, with inverse factors for
forward vectors. Normalization at a fixed exterior point involves
$1-q/x_{\rm out}$, whose convergent logarithm for $|q|<|x_{\rm out}|$
gives its specified contraction. These charts use the same logarithms
and collar cutoff as the split-to-fused frame in step~2.

For a regular old Cartan value $c=w^{-1}t$ write
$U=v(U)c v(U)^{-1}$ in $B$, with $v(U)\in U_B$ unique. The reordered
flag is $v(U)wB$ with the corresponding label convention. Write
$\lambda_w=\kappa(w^{-1}\nu,-)$ for the character in the old Cartan,
clearing a denominator if necessary. In this frame the normalized
tangential action $U^{-a}$ contributes $\lambda_w(c^{-a})$
in multiplicative character notation. The gauge $C^a$ contributes
$\lambda_w(C^a)$ because
$v(C^aUC^{-a})=C^av(U)C^{-a}$. Their product is $\lambda_w(S^a)$
for every regular $n$; the contraction $n\mapsto n_u$ in step~2
commutes with $C^a$. Changing the logarithm lift by $-k\nu$ changes
both fractional powers in $C^ac^{-a}$ by inverse factors. In contrast
the physical ratio $C^a S^{-a}$ retains its integer winding, as in
Proposition~\ref{prop:collar-defect}.

If the open tangent frame has winding $n_p$ relative to this disk
tangent, write its unit transition as
$u(0)\exp(2\pi i n_p b+h(b)-h(0))$, with $h$ periodic. The cap
transition $\exp(-\ell\log u)$ and the boundary gauge $U^{n_p}$
of step~3 give the corresponding character power along $b$, with its
inverse for the opposite frame. The periodic logarithm has its
straight-line contraction. Use the principal $U^{a n_p}$ for the
real action. Changing the auxiliary logarithm lift in this factor gives
the integer tangent rotation in
$\xi'_q=\xi_q-k\nu c_1(T_C)$, precisely the normal term already
retained in the deck calculation. On the unshifted scalar incidence
$U=1$ this factor is the identity for all real $a$.

For a relative multiplier $\exp(aL)$ these formulas give the section
$\exp(bL)$ in the translation model $E\mathbb R=\mathbb R$, with
the inverse section for the inverse multiplier. Restriction along
$E\mathbb Z=\mathbb R\to E\mathbb R$, $b\mapsto b$, gives the
calculated angular path, including its logarithm. Thus the chosen
regular attaching homotopy of the evaluated linear term extends to
the real action on the principal scalar approach. Its characteristic
square is taken with the same section and polarized triangle. This
construction uses the discrete deck calculation to identify $Y$;
it requires a real extension only of this evaluated operation.

\emph{Specialization of the specified map.}
Keep $U$ and $t$ in the ordered base $X$ of
Lemma~\ref{lem:scalar-coefficient-extension}. Let
$K=R\widetilde\mu_*\Q_{\widetilde G}$ and
$D=R(w\widetilde\mu i_L)_*\Q_L$, with the old base identified with
the new one by $w$. The Gysin and deck composite is a map
$f:D\to K[2b]$ (with the usual twist). The regular homotopy just
constructed compares the maps $m_\lambda^2 f$ and
$f[4]m_{p^*\lambda_w}^2$ from $D$ to $K[2b+4]$. Here $\lambda$
is the auxiliary characteristic class in the new flag and
$m_{p^*\lambda_w}$ denotes cup on the triangular source. The affine
unit will identify its pushforward with the disk operation.

Apply rational topological nearby cycles along the chosen final
regular-to-scalar trait to this homotopy before projecting $X$ to
$\mathcal B$. Both $\widetilde\mu$ and
$\widetilde\mu i_L$ are proper, and the spaces $\widetilde G$ and
$L$ are smooth over the ordered parameter. Their constant-sheaf
specializations are therefore isomorphisms, also after these proper
pushforwards. The regular embedding specializes with its oriented
normal Thom map. On the real-action nerve, product with every
$\mathbb R^j$ commutes with nearby cycles by the local product, and
the action maps are homeomorphisms with the parameters retained.
Thus the whole path specializes, together with its power faces.
The scalar target cup is the proper push of the ordinary tautological
operation, and the scalar deck is the proper Springer action. This
specializes sheaf maps, including where $v(U)$ has no limit.

Now project to $\mathcal B$. On the triangle the projection factors
as $i_Cp$ using $n$, so its source is
$i_{C,*}Rp_*\Q=i_{C,*}\Q_{BT}$. Smooth word pullback near the
scalar support commutes with these operations and their normal
orientations by Lemma~\ref{lem:triangular-source}. To identify the
resulting path with the original transport after global cohomology,
use the natural exchange
$\psi Rg_*M_\eta\to Rg_{0*}\psi M_\eta$ for the remaining
projection $g$. Its composite with specialization from
$i^*Rg_*M$ is ordinary fiber base change followed by specialization
on $M$. On the actual weak and fused families the global specialization
and fiber base change are isomorphisms by
Theorem~\ref{thm:weak-betti}; its local product of compactification
pairs also gives constant-sheaf specialization on the total spaces.
Use the affine-coordinate product on the triangle and proper
pushforward for $D,K$. These natural comparisons consequently identify
the transport of the \emph{regular path} with its pushed specialization.
This argument uses the global weak and fused families for fiber base
change and restricts the resulting scalar pasting on its support
before projection.

\emph{The scalar span.}
The scalar embedding, flag trace, affine unit and the transverse
$AG=1$ diagram are equivariant for the real actions above.
The universal degree-zero Springer relations restrict with their
composition faces, as in point~5 of
\S\ref{subsec:coefficient-enhancement}. Their integer restrictions
therefore give the actual $B_C$ and $A_C$ attaching faces. The
unflagged scalar evaluation uses the physical frame $q^{a_p}$, so the
restriction of the universal $BG$ Euler homotopy and the disk $BT$
Weyl homotopy also have this real action.

After composing these specified faces, the two degree-four collar
paths, with the chosen vertex homotopies, form a loop of maps of the
scalar main coefficient. Universally that coefficient is
$i_{S,*}\Q_{BT}$ for the closed inclusion
$i_S:[\mathcal O_S/G]=BT\hookrightarrow[G/G]$.
The base gauge $S^a$ is trivialized on this object by the same $T$-gauge.
By full faithfulness of closed pushforward the group of such loops
before restriction to the integer braids is
$H^3(BT\times B\mathbb R,\Q)=0$, since $\mathbb R$ is contractible.
This also computes maps from the constant source into the homotopy-fixed
target by adjunction. Choose a filling and restrict it to
$\mathbb Z$. Word pullback and cohomology carry it to a filling of the
specified scalar collar pastings. The open-vertex homotopy came from
the single flag Euler identity over $BG$, so these fillings apply
simultaneously to its restrictions at all punctures. The compact
coefficient is a span with one restriction and one attachment at
each collar. The vertex homotopies and these fillings give equality
on its homotopy fiber.
\end{proof}


We have constructed an endomorphism $F_B$ of $H^*(\MB)$ of weight degree
at most $-4$ satisfying~\eqref{eq:whole-commutator}. It remains to prove
that $F_B$, transported by $\Psi$, has perverse degree at most $-2$.

\section{Correspondences and the perverse bound}\label{sec:perverse}\label{app:geometry}

We realize the operator~\eqref{eq:raw-return} on the derived direct image
of the Hitchin map, before taking cohomology. A correspondence whose two
endpoints have the same Hitchin value induces a morphism over the Hitchin
base; integrating the moving point over the compact curve then lowers the
perverse index by two (the trace over $C$ has degree $-2$). The
difficulty is to pass from flat connections to Higgs bundles: the
correspondences are not compact, and their coefficients include gluing
homotopies. We apply the two scalar operations first, place the result on
a correspondence proper over its output, and use an annular estimate to
prove properness through the Higgs endpoint.

The scalar-composition lemmas come first, because the puncture
coefficient must be specialized before ordinary completion in the Hodge
parameter. We then construct the lattice changes, collar fusion and
puncture carrier, and finish with the annular estimate and proper return.

\subsection{Set-up and inputs}\label{subsec:perverse-setup}

The construction applies to reductive $G$, but it is used on the
semisimple factors; the central torus is added in \S\ref{sec:proof}. We
work on a smooth stable fine reference after central rigidification;
finite quotient charts are allowed throughout. The physical parameters
are ordered and separating: character equations of any proper invariant
parabolic are excluded, multiplicatively along the positive parameter
paths and additively at the Higgs endpoint, and all paths
satisfy~\eqref{eq:fuchs}. Let $M_\lambda$, $0\le\lambda\le1$, be the main
Hodge family, whose harmonic product identifies its endpoints by $\Psi$
(Theorem~\ref{thm:harmonic-transport}). At the regular references used
here the residue at a chosen puncture is $a_\lambda\in\ft$ with
\[
 |\operatorname{Im}\alpha(a_\lambda)|\ge\delta>0\qquad(\alpha\in\Phi).
\]

We use the terminology of \S\ref{sec:lowering}: physical and auxiliary
points, caps, collars, carriers and prisms.

The proof uses the family comparisons of
Section~\ref{sec:comparison-families}, the coefficient identities of
\S\ref{subsec:coefficient-enhancement} and
Section~\ref{sec:cup-comparisons}, and the geometric constructions below.
These constructions do not depend on the perverse bound.
\begin{enumerate}[label=(C\arabic*)]
\item\label{C:master} (Propositions~\ref{prop:pg-master}
and~\ref{prop:pg-transport}.) The main and the moving auxiliary fine
masters are smooth over their ordered residue parameters (and over
$q\in\Cc$) and admit smooth proper orbifold compactification pairs
relative to them. Boundary-tangent horizontal lifts give complete
transports along the finite path diagram; a residue path of velocity
$\lambda\rho'(s)$ gives transport equal to the identity at $\lambda=0$,
with displacement $O(\lambda)$ on compact charts.
\item\label{C:graphs} (Proposition~\ref{prop:pg-closing}.) Each integral
closing map is a logarithmic principal lattice change; its graph closes
in a bounded principal Hecke carrier, proper over both endpoints. At
$\lambda=0$ the two Higgs fields are meromorphically isomorphic.
\item\label{C:packet} (Proposition~\ref{prop:finite-packet},
\S\ref{subsec:coefficient-enhancement} and
Section~\ref{sec:cup-comparisons}.) The positive-parameter coefficient diagram is
built from closed triangular embeddings, proper Springer diagrams and
scalar incidences; its scalar, winding, character and two-divisor Thom
identities are given as maps and prisms before integrating out the main
variables.
\item\label{C:fusion} (Proposition~\ref{prop:pg-fusion} and
Lemma~\ref{lem:pg-germ}.) On a physical collar, the triangular-source
main object is the receiving flat object outside the cap, extended with
the ordered nonresonant fused logarithm $a_\lambda+\lambda\rho(s)$, with
$\rho$ bounded. The closed triangular map is proper over the receiving
datum, and the attaching diagram is computed by any sufficiently small
collar.
\end{enumerate}

\subsection{Applying the scalar operations first}\label{subsec:scalar-first}

Let $i:E\hookrightarrow Y$ be the scalar incidence and $\pi:E\to M$ its
flag projection, with $i$ a complex-oriented regular embedding of
codimension $b$ and $\pi$ smooth proper of relative dimension $b$, so that
$B=i_*\pi^*$ and $A=\pi_*i^*$.

\begin{lemma}[Ordered composition with the scalar maps]\label{lem:rc-caps}
Let $Y\xleftarrow{c_1}\Gamma\xrightarrow{c_2}Y$ be a correspondence with
$c_2$ proper and coefficient $c_1^*\Q_Y\to F\to c_2^!\Q_Y[d]$. Put
$\Gamma_1=\Gamma\times_{c_1,Y}E\xrightarrow{e_1}\Gamma$ and
$Z=\Gamma_1\times_{c_2,Y}E\xrightarrow{e_2}\Gamma_1$, with projections
$r_1,r_2:Z\to E$. Then $A\circ u\circ B$ (where $u$ is the operation of
the correspondence) is the operation of
\[
 M\xleftarrow{\pi r_1}Z\xrightarrow{\pi r_2}M
 \qquad\text{with coefficient}\qquad G=e_2^*e_1^!F[2b](b),
\]
of total degree $d$. The construction is natural for coefficient
morphisms and homotopies, and retains excess-intersection terms.
\end{lemma}

\begin{proof}
Proper base change for $i$ identifies $c_1^*Ri_*i^!\Q_Y$ with
$Re_{1*}a_1^*i^!\Q_Y$, $a_1:\Gamma_1\to E$; composing with the counit and
the source map and using adjunction gives $a_1^*i^!\Q_Y\to e_1^!F$, and
purity $i^!\Q_Y=\Q_E[-2b](-b)$ followed by $e_2^*$ gives the source map
to $G$. For the output use the counit $Re_{1*}e_1^!F\to F$ and the output
trace; apply $i^*$ and proper base change to get
$Rr_{2*}e_2^*e_1^!F\to\Q_E[d]$, then $R\pi_*$ and the flag trace. The
triangle identities identify the result with the composite.
\end{proof}

For the identity correspondence $G=\Q_E$ but the output is the
self-intersection map. In the scalar incidence used here,
$N_{E/Y}=T_\pi^*$, giving $AB=(-1)^b\chi(G/B)$; a constant middle
coefficient does not make the operation the identity.

\begin{lemma}[Completion of a proper return]\label{lem:rc-ray}
Let $I=[0,\epsilon]$ and choose product trivializations
$X_I\simeq X\times I$, $Y_I\simeq Y\times I$, together with coefficient
trivializations
\[
 K_I\simeq\operatorname{pr}_X^*K,\qquad
 L_I\simeq\operatorname{pr}_Y^*L.
\]
Here the fibers are locally compact triangulable spaces of finite type
and $K,L$ are bounded constructible complexes. Let
$X_I\xleftarrow{c_1}Z_I\xrightarrow{c_2}Y_I$ be a correspondence over
$I$, with $c_2$ proper, and let $j_Z$ and $i_Z$ denote its positive and
central inclusions. Given
$c_{1+}^*K_+\to F_+\to c_{2+}^!L_+[d]$, put $F_I=Rj_{Z*}F_+$.
Then $F_I$ defines a correspondence over $I$, and its central operation,
with coefficient $i_Z^*F_I$, is the positive operation transported by
the specified products. The assertion is functorial for finite diagrams
of coefficient maps and their homotopies. The middle complexes may be
computed among all rational sheaves.
\end{lemma}

\begin{proof}
The source arrow is the adjunction map
$c_1^*K_I\to Rj_{Z*}c_{1+}^*K_+$ followed by the given source arrow.
For the output, composition of direct images gives
\[
 Rc_{2*}F_I
 =Rj_{Y*}Rc_{2+*}F_+
 \longrightarrow Rj_{Y*}L_+[d]
 \simeq L_I[d].
\]
The last isomorphism uses the coefficient product: locally in $Y$ the
new factor is a contractible interval punctured only at its boundary
endpoint. Since $c_2$ is proper, adjunction gives the required arrow
$F_I\to c_2^!L_I[d]$, and proper base change identifies its fiber
restriction with the output arrow of the restricted correspondence.

The product trivializations also identify the two endpoint direct
images to $I$, with their fiber base-change maps, with
$R\Gamma(X,K)_I$ and $R\Gamma(Y,L)_I$. The operation is a morphism
between these constant complexes on the contractible interval, and
therefore has the same value on every fiber under the specified
identifications. This proves the transport assertion. All the units,
adjunctions and proper base-change maps used here are natural, so they
also extend the arrows and homotopies of a finite diagram.
\end{proof}

In particular the scalar operations are applied \emph{before} the
completion: the central coefficient is
$i_0^*Rj_*\bigl(e_{2,+}^*e_{1,+}^!F_+[2b](b)\bigr)$, and no interchange of
scalar restriction with $\lambda\to0$ is used.

\paragraph{Conventions for the geometric constructions.}
Throughout, cocharacters are in $X_*(T)$ for the
given group $G$, the fine rational weights are those of
Lemma~\ref{lem:rational-chart} and \S\ref{subsec:A-transport} (outside the
walls, with inertia centralizer $T$), all character parabolic degrees
vanish, and the affine residues satisfy~\eqref{eq:fuchs}. The rational
root presentation gives a finite ramified logarithmic lift of each
parahoric chart, with its inertia action and allowed lattice changes; all
constructions are made on the lift and descended. Statements about
smoothness concern the rigidified stacks, which have finite stabilizers;
their coarse spaces may have finite quotient singularities.

\subsection{Integral lattice modifications}\label{subsec:B-modifications}

An integral shift of a logarithm is realized by a change of lattice at
the marked point. A lattice is a point of the affine Grassmannian
$G((z))/G[[z]]$, and a bound on its relative position puts it in a
projective Schubert variety; finite-jet diagonalization defines the
modification where it is nonresonant, and the position bound provides a
proper graph closure where it is not.

\begin{lemma}[Finite-jet lattice change]\label{lem:pg-jet}
Let an ordered logarithmic principal $\lambda$-connection have regular
Cartan raw residue $s$, and let $\nu\in X_*(T)$. Where the finitely many
denominators $\alpha(s)+k\lambda$ below are nonzero, there is a
modification with lattice $g(z)z^{-\nu}G[[z]]$ and new raw residue
$s-\lambda\nu$. The modification depends algebraically on finitely many jets, has an
inverse on the opposite open set, and has relative position bounded by
the dominant representative of $-\nu$. If $s_0$ is regular, it extends
through $\lambda=0$ and preserves the Hitchin section there.
\end{lemma}

\begin{proof}
In an ordered frame, successive root-valued gauges remove root
coefficients in positive powers of $z$; at order $k$ the operator is
$\ad(s)+k\lambda$, with eigenvalues the stated denominators. It suffices
to diagonalize up to order $\max_\alpha|\langle\alpha,\nu\rangle|$: after
conjugation by $z^\nu$ the remaining root terms stay logarithmic and
higher congruence gauges stay holomorphic. Toral coefficients commute with
$z^{-\nu}$, and the derivative of $z^{-\nu}$ is the residue shift
$-\lambda\nu$. Reversing gives the inverse, and multiplication by the
holomorphic $g$ does not change the relative position. At $\lambda=0$
with $s_0$ regular the denominators are units, and the two Higgs fields
agree through the meromorphic identification, so their invariant
polynomials agree. On a ramified chart the recursion is
inertia-equivariant by uniqueness.
\end{proof}

\begin{proposition}[Proper closing graphs]\label{prop:pg-closing}
Let the auxiliary loop be $t_*\exp(-2\pi is\nu)$, so its lifted raw
logarithm ends at $\zeta_*+\nu$. For every power $k$ used below, the
closing map is $j_{-k\nu}$; its positive graph has a closure in a fixed
bounded principal Hecke carrier, proper over both auxiliary endpoints
through $\lambda=0$ and relative to $q\in\Cc$, and its central endpoints
have the same meromorphic Hitchin sections.
\end{proposition}

\begin{proof}
Apply Lemma~\ref{lem:pg-jet} with $s=\lambda\zeta_*$ and cocharacter
$-k\nu$; the denominators $\lambda(\alpha(\zeta_*)+m)$ are nonzero on the
nonresonant anchor for $\lambda>0$ (no claim is made at $\lambda=0$). The
relative position is fixed; enlarge the bound to include the inverse.
Bounded modifications at a moving smooth point form a projective relative
Schubert variety, with a universal identification away from the point.
Logarithmicity of the transported connection, the prescribed residue and
the ordered flag are closed conditions; close the graph in this bounded
incidence. At $\lambda=0$ any new auxiliary residue has zero Cartan part,
hence is nilpotent and contributes zero to every character of an
invariant parabolic, which the separating physical characters exclude; so
the limits are stable. Reversing the modification gives properness over
the other endpoint, and at $\lambda=0$ the Higgs fields agree off $q$.
\end{proof}

For transports $f_a:Y^0\to Y^a$, $f_b:Y^0\to Y^b$, the correspondence
with input $f_b^{-1}\pr_b$ and output $f_a^{-1}\pr_a$ acts by
\begin{equation}\label{eq:pg-loop}
 f_a^*\,j_{-k\nu}^*\,(f_b^{-1})^*=T_\nu^k.
\end{equation}
At a regular physical point the modification itself extends through
$\lambda=0$; with $f_\lambda$ the transport to raw residue
$a_\lambda+\lambda\nu$ put
\begin{equation}\label{eq:pg-physical}
 q_{\nu,\lambda}=f_\lambda^{-1}j_{-\nu,\lambda},\qquad
 q_{\nu,0}=t_{p,-\nu},
\end{equation}
whose pullback is the physical monodromy; the signs follow from
contravariance and the negative exponential in the loop.

\subsection{Fusion along a puncture collar}\label{subsec:B-fusion}

Put $b=\dim G/B$, choose the physical regular $S\in T$ with $Z_G(S)=T$,
and write the remaining generators as $r$ with word $R(r)$, so that the
receiving relation is $R(r)SU=1$. For an ordered auxiliary value $t$, the
fixed-flag triangular source has flag $vB$ and $U\in v(t)N_v$, where $N_v$
is the unipotent radical of $vBv^{-1}$; put $C_v=Sv(t)$ and keep the open
set where each $C_v$ is regular and the full tuple is generic.

\begin{proposition}[The closed triangular source]\label{prop:pg-triangular}
On this open set the triangular source is
$L_v(t)=[R^{-1}(C_v^{-1})\times N_v/T]$; its main projection is an
$\A^b$-bundle and its receiving map is a closed regular embedding of
codimension $b$. The intersection of $L_v(1)$ with the scalar incidence $U=1$ at
$t=1$ is transverse and has the identity as main map.
\end{proposition}

\begin{proof}
For regular $c\in T$, $n\mapsto c^{-1}ncn^{-1}$ is an algebraic
automorphism of $N_v$, solved recursively in the transported height
(the linear coefficient in a root coordinate is $\alpha(c)^{-1}-1$). So
$SU\in C_vN_v$ is uniquely $nC_vn^{-1}$, and $r'=n^{-1}rn$ turns the
relation into $R(r')=C_v^{-1}$. The condition $F=vB$ is closed, and the
word map (with central directions removed by the character equations) is
a submersion along the incidence, by the word-differential calculation
in Lemma~\ref{lem:triangular-source}, after removing the central
character equations, and genericity; so pulling back the point
$vB$ gives a regular embedding of codimension $b$. At $U=1$ we have $n=1$,
and the zero-section/cotangent-fiber chart of the Grothendieck--Springer
resolution gives transversality and the complex orientation.
\end{proof}

\begin{lemma}[From the fixed flag to deck-labeled sources]
\label{lem:pg-deck-transfer}
The fusion, carrier and completion constructions for the fixed flag $B$
realize every packet map $\Phi_w(t)=d_wG_1(w^{-1}t)$ of
Proposition~\ref{prop:finite-packet} ($G_1(s)$ the fixed-$B$ map, with
main source of monodromy $Ss$), with their scalar specialization,
coefficients, composition homotopies and orientations; the completed
central carriers lie over the Hitchin diagonal.
\end{lemma}

\begin{proof}
First work over the regular ordered auxiliary Cartan open. Substitute
$s=w^{-1}t$ in the fixed-$B$ construction. The deck map $d_w$ reorders
the same eigenvalue/logarithm pairs of the auxiliary group element.
The graph of $d_w$ gives the degree-zero coefficient isomorphism, equivalently
inverse pullback along $d_w$. Composing that graph with the source
gives exactly $d_wG_1(w^{-1}t)$, with physical value $Sw^{-1}(t)$.
This composite is not the source with fixed flag $wB$ at parameter $t$.

To specialize, use the proper universal Grothendieck--Springer diagram.
The closure of the regular deck graph in the product of its two flag
resolutions over the unflagged space has proper projections, because
the flag maps are projective. Retain the original graph coefficient
and its proper specialization, rather than replace the closed graph
by a constant coefficient. Concretely, along the chosen punctured
scalar arc, push the regular graph and its coefficient into
this proper closure, then apply nearby cycles together with its
input and output arrows. Proper nearby-cycle compatibility identifies
the resulting arrow with the specialization of $d_w$.

Here is the comparison with the restricted Springer action. Write
$\mu:\widetilde Z\to Z$ for the universal Grothendieck resolution
and $u_w$ for its extended deck operator on $R\mu_*\Q$.
For this comparison on the analytic scalar arc, take rational realization
of the algebraic deck diagram and use ordinary nearby cycles of rational
sheaves. Their exchange with proper maps follows by proper base change
after pulling back the punctured disk to its universal cover.
Pull back by the invariant scalar arc $\chi(t(\tau))$.
The total pullback has closed ordered branches $\widetilde Z_v$
with Cartan value $v^{-1}t(\tau)$. If $\pi$ and $\pi_v$ are
their proper maps to the unflagged arc family, put
\[
 K_\Delta=R\pi_*\Q,\qquad K_v=R\pi_{v*}\Q,
              \qquad r_v:K_\Delta\longrightarrow K_v.
\]
The arrow is restriction to a closed branch. At $\tau=0$ every
branch is the same topological Springer resolution. Proper base
change therefore identifies $i^*K_\Delta=i^*K_v=\mathcal S$
and $i^*r_v=1_{\mathcal S}$. Scheme multiplicities do not change
these rational constant sheaves. Each ordered branch is smooth over
the arc, so smooth local acyclicity upstairs and proper nearby-cycle
compatibility give
$\mathrm{sp}_v:\mathcal S\xrightarrow{\sim}\psi K_v$.
On the punctured arc $r_v$ is projection to the $v$ branch. If
$D_w:K_w\to K_1$ is the regular deck comparison, the deck
convention gives $r_1u_w=D_wr_w$ there. Naturality yields
\[
 \mathrm{sp}_1(i^*u_w)=\psi(D_w)\mathrm{sp}_w,
 \qquad
 \mathrm{sp}_1^{-1}\psi(D_w)\mathrm{sp}_w=i^*u_w.
\]
The right side is precisely the restricted Springer action
\cite[Sections~2.2 and~2.6--2.7]{BM83}, with its original sign.
This description uses a proper correspondence at the singular scalar
fiber, not a deck automorphism of that fiber.

These graphs change only the auxiliary flag and its ordering. The
unflagged bundle, its exterior identification, the physical point and
the ordered node gluing are unchanged. Composing a proper receiving
carrier with the graph remains proper over the receiving datum. Its
main restriction map is the original one composed with projection.
Consequently the annular estimates apply to exactly the same exterior
bundle; the extra flag coordinates lie in proper flag varieties.
All central endpoint pairs still have the same meromorphic Higgs field
off the cap, hence the same invariant polynomials on the whole curve.
This proves the Hitchin-diagonal assertion also on the added strata.

At a fixed node, the bubble-trivial test concerns the underlying
receiving bundle and is unchanged by these graphs. Thus the identity
$pe=\pi r$ on the full scalar inverse image is preserved. The full
nearby coefficient of the composed horizontal source is kept; no
constant-coefficient replacement on its compactifying carrier is made.
Scalar pullback and flag trace are taken at this stage. Only afterward
do we apply ordinary completion in the Hodge parameter, in the order
specified in \S\ref{subsec:scalar-first}.

For a retained finite composable word, close its joint regular graph
with every intermediate flag in the corresponding finite flag fiber
product. Forgetting an intermediate flag is proper. Carry the full
coefficient and its original face arrows by these proper projections
before applying nearby cycles. Thus the same proper composition and
base-change comparisons retain the composition homotopies; no
interchange with a nonproper pullback or an arbitrary graph convolution
is used.

Apply this construction to the original finite graph diagram,
including the character comparisons and the scalar two-face.
The scalar reductions of the resulting diagram remain the
universal Springer relation and the same transverse $AG=1$ square,
with their common composition homotopy. In particular
$A\Phi_w=(Aw)G_1=\operatorname{sgn}(w)AG_1=\operatorname{sgn}(w)$;
only the fixed-$B$ square itself has $AG_1=1$.
A winding by $\nu$ in $t$
is a winding by $w^{-1}\nu$ in $s$; thus the physical endpoint is
the translation $t_{p,-w^{-1}\nu}$. This also retains the original
logarithm convention.
\end{proof}

\begin{proposition}[Exterior fusion]\label{prop:pg-fusion}
On a physical collar and a bounded separated logarithm tube, the main
flat object of the positive triangular source is the receiving flat object
outside the cap, extended across the artificial node by its ordered fused
raw logarithm $a_\lambda+\lambda\rho(s)$. On compact positive prisms the
closed triangular map is proper over the receiving datum (including the
maps of the finite transport word).
\end{proposition}

\begin{proof}
Cut off the two-puncture cap along tangential paths: with main relation
$RC=1$ and cap relation $MUV=1$, the opposite node orientations give
$C=V^{-1}$, and elimination gives $RMU=1$, the same amalgamated groupoid
(including stabilizers and topological sectors). For $\lambda>0$ apply
enriched Riemann--Hilbert \cite[Theorems~A and~D]{Boa11}: the paths keep
the regular imaginary root margins, so the fused logarithms are
nonresonant and the ordered extension at the node is unique; gluing and
restriction are inverse, analytically in the parameters. A horizontal
Borel on the cap gives the triangular object of
Proposition~\ref{prop:pg-triangular}; conversely the toral residues sum
to zero, so the cap torus bundle has degree zero, every root-line quotient
has degree zero on $\PP^1$, and $H^1(\PP^1,\mathcal O)=0$ makes the Borel
bundle constant. The closure satisfies the closed equation
$\chi(MU)=\chi(\exp[-2\pi i(a_\lambda/\lambda+\rho(s))])$ for invariant
functions $\chi$, whose right side is regular semisimple, so no support
escapes the logarithm chart; the fixed-flag source map is closed, hence
proper over compact prisms.
\end{proof}

\begin{lemma}[Small collars]\label{lem:pg-germ}
On each compact receiving chart the puncture attachment of the Betti
coefficient diagram (with its restriction, character and two-divisor Thom
prisms) is computed by any sufficiently small collar, and two such
collars are compared by the radial prisms; these can be completed
separately after shortening the Hodge interval.
\end{lemma}

\begin{proof}
The attaching object is the restriction germ of the punctured-disk
diagram; removing a smaller closed disk changes it by the radial annulus
cylinder, and nested collars are cofinal. By
Proposition~\ref{prop:pg-fusion} these cylinders are the same maps of the
coefficient diagram, including the character line and its divisor
section. The radial prism between two radii has compact spatial parameter
disjoint from $D$, so the bounded-return construction applies after
shortening the interval, and finitely many faces use a common interval.
Over a relatively compact Hitchin chart, properness allows finitely many
charts and one collar. No joint limit of radius and $\lambda$ is taken.
\end{proof}

\subsection{The coefficient at a puncture}\label{subsec:B-node}

We construct the puncture vertex in a nodal receiving neighborhood
containing its entire triangular source and scalar incidence. This
suffices for the scalar-capped coefficient: a supported Gysin map
followed by restriction to the scalar incidence can be computed in
such a neighborhood. We do not identify the cohomology of the whole
receiving space with that of this open set.

\emph{Plan of the construction.} (1)~The \emph{receiving family}
$Y$ of~\eqref{eq:pg-nodal-receiving} glues a main object, at a node, to a
bubble $\PP^1$ carrying the physical residue $A$ and an auxiliary
residue $B$ with a flag $F$. (2)~Inside $Y$, the \emph{horizontal
source} $H$ (on which the flag is forced to be the flag of $A$) and the
\emph{scalar incidence} $E$ (on which $B=0$) meet transversely along the
main space (Lemma~\ref{lem:pg-nodal-carrier}); this is the algebraic
form of the identity $A\circ G(1)=\id$ of~\eqref{eq:triangle-scalar}.
(3)~For $\lambda>0$ the family $Y$ is an open neighborhood of the
relevant part of the Betti collar diagram, so the scalar-capped
coefficient may be computed in $Y$ (Lemma~\ref{lem:pg-nodal-excision}).
(4)~Deck-labeled carriers realize the packet maps $\Phi_w$ properly, and
completion in $\lambda$ gives a correspondence over the diagonal of the
Hitchin base (Proposition~\ref{prop:pg-fixed-return}).

\paragraph{The receiving functor.}
Fix an ordinary extension on each ordered fine parahoric chart as
follows. In the apartment modulo the central directions, the affine
Weyl group takes the standard alcove to the chosen alcove. Use the
image of the origin under this element as the hyperspecial vertex.
Its translation lies in the coroot lattice, hence in $X_*(T)$ for
every global form of $G$. Lifting the translation and finite Weyl
element to the loop group identifies this Iwahori with the inverse
image of a Borel under evaluation at zero
\cite[\S5]{Boa11}. Its torsor is thus an
ordinary principal bundle with a full flag. Hold this extension fixed. Ordinary frame changes of this extension preserve
the construction below; a different lattice need not preserve bubble
triviality and will be compared by the capped excision argument instead.
All root inertia, central lifting and topological sectors are retained,
and the raw-to-affine residue correction is always $\lambda q$.

Attach a bubble $\PP^1$ to $C$ at $p$, with physical, auxiliary and node
points $0,1,\infty$. Bundles on this nodal curve are pairs of bundles
on its normalization with identified node fibers. This follows locally
from the node ring $\C[[x]]\times_\C\C[[y]]$; gluing locally free
modules and their defining tensors gives the assertion for principal
$G$-bundles. Connections valued in the dualizing differential glue
precisely when their node residues are opposite. Thus the gluing of
connections is over the residue stack $[\fg/G]$, with the sign reversed
on one branch. On a root presentation the node bands are balanced and
their inertia homomorphisms are inverse.

Write $a$ for the regular physical raw residue, let $\rho=t\rho_*$
be the small scalar arc, with $\rho_*$ regular in the root directions,
and put $c=a+\lambda\rho$. The parameters
satisfy the separating and character equations of
Section~\ref{sec:perverse}. For reductive $G$, all residue evaluations
below are over the compatible central affine slice: with the other
residues fixed, the central component is prescribed by the Fuchs
equations, and $\rho$ has zero central component. If other central
residues vary, eliminate the character equations first and use the
resulting smooth parameter space.

Let $X_C^{\rm fr}$ be the stable main family with an ordinary node
frame, variable regular residue $D$, and an ordering of $D$. The residue
evaluation of $X_C^{\rm fr}$ is smooth over that compatible slice by
\S\ref{subsec:A-restrictions}: it is the fine master over the
ordered Cartan together with its smooth orbit framing. Let
\[
 \mathcal C_{a,\lambda,\rho}
 =\{(A,B,F):A\in\mathcal O_a,\ F\in G/B,\
       B\in\mathfrak b_F,\ \overline B=\lambda\rho\},
\]
where $\mathcal O_a\subset\fg$ is the adjoint orbit of the regular
semisimple element $a$ (so $\mathcal O_a\cong G/T$), $\mathfrak b_F$ is
the Borel subalgebra of the flag $F$, and $\overline B\in\ft$ is the
ordered Cartan image of $B$. Pull this space, by
$D=A+B$, to the regular ordered node cover, and restrict to the locus
where the main object is stable. Define
\begin{equation}\label{eq:pg-nodal-receiving}
 Y_{a,\lambda,\rho}
 =\left[
 X_C^{\rm fr}\mathop{\times}_{D=A+B}
 \mathcal C_{a,\lambda,\rho}\;\middle/\;G
 \right].
\end{equation}
The quotient changes the node frame. Off the horizontal source the
ordered eigenvalues of $D$ vary; no condition $D\sim c$ is imposed
on the entire receiving family.

The receiving functor and the affine- and flag-bundle assertions below
are first formed on the centrally rigidified stacks. Coarse fibers may
be finite quotients; rational units, Gysin maps and traces descend by
finite invariants, as in Lemma~\ref{lem:finite-quotient}.
This is the receiving functor with trivial ordinary bubble extension.
Indeed a node frame trivializes such a bubble bundle, and its connection
is uniquely
\[
 \lambda d+A\frac{dz}{z}+B\frac{dz}{z-1},
\]
with residues $(A,B,-A-B)$. This is also a relative assertion. A family
of fiberwise trivial bundles on $\PP^1$ is pulled back from its node
fiber: after a faithful representation, evaluation and cohomology base
change prove it for vector bundles, and the reduction to $G$ is a map
to the affine variety $\GL(V)/G$, hence factors through the base.
The node frame fixes the resulting trivialization.

The stack in~\eqref{eq:pg-nodal-receiving} is of finite type.
A bubble automorphism is constant and determined by its node value,
so its automorphisms inject into those of the stable main object.
Constructing with principal frames and then rigidifying the common
center therefore gives finite stabilizers. No stability condition on a
singular curve is needed, and the construction is independent of the
genus of $C$.

\begin{lemma}[A proper horizontal carrier]\label{lem:pg-nodal-carrier}
Let $F_a(A)$ be the unique Borel flag of regular $A$ with ordered Cartan
image $a$. Define $H\subset Y$ by $F=F_a(A)$, selecting the node order
induced by this same flag, with Cartan value $c=a+\lambda\rho$.
Then $k:H\to Y$ is a closed regular immersion of codimension
$b=\dim G/B$, and main restriction gives an affine-space bundle
$p:H\to M_{c,\lambda}$ of relative dimension $b$.
At $\rho=0$, let $E\subset Y^{\rm sc}$ be $B=0$ with node order $a$.
Then $\pi:E\to M_{a,\lambda}$ is the flag bundle,
$N_{E/Y^{\rm sc}}=T_\pi^*$, and the scalar inverse image $Z=H\times_YE$
is the main space. In particular
\begin{equation}\label{eq:pg-nodal-diagonal}
 p\circ e=\pi\circ r
\end{equation}
for the maps $e:Z\to H$, $r:Z\to E$, throughout the Hodge interval.
The map $\pi r$ is proper, and $Y$ is smooth over the scalar parameter.
\end{lemma}

\begin{proof}
Since $\mathcal O_a=G/T$, the map $F_a:gT\mapsto gB$ is algebraic.
The flag incidence is regular of codimension $b$ in the cap variables;
smooth framed residue evaluation preserves it under the main gluing.
On this incidence the ordered node cover splits into the separated
Weyl translates of the regular section $c$. The selected lift is a
closed component, so $k$ is a closed immersion and hence proper.

Both $A$ and $B$ preserve $F$, and $A+B$ has Cartan value $c$.
Since $c$ is regular, $F$ is the selected node flag. In a main
node frame with $D=c$, the remaining choices are exactly
\[
 B=\lambda\rho+n,\qquad A=a-n,\qquad n\in\mathfrak n.
\]
Every $a-n$ is regular semisimple with Cartan value $a$: conjugation
by the unipotent radical is solved recursively in root height, with
nonzero linear coefficients $\alpha(a)$. Thus $p$ is the affine bundle
associated to the node $T$-torsor with fiber $\mathfrak n$. The argument
uses only the root spaces of $G$, so it holds for every reductive global
form.
To identify its varying main target with the designated reference,
apply the stationary transport along $a+s\lambda\rho$; it is the
identity at $\lambda=0$.

At $B=0$ the bubble is the trivial two-pole extension and its auxiliary
flag is arbitrary, giving $E$. The auxiliary Grothendieck fiber at
Cartan value zero is $T^*(G/B)$; the condition $B=0$ is its zero
section. Therefore $E\subset Y^{\rm sc}$ has codimension $b$ and
normal bundle $T_\pi^*$. The horizontal incidence fixes the flag,
while the scalar incidence fixes its cotangent coordinate; they meet
transversely. Their intersection is the unique common-flag lift of
each main object, proving~\eqref{eq:pg-nodal-diagonal} on the whole
intersection, not merely on a dense open.

Finally $\widetilde{\fg}\to\ft$ is smooth: it is an affine bundle
with nilradical fiber over the flag variety and the Cartan. The base
change of $\widetilde\fg\to\ft$ by $(\lambda,\rho)\mapsto\lambda\rho$ remains smooth even
at their common zero. The physical orbit, ordered regular node cover,
smooth framed gluing and node-frame quotient preserve this relative
smoothness. After rigidification finite quotient charts give rational
local acyclicity, hence $\psi_t\Q_Y=\Q_{Y^{\rm sc}}$.
\end{proof}

\begin{lemma}[Comparison and capped excision]
\label{lem:pg-nodal-excision}
For $\lambda>0$, the receiving family above has an analytic branch
identified with an open neighborhood of the triangular source and
the scalar incidence in the original collar diagram. The identification
retains the source, scalar and character arrows and their finite
continuation faces. The scalar-capped puncture coefficient and its
nearby cycles can be computed in this neighborhood.
\end{lemma}

\begin{proof}
Use the same collar paths as in Proposition~\ref{prop:pg-fusion}.
The main and cap groupoids have relations $RC=1$ and $MUV=1$, with
$C=V^{-1}$, and gluing gives $RMU=1$. Retain node frames and central
sectors. On each positive parameter chart choose the physical and node
logarithms with their specified order. The imaginary root margins
exclude nonzero integral root differences. At the auxiliary point the
scalar arc is chosen small enough to exclude nonzero integral
differences; nilpotent residues and zero differences at $t=0$ are
allowed.

Here the local enriched Riemann--Hilbert equivalence
\cite[Theorems~A and~D]{Boa11} is valid in families. To see the needed
parameter assertion, construct the horizontal local system from its
representation matrices on the complement of the three disks, use
the chosen logarithmic models on the disks, and glue on annuli.
Here the ordinary disk normal form must be computed in $\fg$. Write
the connection as $\lambda d+(R+B(z))dz/z$, $B(0)=0$.
On the series without constant term let
$L=z\partial_z+\ad(R/\lambda)$. On the $n$th coefficient its inverse
is $(n+\ad(R/\lambda))^{-1}$: the root conditions above exclude all
poles for $n\ge1$, including at a nilpotent residue. On compact
positive parameter charts these inverses are $O(1/n)$ uniformly.
Thus both $L^{-1}$ and $z\partial_z L^{-1}$ are bounded in the norm
$\sum_{n\ge1}\|u_n\|r^n$. Write the gauge as $\exp(L^{-1}y)$.
After dividing by $\lambda$, the transformed nonconstant coefficient
is $y+B/\lambda$ plus terms of order
$O(\|y\|^2+\|B/\lambda\|\,\|y\|)$, with the corresponding Lipschitz
bound by the exponential series. Shrinking $r$ makes $B/\lambda$
small and contraction gives the unique based normalizing gauge,
holomorphically in parameters for fixed $\lambda>0$ and locally also
in $\lambda$ away from zero. The inverse monodromy branch is
holomorphic on the group: the differential of the exponential at
$X=-2\pi i R/\lambda$, after left translation, is
$(1-\exp(-\ad X))/\ad X$, with value $1$ at zero. This differential
is invertible on the chosen branch even for nilpotent $R$. Conjugation
equivariance and restriction to a preserving Borel retain the
ordered flags. A faithful representation can have integral differences
between weights whose difference is not a root; the criterion used
here concerns the adjoint action on $\fg$. Nonsingular parallel transport on
the complement has the same dependence. These constructions are
inverse and preserve morphisms, flags and defining tensors, so give
an equivalence of analytic groupoids, not just a bijection of points.
No uniform analytic inverse at $\lambda=0$ is asserted.

Triviality of the ordinary bubble extension is open. For a faithful
associated vector bundle this follows from semicontinuity and the
criterion for triviality near $\mathcal O^{\oplus r}$; the affine
quotient $\GL(V)/G$ then gives the principal assertion. Main stability
is open as well. Choose a small invariant spectral neighborhood of
the fused node values, with disjoint Weyl and logarithm translates.
After shortening the scalar arc, this neighborhood contains all
$a+\lambda w^{-1}\rho$ for the finite packet, and its selected ordered
branch gives an open
receiving neighborhood $U$ in the Betti collar family. This is a
positive-parameter choice and need not extend as a logarithm chart to
$\lambda=0$.

The whole triangular source lies in $U$: the horizontal Borel reduction of a source point has
toral raw residues $(a,\lambda\rho,-a-\lambda\rho)$; since
$\lambda\ne0$, the residue theorem gives degree zero for the line bundle
of every character of $T$. The torus bundle on $\PP^1$ is therefore trivial,
and a root-height filtration of the unipotent radical, together with
$H^1(\PP^1,\mathcal O)=0$, makes the Borel bundle trivial. The reduction
to $B$ is then constant, since a map $\PP^1\to G$ is constant. The scalar incidence
also lies in $U$, with its trivial two-pole bubble. The selected
physical ordering identifies the horizontal condition with $F=F_a(A)$.
At $\lambda=0$ no degree conclusion is deduced from a residue sum;
the algebraic source has already been defined in the degree-zero
ordinary extension.

For the precise excision, write $j:U\hookrightarrow Y_B$,
$k_B:H_B\to Y_B$, and $i_B:E_B\to Y_B$ for the positive Betti diagram.
Both $k_B$ and $i_B$ factor through $U$. For a retained coefficient
$\mathcal F$ on the proper source, open restriction and proper base
change give
\[
 j^*Rk_{B*}\mathcal F=Rk_{U*}\mathcal F,\qquad
 i_B^*Rk_{B*}\mathcal F=i_U^*Rk_{U*}\mathcal F.
\]
The Gysin arrow restricts to the same arrow on $U$, with the same
complex orientation; the entire affine-bundle source and its unit
are retained. Applying $R\pi_*$ and the flag trace gives the original
capped map. For the scalar family, nearby cycles commute with open
restriction, and $k_B$ is proper. Consequently the same identities
retain its specialized input and output arrows. For the nonproper
main projection we use the natural pullback-to-nearby-cycles arrow,
not a base-change isomorphism.

The associated character line and its section are tensor constructions
on the flag torsor and are preserved by this functor. With the paths
fixed, groupoid gluing respects radial continuation and the retained
finite braid words. Restrict the same finite diagram to common
refinements of these relative opens; excision is natural for its maps
and homotopies. Thus all its faces are retained. Ordinary frame changes
of the fixed extension descend the construction. If another ordinary
extension is chosen, its capped diagram is compared through this same
Betti diagram and the displayed excision identities; its bubble-trivial
open need not be identical to $U$.
\end{proof}

\paragraph{Deck-labeled carriers.}
For $w\in\WG$, use the horizontal source with parameter $w^{-1}\rho$
and compose with the closed regular deck graph in the two auxiliary
flag resolutions over their common unflagged residue. The projections
of its closure are proper because flag varieties are projective.
Retain the proper specialization of the original graph coefficient,
as in Lemma~\ref{lem:pg-deck-transfer}. Denote the resulting carrier
by $H_w$ and its composed source over the punctured scalar
arc by $j_w:L_w^\circ\hookrightarrow H_w^\circ$.

The graph changes the auxiliary flag and its ordering but not the
unflagged $B$, $D=A+B$, bubble bundle or main object. The selected node
order on the source is $a+\lambda w^{-1}\rho$, and all these orders
lie in the branch just chosen. Hence the composite carrier is proper
over the receiving family. Close its positive graph with all retained
intermediate flags over the whole Hodge interval; this remains a
closed subspace of the same proper flag carrier. At the scalar
inverse image $B=0$, every added flag stratum has the same main object.
Thus~\eqref{eq:pg-nodal-diagonal} persists on the whole scalar inverse
image, and its output is proper over the main family. This argument
also treats finite words and their proper face projections.

\begin{proposition}[The return at a puncture]\label{prop:pg-fixed-return}
Pull back the generic coefficient diagram to the scalar arc
$\rho=t\rho_*$ at $\lambda>0$. With $j_w$ as above put
\begin{equation}\label{eq:pg-fixed-coefficient}
 K_w^{\rm sc}=\psi_tRj_{w*}\Q,\qquad
 F_{w,+}=e_{w,+}^*K_w^{\rm sc},\qquad
 F_{w,0}=i_0^*Rj_{Z_w*}F_{w,+}.
\end{equation}
With its canceling shifts this is a proper cohomological realization
of $A\Phi_w$ and its scalar homotopy. The central operation of this realization is the
positive one transported by the harmonic product. After a physical
winding its central endpoints differ by $t_{p,-w^{-1}\nu}$, so the
returned coefficient lies over the Hitchin diagonal; the same holds
for powers and finite logarithm polynomials. The operation has degree
zero without a character insertion and degree two with a character
insertion, before the trace over $C$.
\end{proposition}

\begin{proof}
Lemma~\ref{lem:pg-nodal-excision} identifies the positive capped
diagram with the original one, including its coefficient and faces.
The receiving smoothness in Lemma~\ref{lem:pg-nodal-carrier} gives
$\psi_t\Q_Y=\Q_{Y^{\rm sc}}$. Apply nearby cycles to the generic
Gysin map, using proper compatibility for $k_w$. The source arrow is
\[
 p_{w,0}^*\psi_t\Q_M\longrightarrow
 \psi_t p_{w,\circ}^*\Q_M\longrightarrow
 \psi_tRj_{w*}j_w^*p_{w,\circ}^*\Q_M;
\]
the first arrow is natural, with no claim of nonproper base change.
Apply scalar pullback to the proper output and then the flag trace.
Proper base change gives $F_{w,+}$ and the canceling shifts. The
proper deck comparison of Lemma~\ref{lem:pg-deck-transfer} identifies
the resulting map and homotopy with $A\Phi_w$.

The scalar comparison on analytic logarithm charts and the following
completion in the real parameter $\lambda$ are taken in rational sheaves.
The weight computation at the scalar point uses the ordinary algebraic
fiber of Lemma~\ref{lem:scalar-coefficient-extension}.
The main harmonic product trivializes
the endpoint constant complexes as product coefficients; stationary
main transports and endpoint graph pullbacks carry these same
trivializations. The scalar carrier is proper over its output through
$\lambda=0$. Thus Lemma~\ref{lem:rc-ray}, with its coefficient product
hypotheses, gives the last term of~\eqref{eq:pg-fixed-coefficient} and
identifies its operation. Scalar pullback has not been interchanged
with this ordinary completion, and the middle graph coefficient has
not been replaced by a constant sheaf on its singular closure.

Insert the source endpoint map
$q_{w^{-1}\nu,\lambda}p_we_w$ of~\eqref{eq:pg-physical}.
Proposition~\ref{prop:pg-transport} compares it with continuation;
its central value is $t_{p,-w^{-1}\nu}p_we_w$.
Equation~\eqref{eq:pg-nodal-diagonal} then identifies the two central
main objects up to this specified principal Hecke modification, which
preserves the Hitchin sections. This is a modification of the main
endpoint, not an assertion that arbitrary lattice changes preserve
bubble triviality. Character insertion is a degree-two coefficient
morphism on the same diagram. Proper joint flag carriers and the
functorial excision and completion identities retain the finite
composition faces, hence powers and logarithm polynomials.
\end{proof}

\subsection{An annular estimate}\label{subsec:annular}

To prove properness we control extensions across the cap as
$\lambda\to0$. The imaginary root gaps bound the root-space gauge
equations uniformly. The Cartan-valued holomorphic part must be kept as
it is; dividing it by $\lambda$ would destroy the estimate. Normalize an
annular coordinate so that $|z|=1$ is the middle circle, and let
\[
 \mathscr A_\sigma=\Bigl\{u=\sum_{n\in\Z}u_nz^n:\ \|u\|_\sigma=
 \sum_n\|u_n\|e^{\sigma|n|}<\infty\Bigr\},
\]
a Banach algebra of holomorphic functions on an annulus. Norms are taken
through a faithful representation; all parameters other than $\lambda$
range over a compact set.

\begin{lemma}[Normalization keeping the Cartan jet]\label{lem:rc-annulus}
Consider an annular $\lambda$-connection
\[
 \lambda d+\bigl(a_\lambda+H_{\lambda,k}(z)+\lambda B_{\lambda,k}(z)
 \bigr)\frac{dz}{z},
\]
where $a_\lambda$ is bounded and continuous with imaginary root gaps
$\ge\delta>0$, $H$ is $\ft$-valued, holomorphic on the inner disk and
vanishes at its center, and $H,B$ are continuous in $\mathscr A_\sigma$
with $B$ bounded and $H$ sufficiently small. For small $\epsilon$ there
is a continuous $G$-valued gauge, bounded with bounded inverse through
$\lambda=0$, bringing the connection to
\[
 \lambda d+\bigl(a_\lambda+H_{\lambda,k}(z)+\lambda\beta_{\lambda,k}
 \bigr)\frac{dz}{z},\qquad\beta_{\lambda,k}\in\ft
\]
bounded and continuous. The root factor $\exp X$ of the gauge satisfies
$\|X\|_\sigma+\|\lambda z\partial_zX\|_\sigma=O(\lambda)$.
\end{lemma}

\begin{proof}
Put $\mathsf D=z\partial_z$. On a root space,
$\mathsf L_\lambda=\lambda\mathsf D+\ad(a_\lambda)$ has inverse with
Laurent multipliers $(\lambda n+\alpha(a_\lambda))^{-1}$. Since $\lambda$
is real, the imaginary gap gives
$\|\mathsf L_\lambda^{-1}y\|_\sigma\le C\|y\|_\sigma$ and, from
$\lambda\mathsf D\mathsf L_\lambda^{-1}=1-\ad(a_\lambda)\mathsf
L_\lambda^{-1}$, also $\|\lambda\mathsf D\mathsf L_\lambda^{-1}y\|_\sigma
\le C'\|y\|_\sigma$ (differentiation alone is unbounded on
$\mathscr A_\sigma$). With $X=\mathsf L_\lambda^{-1}y$ in the root
spaces, the root part of the gauge equation for $\exp X$ is
$y+\pr_{\rm rt}[H,\mathsf L_\lambda^{-1}y]+\lambda B_{\rm rt}+R_\lambda(y)
=0$, where on a small ball $\|R_\lambda(y)\|_\sigma\le C(\|y\|_\sigma^2+
\lambda\|B\|_\sigma\|y\|_\sigma)$ with a matching Lipschitz bound; all
derivatives occur as $\lambda\mathsf DX$. If the linear term in $H$ has
norm $<\tfrac14$, the contraction principle on a ball of radius
$C_0\lambda$ gives the root factor.

The remaining Cartan coefficient is $a_\lambda+H+\lambda C(z)$ with $C$
bounded (linear root terms have no Cartan projection; later root terms
are $O(\lambda^2)$). With $\beta=C_0$ and $\xi=\sum_{n\ne0}C_nz^n/n$, we
have $\mathsf D\xi=C-\beta$ and $\|\xi\|_\sigma\le\|C\|_\sigma$, and the
Cartan gauge $\exp(-\xi)$ removes the nonconstant part of $\lambda C$ and
commutes with $a_\lambda+H$. Continuity at $\lambda=0$ follows because
the multipliers of $\mathsf L_\lambda^{-1}$ and
$\lambda\mathsf D\mathsf L_\lambda^{-1}$ converge with uniform majorants;
applying this to the contractive equation for $y/\lambda$ gives
continuity of the scaled solution. At $\lambda=0$ the root factor is the
identity and the remaining factor is $\ft$-valued.
\end{proof}

\begin{corollary}[Finitely many extension types]\label{cor:rc-types}
Suppose the fused main object has the same exterior flat object as in
Lemma~\ref{lem:rc-annulus} and the ordered logarithm
$a_\lambda+\lambda\rho_k$ with $\rho_k$ bounded. Then
$\mu=\beta_{\lambda,k}-\rho_k\in X_*(T)$ belongs to a finite set; the
extension is obtained by the lattice $z^{-\mu}$ and has disk coefficient
$a_\lambda+H_{\lambda,k}(z)+\lambda\rho_k$. Hence the fused main objects
lie in finitely many continuous families through $\lambda=0$.
\end{corollary}

\begin{proof}
The holomorphic toral term $H$ has no monodromy period, so the two
ordered monodromies are $\exp(-2\pi i(a_\lambda/\lambda+\beta))$ and
$\exp(-2\pi i(a_\lambda/\lambda+\rho_k))$; their equality in $T$ gives
$\mu\in X_*(T)$, and boundedness gives finiteness. The change $z^{-\mu}$
subtracts $\lambda\mu$ and commutes with $H$; nonresonance and the
ordered residue determine the logarithmic extension uniquely. If one
started from an unordered conjugacy class, a Weyl element $w$ would give
$a_\lambda+\lambda\beta-w(a_\lambda+\lambda\rho_k)\in\lambda X_*(T)$;
taking imaginary parts and using regularity forces $w=1$ for small
$\lambda$.
\end{proof}

\subsection{Proper completion and the trace over the curve}
\label{subsec:completion}

\begin{proposition}[Compactness]\label{prop:rc-compactness}
Under~\ref{C:master}--\ref{C:fusion}, over a compact set of receiving
main data and on every sufficiently small physical collar, the set of
input main objects of the output-capped coefficient cylinder is
relatively compact, uniformly over the finitely many scalar, character
and composition prisms.
\end{proposition}

\begin{proof}
On a compact receiving chart, a common small disk and the ordered frame
give $\lambda d+(a_\lambda+zR_\lambda(z))dz/z$ with uniform bounds.
Here one uses the root contraction in the proof of
Lemma~\ref{lem:rc-annulus} with small forcing $zR_\lambda$, rather than
its $O(\lambda)$ conclusion. On the positive-power subalgebra the same
multipliers are uniformly bounded, preserve vanishing at the center,
and give a root gauge of size $O(\|zR_\lambda\|)$, with scaled
derivative of the same size. The Cartan part left after this gauge is
$a_\lambda+H_\lambda(z)$, with $H_\lambda(0)=0$; it is retained in
full. Thus the disk gauge and its inverse are uniformly bounded, and
shrinking the disk makes $H$ small. Read each term of the capped word from its receiving endpoint.
As the block formula in the proof of Theorem~\ref{thm:rc-perverse}
shows, at most one face requires triangular fusion. Up to that face,
a proper closing or flag graph takes a compact set of outputs to a
compact set of input choices. Transport over a compact parameter
prism also preserves compactness. On the resulting compact sets use
finitely many framed charts with universal families. Off the moving
point the modifications identify the two exteriors. A stationary
master transport has velocity $\lambda V$ in these charts
(\ref{C:master}); integrating the derivative of the family
coefficients gives an exterior change $\lambda B$ with $B$ bounded.
Trivialize on a slightly larger annulus. The transition maps and
their inverses are bounded there, and Cauchy's estimate on a smaller
annulus controls their derivatives in the connection transformation.
The finite word and its faces require only finitely many shrinkings.

At the fusion face use the identification of the exterior flat
object in~\ref{C:fusion} before choosing its annular frame. Apply
Lemma~\ref{lem:rc-annulus} to that exterior, retaining the entire
Cartan jet $H$. Even in rank one, removing a term $H(z)=bz$ by the
flat gauge $\exp(-bz/\lambda)$ need not give a bounded gauge.
Corollary~\ref{cor:rc-types} now places the fused inputs in finitely
many extending families; the additive character test gives stable
central limits. Only at this point use compact charts for the fused
inputs and finish the remaining proper graphs and transports in the
same order. The estimate applies separately to every face.
\end{proof}

\begin{proposition}[Proper completion]\label{prop:rc-completion}
The output-capped cylinder and its faces admit a completion proper over
the receiving main space through $\lambda=0$, with central support over
the diagonal of the Hitchin base. The positive coefficient, with its
homotopies, extends by Lemma~\ref{lem:rc-ray}.
\end{proposition}

\begin{proof}
Close the graph of the positive input map in the product of a proper
compactification $\overline M_0$ of the input main family with the
output cap datum. The projection of this closure to the cap datum is proper
by~\ref{C:fusion}, and Proposition~\ref{prop:rc-compactness} excludes
points of $\overline M_0\setminus M_0$ over compact receiving sets; so
the closure is proper over the receiving main family. Every central point
is a limit in one of the finitely many families of
Corollary~\ref{cor:rc-types}; its two endpoint Higgs fields agree outside
the cap, so their invariant polynomials agree on all of $C$. Push the
positive coefficient diagram to these carriers and apply $i_0^*Rj_*$;
for a closed face $v_F$ use $Rj_{I*}Rv_{F,+*}=Rv_{F,I*}Rj_{F*}$ and proper
base change, so that boundary maps and composition homotopies extend
without commuting a closed pullback with $Rj_*$.
\end{proof}

\begin{proposition}[Independence of the collar]\label{prop:rc-collars}
Completions made on different sufficiently small collars, and on shorter
Hodge intervals, are identified by their radial comparison homotopies,
and they glue over the whole receiving space.
\end{proposition}

\begin{proof}
Two choices of circle over a relatively compact chart are compared by
the radial prism in the punctured disk, which is compact; the previous
two propositions apply to it. For the graph presentations here, retain the same positive capped
carrier $Z_+$ and its original coefficient diagram $\mathcal F_+$.
Close its joint positive graph in the fiber product of the two proper
presentations over the receiving datum. Write $j_{12}:Z_+\to Z_{12}$
and $r_i:Z_{12}\to Z_i$ for its graph inclusion and proper projections,
with $j_i=r_i j_{12}$. Give this refinement the coefficient
$Rj_{12*}\mathcal F_+$, so composition gives
$Rr_{i*}Rj_{12*}\mathcal F_+=Rj_{i*}\mathcal F_+$, including the original
source and output arrows. Proper base change then identifies the
central arrows. For a radial comparison use the same original positive
prism and retain its face maps before completion; finite joint graphs
give the higher comparisons.
Restricting to a shorter interval does not change $i_0^*Rj_*$. For a
relatively compact open $V$ of the Hitchin base, $h^{-1}(\overline V)$ is
compact, so finitely many charts and one common circle suffice, and the
comparisons glue.
\end{proof}

\begin{theorem}[Perverse decrease]\label{thm:rc-perverse}
Transported by $\Psi$, the operator $F_\nu^{\rm raw}$ of~\eqref{eq:raw-return}
has perverse degree at most $-2$.
\end{theorem}

\begin{proof}
\emph{Idea.} A morphism $K_h\boxtimes\Q_C\to K_h\boxtimes\Q_C$ of
degree zero over $\mathcal A\times C$ preserves perverse truncations;
the trace over $C$ then has degree $-2$. So it suffices to realize the
whole return, including its attachments at the punctures, as such a
morphism, i.e.\ by correspondences that are proper and lie over the
diagonal of the Hitchin base. Propositions~\ref{prop:rc-compactness}--%
\ref{prop:rc-collars} provide these correspondences.

Over $\Cc$ the closing graphs and the two scalar caps give a
correspondence proper over both main endpoints, with coefficient
extended by Lemma~\ref{lem:rc-ray}; at $\lambda=0$ its endpoints are
meromorphically isomorphic and scalar at $q$, so it lies over the
diagonal of the Hitchin base. At each puncture the vertex map is the
return of Proposition~\ref{prop:pg-fixed-return}, and
Propositions~\ref{prop:rc-completion} and~\ref{prop:rc-collars} supply
the attaching cylinder and its faces. To see that every term of the
logarithm polynomial is covered, write the homotopy-fiber model as
$\mathsf D\oplus\mathsf R[-1]$ ($\mathsf D$ open and disk coefficients,
$\mathsf R$ circle coefficients). A compatible map has block form
\[
 u=\begin{pmatrix}u_{\mathsf D}&0\\ h_u&u_{\mathsf R}\end{pmatrix},\qquad
 (u^m)_{\rm off}=\sum_{j=0}^{m-1}u_{\mathsf R}^{m-1-j}h_uu_{\mathsf D}^j,
\]
and two off-diagonal blocks multiply to zero. The scalar maps have the
same block form. In particular,
\begin{align*}
 (Au^mB)_{\rm off}={}&h_Au_{\mathsf D}^mB_{\mathsf D}
   +A_{\mathsf R}u_{\mathsf R}^mh_B\\
 &+A_{\mathsf R}\sum_{j=0}^{m-1}
   u_{\mathsf R}^{m-1-j}h_uu_{\mathsf D}^jB_{\mathsf D}.
\end{align*}
Here the packet component of $h_A$ is the specialization of the transverse
$AG=1$ square
with the proper deck flags retained. The face $h_B$ first applies
$ce_-=BA$ over the unflagged datum and then the same $AG=1$ square.
The first comparison changes the coefficient on the proper flag
correspondence; it introduces no second fusion. The packet component of $h_u$ is the
continuation rectangle for $\Phi_w$, with the two physical-circle
presentations identified by the gauge $C_w(t)^\ell$ of
\S\ref{subsec:coefficient-enhancement}. The open-restriction components need no fusion.
Thus each lower-left term has at most one triangular incidence, flanked by the indicated finite words of
transports and closing graphs. On its fusion face identify the exterior
by restriction before choosing annular frames; the change of flat
presentation by $C_w(t)^\ell$ does not change this exterior. In the
remaining word, Hecke identifications agree off the moving point and
flag projections leave the exterior fixed. A stationary transport adds
$O(\lambda)$ to the exterior coefficient on a slightly larger annulus.
On compact graph charts a change of frame $g$ and its inverse are
bounded. Cauchy's estimate bounds the derivative on a smaller annulus,
so the additional connection term $\lambda g^{-1}dg$ is $O(\lambda)$.
A finite iteration therefore retains the bound used in
Proposition~\ref{prop:rc-compactness}. The common scalar face of
Proposition~\ref{prop:finite-packet} and the functorial completions
therefore retain all attaching faces, with their original coefficients
and homotopies.

Hence the unintegrated return is a degree-zero morphism
\[
 K_h\boxtimes\Q_C\longrightarrow K_h\boxtimes\Q_C\quad\text{over }
 \mathcal A\times C,\qquad K_h=Rh_*\Q_{M_0}[N],
\]
with $2N=\dim M_0$. The middle completion complexes may be computed
among all rational sheaves; the two endpoint complexes are bounded
constructible, so the resulting arrow belongs to their full
constructible derived category. The curve unit followed by the oriented
trace over $C$ gives a morphism $\varphi:K_h\to K_h[-2]$.
Applying perverse truncation gives
\[
 \ptau_{\le k}K_h\longrightarrow
 \ptau_{\le k}(K_h[-2])
 = (\ptau_{\le k-2}K_h)[-2].
\]
After taking cohomology this sends $P_kH^j(M_0)$ into
$P_{k-2}H^{j-2}(M_0)$, so its perverse degree is at most $-2$.
Lemma~\ref{lem:rc-ray} identifies its action on cohomology
with the Betti operator.
\end{proof}


\section{Proof of the main results}\label{sec:proof}

We first combine the previous sections on a fine reference, then recover
an arbitrary closed Jordan boundary condition by Springer extraction, and
finally return from $Z(G)^0\times G_\der^{\rm sc}$ to the given group $G$
by a finite central quotient. Both steps keep all compatible topological
types and lifting sectors.

\begin{theorem}[The fine comparison]\label{thm:fine-comparison}
On each regular separating fine reference used in the proof of
Theorem~\ref{thm:A},
\[
 P_kH^j(\MD,\Q)=\Psi^*W_{2k}H^j(\MB,\Q)=\Psi^*W_{2k+1}H^j(\MB,\Q),
\]
and the same holds on the image of every finite projector of
\S\ref{subsec:finite}.
\end{theorem}

\begin{proof}
Let first $G$ be semisimple. Theorem~\ref{thm:common-leading} gives a
degree-two operator $E$ whose leading parts make
$(P_\bullet,E)$ and $(W_{2\bullet},E)$ filtered Lefschetz modules with the
same center $N$. Proposition~\ref{prop:betti-weight} and
Lemma~\ref{lem:whole-quadratic-return} give $F_B$ of weight degree at most
$-4$ with $[[F_B,E],E]=-2E$, and Theorem~\ref{thm:rc-perverse} shows that
$\Psi^*F_B(\Psi^*)^{-1}$ has perverse degree at most $-2$.
Theorem~\ref{thm:common-lowering} gives $P_k=\Psi^*W_{2k}$, and
$W_{2k}=W_{2k+1}$ because the cohomology is Hodge--Tate. If the Goldman
class has complex coefficients we argue over $\C$; equality of the
complexifications of rational subspaces implies equality over $\Q$. For any compatible finite idempotent $e$, both filtrations are
preserved. Applying $e$ to the full-space equality gives the equality
on its image, since $eP_k=P_k\cap\operatorname{im}e$ and likewise for $W$.

For general $G$ take the product with the torus factor
(Example~\ref{ex:intro-rank-one}): weight filtrations of products are
tensor-product filtrations, perverse filtrations of products of proper
maps are the convolution filtrations by the perverse K\"unneth formula,
and defects add. The finite central quotient of
Proposition~\ref{prop:central-quotient} then descends the equality by
Lemma~\ref{lem:finite-quotient}.
\end{proof}

\subsection{Extracting a closed Jordan boundary condition}
\label{subsec:extraction}

The Springer resolution separates the cohomology of a full spectral
problem into intersection complexes of nilpotent orbit closures, with
multiplicity spaces that are $\WG$-representations. Choosing a
representation therefore selects a geometric boundary condition;
choosing an arbitrary vector in a multiplicity space would not. We first
identify the summand of the full nilpotent cone and then compute the
shifts for smaller orbit closures.

Write $\mathcal I_Z=\IC_Z[-\dim Z]$ for the unshifted intersection
complex. For a semisimple boundary value $S$, the nilpotent directions lie
in $\mathfrak l=\Lie Z_G(S)^0$, and the full spectral resolution uses the
Springer resolution of each such $\mathfrak l$. When $Z_G(S)$ is
disconnected, its component group permutes these data; all translates are
kept before taking invariants.

\begin{lemma}[The degree-zero unit]\label{lem:degree-zero-unit}
Let $e_\WG$ be the projector onto $\WG$-invariants for the deck action of
the Grothendieck--Springer resolution. On the nilpotent cone,
$e_\WG R\mu_*\Q_{\widetilde\Nil}=\Q_\Nil$, the inclusion being the
pullback unit, with no shift and no Tate twist; this commutes with
arbitrary base change.
\end{lemma}

\begin{proof}
The Grothendieck--Springer resolution is small, so its shifted direct
image is the intermediate extension of the regular representation local
system of the $\WG$-cover; averaging selects the constant local system,
whose intermediate extension on the smooth Lie algebra is constant.
Removing the common shift and applying proper base change gives the
statement on $\Nil$; the inclusion is the unit on the dense regular
locus, hence everywhere \cite[\S\S2.2, 2.6--2.7]{BM83},
\cite[\S2.5, Lemma~2.1]{AHJR14}.
\end{proof}

This full spectral unit (trivial representation) differs from the sign
summand used in~\eqref{eq:scalar-sign}; confusing the two would introduce
a spurious shift or a factor $|\WG|$.

\begin{proposition}[Closed Jordan multiplicities]\label{prop:jordan-multiplicity}
Let $\pi:Y\to X$ be a resolved full spectral boundary problem of
dimension $2N$, and let $i:Z\hookrightarrow X$ be the locus of a
closed Jordan condition, with constant intersection-complex coefficients. Let $b$ be the sum over the punctures of the
dimensions of the corresponding Springer fibers, and $V_\rho$ the
irreducible representation of the local Weyl group $\WG_{\rm loc}$
attached to the orbits and the constant local system. In the
multiplicity functor we give this rational representation Hodge type
$(0,0)$, so the Tate twist from the fiber is displayed separately.
Then, on the
ordered connected-centralizer cover,
\begin{equation}\label{eq:jordan-multiplicity}
 \Hom_{\WG_{\rm loc}}\bigl(V_\rho,R\pi_*\Q_Y\bigr)
 \cong i_*\mathcal I_Z[-2b](-b),
\end{equation}
as mixed Hodge module complexes, with zero right-hand side if $Z$ is
empty. If $Z$ is nonempty, its specified dense orbit type occurs and
$Z$ has pure dimension $2(N-b)$. For a nonconstant Springer local
system the corresponding summand uses its intersection complex.
\end{proposition}

\begin{proof}
By Proposition~\ref{prop:springer}\ref{it:semismall}, the local
Springer map is semismall, so the decomposition theorem involves no
further shifts and the multiplicity of an orbit closure is read from the
top cohomology of the Springer fiber, of Tate weight $2b$. Every
\emph{local} orbit with constant coefficients occurs by
Proposition~\ref{prop:springer}\ref{it:semismall}. The framed residue
evaluation is smooth over the product of nilpotent cones
(\eqref{eq:framed-complex} in Section~\ref{sec:comparison-families}).
After base change to the chosen product of orbit closures it is still
smooth, hence open. Every nonempty open of the inverse image meets
the preimage of the dense product orbit. Smooth pullback thus gives
the intersection complex with its constant top-stratum line wherever
the inverse image is nonempty, preserving codimension $2b$.
Removing the dimension shifts gives $[-2b](-b)$; products add the
shifts and twists. Descending through frames and finite component
actions gives~\eqref{eq:jordan-multiplicity}. For $\GL_n$ the
character-variety form of this decomposition is
\cite[Theorem~3.8]{Bal23}; the mixed Hodge module decomposition and
its Tate twists are explicit in \cite[Theorem~5.4]{Let15}.
\end{proof}

In type $A$, if $\lambda_{pa}$ is the Jordan type in the block of size
$m_{pa}$, then $b=\sum_{p,a,j}\binom{(\lambda_{pa})'_j}{2}$. For a block
of size $2$, the semisimple condition (zero nilpotent orbit) has $b=1$:
its cohomology appears in degree $j+2$ of the resolved space, with weights
shifted by two.

Local occurrence need not give a nonempty global locus. On the
three-punctured sphere take $G=\SL_2$, $\zeta=1$ and the tuple
$(U,A,(UA)^{-1})$, where
\[
 U=\begin{pmatrix}1&1\\0&1\end{pmatrix},\qquad
 A=\begin{pmatrix}0&-1\\1&0\end{pmatrix}.
\]
The last two traces are $0$ and $1$. Their eigenvalues cannot be
mutual inverses, so the semisimple tuple is intrinsically generic.
The full spectral condition at the first puncture is realized by $U$.
Prescribing the identity there would force the last two matrices to
be inverses and to have equal traces. Its global summand is zero.

When a component group permutes local Springer pairs, first form
the direct sum over the orbit of the pair, with the equivariant structure
inherited from the Springer incidence, and then take invariants.
A single multiplicity representation descends only when its isomorphism
class and its coefficient-system descent data are preserved.
For the comparison below we use the constant coefficient on each
target orbit closure. Extending the conclusion to nonconstant Springer
systems would additionally require transport of those systems under the
equivariant harmonic product.

\begin{lemma}[Matching shifts]\label{lem:matching-shifts}
With the defect $N-b$ on $Z$, the multiplicity~\eqref{eq:jordan-multiplicity}
gives
\begin{gather*}
 \Hom_{\WG_{\rm loc}}\bigl(V_\rho,P_{k+b}H^{j+2b}(Y_\Dol)\bigr)
 =P_k\IH^j(Z_\Dol)(-b),\\
 \Hom_{\WG_{\rm loc}}\bigl(V_\rho,W_{2k+2b}H^{j+2b}(Y_\B)\bigr)
 =W_{2k}\IH^j(Z_\B)(-b),
\end{gather*}
and similarly for odd weight steps.
\end{lemma}

\begin{proof}
The normalized source summand is
$\mathcal I_Z[N-2b](-b)=\bigl(\mathcal I_Z[N-b]\bigr)[-b](-b)$. The extra
shift $[-b]$ raises the perverse index by $b$ and the twist $(-b)$ raises
weights by $2b$. Proper direct image to the common Hitchin base preserves
the direct summand, perverse truncation commutes with the projector, and
morphisms of mixed Hodge structures are strict.
\end{proof}

\begin{theorem}[Propagation to closed Jordan conditions]
\label{thm:boundary-descent}
Theorem~\ref{thm:fine-comparison} implies~\eqref{eq:thmA} for every closed
Jordan boundary condition with constant intersection-complex
coefficients, and on every compatible finite-projector summand. Together with
\S\ref{subsec:central-quotient}, this proves Theorem~\ref{thm:A}.
\end{theorem}

\begin{proof}
First transport the equality on the entire cohomology of the fine
resolution. Work on the ordered cover before finite central descent.
Choose the rational fine refinement $q_{\rm f}$ sufficiently close
to the target growth weight $a$ and generic in the derived directions,
so $a-q_{\rm f}$ lies in an open Weyl chamber at each puncture.
Use its Borel ordering. Choose a nearby regular growth weight
$a_s=(1-s)a+s q_{\rm f}$ with $0<s<1$ small and perturb the imaginary
residues in the derived Cartans off the finitely many root and
testing-character hyperplanes. The perturbation preserves the central
equations. The target is intrinsically generic, an open condition,
so both endpoints and the small Betti path between them remain generic.
The fine weight and $a_s$ can be taken in the same facet and stability
chamber: in a sufficiently small neighborhood only walls through $a$
need be considered, and both points are on the same side of each wall
they leave. The imaginary perturbation gives a separating reference.
At $\lambda=1$, the real part of the horizontal exponent minus the fine
weight is $a-q_{\rm f}$ at the target and $(1-s)(a-q_{\rm f})$ at
the reference. One logarithm chart and ordering thus give the enriched
Riemann--Hilbert comparison along the short path
(Lemma~\ref{lem:enriched-rh}).

For $u=s,c$ (reference and target), let $r_{D,u}$ and $r_{B,u}$ denote
restriction from the cohomology of this one fine Hodge master to its
Higgs and enriched Betti endpoints. They are isomorphisms by
Lemma~\ref{lem:basechange}. Write $J_u=r_{D,u}r_{B,u}^{-1}$ and
$T_*=r_{*,c}r_{*,s}^{-1}$. Then $J_cT_B=T_DJ_s$.
At the regular reference $J_s=\Psi_s^*$ by
Theorem~\ref{thm:harmonic-transport}, and
Theorem~\ref{thm:fine-comparison} identifies the two filtrations there.
Theorem~\ref{thm:fine-master} transports $P$ by $T_D$, and
Theorem~\ref{thm:weak-betti} transports $W$ by $T_B$. Hence $J_c$
identifies the two filtrations on the entire target fine cohomology.

Now take the Springer multiplicity along the target's own real Hodge
path. Its graded centralizer is fixed, so on framed charts the
incidence is the pullback of the Springer incidence for that
centralizer. Proper base change transports its restriction action
along this path. Theorem~\ref{thm:harmonic-transport} provides a
stratified product structure preserving the tuples of nilpotent orbit
types. The unshifted intersection complex along the path is therefore the
pullback of the initial one (both are the intermediate extension of the
same top-stratum coefficient), and the ratio of the endpoint restrictions
is the map on intersection cohomology induced by $\Psi$. Applying the exact
multiplicity functor and Lemma~\ref{lem:matching-shifts}, and canceling
the common shift and twist, gives~\eqref{eq:thmA} with the defect of each
target. Component-group translates are included before taking invariants,
and the full-space equality restricts to each compatible finite
projector because its image inherits both filtrations by intersection.
Combined with
\S\ref{subsec:central-quotient} this proves Theorem~\ref{thm:A}.
\end{proof}

\subsection{Central lifting on root curves}\label{subsec:central-lifts}\label{app:central-lift}

The reduction from $G$ to $\widehat G=Z(G)^0\times G_\der^{\rm sc}$ must
preserve objects and automorphisms. The obstruction to lifting a bundle
along a central isogeny is a degree-two cohomology class; the nonempty
puncture set gives a uniform way to kill it without changing the
problem downstairs. We then identify the finite group of changes of lift,
including changes at the punctures. All compatible lifted local types
and degree components are retained: the full group of changes of lift
acts on their union, and may permute its components.

Fix a central isogeny $1\to K\to\widehat G\xrightarrow{q}G\to1$ and let
$e$ be the exponent of $K$. A \emph{lift} of a $G$-object always includes
an identification of its pushout with the given object. A \emph{root
curve} $C[\sqrt[N_1]{p_1},\ldots,\sqrt[N_s]{p_s}]$ replaces a disk at
$p_i$ by $[\Spec\C[u]/\mu_{N_i}]$ with $u^{N_i}=z$; it has coarse space
$C$ and tautological root lines $L_i$ with $L_i^{N_i}=\mathcal O_C(p_i)$
\cite[Definitions~2.2.1, 2.2.4 and Example~2.4.1]{Cad07}. Rational parabolic bundles are bundles on root curves with
prescribed stabilizer actions \cite[Theorem~6.1]{BV12}, and the principal
version is \cite[Proposition~6.1]{BMW12} (arXiv v1 numbering).

\begin{lemma}[Uniform root refinement]\label{lem:central-root-H2}
For $\mathcal C_{\mathbf N}=C[\sqrt[N_1]{p_1},\ldots,\sqrt[N_s]{p_s}]$
and the refinement $r:\mathcal C_{e\mathbf N}\to\mathcal C_{\mathbf N}$,
\[
 r^*:H^2_{\rm et}(\mathcal C_{\mathbf N},K)\to
 H^2_{\rm et}(\mathcal C_{e\mathbf N},K)\ \text{is zero},\qquad
 H^1_{\rm et}(\mathcal C_{e\mathbf N},K)=\Hom(\pi_1(\Cc),K)\cong
 K^{2g+s-1}.
\]
Consequently every $G$-bundle on $\mathcal C_{\mathbf N}$ lifts after
pullback to $\mathcal C_{e\mathbf N}$, families lift \'etale locally on
the parameter space, and the lifts form a torsor under the Picard groupoid
of $K$-torsors.
\end{lemma}

\begin{proof}
Write $K=\prod\mu_m$ with $m\mid e$. Set $X=\mathcal C_{\mathbf N}$
and $U=C\setminus D$. The open curve has $H^2_{\rm et}(U,\mu_m)=0$
\cite[Tag~03RR]{Stacks}; over $\C$ this also follows by finite-coefficient
comparison since $U$ retracts onto a graph. Localization therefore
generates $H^2_{\rm et}(X,\mu_m)$ by the classes of the root divisors
$\mathcal D_i$ above the punctures. Here is the local calculation.
On a smooth \'etale chart at $\mathcal D_i$ the strict henselian local
ring is a discrete valuation ring. Its units have $m$th roots, so by
Kummer theory the degree-one cohomology of its punctured trait with
coefficients $\mu_m$ is $\Z/m$, generated by the valuation of a
uniformizer. The support cohomology sheaves vanish below degree two
and in degree two give this constant $\Z/m$ on the divisor. A change
of uniformizer multiplies by a unit, so the identification descends
to the residual gerbe. In total degree two the local-to-global sequence
for supports therefore uses only $H^0(\mathcal D_i,\Z/m)$. The support-to-global map sends $1$ to
$c_1(\mathcal O_X(\mathcal D_i))$ by the Kummer boundary for the local
equations of the divisor. Refinement pulls that line back to
$\mathcal O(e\mathcal D'_i)$; its class vanishes modulo $m$. This proves
the assertion about $H^2$. Equivalently, one can use the vanishing of
the Brauer group and the root-curve Picard presentation
\cite[Proposition~5.1]{San23}, \cite[Theorem~1.1]{Lop23}, together
with Tsen's theorem for $\C(C)$ \cite[Hauptsatz]{Tse33}.

The obstruction to a central lift is a $K$-gerbe whose class lies in this
$H^2$; its vanishing gives a lift, two lifts differ by a $K$-torsor, and
automorphisms of a lift over the given object form $K$. For a family over $S$, put $X=\mathcal C_{e\mathbf N}$.
The stack of lifts of the pulled-back bundle, with its pushout
identification, is locally of finite presentation. We give the curve-chart
argument, which also applies to nonreduced parameter schemes.
Choose $\widehat G\subset\GL(V)$. The homogeneous space
$\GL(V)/\widehat G$ is affine \cite[Theorem~12.15(ii)]{Alp13}, and finite
stable coordinate generators embed this space as a closed orbit in a
representation $W$. On a curve Quot chart $T$, retain the prescribed
deck and inertia data of a fixed covering curve $Y$ and let
$\mathcal F$ be the associated $W$-bundle on $Y\times T$.
Write $p:Y\times T\to T$. For an ordinary point $y\in Y$ and $m\gg0$,
after shrinking $T$, the sequence
\[
 0\longrightarrow\mathcal F\longrightarrow\mathcal F(my)
 \longrightarrow\mathcal F(my)|_{my}\longrightarrow0
\]
computes sections after every base change as the kernel of
$U\to V_1$, where $U=p_*\mathcal F(my)$ and
$V_1=p_*(\mathcal F(my)|_{my})$ are finite locally free.
Thus the section functor is a closed affine kernel, without assuming
base change for $p_*\mathcal F$. Finite equivariant polynomial equations
for the closed orbit, applied in the same section models for their
target bundles, cut out the reductions scheme-theoretically.
Isomorphisms form the open invertible locus in Hom-section spaces,
with closed preservation equations. The pushout identification and
finitely many deck arrows and cocycle equations therefore give the
required finite-presentation charts.
Over a nilpotent thickening, lift local trivializations of the given
$G$-bundle compatibly with a specified lift on the reduction. Its
transition functions then lift uniquely through the \'etale $q$, so
their cocycle and the prescribed comparison lift as well. A smooth
scheme atlas of the lift stack is therefore smooth over $S$ by the
infinitesimal criterion \cite[Tag~0DP0]{Stacks}. Each geometric fiber has a lift by the preceding
calculation, giving sections \'etale locally on $S$.
On this fixed root curve, rigidify a $\mu_m$-torsor at an ordinary point.
Its inertia weights are $(eN_i/m)b_i$ with $b_i\in\Z/m$; for representatives
$0\le b_i<m$, its coarse line factor satisfies
$M^m=\mathcal O_C(-\sum_i b_i p_i)$ in the rigidified Picard scheme.
Thus $\sum_i b_i=0\bmod m$, and each permitted tuple gives a fiber of
multiplication by $m$ on a Picard component, finite \'etale and constant.
Restoring the framing gives the relative automorphisms $\mu_m$.
Taking products over the cyclic factors gives the constant group of
changes $H^1(X,K)$ after rigidification.
For Higgs fields and $\lambda$-connections, the isomorphism
$dq:\widehat{\mathfrak g}\to\mathfrak g$ lifts the local forms
uniquely on a lifted bundle; the transition formulas make these
lifts compatible. Thus the same assertion holds for the fields. Finally the orbifold fundamental group
is $\pi_1(\Cc)$ with the extra relations $c_i^{eN_i}=1$, which every
homomorphism to $K$ satisfies; hence
$H^1(\mathcal C_{e\mathbf N},K)=\Hom(\pi_1(\Cc),K)$.
\end{proof}

\begin{lemma}[The full local groupoid]\label{lem:central-full-local}
The root presentation keeps the full centralizer of the local cyclic
stabilizer action; its arrows include all component identifications of
the tame problem.
\end{lemma}

\begin{proof}
In a chart $z=u^N$ with rational weight $\theta$ ($N\theta$ integral) put
$s_\theta=\exp(2\pi i\theta)$. An automorphism of the root bundle is an
equivariant holomorphic $a(u)$ with $a(0)\in Z_G(s_\theta)$, the full
group; in the punctured coarse trivialization it is
$g(z)=u^{-N\theta}a(u)u^{N\theta}$. Conversely, a single-valued $g$
defines an automorphism precisely when its conjugate extends
holomorphically and invertibly over the root chart. With a prescribed
graded residue one takes its stabilizer in this full group, retaining
all components. Replacing $Z_G(s_\theta)$ by its identity
component would change the groupoid. A change of lift at one puncture
need not extend globally; global compatibility is exactly the peripheral
product condition in $H^1(\Cc,K)$.
\end{proof}

\begin{example}[A disconnected local automorphism group]
\label{ex:full-root-PGL2}
Let $G=\PGL_2$ and let $\alpha^\vee(t)$ be the class of
$\operatorname{diag}(t,t^{-1})$. Take $\theta=\alpha^\vee/4$ and the
root chart $z=u^4$. The class $s_\theta=[\operatorname{diag}(i,-i)]$
is centralized by $w=[\left(\begin{smallmatrix}0&-1\\1&0\end{smallmatrix}\right)]$,
so the constant root-chart gauge $a(u)=w$ is an automorphism. In the
coarse punctured frame this gauge is
\[
 u^{-\alpha^\vee}wu^{\alpha^\vee}
 =\left[\begin{pmatrix}0&-z^{-1}\\1&0\end{pmatrix}\right].
\]
The determinant valuation is odd, whereas every element of the connected
parahoric has even determinant valuation modulo scalars. Thus this
automorphism belongs to the full local group but not its connected
parahoric. Discarding the extra component would change the boundary
problem.
\end{example}

\begin{proposition}[The Dolbeault central quotient]
\label{prop:central-Dolbeault-quotient}
For a rational parahoric Higgs problem with full local automorphism groups
and closed invariant residue conditions, let $\widehat\MD$ be the disjoint
union of all compatible $\widehat G$-lift spaces. Extension of structure
group induces an algebraic finite quotient
$\MD=\widehat\MD/\Delta$, $\Delta=H^1(\Cc,K)$, over the same Hitchin
base, preserving (semi)stability and commuting with $\Psi$.
\end{proposition}

\begin{proof}
Refine the root denominators by $e$; the $G$-objects are unchanged (their
inertia factors through the old inertia). Lemma~\ref{lem:central-root-H2}
gives lifts \'etale locally in families, of exactly the permitted local
types. Since $dq$ is an isomorphism, adjoint bundles agree and a
logarithmic Higgs field lifts uniquely, with the same invariant
polynomials, graded residues and nilpotent orbit conditions; the
refinement adds no flag inside an equal-weight block, since its extra
inertia is central. Rigidifying the relative automorphisms $K$ gives a
finite \'etale $\Delta$-torsor over the stack downstairs; the residual
$BK$-gerbe does not change coarse spaces, so $\MD=\widehat\MD/\Delta$
algebraically (finite geometric fibers alone would not prove this).
Parabolic reductions correspond, and a positive power of each character of
a lifted parabolic descends, so stability inequalities keep their signs;
finite central torsors have parabolic degree zero. Finally the isogeny
identifies the symmetric spaces of compatible maximal compact subgroups,
so harmonic metrics correspond and $\Psi$ commutes with $\Delta$.
\end{proof}

For real weights, the rational presentation must keep the facet and the
full local group: in a bounded chart keep the finitely many equations
$w\theta-\theta=\mu$ ($\mu\in X_*(T)$) holding at the given weight and
acquire no new one; with~\eqref{eq:fuchs} and the stability chamber these
are the charts of Lemma~\ref{lem:rational-chart}. The loop representatives
$z^{w\theta-\theta}\dot w$ of the remaining components are then constant,
so the local groupoid is preserved; the harmonic square uses the real
weights.


\subsubsection{The quotient for the given reductive group}
\label{subsec:central-quotient}

A connected reductive group is a finite central quotient of a torus times
a simply connected semisimple group. To use this, every object downstairs
must lift, and the full finite group of changes of lift must be kept; on a
punctured curve this group contains peripheral changes as well as handle
changes. Choose a central isogeny
\[
 1\longrightarrow K\longrightarrow\widehat G=Z(G)^0\times G_\der^{\rm sc}
 \longrightarrow G\longrightarrow1,
\]
let the moduli spaces upstairs be the disjoint unions over all compatible
peripheral lifts and topological types, and put
\begin{equation}\label{eq:central-group}
 \Delta=H^1(\Cc,K)=\Hom(\pi_1(\Cc),K)\cong K^{2g+s-1}.
\end{equation}
Here $\Delta$ acts on this complete union of lifted sectors. A fixed
tuple of lifted boundary classes is preserved only by the subgroup
whose peripheral multipliers preserve that tuple.

\begin{proposition}[The central quotient]\label{prop:central-quotient}
The Betti and Dolbeault spaces for $G$ are the quotients by $\Delta$ of
the corresponding spaces for $\widehat G$. The Dolbeault quotient map
lies over the common Hitchin base, and both quotient maps commute with
the tame nonabelian Hodge maps.
\end{proposition}

\begin{proof}
On the Betti side, $\pi_1(\Cc)$ is free; lift all generators to
$\widehat G$, and absorb the discrepancy in $K$ of the lifted surface word
into one puncture. Keeping all boundary lifts gives exactly a torsor under
$\Hom(\pi_1(\Cc),K)=\Delta$ before conjugation, and the invariant rings
commute with this finite action. (Using only the handle subgroup $K^{2g}$
would omit peripheral sectors.) The Dolbeault statement is
Proposition~\ref{prop:central-Dolbeault-quotient}. A central isogeny does
not change $\ad(E)$ or the invariant polynomials, so the Hitchin bases
agree. Finally, the isogeny identifies the noncompact symmetric spaces of
compatible maximal compact subgroups, so a lifted harmonic metric is the
same harmonic map; finite central twists are unitary with zero Higgs
field. Hence the construction commutes with $\Psi$ and with $\Delta$.
\end{proof}

Products have tensor-product weight filtrations and convolution perverse
filtrations with additive defects (perverse K\"unneth). By
Lemma~\ref{lem:finite-quotient}, $\Delta$-invariants descend the comparison
from $\widehat G$ to $G$ with no shift. On the complete union upstairs,
isotypic components for characters of $\Delta$ are taken over a splitting
field, or summed over Galois orbits over $\Q$. The invariant component
is the comparison on the prescribed downstairs problem; an isotypic
statement on a single fixed boundary sector uses the stabilizer
subgroup of that sector. This completes the proof of Theorem~\ref{thm:A}.

\subsection{The usual fixed-determinant problem}
\label{subsec:fixed-determinant-root}

For $\SL_n$, the customary degree-$d$ problem consists of rank-$n$
parabolic Higgs bundles with a specified determinant line of degree $d$;
when that line is nontrivial these are not ordinary principal
$\SL_n$-bundles on $C$. The following construction identifies that
problem with a specified central-inertia sector on a root curve and
is compatible with the nonabelian Hodge correspondence.

Work first with rational weights on a root curve $\mathcal C$.
Fix the \emph{full parabolic determinant} $(\mathcal L,\psi)$: this
includes its stabilizer characters, not just its coarse degree, and
$\psi=\operatorname{tr}\varphi$ is the fixed determinant Higgs field.
Apply Lemma~\ref{lem:central-root-H2} to $\mu_n$. After a refinement
$r:\mathcal C'\to\mathcal C$, the Kummer obstruction to taking an
$n$th root of $\mathcal L$ vanishes. Choose once and for all a line
$Q$ and an isomorphism $Q^{\otimes n}\cong r^*\mathcal L$, and give $Q$
the Higgs field $r^*\psi/n$. Then
\begin{equation}\label{eq:fixed-determinant-root}
 (E,\varphi)\longmapsto
 \left(F=r^*E\otimes Q^{-1},\quad
 \varphi_F=r^*\varphi-\tfrac1n r^*\psi\,\mathrm{id}\right)
\end{equation}
has trivial determinant and trace-free Higgs field. At each stacky
point its inertia representation is the original one multiplied by the
inverse scalar character of $Q$. In particular, the extra inertia of
$r$ acts centrally. We retain exactly these local types, with the
original flags and their component identifications.

This is an algebraic equivalence of moduli groupoids, including families.
The inverse equivalence tensors $F$ by $Q$ and adds $r^*\psi/n$ to the Higgs field;
the extra inertia then acts trivially, so the resulting bundle and
Higgs field descend to $\mathcal C$. The determinant identification
descends to the prescribed one with $\mathcal L$. Tensoring changes
the parabolic degree of a rank-$a$ subbundle by
$-a\,\operatorname{pardeg}(Q)$, so every subbundle slope and the total
slope change by the same amount. Thus stability and semistability are
preserved. The nilpotent parts of the residues and their flags are
unchanged; semisimple residues are shifted only by the prescribed
scalar $\psi/n$. The endomorphism bundles agree, and the traceless
Hitchin coefficients agree. Hence this equivalence is over the same
traceless Hitchin base. In the fixed-trace $\GL_n$ notation the change
of characteristic-polynomial coefficients is the fixed translation
$T\mapsto T+\psi/n$, an isomorphism of Hitchin bases.

The equivalence~\eqref{eq:fixed-determinant-root} also identifies the
harmonic constructions. Fix an adapted
determinant metric $h_{\mathcal L}$ and give $Q$ its root metric
$h_Q^{\otimes n}=r^*h_{\mathcal L}$. If $h_E$ is a projectively
harmonic metric with this determinant, then
$h_F=r^*h_E\otimes h_Q^{-1}$ satisfies
\[
 F_{h_F}+[\varphi_F,\varphi_F^{\dagger}]
 =r^*\left(F_{h_E}+[\varphi,\varphi^{\dagger}]
       -\tfrac1n F_{h_{\mathcal L}}\,\mathrm{id}\right)=0.
\]
Conversely, tensoring a harmonic metric on $F$ by $h_Q$ gives the
projectively harmonic metric on $E$. The scalar root shifts the
prescribed growth weights, and subtracting its connection subtracts
exactly one $n$th of the determinant connection. Thus the ordinary
flat connection on the root curve is the one defining the usual
fixed-determinant twisted monodromy, with these choices; this comparison
is not an inference from an abstract homeomorphism of moduli spaces.

For an explicit degree example, suppose the full parabolic determinant
has trivial inertia and is the pullback of an ordinary line $L$ of
degree $d$ (in particular this holds with trace-free local weights).
Choose $p\in D$, a degree-zero line $A$ with
$A^{\otimes n}\cong L(-dp)$, and, on a common root refinement, a line
$R$ with $R^{\otimes n}\cong\mathcal O_C(p)$. Such $A$ exists because
multiplication by $n$ on $\Pic^0(C)$ is surjective. One can take
\[
 Q=\pi^*A\otimes R^d.
\]
Thus the original ordinary degree $d$ is recovered by the inverse
construction. The root twist contributes scalar inertia $\epsilon^{-d}$
in addition to the inherited inertia at $p$, where $\epsilon=e^{2\pi i/n}$. With the convention
$c=\exp(-2\pi i\vartheta)$ of~\eqref{eq:nah-dictionary}, the scalar
weight shift $-d/n$ multiplies that peripheral monodromy by
$\exp(2\pi i d/n)$. Accordingly the surface-relation convention
$\prod_i[A_i,B_i]\prod_jM_j=\zeta$ uses
$\zeta=\exp(-2\pi i d/n)$ here, and absorption replaces
$M_p$ by $\zeta^{-1}M_p$. Reversing the surface-loop convention reverses
both signs. If the determinant has nontrivial parabolic characters,
one must use its full root $Q$ above; the formula using $L$ alone does
not apply.

The choices are not canonical: two roots differ by a $\mu_n$-torsor,
and their two identifications differ by the corresponding finite
central twist of $F$. All compatible sectors are retained, as in
Proposition~\ref{prop:central-Dolbeault-quotient}; on a single prescribed
sector only its stabilizer subgroup acts. The same construction applies
to real weights in the compatible rational presentations of
Lemma~\ref{lem:rational-chart}, retaining the original real growth data
in the harmonic correspondence. This gives the fixed-determinant
interpretation used in Corollary~\ref{cor:B}, without replacing its
ordinary degree by an arbitrarily chosen degree-zero extension.

\subsection{Ordinary cohomology and finite symmetry summands}
\label{subsec:rational-smooth}

\begin{proof}[Proof of Corollary~\ref{cor:ordinary}]
For a rationally smooth variety $X$, the identification
$\Q_X[\dim X]=\IC_X$ is canonical and compatible with mixed Hodge
structures. Theorem~\ref{thm:A} therefore gives \eqref{eq:PW}.

For the stated sufficient condition, the nilpotent cone is rationally
smooth (Proposition~\ref{prop:springer}\ref{it:rational-smooth}), and so
are its products with smooth varieties and its smooth pullbacks. The
framed residue evaluation is smooth (\eqref{eq:framed-complex}), so the
generic full spectral space is, locally, a finite quotient of a smooth
pullback of a product of nilpotent cones; at semisimple punctures one
replaces the nilpotent cone by its vertex. Finite quotients preserve
rational smoothness (Lemma~\ref{lem:finite-quotient}). A smaller orbit
closure need not be rationally smooth, and no replacement is made there.
\end{proof}

\begin{proof}[Proof of Corollary~\ref{cor:equivariant}]
The Betti action preserves the weight filtration, and the Dolbeault
action preserves the perverse filtration because it is over the
Hitchin base (\S\ref{subsec:finite}). Compatibility with $\Psi$ makes
the two actions agree under the cohomological identification.
Every $e\in\Q[\Gamma]$ therefore preserves both filtrations. Applying
$e$ to the equality of subspaces in \eqref{eq:thmA} gives that equality
on its image, with the induced filtrations. Corollary~\ref{cor:ordinary}
gives the ordinary-cohomology version. Extending scalars gives the
individual character projectors described in \S\ref{subsec:finite}.
\end{proof}

\begin{proof}[Proof of Corollary~\ref{cor:B}]
The determinant-root equivalence of
\S\ref{subsec:fixed-determinant-root} identifies the usual
fixed-determinant degree problem with its specified central-inertia
sector, over the same traceless Hitchin base and compatibly with the
harmonic correspondence. For $\GL_n$ and fixed-determinant
$\SL_n$, generic semisimple tuples are
stable with stabilizer the center (Lemma~\ref{lem:generic}). Choose a
small full-flag refinement of every scalar block; the invariant-subspace
test is unchanged on each block, and the generic character margins keep the
refinement inside the original stable problem. The Springer family and its
scalar extraction recover the original multiplicities and degree, so
Theorem~\ref{thm:boundary-descent} and Corollary~\ref{cor:ordinary}
apply, and Corollary~\ref{cor:equivariant} gives the finite symmetry
summands. No assumption on primality of $n$, coprimality of the degree, or
generation by tautological classes is needed.
\end{proof}

\appendix
\addtocontents{toc}{\protect\setcounter{tocdepth}{2}}

\section{Analytic foundations of harmonic transport}\label{app:analytic}

Section~\ref{sec:comparison-families} identifies the comparison maps
using two analytic inputs: the harmonic inverse at a nonzero real Hodge
parameter and continuity of preferred evaluation. We give their proofs
here, keeping the local models, prolongation estimates and tensor
compactness together. The algebraic residue families and the transport of
intersection complexes remain in the main argument.

\subsection{Harmonic metrics at nonzero real Hodge parameter}
\label{app:harmonic}

The preferred evaluation defined in \S\ref{subsec:A-transport} is
the map studied here.

The proof of Theorem~\ref{thm:harmonic-transport} needs the inverse of
preferred evaluation at each fixed real $\tau\ne0$. A filtered local
system has the \emph{zero filtration} if all its filtration weights are
zero; its local monodromy need not be trivial. A metric is
\emph{adapted} if it induces the prescribed growth lattices, and
\emph{tame} if its Higgs characteristic coefficients have logarithmic
pole orders. The fixed-parameter existence theorem does not supply the
continuity as $\tau\to0$; that is the separate role of
Proposition~\ref{prop:preferred-continuity}.

\begin{theorem}\label{thm:curve-inverse}
Fix $\tau\in\R^\times$. Let $(\mathcal E,\mathbb D^\tau)$ be the regular
filtered $\tau$-flat bundle on $(C,D)$ associated with a semisimple local
system on $\Cc$ with the zero filtration. Then it admits an adapted
 tame $\tau$-harmonic metric. Such metrics are unique up to parallel
positive automorphisms, and the resulting correspondence is compatible
with direct sums, tensor products, duals and their morphisms.
\end{theorem}

\begin{proof}
Decompose the semisimple local system into irreducible summands. Every
summand inherits the zero filtration. To see the local dictionary,
choose a logarithm $A$ of the monodromy $T$, with
$\exp(-2\pi i A)=T$ and the real parts of its eigenvalues in one
half-open interval of length one. If a multivalued horizontal frame
$u$ continues to $uT$, then $v=u\exp(A\log z)$ is single valued
and has $\tau$-connection matrix $\tau A\,dz/z$. An eigenvalue
$e/\tau$ of $A$ gives the growth weight $a=-\operatorname{Re}(e/\tau)$
for the zero filtration on horizontal sections. The nilpotent part adds
only powers of $\log z$. Replacing $v$ by $z^m v$ shifts $a$ by $-m$
and $e$ by $\tau m$. On the chosen extension the determinant residue
theorem gives $\deg E=-\sum_p\operatorname{tr} A_p$. Thus
$\deg E-\sum_{p,j}a_{pj}=0$, counting weights with multiplicity.
This is parabolic degree zero in the growth convention. Since $\tau\ne0$, an invariant holomorphic
subbundle on the punctured curve is a local subsystem; irreducibility
therefore gives stability of each summand, and semisimplicity gives
polystability. On a curve there is no parabolic second-character-number
condition to check.

The curve existence result is
\cite[\S2.6.2, Proposition~2.54]{Moc09}. Mochizuki states the
arbitrary-parameter extension of Simpson's curve correspondence there.
The proof of \cite[Theorem~5.16]{Moc09} uses the same proposition in
dimension one. To indicate the analytic input, take the explicit
logarithmic models below at the punctures and extend their metrics
across the compact complement. Their pseudo-curvature is supported
away from the punctures. By \cite[equation~(11)]{Moc09} the trace is
$1+\tau^2$ times the determinant Chern curvature. The normalized integral
of that determinant curvature is the parabolic degree by \cite[Lemma~6.1]{Sim90}, hence vanishes.
Solving the scalar Poisson equation on $C$ rescales the metric by a
smooth positive factor, harmonic in logarithm near the punctures,
and kills the trace. The fixed-parameter heat equation then solves the
trace-free equation on each stable summand, with a metric mutually
bounded with the initial one; see \cite[\S\S2.4--2.5]{Moc09} and the
curve estimates in \cite[Proposition~2.2 and \S\S5--7]{Sim88}.
In the stability estimate a putative proper limiting projection
$\pi$ satisfies $(1-\pi)\mathbb D^\tau\pi=0$. In a horizontal frame,
separating types gives $(1-\pi)d\pi=0$ since $\tau\ne0$.
The map $\pi\mapsto\operatorname{im}\pi$ to the Grassmannian has
this normal derivative on the image. The Sobolev chain rule therefore
makes the image locally constant, contrary to irreducibility.
On a curve contraction detects every $(1,1)$-form. The resulting
pseudo-curvature therefore vanishes. The uniqueness and tensor
assertions follow from \cite[Propositions~2.55--2.56 and
Corollary~5.18]{Moc09}. This proves the theorem.
\end{proof}

\begin{remark}[The fixed-parameter analytic identities]
\label{rem:lambda-analytic-identities}
For intermediate metrics, the second identity in
\cite[Lemma~2.29]{Moc09} needs correction when $\tau^2\ne1$.
Put $q=1+\tau^2$ and write $\mathbb D^\tau=d''+d'$. If
$\delta'_h,\delta''_h$ are the metric-dual operators of $d'',d'$ in
\cite[\S2.2.1, equation~(4)]{Moc09}, set
$\mathbb D_h^{\tau\star}=\delta'_h-\delta''_h$ and
$B_h=\bar\partial_h+\theta_h$. The definitions give
\[
 B_h=q^{-1}(\mathbb D^\tau-\tau\mathbb D_h^{\tau\star}),\qquad
 (\mathbb D^\tau)^2=(\mathbb D_h^{\tau\star})^2=0,
 \qquad G_h=[\mathbb D^\tau,\mathbb D_h^{\tau\star}]
       =-q^2\tau^{-1}B_h^2.
\]
On a curve $B_h^2=\bar\partial_h\theta_h$. With
$Y=\bar\partial_h\theta_h$, $Z=\partial_h\theta_h^\dagger$ and
$X=[\partial_h,\bar\partial_h]+[\theta_h,\theta_h^\dagger]$, flatness
and its metric adjoint imply
\[
 \tau^{-1}(Y+Z)=0,\qquad X=(\tau^{-1}-\tau)Z.
\]
Thus $G_h=0$ gives $Y=Z=X=0$, the required harmonic equations. In
contrast, for the trivial line connection $d''=\bar\partial$,
$d'=\tau\partial$, metric $h=e^{|z|^2}$ and $\tau=2$, one has
$X=\tfrac35\,dz\wedge d\bar z\ne0$. Hence $X=0$ cannot be imposed on
arbitrary intermediate metrics.

There is a simple test for why the parameter matters. Write an ordinary
harmonic flat connection as $D=\partial_h+\bar\partial_h+\psi+\psi^\dagger$.
For the $\tau$-connection $D^{0,1}+\tau D^{1,0}$, using the same
metric gives
\[
 \bar\partial_{h,\tau}=\bar\partial_h+(1-\tau^2)\psi^\dagger/q,
 \qquad \theta_{h,\tau}=2\tau\psi/q.
\]
Consequently $\bar\partial_{h,\tau}\theta_{h,\tau}
=2\tau(1-\tau^2)[\psi^\dagger,\psi]/q^2$. The rank-two model below
has nonzero commutator. For $\tau^2\ne1$ ordinary harmonicity thus
does not supply the metric in Theorem~\ref{thm:curve-inverse};
compare \cite[Remark~1.4]{Moc09}. In a horizontal frame of the
normalized flat connection the actual fixed-parameter equation for
the Hermitian matrix $H$ is $G_{z\bar z}=0$, where
\[
 G_{z\bar z}=-\partial_{\bar z}(H^{-1}\partial_z H)
              -\tau^2\partial_z(H^{-1}\partial_{\bar z}H).
\]
The principal term is $-qH^{-1}\partial_z\partial_{\bar z}H$.

The heat-flow calculation uses the K\"ahler adjoint, trace-curvature
and metric-change identities
\cite[Lemma~2.26, equation~(11) and \S2.2.5]{Moc09}, which follow from
the definitions above without that extra assertion. Its noncompact
estimates use the bounded-Laplacian exhaustion of
\cite[Proposition~2.2]{Sim88}: take $-\log|z|$ on the ends with
a smooth cutoff on the compact complement. The smooth compact
background metric gives finite distance to a puncture, so a proper
exhaustion cannot have bounded gradient there. This explains the use
of Simpson's curve hypotheses here. On the Neumann exhaustions of
\cite[\S6 and Lemma~7.1]{Sim88}, the puncture boundary terms vanish by
logarithmic cutoffs: for $t=-\log|z|$, take
$\chi_R=\eta(t/R)$, with $\eta=1$ on $[0,1]$ and $\eta=0$ on
$[2,\infty)$. Then $\int|d\chi_R|^2=O(R^{-1})$, so Cauchy--Schwarz
removes every boundary term with an $L^2$ transgression form.
For the weak metric limit, use quadratic lower semicontinuity. Let
$\Phi$ be a continuous nonnegative self-adjoint coefficient function
and set $T=\Phi(u)^{1/2}$. If
$u_j\to u$ almost everywhere, $D u_j\rightharpoonup D u$ in $L^2$, and
$T_j=\Phi(u_j)^{1/2}$ are uniformly bounded, then
$T_jv\to Tv$ strongly for each $L^2$ test section $v$. Thus
$T_jD u_j\rightharpoonup TD u$ and
$\|TDu\|_2\le\liminf\|T_jD u_j\|_2$. This is the needed conclusion in
the limit argument of \cite[Lemma~2.42]{Moc09}; pointwise convergence of
the derivatives is not required.

For the degree-zero principal existence input, use the degree identity
in \cite[\S5.3]{BGPM20} with the full derivative
$D''u=\bar\partial_hu+[\varphi,u]$. Here $u$ is the self-adjoint
endomorphism defining the limiting reduction. In the first expression
of equation~(5.27) of that version, the curvature must be the Chern
curvature $R(h)$ if only the $\bar\partial_hu$ energy is retained.
Equivalently, when $F(h)=R(h)-[\varphi,\sigma_h\varphi]$ is used, retain
the full $D''u$ energy; $\sigma_h$ denotes compact conjugation determined
by the reduction $h$. Its Higgs term cancels the curvature commutator,
giving the stated parabolic degree.
\end{remark}

We record the local models to fix the normalizations used in
\eqref{eq:kms}. Put $t=-\log|z|$. For a residue eigenvalue $e$ and growth
weight $a$, the zero-filtration relation makes
$\beta=(e+\tau a)/(1+\tau^2)$ purely imaginary. The rank-one Higgs model
\[
                      h=|z|^{-2a},\qquad \theta=\beta\,dz/z
\]
is harmonic. The preferred holomorphic frame of this model differs from the Higgs
frame by $\exp(-\tau\bar\beta\log|z|^2)$, of modulus one, and its
$\tau$-residue is $\beta-\tau a-\tau^2\bar\beta=e$.
For the nilpotent part take
\[
 N=\begin{pmatrix}0&1\\0&0\end{pmatrix},\qquad
 \theta=\frac N2\frac{dz}{z},\qquad
 h=\begin{pmatrix}t^{-1}&0\\0&t\end{pmatrix}.
\]
Direct calculation gives $R(h)+[\theta,\theta^\dagger]=0$.
In the preferred holomorphic frame
$g=\left(\begin{smallmatrix}1&0\\-\tau/t&1\end{smallmatrix}\right)$
the connection is $\tau\partial+(N/2)dz/z$. Symmetric powers give the
nilpotent Jordan blocks; tensor products with the rank-one models and
direct sums give all regular local types.

The uniqueness assertion applies to all adapted tame metrics inducing
the same filtered bundle and its dual, not only to a preselected bounded
class. One can see the boundary point directly. In a common filtered
frame, adaptation for the bundle and its dual makes both relative
traces subpolynomial in $|z|^{-1}$. The logarithm of the determinant
ratio is harmonic and $o(\log|z|^{-1})$ at each puncture, hence
extends harmonically across $D$ and is constant on $C$. After scaling
one metric by a constant, the determinants agree. The log-trace comparison
\[
 u=\log\frac{\operatorname{tr}(h_1^{-1}h_2)}{\operatorname{rank}\mathcal E}
 \ge0,\qquad \Delta^\tau u\le0
\]
uses the positive-Laplacian convention and satisfies
$u=o(\log|z|^{-1})$ at each puncture. Comparing it on a punctured disk
with $\varepsilon\log(r_0/|z|)$ and then letting
$\varepsilon\downarrow0$ bounds $u$ by its value on the outer boundary
circle. Hence $u$ extends as a bounded weak subsolution across the puncture.
The compact-curve maximum principle and the equality case of the
metric-change identity make the relative endomorphism parallel; see
\cite[\S\S2.2, 2.4 and Propositions~2.55--2.56]{Moc09}.
This also explains why applying the inverse in different representations
of $G$ preserves the defining tensors and their relations.


\subsection{Harmonic compactness and continuity}\label{subsec:A-harmonic}

A metric on a bundle over a punctured disk determines \emph{growth
lattices}: holomorphic sections $s$ with $|s|=O(|z|^{-c-\varepsilon})$
for all $\varepsilon>0$. A \emph{growth cut} is a value of $c$ between
consecutive weights; it specifies one holomorphic extension across the
puncture. Convergence on the punctured curve must therefore be supplemented
by convergence of full-rank frames of these lattices. We use the numbering
of the Kyoto preprint \cite{Moc05} and of arXiv v2 of \cite{Moc07}.

\begin{lemma}[Pointwise comparison on the punctured disk]
\label{lem:pointwise-kms}
Fix a tame harmonic bundle and $\tau\in\R$. For every joint Higgs
label $(a,\beta)$, the centered nilpotent weight-graded dimensions
agree with those of the corresponding joint label
$(p_\tau,e_\tau)$ of its preferred $\tau$-connection.
\end{lemma}

This is \cite[Corollary~7.16]{Moc07}. Within each joint label the
number of blocks of length $s$ is $d_{s-1}-d_{s+1}$, where $d_j$ is
the dimension in degree $j$.

\begin{proof}
Assume $\tau\ne0$, the case $\tau=0$ being immediate. Build the model
of the Higgs labels and Jordan blocks from the lines and symmetric
powers above. For an arbitrary line residue
$\beta$, the preferred frame $\exp(-\tau\bar\beta\log|z|^2)v$
has the general label $(p_\tau,e_\tau)$ given before Theorem~\ref{thm:NAH}; a change of strip by $z^m$
sends $(b,e)$ to $(b-m,e+\tau m)$. It suffices to compare
this model with the actual bundle at $\tau$.

Put $r=|z|$, $t=-\log r$. At each parameter choose holomorphic
columns $v_j$ adapted to the growth weights $b_j\in(-1,0]$, graded
residue eigenvalues $e_j$ and nilpotent weights $2k_j$ measured there.
The normalized columns $r^{b_j}t^{-k_j}v_j$ have bounded Gram matrix
and bounded inverse. The single-bundle construction in
Lemma~\ref{lem:frames} gives the growth lattices and their duals from
the Poincar\'e curvature bound and flat determinant. For the norm
estimate recall \cite[\S4.3.3, Corollaries~4.4--4.5]{Moc02} and
\cite[\S\S7.2.1, 7.3.1]{Moc07}.
For those measured data take a line with Higgs residue
$\beta_j=(e_j+\tau b_j)/(1+\tau^2)$ and Higgs growth
$b_j-2\tau\operatorname{Re}\beta_j$, and a nilpotent block of the
measured length. In a filtered identification the difference of the logarithmic
operators vanishes on each growth-graded space. In the model norm
the difference is $O(r^{-1+\epsilon})$ for some $\epsilon>0$;
one derivative costs at most a factor $r^{-1}$ after shrinking
$\epsilon$. The model metric for the actual operator consequently
has pseudo-curvature error $O(r^{-2+\epsilon})$. The log-trace comparison
in both directions, with barrier $C r^\epsilon$ for this error and
$\eta\log(1/r)$ for the inner boundary, gives mutual boundedness
as $\eta\downarrow0$. The inner boundary estimate uses adaptation
for the bundle and its dual. This obtains the frame bounds before
comparing the data at different parameters.

Choose a Higgs-side filtered holomorphic isomorphism $I$ from the
first model to the actual bundle intertwining their graded residues.
The frame bounds at zero give bounded $I,I^{-1}$. Write
$Q=\theta I-I\theta_0$, $M=\theta^\dagger I-I\theta_0^\dagger$;
subscript zero here denotes the model. Since $zQ/dz$ strictly lowers
the growth filtration, its finitely many jumps give
$|Q|=O(r^{-1+\epsilon})$. With $D'_h$ the Chern $(1,0)$ derivative
on the Higgs Hom bundle, the Bochner identity for the holomorphic
section $I$ gives, for the positive Euclidean Laplacian,
\[
 \Delta|I|^2\le C|Q|^2-c\bigl(|D'_h I|^2+|M|^2\bigr),\qquad c>0 .
\]
Indeed the harmonic curvature contributes the difference of the two
Higgs squares. Adding $C' r^{2\epsilon}$ absorbs the positive error.
The resulting bounded subharmonic function extends across zero, and
the disk Green kernel gives finite integral of both negative-square
terms against $t\,dA$. This supplies the estimate of
\cite[Lemma~7.7]{Sim90} while retaining the graded-residue error $Q$.
The two deformed Hom derivatives $E''=\tau M$ and
$E'=\tau D'_h I+Q$ therefore both belong to $L^2(t\,dA)$.

Use adapted holomorphic frames at $\tau$ on both sides. For the
entry from model column $i$ to actual column $j$ put
\[
 u=I_{ji}\,e^{-\gamma t}t^\delta,\qquad
 \gamma=b_{0i}-b_j,\qquad \delta=k_j-k_{0i}.
\]
These are entries of a bounded matrix with bounded inverse.
In polar coordinates $(r,\phi)$ the equation for $E''$ gives
\[
 (-\partial_t+i\partial_\phi-\gamma+\delta/t)u
       =f,\qquad \int_{t>T}t|f|^2\,dt\,d\phi<\infty .
\]
On the Fourier modes with $n+\gamma\ne0$, the self-adjoint operator
$A=i\partial_\phi-\gamma$ has a spectral gap at zero. For compactly
supported $v$, integration by parts gives
$\|(-\partial_t+A)v\|_2^2=\|\partial_t v\|_2^2+\|Av\|_2^2$.
Apply this to the projection of $u$ times a cutoff, equal to one
on $[T+1,R]$ and vanishing past $2R$, with outer derivative $O(R^{-1})$.
The right-hand forcing, $\delta u/t$ and the cutoff errors are
bounded in $L^2(dt\,d\phi)$ independently of $R$. The spectral gap
therefore puts these modes in $L^2(dt\,d\phi)$ on the tail.
Since $|\gamma|<1$, the only remaining mode has $n=\gamma=0$.
For its circular mean $u_0$, put $y=\log t$. If $\delta\ne0$,
$(-\partial_y+\delta)u_0=t f_0\in L^2(dy)$. The same cutoff identity
gives $u_0\in L^2(dt/t)$.

It remains to separate unequal eigenvalues when $\gamma=\delta=0$.
In the normalized columns the matrix of the logarithmic connection
coefficient is $\operatorname{diag}(e_j)+O(t^{-1})$: on a single
growth grade the nilpotent lowers $2k$ by two, strictly lower growths
gain a positive power of $r$, and a holomorphic term vanishing at zero
from column $i$ to row $j$ gains $r^{1+b_i-b_j}$ up to logarithms. Here only the
coefficient matrix is conjugated; derivatives still act in the
holomorphic frames. In the normalized entry of $zE'_z-\tau\bar z
E''_{\bar z}$ the radial derivatives cancel, leaving
\[
 -i\tau\partial_\phi u+(e_j-e_{0i})u+O(t^{-1})
                      \ \in L^2(dt\,d\phi/t).
\]
Taking the circular mean shows $u_0\in L^2(dt/t)$ whenever
$e_j\ne e_{0i}$. Thus every entry between unequal triples has finite
squared integral against $dt\,d\phi/t$. This measure has infinite mass
on every tail, so one sequence approaching the puncture makes all
such entries tend to zero. A subsequential limit of the matrices is
invertible and block diagonal in the triples $(b,e,k)$. Each block is
square, giving the asserted equality of dimensions.
\end{proof}

In each associated vector bundle, for $\tau\ne0$ this also identifies
the unipotent monodromy blocks. Choose a logarithmic lattice in this
representation whose growth representatives $a$ lie in one half-open
interval of length one. Since $\beta\in i\R$, the connection exponents
$-a+(1+\tau^2)\beta/\tau$ have no differences in $\Z\setminus\{0\}$.
Write $\tau^{-1}\nabla_\tau=d+A(z)\,dz/z$, with
$A(z)=A_0+\sum_{j>0}A_jz^j$. A holomorphic gauge $P(0)=1$ satisfies
\[
 (j+\ad A_0)P_j=-\sum_{i=1}^j A_iP_{j-i}.
\]
These operators are invertible with inverse norm $O(j^{-1})$;
the resulting scalar majorant gives a convergent gauge.
The monodromy is therefore conjugate to $\exp(-2\pi iA_0)$.
For a nilpotent block $N$ and $c\ne0$, $\exp(cN)-1$ is $N$ times
an invertible polynomial in $N$, so all kernel dimensions and Jordan
lengths agree. Nonresonance is essential for this conclusion: an
arbitrary logarithmic lattice need not have monodromy determined by
its residue alone.

At the Betti endpoint the principal orbit has a direct description.
For the horizontal equation Boalch's parameters are
$\theta=-a$, $\tau_{\rm B}=a$, $\sigma=-2\beta$ and
$\phi=\theta+\tau_{\rm B}=0$. Thus the parabolic for the Betti
filtration is $G$ itself. If $N_1$ is the nilpotent part of the
graded affine residue of the connection operator, the semisimple and
unipotent parts of monodromy are conjugate as a pair to
$\exp(2\pi i(a-2\beta))$ and $\exp(-2\pi iN_1)$
\cite[\S3, Lemma~5 and Theorem~6]{Boa11}. A Jacobson--Morozov
cocharacter in the connected centralizer conjugates $N_1$ to every
nonzero scalar multiple of $N_1$. This identifies the principal
nilpotent orbit even when a common strip for the vector exponents
requires a lattice change outside $G$.

The passage from vector bundles to principal bundles requires tensor
compatibility at the actual real growth weights. The following local
calculation supplies that compatibility without a restriction on the
global form of $G$.

\begin{lemma}[Tensor products of growth lattices]\label{lem:real-growth-tensor}
Let $V,W$ be tame Hermitian holomorphic bundles on a punctured disk, with
adapted holomorphic growth frames and the corresponding dual growth
estimates. Define
\[
 \mathcal P_c(V)=\{s:|s|=O(|z|^{-c-\epsilon})
                    \text{ for every }\epsilon>0\}.
\]
Then the image of the parabolic tensor product is
\[
 \mathcal P_c(V\otimes W)
 =\sum_{u+v=c}\operatorname{im}
       \bigl(\mathcal P_u(V)\otimes\mathcal P_v(W)
                      \longrightarrow j_*(V\otimes W)\bigr),
\]
where $j$ is the punctured-disk inclusion. The prolongation of $V^*$ is
the filtered dual: a section $f$ belongs to $\mathcal P_c(V^*)$ precisely
when $f(\mathcal P_d(V))\subseteq\mathcal P_{c+d}(\mathcal O)$ for every
$d\in\R$.
\end{lemma}

\begin{proof}
For a frame $v_i$ of weights $a_i$, the adapted and dual estimates give
\[
 \mathcal P_c(V)=\bigoplus_i z^{\lceil a_i-c\rceil}\mathcal O\,v_i.
\]
The necessary frames and inverse estimates are supplied by
\cite[Proposition~8.20 and Lemmas~8.28--8.29]{Moc05}; logarithmic factors
are absorbed by the arbitrary $\epsilon>0$ in the definition. If $b_j$
are the weights of $W$, the least exponent in the tensor convolution is
\[
 \min_{u+v=c}\bigl(\lceil a_i-u\rceil+\lceil b_j-v\rceil\bigr)
       =\lceil a_i+b_j-c\rceil,
\]
with equality at $u=a_i$. This proves the tensor formula. Applying the
same calculation to the dual frame gives weights $-a_i$ and the stated
filtered-dual formula. To track residues as well, take a representative
$\kappa(x)\in(-1,0]$ of $x$ modulo $\Z$ and write $n=x-\kappa(x)$.
For raw tensor weight $x=b_1+b_2$ and raw residue $e=e_1+e_2$, the
frame of representative weight $\kappa(x)$ is $z^n$ times the raw
frame, so the residue is $e+\tau n$. For a dual frame take $x=-b$
and $e=-e_1$, for an original label $(b,e_1)$.
Both rules follow from $\nabla_\tau(z^n v)=z^n\nabla_\tau v+
\tau n z^n v\,dz/z$. Every defining tensor and its dual are therefore
retained in the prolongation. The argument uses the weights in actual
representations of $G$ and is independent of its global form. This is
the representationwise version of \cite[Proposition~5.1]{BBP17}, whose
stated hypotheses are semisimple simply connected $G$; rational
reconstruction is used only after passing to a rational weight chart.
\end{proof}

\begin{lemma}[Uniform prolongation frames]\label{lem:frames}
Let Hermitian holomorphic bundles of fixed rank on a punctured disk
have Chern connections converging smoothly in unitary gauges on compact annuli, a
uniform Poincar\'e curvature bound, flat determinant metrics, and a fixed
finite set of growth weights. For each selected growth cut, choose the
holomorphic determinant frame of the corresponding exponent whose norm
is a pure power of $|z|$, normalized in norm and phase at one fixed
outer point. For finitely many cuts at positive distance from the
weights modulo $\Z$, a common shrinking and a subsequence admit
holomorphic frames of the growth lattices converging on compact annuli.
Their wedge coefficients relative to the chosen determinant frames are
uniformly bounded above and below. The limits frame the full limiting
lattices.
\end{lemma}

\begin{proof}
We give the family estimate explicitly. For the corresponding
single-bundle construction see \cite[Proposition~8.10 and
\S\S8.3.4--8.3.5]{Moc05}. Write $z=\exp(-t+i\theta)$. In a radial
unitary frame the angular connection coefficient satisfies
$\partial_t A_\theta=F_{t\theta}=O(t^{-2})$. Its limit on the circle
therefore has error $O(t^{-1})$, uniformly in the family. By parallel
transport and diagonalization of its unitary holonomy, a periodic
angular gauge puts the limiting connection in the form
$i\Gamma\,d\theta$, with $\Gamma=\operatorname{diag}(\gamma_1,\ldots,\gamma_r)$
real. Choose the representatives in a strip of width one. The
holomorphic equation is $(D+B)F=0$, where
\[
 D=\partial_t-i\partial_\theta+\Gamma,\qquad |B|\le C/t.
\]
Flatness of the determinant gives $\tr B=0$. In particular the
wedge of the model columns $f_j=e^{-\gamma_j(t-T)}v_j$, in this
unitary basis $v_j$, is an exact holomorphic determinant frame.

Choose $\eta>0$ with $r\eta<1$ and
$\gamma_l-\gamma_j+\eta\notin\Z$ for every $j,l$. On the full cylinder
$D-\gamma_j+\eta$ has Fourier multipliers
$i\xi+n+\gamma_l-\gamma_j+\eta$, $n\in\Z$. Division by these
multipliers gives a bounded inverse from $L^2$ to $W^{1,2}$.
Extend $B$ by zero for $t<T$. Solve
\[
 (D-\gamma_j+\eta+B)Y_j
       =-\mathbf 1_{t\ge T}e^{-\eta(t-T)}Bv_j .
\]
The right side has $L^2$ norm $O_\eta(T^{-1})$. For large $T$ the
inverse and a Neumann series give the same bound in $W^{1,2}$.
Sobolev embedding on unit strips gives uniform $L^p$ bounds, for
each finite $p$. The interior $W^{1,p}$ estimate for the
Cauchy--Riemann operator, with $p>2$, then gives
$\|Y_j\|_\infty=O_\eta(T^{-1})$ on $t\ge T+1$.
Consequently $F_j=f_j+e^{-(\gamma_j-\eta)(t-T)}Y_j$ is holomorphic
there. The ratio of $\bigwedge_j F_j$ to the model determinant is
holomorphic with growth $O(|z|^{-r\eta})$ and so extends across zero.
On the circle $t=T+1$ the ratio differs from $1$ by $O_\eta(T^{-1})$.
The maximum principle gives upper and lower wedge bounds throughout
the smaller disk. Cofactors now bound the dual basis $F^j$:
\[
 |F_j|\le C|z|^{\gamma_j-\eta},\qquad
 |F^j|\le C|z|^{-\gamma_j-(r-1)\eta}.
\]
For a cut $c$ away from $-\gamma_j+\Z$, take $r\eta$ smaller than
the distance to every such jump. The dual bound forces the $j$th
coefficient of a section of $\mathcal P_c$ to have Laurent order at
least $\lceil-\gamma_j-c\rceil$. The column bound gives the reverse
inclusion. Hence the columns $z^{\lceil-\gamma_j-c\rceil}F_j$ form a
basis. Applied first to a single bundle on both sides of each possible
jump, this calculation identifies all growth weights and their
multiplicities with $-\gamma_j$ modulo $\Z$. The fixed growth data
now give a common spectral margin and a common $T$ for the sequence,
after fixing multiplicities by a subsequence. One choice of $\eta$
works for all the selected cuts.

In the original convergent gauges the columns have a subsequence
converging on compact annuli by interior elliptic estimates.
The same column and dual bounds prove that the limits generate each
full limiting lattice. The integer powers of $z$ shift the determinant
exponent by their sum; the pure-power normalization at the outer point
only changes the bounds by a fixed factor. Compare
\cite[Lemmas~8.26--8.29 and~9.8--9.9]{Moc05}.
\end{proof}

\begin{proposition}[Principal harmonic compactness]\label{prop:tensor-compactness}
Let polystable tame principal Higgs objects $x_i$ converge to $x$ in a
separated coarse moduli space, with fixed growth weights and semisimple
residue labels (nilpotent residues may degenerate). Choose their adapted
harmonic reductions modulo parallel positive automorphisms. After passing
to a subsequence and choosing unitary identifications of their underlying
Hermitian bundles on $\Cc$, the harmonic bundles converge there, and
their principal prolongations converge to an object representing $x$.
\end{proposition}

\begin{proof}
Realize $G\subset\GL(V)$ as the stabilizer of a finite tuple $s$ of mixed
tensors with closed $\GL(V)$-orbit (embed $\GL(V)/G$ equivariantly in a
representation). The characteristic polynomials of the associated Higgs
fields converge. Here is the curve compactness argument. In a unitary
Schur basis $T=\operatorname{diag}(\lambda_i)+N$, with $N$ strictly upper
triangular and $H=\operatorname{diag}(r-1,\ldots,0)$,
\[
 \tr(H[T,T^*])=\sum_{i<j}(j-i)|N_{ij}|^2\ge\|N\|^2.
\]
The Bochner inequality for the holomorphic Higgs coefficient, combined
with this estimate and the barrier $c/(R^2-|z|^2)^2$, bounds its norm on
each smaller coordinate disk in terms of its characteristic coefficients.
The Hitchin equation then bounds the curvature. In unitary radial gauge,
\[
 A_k(x)=\int_0^1t x^jF_{jk}(tx)\,dt,
\]
so $A$ is bounded. Fix $p>2$. The equation
$\bar\partial\phi=-[A^{0,1},\phi]$ gives interior $W^{1,p}$ bounds for
$\phi$. The curvature equation and the radial formula then bound $A$
in $W^{1,p}$: the differentiated integral has factor $t^{2-2/p}$, which
is integrable for $p>1$. Iteration gives smooth compactness. Unitary
transition functions satisfy $dg=gA'-Ag$ and have values in the compact
group $U(V)$, so the local limits patch on $\Cc$.

Near the punctures apply \cite[Remark~7.7 and Proposition~7.8]{Moc05}.
The curvature appearing in the hypothesis there is
$R(h)+[\theta,\theta^{\dagger_h}]$, which vanishes by the harmonic
equation. The fixed spectral labels and bounded characteristic
coefficients then give the following uniform estimates as the
nilpotent residue varies.
Write $\theta=f_0\,dz/z$ and $t=-\log|z|$. For each root $a$ of the
limiting characteristic polynomial $\det(s-f_0)|_{z=0}$, let $R$ act by
the constant $a$ on its spectral cluster $E_a$. Let $R^\perp$ be the normal
endomorphism obtained by orthogonalizing the ordered cluster filtration.
The cited estimates say
\[
 f_0=R^\perp+Q,\qquad R=R^\perp+S,\qquad
 Q=O(t^{-1}),\quad S=O(t^{-2}).
\]
Since $[f_0,R]=0$, we have
$[Q,R^\perp]=-[f_0,S]=O(t^{-2})$.
In a unitary eigenbasis of $R^\perp$, conjugating the eigenvalue gaps
preserves their absolute values, so
$[Q,(R^\perp)^*]=O(t^{-2})$ as well. Consequently
$[f_0,f_0^*]=O(t^{-2})$, which is the required Poincar\'e curvature bound.
If all residue eigenvalues are zero, first tensor the associated vector
bundle locally by the harmonic line with Higgs field $dz/z$ and constant
metric. This makes the spectral constant in that proposition positive
and leaves the curvature and commutator unchanged.
Choose the pure-power determinant frames and common growth cuts of
Lemma~\ref{lem:frames}. That lemma gives convergence of the full
prolongations. The characteristic coefficients retain their logarithmic
pole bounds in the limit, so the limiting harmonic bundle is tame. Each defining tensor is parallel and annihilated
by the Higgs field and its adjoint. Choose the Hermitian form on $V$
invariant under a maximal compact subgroup of $G$, and use the induced
Hermitian reductions. The tensor values in these unitary frames then lie
in the fixed compact orbit $U(V)s$. They converge on compact
annuli after a subsequence. To see what happens at the puncture,
take the simultaneous frames $F_j$ above and put $b_j=-\gamma_j$.
For a basis monomial in a mixed tensor let $A$ be the sum of the
$b_j$ in its vector factors minus the sum in its covector factors.
Pairing the bounded tensor with the dual monomial bounds its scalar
coefficient by $C|z|^{A-dr\eta}$, where $d$ bounds the tensor
length and $r=\dim V$. Choose $dr\eta$ less than the distance from
each of the finitely many $A$ to the largest integer strictly below
$A$. The Laurent order of the coefficient is then at least
$\lceil A\rceil$, including at an integral $A$.

Choose the coarse rational polarization $q$ on the face retaining
the integral coincidences and the order of fractional parts in the
finite tensor presentation and its dual, including the trivial
weight. Impose also $\lceil Q\rceil=\lceil A\rceil$ for each signed
sum $Q$ of the corresponding representatives $q_j$. These conditions
hold on a rational affine face through $a$ as in
Lemma~\ref{lem:rational-chart}. On $z=u^N$, with $Nq_j\in\Z$,
the regularized columns are $u^{Nq_j}F_j$. A tensor coefficient
$c(z)$ becomes $u^{-NQ}c(u^N)$, of order at least
$N(\lceil A\rceil-Q)\ge0$. Cauchy's formula thus gives convergence
on the disk in these regular frames. The limiting tuple at the
puncture belongs to the closed orbit $\GL(V)s$. Local holomorphic
sections of $\GL(V)\to\GL(V)/G$ recover principal frames, and the
regularized columns retain the specified inertia characters. This
choice of $q$ concerns only the coarse chart.

On the compact covering curve, patch the local smooth near-identifications
of these bundles with the limit using a partition of unity, and
average over the finite deck action. They are invertible for all
large indices, since in the finite cover of smaller charts every
summand tends to the same local identification. In these
identifications the holomorphic structures converge smoothly.
Twist by a line sufficiently positive that the limiting faithful
bundle is generated with vanishing first cohomology. If $D_i\to D$
are the twisted Dolbeault operators, choose a bounded right inverse
$T$ to $D:W^{1,p}\to L^p$, $p>2$. Then
$T_i=T(1+(D_i-D)T)^{-1}$ is a bounded right inverse for all large
$i$, and the projections $1-T_iD_i$ give bases of holomorphic
sections converging to a chosen limiting basis, by constancy of
the index and elliptic regularity. Choose the twist sufficiently positive for the finite multiplication
targets as well. In the resulting Grassmannian presentation of Quot, the kernels of the finite
section multiplication maps therefore converge, because their
ranks are locally constant for that twist. The logarithmic
operators, tensors, flags and deck arrows converge in these frames.
Thus the framed points, and hence their good-moduli points, converge
in the analytic topology.
The limit represents a point which is also a limit of $x_i$, hence equals $x$ by
separatedness. Stability of associated vector bundles is never used.
\end{proof}

\begin{proposition}[Continuity of preferred evaluation]
\label{prop:preferred-continuity}
On a fixed generic full spectral principal Higgs space, preferred
evaluation is jointly continuous in the object and in $\tau$ in a bounded
interval (including $\tau=0$), with values in the logarithmic coarse
moduli space with the given parahoric prolongations.
\end{proposition}

\begin{proof}
Let $x_i\to x$, $\tau_i\to\tau$. By Proposition~\ref{prop:tensor-compactness}
the operators converge on compact sets. The Chern connection of
$\bar\partial_{i,\tau_i}$ with metric $h_i$ has curvature
$-(1+\tau_i^2)[\Phi_i,\Phi_i^\dagger]$, uniformly bounded near the
punctures by \cite[Proposition~7.8]{Moc05}, with trace zero; normalizing
the flat determinant frame at one point gives the determinant bounds, and
\eqref{eq:kms} fixes the growth cuts. Lemma~\ref{lem:frames} gives
convergent full-rank frames in which $\nabla_{i,\tau_i}=\tau_id+M_i(z)dz/z$;
the matrices converge on an annulus and, by Cauchy's formula, on the disk,
including residues. The tensor argument keeps the principal group and
inertia. Distinct semisimple labels in one angular block stay separated,
since their differences are multiplied by $1+\tau_i^2$. Every subsequence
thus has the prescribed limit, and the moduli spaces are metrizable.
\end{proof}

\section{Reductive-group and Bruhat calculations}\label{app:calculations}
\label{app:auxiliary}

This appendix collects the algebraic calculations used in the main
proof: structural facts, nilpotent-cone normality and low-degree Springer
cohomology, followed by stabilized Bruhat charts and their Goldman-form
comparison. The main text states the consequences at their points of use.

\subsection{Two structural facts about reductive groups}\label{app:structural}

The two facts used in \S\ref{sec:background} admit short proofs over
$\C$. We include the arguments to make the closed-orbit and
connected-centralizer inputs explicit. The first proof uses Cartan
decomposition; the second uses the finite and affine reflection chambers
of the root system.

\begin{lemma}[The irreducible case of the closed-orbit criterion]
\label{lem:irreducible-closed-direct}
Let $G$ be connected complex reductive and $x=(x_1,\ldots,x_m)\in G^m$.
If the $x_i$ lie in no common proper parabolic subgroup, their simultaneous
conjugacy orbit is closed and their common centralizer is finite modulo
$Z(G)$. This is the irreducible case of Richardson's criterion
\cite[Theorem~16.4]{Ric88}.
\end{lemma}

\begin{proof}
We give the complex-group argument, using the compact-form and Cartan
decomposition theorems \cite[Theorems~6.11, 6.31 and~7.39]{Kna23},
with the central torus treated separately.
Choose a maximal compact subgroup of $G_\der$ and a faithful representation
of $G$ whose Hermitian form is invariant under it. Modulo the center,
every sequence $g_n\in G$ can be written
$g_n=k_n\exp(H_n)\ell_n$, with $k_n,\ell_n$ in that compact subgroup
and $H_n$ in a fixed real positive chamber for $G_\der$. Suppose that
$g_nx_ig_n^{-1}$ converges for every $i$. Pass to a subsequence on which
$k_n,\ell_n$ converge, say $\ell_n\to\ell$.

If $H_n$ is unbounded, pass again to a subsequence with
$H_n/\lVert H_n\rVert\to H\ne0$. In the weight decomposition of the
chosen representation, the block from weight $\nu$ to weight $\mu$
is multiplied by $\exp((\mu-\nu)(H_n))$. Boundedness therefore forces
the corresponding block of $\ell x_i\ell^{-1}$ to vanish whenever
$(\mu-\nu)(H)>0$. Choose an integral cocharacter $\lambda$ having
the sign pattern of $-H$ on all these weight differences and on the roots.
Such a choice exists: the equalities and strict inequalities define a
rational face, and a rational point in it can be multiplied into
$X_*(T)$. Thus every limit
$\lim_{t\to0}\lambda(t)\ell x_i\ell^{-1}\lambda(t)^{-1}$ exists.
The limits remain in $G$: their determinants are unchanged, and $G$ is
closed in the chosen general linear group. Hence all $\ell x_i\ell^{-1}$
lie in $P_\lambda$. This parabolic is proper, since $H$ lies in the
semisimple direction and the preserved root signs make $\lambda$
noncentral. This contradicts the hypothesis.

Consequently $H_n$ is bounded, and $g_n$ has a convergent subsequence
modulo $Z(G)$. The limit of this subsequence carries $x$ to the given limiting tuple. The
orbit is therefore closed in the complex topology, hence Zariski closed:
an algebraic orbit is locally closed and dense in its Zariski closure,
so the two closures agree. Apply the same argument to sequences in the
common centralizer. The quotient of the centralizer by $Z(G)$ is a compact affine complex
algebraic group and is therefore finite, since regular functions are
constant on each compact connected complex component.
\end{proof}

\begin{lemma}[Closed monodromy orbits]
\label{lem:general-closed-direct}
For $G$ connected complex reductive, the orbit of $x\in G^m$ is closed
if and only if every parabolic containing its entries has a Levi subgroup
containing them.
\end{lemma}
\begin{proof}
Choose a faithful matrix representation unitary on a compact real form $K$.
A closed orbit contains a minimum $y$ of the sum of the squared
Hilbert--Schmidt norms of the entries: determinants are constant and
nonzero, so the orbit is closed even in the ambient matrix space.
Given a parabolic $P$ containing $y$, choose a cocharacter $\lambda$
of a $K$-compatible torus defining $P=P(\lambda)$. Such a choice is
possible since $K$ is transitive on each $G/P$. The real positive
part of $\lambda$ acts by Hermitian matrices. In the induced orthogonal
weight decomposition write $y=\sum_d y_d$. Membership in $P^m$ gives
$y_d=0$ for $d<0$. The derivative at zero of
$\sum_d e^{2td}\|y_d\|^2$ is zero by minimality, so all $y_d$ with
$d>0$ vanish as well. Thus every entry belongs to the Levi
$Z_G(\lambda)$.

Conversely take a minimal parabolic among those containing $x$ (allow
$G$ itself), and write it as $P=L U$ with $x\in L^m$, by hypothesis.
A proper parabolic $Q$ of $L$ containing $x$ would give the smaller
parabolic $Q U$ in $G$. Hence $Lx$ is closed by
Lemma~\ref{lem:irreducible-closed-direct}.
The unipotent group $U$ has closed orbits on the matrix tuple space.
Indeed, in a unitriangular basis ordered so each increment uses earlier
coordinates, let $j$ be the first coordinate of $uv$ varying with $u$.
Then $c(u)=(uv-v)_j$ is an additive homomorphism, since all earlier
coordinates of $uv-v$ vanish. Choose $X\in\operatorname{Lie}U$ with
$dc(X)=1$. For a convergent sequence $u_nv$, the values $c(u_n)$
converge, and multiplying by $\exp(-c(u_n)X)$ reduces to a convergent
sequence in $(\ker c)v$. Dimension induction gives a limit in the same
orbit. If every coordinate is constant there is nothing to prove.

The Levi projection $P^m\to L^m$ sends $Px$ onto $Lx$. For a convergent
sequence in $Px$ the projected limit lies in $Lx$. The smooth orbit
map $L\to Lx$ has analytic local sections, so one may conjugate by
convergent Levi representatives of the projections to reduce the sequence
to $Ux$, which is closed. Thus $Px$ is closed in $P^m$, hence in
$G^m$. Finally $G=KP$ and $K$ is compact, so $Gx=K(Px)$ is closed
in the complex topology. The closure of a complex algebraic orbit
agrees with its Zariski closure.
\end{proof}

\begin{lemma}[The complex semisimple-element connectedness theorem]
\label{lem:semisimple-centralizer-direct}
If $G_\der$ is simply connected, the centralizer in $G$ of every
semisimple element is connected.
\end{lemma}

\begin{proof}
We give a root-system proof of the complex case of Steinberg's theorem
used here \cite{Ste75}. For the earlier compact-group proof, see
\cite[Theorem~3.4 and Corollary~3.5]{Bor61}.
First let $G$ be semisimple and simply connected, choose a
maximal torus $T$ containing the element $s$, and put
$V=X_*(T)\otimes_\Z\R$. Write $s=\exp(2\pi i(a+ib))$ with $a,b\in V$.
Simple connectedness says $X_*(T)=Q^\vee$, the coroot lattice.
Put $C=Z_G(s)$ and $\Phi_s=\{\alpha:\alpha(s)=1\}$. The root-space
calculation gives $C^0$ reductive with maximal torus $T$ and root system
$\Phi_s$. Conjugacy of maximal tori in $C^0$ then gives
\[
 \pi_0(C)=\{w\in\WG:w(s)=s\}/W(\Phi_s).
\]
Indeed, every component of $C$ has a representative normalizing $T$.

Consider the affine Weyl group $W_{\rm aff}=Q^\vee\rtimes\WG$ on $V$.
An \emph{alcove} is a connected component of the complement of the
affine root hyperplanes $\alpha(v)=k$, $k\in\Z$. The coroot affine
Weyl group acts simply transitively on alcoves
\cite[\S4, no.~5, Satz~15, p.~375]{Sti42}; this freeness is not asserted
for the extended group formed with the larger coweight lattice.
The condition $w(s)=s$ means $wb=b$ and $wa-a\in Q^\vee$; hence $w$
corresponds uniquely to the affine transformation $v\mapsto wv+a-wa$
which fixes $a$ and whose linear part fixes $b$. The stabilizer of $a$
in $W_{\rm aff}$ is generated by the affine reflections through $a$.
After translating $a$ to the origin it is a finite reflection group,
and its stabilizer of $b$ is generated by the reflections whose
linear hyperplanes contain $b$. These assertions follow from the
chamber decomposition: chambers incident to the point are connected
by galleries crossing only walls through that point, and the chamber
action is free. The same argument applies to affine alcoves, using a
small ball about the point and local finiteness of the reflecting walls.
It follows that $\{w:w(s)=s\}$ is generated by the reflections $s_\alpha$
with $\alpha(a)\in\Z$ and $\alpha(b)=0$. These conditions are precisely
$\alpha(s)=1$, so this group is $W(\Phi_s)$ and $C$ is connected.

For connected reductive $G$ with simply connected derived group, write
$s=zh$ with $z\in Z(G)^0$ and $h\in G_\der$. The commuting factors are
semisimple, and $Z_G(s)=Z(G)^0Z_{G_\der}(h)$. This is the image of a
product of connected groups, so is connected.
\end{proof}

\subsection{Finite deck actions}\label{subsec:finite-deck-actions}

The following elementary fact about affine GIT quotients supplies the
proper forgetful map needed in Proposition~\ref{prop:pg-master}.

\begin{lemma}[Forgetting a finite action]\label{lem:finite-deck-forgetting}
Let $H=\mathrm{GL}(W)$ act on an affine finite-type $\C$-scheme $R$,
and let a finite group $\Gamma$ act commuting with $H$.
For a homomorphism $\sigma:\Gamma\to H$, let $R_\sigma$ be the
fixed scheme of $\gamma\mapsto(\sigma(\gamma),\gamma)$ and
$J=Z_H(\sigma(\Gamma))$. Then
$R_\sigma/\!/J\to R/\!/H$ is finite.
\end{lemma}
\begin{proof}
Embed $R$ equivariantly as a closed subscheme of an
$H\times\Gamma$-module $V$ by choosing a finite-dimensional invariant
space of algebra generators. Choose metrics on $W,V$ with
$\sigma(\Gamma)\subset U(W)$ and $U(W)\times\Gamma$ unitary on $V$.
In the fixed linear subspace $V_\sigma$, the $J$-nullcone is the
intersection with the $H$-nullcone. To prove the nontrivial inclusion,
take $v$ in this intersection and minimize norm in the closure of
$Jv$, obtaining $y\in V_\sigma$. For Hermitian $A\in\mathrm{End}(W)$
the derivative at zero of $\|\exp(tA)y\|^2$ is real-linear in $A$
and invariant under $A\mapsto\sigma(\gamma)A\sigma(\gamma)^{-1}$.
It vanishes on the Hermitian commutant by minimality. Averaging gives
its vanishing for every $A$. By the spectral theorem the same function
of $t$ is convex. Polar decomposition now gives
$\|hy\|\ge\|y\|$ for all $h\in H$. All positive-degree homogeneous
$H$-invariants vanish at $y$ since they vanish at $v$. The nullcone
criterion gives $0\in\overline{Hy}$, hence $y=0$. Here Euclidean
and Zariski orbit closures have the same underlying complex points,
since an algebraic orbit is locally closed. Conversely, the
$J$-nullcone is contained in the $H$-nullcone since $J\subset H$.

Put $A_0=\C[V]^H$ and $B_0=\C[V_\sigma]^J$.
The nullcone equality and surjectivity of the affine GIT quotient by $J$
imply
$\sqrt{A_{0,+}B_0}=B_{0,+}$. Finite generation gives
$\dim_\C B_0/A_{0,+}B_0<\infty$; homogeneous lifts of a basis generate
$B_0$ over $A_0$ by induction on degree. Since $J$ is a product of general
linear groups, both actions have exact invariants. Restricting the two
rings to $R$ and its closed fixed scheme $R_\sigma\subset V_\sigma$
preserves the surjections, proving the asserted finiteness also for
nonhomogeneous equations.
\end{proof}

\subsection{Nilpotent cones and Springer fibers}\label{subsec:lie-springer}

The following proof of Proposition~\ref{prop:springer}\ref{it:kostant}
deduces the complex Lie-algebra statement from the group adjoint
quotient. The scheme structure matters: an irreducible set cut out by
the expected number of equations need not be reduced.

\begin{lemma}[The nilpotent cone over $\C$]\label{lem:nilpotent-cone}
Let $L$ be a connected complex reductive group. The zero fiber of
$\mathfrak l\to\mathfrak l\sslash L$, with its scheme structure, is
an irreducible normal complete intersection.
\end{lemma}
\begin{proof}
The central coordinates vanish on this fiber. Thus it suffices to treat
$\mathfrak h=\Lie H$ for the simply connected semisimple group $H$ with
Lie algebra $[\mathfrak l,\mathfrak l]$; changing the isogeny type does
not change this adjoint quotient. Put $r=\operatorname{rk}H$ and let
$\chi:H\to\mathbb A^r$ be the map given by the fundamental characters.
Steinberg's Theorem~6.1, Corollary~6.7 and Theorem~6.11 identify
$U=\chi^{-1}(\chi(1))$ set-theoretically with the unipotent variety:
it is irreducible of codimension $r$, and its regular conjugacy class
has complement of codimension at least two. His Theorems~7.9 and~8.1
show that the differentials of these $r$ characters are independent
on that class \cite[Theorem~6.1, Corollary~6.7 and Theorems~6.11, 7.9 and 8.1,
pp.~62--71]{Ste65}. Consequently the scheme fiber is a complete
intersection and is generically reduced. A complete intersection has
no embedded associated primes, so it is reduced; it is Cohen--Macaulay
and regular in codimension one, and hence normal by Serre's criterion.

For clarity, the transfer to the Lie-algebra zero fiber includes its
scheme structure. Chevalley restriction identifies its adjoint quotient
with $q:\mathfrak h\to\mathfrak t/W$. Near $0\in\mathfrak t$, the
exponential map is a $W$-equivariant analytic isomorphism onto a
neighborhood of $1\in T$ and induces an analytic isomorphism
$\phi:(\mathfrak t/W,0)\simeq(T/W,[1])$. Near any nilpotent $x$ we have
$\chi\circ\exp=\phi\circ q$, after identifying $T/W$ with the
character quotient. This follows first on semisimple elements and then
by density; their Cartan representatives may be chosen near $0$ when
$q$ is near $0$. Moreover the differential of $\exp$ at $x$, after
translation to the identity, is
\[
 \frac{1-e^{-\operatorname{ad}x}}{\operatorname{ad}x},
\]
an invertible polynomial in the nilpotent operator
$\operatorname{ad}x$. Thus the two scheme fibers are locally analytically
isomorphic. Faithful flatness of the algebraic local rings in their
analytifications therefore shows that $q^{-1}(0)$ is reduced. The fiber $q^{-1}(0)$ is cut out by
$r$ basic homogeneous invariants and has codimension $r$, so these
invariants form a regular sequence. The polynomial maps $\exp$ and
$\log$ on nilpotent and unipotent matrices restrict to mutually inverse
algebraic isomorphisms of the reduced varieties (choose a faithful
representation of $H$); hence irreducibility and normality also transfer.
Adjoining the vanishing central linear coordinates proves the assertion
for $\mathfrak l$.
\end{proof}

\subsubsection{Dimensions and low-degree cohomology of Springer fibers}\label{app:springer}

We first explain the dimension equality in
Proposition~\ref{prop:springer}\ref{it:semismall}. Write
$n=\dim G/B$, $\mathfrak n=\Lie R_u(B)$, and let
$Z=\widetilde\Nil\times_{\Nil}\widetilde\Nil$. Above a pair of
flags in relative position $w$, the fiber of $Z$ is a conjugate of
$\mathfrak n\cap w\mathfrak n$, of dimension $n-\ell(w)$.
The orbit of the flag pair has dimension $n+\ell(w)$. Thus every
such part of $Z$ has dimension $2n$.

For $e=0$ the fiber is $G/B$, which gives the asserted equality.
For a nonzero nilpotent $e$, choose a maximal torus in the centralizer of
an $\mathfrak{sl}_2$-triple for $e$. Its centralizer in $G$ is a
Levi $L$ in which $e$ is distinguished: the reductive centralizer
of the triple has no further torus modulo the center of $L$.
Choose the positive roots so $L$ is a standard Levi of type $I$
and the triple's grading is dominant in $L$. This distinguished
grading on $\Lie L$ is even on each simple factor, by the classical block descriptions
and the exceptional root calculation below. Its degree-zero roots
form a Levi of type $J\subset I$ and its positive part
$\mathfrak u$ is the nilradical of a parabolic $Q\subset L$.
The raising map is onto every positive degree from the degree two
steps below. Hence $[\Lie Q,e]=\mathfrak u$ and $Qe$ is open dense in
$\mathfrak u$. For the longest elements of the indicated Weyl groups,
take $w=w_Jw_Iw_0$. The positive roots that $w^{-1}$ keeps positive
are precisely $\Phi_I^+\setminus\Phi_J^+$, so
$\mathfrak n\cap w\mathfrak n=\mathfrak u$. The inverse image of
$G e$ is dense in this $2n$-dimensional flag-pair part and has
dimension at most $2n$ in all of $Z$. Its dimension is also
$\dim G e+2\dim\Bfl_e$ by equivariance. Since
$\dim\Nil=2n$, this proves relevance for every orbit.

The same flag-pair calculation explains why a Weyl multiplicity
selects a single orbit and local system. Put $d=2n$ and
$K=R\mu_*\Q_{\widetilde\Nil}[d]$. Adjunction and proper base change
identify $\operatorname{End}(K)$ with the top Borel--Moore homology
$H^{\rm BM}_{2d}(Z,\Q)$: in the correspondence square the Hom complex
is the cohomology of the dualizing complex on $Z$ shifted by $-2d$,
since $\widetilde\Nil$ is smooth of dimension $d$. Top homology has
one fundamental class per component of dimension $d$. The flag-pair
parts above are irreducible of that dimension, one for each $w$.
Hence $\dim\operatorname{End}(K)=|\WG|$.

For the simultaneous resolution replace $\mathfrak n$ by
$\mathfrak b=\Lie B$ in the double-flag calculation. Each part
has dimension $\dim\fg$, and its flag-pair fiber contains a common
Cartan. Regular semisimple elements are therefore dense in every
part. The complement in the double incidence has dimension strictly
less than $\dim\fg$, which gives the strict fiber-dimension
inequalities for smallness. The deck action on the regular cover
thus extends by intermediate extension, and restriction gives an
action on $K$.

The restriction action of $\Q[\WG]$ on $K$ is faithful already on
the stalk cohomology at zero, $H^*(G/B)$. To see the action, restrict
the cohomology of the simultaneous resolution to one ordered sheet
over the orbit of a regular $t_0\in\ft$. The sheet is $G/T$ by
$gT\mapsto(\operatorname{Ad}(g)t_0,gB)$. Projection onto $G/B$
induces a cohomology isomorphism both from this sheet (fiber $B/T$
is contractible) and from the whole resolution (a vector bundle).
Relabeling the sheet changes $gB$ to $gwB$. Naturality of restriction,
also at the zero fiber, identifies the action with that of the right
normalizer on the compact flag manifold $K_0/T_0$, for a maximal
compact group $K_0$ and its maximal torus $T_0$. Each nonidentity Weyl
element acts freely there, so has Lefschetz number zero. The Bruhat
cells give even cohomology of total rank $|\WG|$. Its character is
therefore the regular character, proving faithfulness.
It follows by dimension that $\Q[\WG]\simeq\operatorname{End}(K)$.
Over $\C$ the decomposition theorem identifies this last algebra
with the product of the matrix algebras of the nonzero multiplicity
spaces. Each such space is consequently irreducible, and distinct
factors give inequivalent representations. For a constant
coefficient the invariant top cohomology in
Proposition~\ref{prop:springer}\ref{it:semismall} gives the rational
model, irreducible also after extension of scalars.

The descent of the Goldman class in
Lemma~\ref{lem:whole-goldman-descent} uses connectedness, vanishing
$H^1$, and Tate $H^2$ of Springer fibers. These follow from the stronger
cycle-class theorem of De~Concini, Lusztig and Procesi
\cite[Theorem~3.9]{DLP88}. We spell out their low-degree argument and
supply the exceptional incidence calculation, which their proof describes
as mechanical and omits. The argument combines smooth orbit strata with restriction of
holomorphic forms to the fixed loci of smooth projective varieties.
The finite calculation below has independently checked exact certificates.

\begin{proposition}\label{prop:springer-low-degree}
For every nilpotent $N$ in the Lie algebra of a connected complex
reductive group, the Springer fiber $\Bfl_N$ is connected,
$H^1(\Bfl_N,\Q)=0$, and $H^2(\Bfl_N,\Q)\cong\Q(-1)^{\oplus m}$ for some
$m$.
\end{proposition}

\begin{proof}
Connectedness is Proposition~\ref{prop:springer}\ref{it:kostant}.
Central isogenies preserve flag varieties and Springer fibers, and
products are handled by K\"unneth. Thus it suffices to treat adjoint
simple groups. All cohomological comparisons below respect mixed Hodge
structures.

\emph{Fixed varieties and the incidence relation.}
Let $M$ be connected reductive with Borel $B_0$, and let $V$ be an
$M$-module with an open orbit $Mv$. For a $B_0$-stable subspace $U$, set
\[
 X_U=\{gB_0:g^{-1}v\in U\},\qquad
 \delta(U)=\dim(M/P_U)-\operatorname{codim}_V U,
\]
where $P_U$ is the full stabilizer of $U$. The group $P_U$ is parabolic: its identity
component contains $B_0$, hence is parabolic and self-normalizing.
The smooth orbit map $M\to Mv$ shows that $X_U$ is smooth; it is also
closed in $M/B_0$, hence projective. The same argument on $M/P_U$ shows
that $X_U$ is empty for $\delta(U)<0$, and is a finite disjoint union of
copies of $P_U/B_0$ for $\delta(U)=0$. These are the seed cases.

Suppose $U''\subset U\subset U'$ have one-dimensional successive
quotients, and a minimal parabolic $P_i$ stabilizes $U''$ and $U'$ but
not $U$. The incidence construction of
\cite[\S\S2.7--2.11]{DLP88} identifies the same variety as a
$\PP^1$-bundle over $X_U$ and as the blow-up of $X_{U'}$ along
$X_{U''}$. Here is the local geometry behind both descriptions. In the
rank-two bundle with fiber $gU'/gU''$, the image of $v$ is a section
whose zero scheme is the smooth codimension-two locus $X_{U''}$.
The closure of its graph in the projectivization is the blow-up.
The rank-one semisimple factor of $P_i$ acts irreducibly on
$U'/U''$: otherwise its two one-dimensional constituents would both
be trivial, contradicting the failure to stabilize $U/U''$.
The unipotent radical acts trivially, since its invariant subspace is
nonzero and stable under this rank-one factor. Thus
$P_i/B_0\cong\PP(U'/U'')$, giving the other description.
The projective-bundle and blow-up formulas give
\[
 H^1(X_U)\cong H^1(X_{U'}),\qquad
 H^2(X_U)\oplus H^0(X_U)(-1)
 \cong H^2(X_{U'})\oplus H^0(X_{U''})(-1).
\]
Consequently ``$H^1=0$ and $H^2$ Tate'' propagates in both directions
along an incidence edge.

\emph{The reduction and its induction hypothesis.}
If $N=0$, then $\Bfl_N=G/B$. For semisimple $s$, the components of
$(G/B)^s$ are flag varieties of $Z_G(s)^0$; their Bruhat decompositions
give the asserted vanishing and Tate property, including the stronger
induction assertion below. Hence assume $N\ne0$.
Choose an $\mathfrak{sl}_2$-triple $(N,h,f)$ (Jacobson--Morozov
\cite[Theorem~3.3.1]{CM93}) with cocharacter $\tau$,
write $\fg=\bigoplus_j\fg_j$, and set $G_0=Z_G(\tau)$ and
$P=P(\tau)$, with Lie algebra $\fg_{\geq0}$. The raising maps in the
$\mathfrak{sl}_2$ decomposition show that $G_0N$ is open in $\fg_2$
and $PN$ is open in $\fg_{\geq2}$. Intersecting the finitely many
$P$-orbits in the flag variety with $\Bfl_N$ gives strata ordered
by a closed filtration. On an orbit through a $\tau$-fixed Borel
$B'$ the inverse image of the stratum in $P$ is cut out by
$\operatorname{Ad}(p^{-1})N\in\Lie B'\cap\fg_{\geq2}$.
The orbit map is smooth onto an open of $\fg_{\geq2}$, so this
inverse image and its quotient by the flag stabilizer are smooth.
Write $p=nm$ with $n\in R_u(P)$ and $m\in G_0$. Acting by
$\tau(t)$ contracts $nmB'$ to $mB'$ as $t\to0$, since the
unipotent coordinates have positive weights. The contraction stays
in $\Bfl_N$ and gives a deformation retraction onto the smooth
projective $X_U$ with $M=G_0$, $V=\fg_2$ and
$U=\Lie B'\cap\fg_2$ \cite[\S\S3.2--3.4]{DLP88}.

For the subsequent dimension induction one must prove a slightly
stronger assertion: $(\Bfl_N)^s$ has $H^1=0$ and Tate $H^2$ whenever
$s$ is semisimple and $\operatorname{Ad}(s)N=qN$ for some $q\ne0$.
This includes semisimple actions induced on smaller parabolic quotients;
those actions need not be the quotients' own canonical grading actions.
The fixed variety is projective, but need not be smooth.
An $s$-adapted triple can be obtained directly from an ordinary triple
$(N,h_0,f_0)$. The set
\[
 A=\{h:[h,N]=2N,\ h\in\operatorname{im}(\operatorname{ad}N)\}
 =h_0+(\ker\operatorname{ad}N\cap\operatorname{im}\operatorname{ad}N)
\]
is affine and is preserved by $\operatorname{Ad}(s)$. Its direction
space consists of the positive-weight highest vectors in the
$\mathfrak{sl}_2$-decomposition for the ordinary triple. Each
$h\in A$ is conjugate to $h_0$ by elements centralizing $N$:
if the lowest remaining weight in $h-h_0$ is $j>0$ with component
$v_j$, conjugation by $\exp(v_j/j)$ kills that component and changes
only higher weights. The process terminates. As $\operatorname{Ad}(s)$
is semisimple, its projection to the $1$-eigenspace can be written as
$p(\operatorname{Ad}(s))$ with $p(1)=1$ (interpolate on the spectrum
together with $1$). This is an affine combination of powers of
$\operatorname{Ad}(s)$ and sends $h_0$ to a fixed $h\in A$. Transporting
$f_0$ gives a completing $f$, unique since $\ker\operatorname{ad}N$
contains no vectors of $h$-weight $-2$. Thus $q\operatorname{Ad}(s)f=f$.
The triple integrates to $\mathrm{SL}_2\to G$, with diagonal
cocharacter $\tau$. The element $s$ commutes with $\tau$, and
$s_1=s\tau(q^{-1/2})$ centralizes the triple.
Set $G_1=Z_G(s_1)^0$, which is reductive, and take the preceding
$P(\tau)$-orbit pieces for this adapted triple.
The $s$-fixed parts of the preceding smooth orbit pieces are smooth
and contract under $\tau$ to their simultaneous $s_1,\tau$-fixed loci.
Components of the ambient $s_1$-fixed flag variety are flag varieties
of $G_1$, so these bases are unions of connected components of
copies of $(\Bfl_N^{G_1})^\tau$.
The contraction is equivariant for $s$ and hence also retracts
these smooth fixed strata onto their bases
\cite[\S\S3.4--3.5]{DLP88}. If $G_1\ne G$, dimension induction applies,
using an element of $\tau(\C^*)$ with Zariski-dense cyclic subgroup.
If $G_1=G$, the remaining base is $\Bfl_N^\tau$.

This last variety is a disjoint union of the smooth projective $X_U$.
Choose a maximal torus $T$ in the connected centralizer of $N$, chosen
to centralize the triple, and set $L=Z_G(T)$. The $T$-fixed locus is a
finite disjoint union of copies of $(\Bfl_N^L)^\tau$. If $L\ne G$,
dimension induction gives vanishing $H^1$ and Tate $H^2$ on these
fixed components. To deduce the same on a connected component $X$
of $X_U$, use a $T$-equivariant projective embedding of the ambient
flag variety and choose a cocharacter separating the finitely many
torus weights of the coordinates. On the dense open of $X$ where
the lowest occurring weight coordinate is nonzero, projection to the
whole lowest-weight space is the limit at zero of the cocharacter
action. This gives a morphism to $X^T$: after division by the lowest
power of the parameter, all remaining coordinates have nonnegative
exponents. A holomorphic $j$-form on $X$ is invariant under the
connected torus since $H^0(X,\Omega_X^j)$ injects into $H^j(X,\C)$
and each torus element is homotopic to the identity. Passing to
the limit in its pullback on this open gives the pullback of its
restriction to $X^T$. For $j=1,2$ that restriction is zero by
induction. The form therefore vanishes on $X$ by density. Hodge
decomposition and complex conjugation now give $H^1(X)=0$ and
$H^2(X)$ of type $(1,1)$, as required. If $L=G$, then $N$ is distinguished, meaning
that its connected centralizer has no torus beyond the center
\cite[\S\S3.5--3.7]{DLP88}. Compact-support assembly, verified below,
justifies each of these reductions.

For distinguished classical nilpotents the standard module has one
Jordan block in type $A$, distinct even blocks in type $C$, and
distinct odd blocks in the orthogonal types. These restrictions follow
from the absence of positive-dimensional tori in the orthogonal or
symplectic groups of the multiplicity spaces. In particular, the
$\tau$-weights on $\ker N$ are distinct. Sending a $\tau$-fixed
invariant flag to its first $N$-killed line therefore has finite image
(in the orthogonal case retain isotropic lines). The fibers are open
and closed in $\Bfl_N^\tau$ and are semisimple-fixed Springer fibers
for a smaller reductive quotient of a proper parabolic. The induced
semisimple element scales the induced nilpotent, so the stronger
induction hypothesis applies \cite[\S3.8]{DLP88}.

\emph{The exceptional calculation.}
For the five exceptional root systems, we check that every $X_U$ is
joined by incidence edges to a seed with $\delta\leq0$, supplying the
finite step in \cite[Theorem~3.9]{DLP88}. Choose $h$ dominant and put
$d_i=\alpha_i(h)$. Then $d_i\in\{0,1,2\}$: if $d_i>2$, surjectivity
of $\operatorname{ad}N:\fg_{d_i-2}\to\fg_{d_i}$ would express the simple
root $\alpha_i$ as a sum of two positive roots. For a distinguished
nilpotent in an adjoint simple group, the triple's reductive centralizer
has zero Lie algebra. The $\mathfrak{sl}_2$ decomposition of $\fg$ therefore gives
$\dim\fg_0=\dim\fg_2$. Only this necessary condition is used; no converse
classification claim is needed.

The exact calculation enumerates all $9{,}567$ labelings in
$\{0,1,2\}^r$ across these five root systems. The dimension equality
retains $2,4,3,6,11$ candidates in $G_2,F_4,E_6,E_7,E_8$, respectively.
None has an odd label; evenness is thus an output, not an assumption.
For each candidate, $V=\fg_2$ has distinct root weights, so its
$B_0$-stable subspaces are precisely the subsets closed under addition
of zero-labeled simple roots. A minimal parabolic stabilizes such a
subset exactly when it is also closed under lowering by that simple
root, equivalently invariant under its simple reflection. If $I_U$
is the set of those simple roots, its full stabilizer is $P_{I_U}$ and
\[
 \delta(U)=|\Phi_0^+|-|\Phi_{I_U}^+|
            -\bigl(|\Phi_2|-\dim U\bigr).
\]
The verifier enumerates every stable subset, checks every stored edge
against the three stabilizer conditions of the incidence construction,
and checks that every graph component meets a seed. The totals are
\[
\begin{array}{c|r|r|r}
\text{type}&\text{candidates}&\text{subspaces}&\text{verified edges}\\\hline
G_2&2&9&0\\
F_4&4&105&30\\
E_6&3&352&182\\
E_7&6&2328&1938\\
E_8&11&14208&16183\\\hline
\text{total}&26&17002&18333
\end{array}
\]
All stable subspaces must be present; only enough valid edges to connect
them to seeds are required.

The accompanying directory \texttt{springer\_verification} contains the
26 certificates, their hashes and \texttt{FINITE\_PROOF.md}.
The generator uses integral Weyl-reflection closure and descending-height
subset enumeration. The independent verifier uses Euclidean root sets
in exact rational coordinates, upward-closure enumeration and reflection
invariance. The verifier also reconstructs the Cartan matrices from root strings
and exhausts the ternary labelings. An additional verifier tests every
negative $\fg_0$-root direction on every subspace and recovers the same
full stabilizers. From that directory run, without Python's optimization
flag,
\begin{verbatim}
python3 check_certificates.py
python3 reaudit_finite_proof.py
\end{verbatim}
Both programs use only exact integer or rational arithmetic and the
Python standard library. The first finishes with
\texttt{TOTAL 17002 18333}; the second reports all 26 stabilizer checks
as passed. The completeness argument above links this exhaustive
calculation to the geometric induction. The certificates verify this
finite exceptional-goodness assertion. For the original
geometric reduction by attracting strata, see
\cite[\S\S1.2--1.5 and 3.6]{DLP88}.

\emph{Assembly.}
For a smooth connected stratum $E$ retracting onto the smooth
projective fixed component $Y$, the algebraic inclusion induces an
isomorphism $H^*(E)\simeq H^*(Y)$ of mixed Hodge structures.
Poincar\'e duality on both sides gives
$H_c^j(E)\simeq H^{j-2r}(Y)(-r)$, where $r=\dim E-\dim Y$.
The integer $r$ is nonnegative because $Y$ is a fixed subvariety of $E$. Thus $H^1(Y)=0$ and Tate $H^2(Y)$ imply
$H_c^1(E)=0$ and Tate $H_c^2(E)$: only $r=0,1$ can contribute.
For $Z\subset X$ closed with complement $O$, the same properties for
$Z$ and $O$ give $H_c^1(X)=0$ and the exact sequence
\[
 0\longrightarrow H_c^2(O)\longrightarrow H_c^2(X)
 \longrightarrow\ker\bigl(H_c^2(Z)\to H_c^3(O)\bigr)
 \longrightarrow0.
\]
Strictness of mixed Hodge structures makes the middle term Tate of
weight two. Induction along the stated closed filtrations proves the
claim, and properness of $\Bfl_N$ identifies compactly supported and
ordinary cohomology.
\end{proof}


\subsection{Stabilized Bruhat coordinates and the Goldman form}\label{subsec:stabilized-bruhat}
\label{app:stabilized-bruhat}

This subsection proves the geometric input used in
\S\ref{subsec:torus-coordinates}. As there, $G_{\der}$ is simply
connected and all boundary classes are regular semisimple; general global
forms are treated by the finite central descent of
\S\ref{subsec:central-quotient}. We work throughout on the stabilized representation variety
$R\times U_B$ of Definition~\ref{def:stabilization}, in its description
as the space of tuples whose surface product lies in $U_BC_s^{-1}$; the
leading factor $u_0\in U_B$ below records that product. The coordinate
constructions extend Mellit's type-$A$ construction
\cite[\S\S5.4, 7.3--7.4]{Mel25} and follow \cite[\S\S3--8]{dA25}. Only
the action of the torus $T$ is used: no formula for an action of $B$ and
no cancellation after a quotient by its unipotent radical enters the
proof.

The calculation has two parts. Subsection~\ref{subsec:H-coordinates} builds
$T$-equivariant product charts from two elementary moves (the Bruhat
switches of Lemma~\ref{lem:H-switch}), assembles them into charts for
handles and punctures (Lemma~\ref{lem:H-charts}) and for the whole
relation space (Corollary~\ref{cor:H-charts}), and proves the closed
filtration (Lemma~\ref{lem:H-filtration}).
Subsections~\ref{subsec:H-bars}--\ref{subsec:H-word} track the Goldman form
through each of these coordinate changes; an abstract isomorphism of
varieties would not justify this comparison. A rank-one computation at
the end checks all signs.

Write $U=U_B$, and choose the opposite unipotent group $U^-$.
For $w\in\WG$, choose a lift $\dot w\in N_G(T)$ and put
\[
 U_w^-=U\cap\dot w^{-1}U^-\dot w,\qquad
 U_w^+=U\cap\dot w^{-1}U\dot w.
\]
Ordered root multiplication gives $U=U_w^-U_w^+$ and
$\dim U_w^-=\ell(w)$. Polynomial logarithm identifies every unused
unipotent factor below with its Lie algebra, equivariantly for torus
conjugation. A \emph{root chart} is a map
$F_w(z)=\dot w\,r(z)$ with $r:\A^{\ell(w)}\xrightarrow{\sim}U_w^-$
given by ordered root coordinates; the torus acts on its coordinates by
$F_w(t\cdot z)=w(t)F_w(z)t^{-1}$.

\Needspace{5\baselineskip}
\subsubsection{Equivariant coordinates and the closed filtration}
\label{subsec:H-coordinates}

\begin{lemma}[Equivariant Bruhat switches]\label{lem:H-switch}
Fix $w\in\WG$ and a lift $a\in N_G(T)$ of $w$. Choose root-group charts
$P(z)=a r_w(z)$ and $\bar P(z)=a^{-1}r_{w^{-1}}(z)$, where
$r_v:\A^{\ell(v)}\simeq U_v^-$. For $d\in T$ the identities
\[
 \bar P(z)^{-1}du=u'w(d)P(z'),\qquad
 P(z)^{-1}w(d)u=u'd\bar P(z')
\]
define algebraic isomorphisms $(z,u)\leftrightarrow(u',z')$, jointly
algebraic in $d$. They intertwine left and right multiplication by $T$.
\end{lemma}
\begin{proof}
These are the two unique Bruhat parametrizations of the same Cartan
fiber of $BwB$ and $Bw^{-1}B$. Their inverse maps are the Bruhat
coordinate projections. Ordered root-subgroup multiplication proves
that these maps and their inverses are regular. The root charts satisfy
$F_v(t\cdot z)=v(t)F_v(z)t^{-1}$ with a diagonal action on $z$;
uniqueness proves the stated equivariance.
\end{proof}

\begin{lemma}[The $T$-equivariant handle and puncture charts]\label{lem:H-charts}
Let $b=\dim U$. After adjoining the leading $U$-coordinate, a handle
Bruhat pair $(w,v)$ has a product chart
\[
 U\times T^2\times\A^{2\ell(w)+2\ell(v)}
   \times U_w^+\times U_{v^{-1}}^+,
\]
whose map to $G$ is
\[
 (u_0,t_1,t_2,z)\longmapsto
 u_0D_{w,v}(t_1,t_2)P(z_1)Q(z_2)\bar P(z_3)\bar Q(z_4),
\]
where $P,\bar P$ and $Q,\bar Q$ use inverse lifts of $w,v$ and
\[
 D_{w,v}=w(t_1)wv(t_2t_1^{-1})wvw^{-1}(t_2^{-1}).
\]
The last two unipotent factors do not occur in this map and carry
conjugation by $T$. For a regular puncture $C\in T$, its Schubert
stratum $w$ similarly gives
\[
 U\times\A^{2\ell(w)}\times U_w^+,
 \qquad (u_0,z,z')\longmapsto u_0P(z)C\bar P(z'),
\]
with $T$ acting on the unused $U_w^+$ by conjugation through $w^{-1}$.
All these isomorphisms are $T$-equivariant for the original simultaneous
conjugation action.
\end{lemma}
\begin{proof}
Here are explicit reversible constructions. Write
$x=\bar P(\xi)^{-1}t_1r_1$ and $y=r_2v(t_2)Q(\eta)$, and replace
$(r_1,r_2,u)$ by
\[
 h=r_1Q(\eta)^{-1}b_v,\quad k=r_1r_2,\quad l=r_2^{-1}u,
\]
where $a$ and $b_v$ are the chosen lifts of $w$ and $v$. Apply the second switch to
$Q(\eta)^{-1}v(t_2^{-1})l=l_1t_2^{-1}\bar Q(\eta_1)$ and put
$u_1=a\bar P(\xi)l_1$. Split
$h^{-1}=P(\zeta)^{-1}a h_+$, with $h_+\in U_w^+$, put $u_2=\operatorname{Ad}_{at_1}(h_+)u_1$, and use the second switch
\[
 P(\zeta)^{-1}w(t_1^{-1})u_2=u_3t_1^{-1}\bar P(\zeta_1).
\]
Put $\zeta_2=t_2\cdot\zeta_1$. Moving $t_2^{-1}$ left through
$\bar P(\zeta_1)$ gives
\[
 \bar P(\xi)^{-1}t_1k b_vt_2u_3
 t_1^{-1}w^{-1}(t_2^{-1})\bar P(\zeta_2)\bar Q(\eta_1).
\]
Split $k=\bar Q(\theta)^{-1}b_v^{-1}k_+$ with
$k_+\in U_{v^{-1}}^+$; put $u_4=\operatorname{Ad}_{t_2^{-1}b_v^{-1}}(k_+)u_3$, and apply
the first switch for $v$:
\[
 \bar Q(\theta)^{-1}t_2u_4=u_5v(t_2)Q(\theta_1).
\]
Put $d=t_1^{-1}w^{-1}(t_2^{-1})$,
$\theta_2=d^{-1}\cdot\theta_1$, and $D_0=v(t_2d)$.
Move $d$ left through $Q$ and use the first switch for $w$:
\[
 \bar P(\xi)^{-1}t_1u_5=u_0w(t_1)P(\xi_1).
\]
Finally put $\xi_2=D_0^{-1}\cdot\xi_1$. The displayed handle formula,
with word coordinates $(\xi_2,\theta_2,\zeta_2,\eta_1)$, follows. Every step is a Bruhat
switch, the product decomposition $U=U_w^-U_w^+$, an invertible
unipotent multiplication, or a rescaling by a nonzero torus character.
This also constructs the inverse.

Under conjugation by $s\in T$, the initial torus parameters transform
as $t_1\mapsto w^{-1}(s)s^{-1}t_1$ and
$t_2\mapsto v^{-1}(s)s^{-1}t_2$. The variables $h,k$ transform by
conjugation, so both unused factors do too. The right torus arguments
on the four word blocks are respectively
\[
 vw^{-1}v^{-1}(s),\quad w^{-1}v^{-1}(s),\quad v^{-1}(s),\quad s.
\]
These follow directly from the two-sided equivariance of each switch.
They are the ordinary convolution weights. The leading $u_0$ transforms
by conjugation.

For the puncture, the automorphism $n\mapsto C^{-1}nCn^{-1}$ of $U$
puts its conjugacy-class chart in the form
$\bar P(\xi)^{-1}Cn\bar P(\xi)$. Split
$n=P(\eta)^{-1}a n_+$ and absorb $an_+\bar P(\xi)$ into the following
$u\in U$, setting $u'=an_+\bar P(\xi)u$. Apply the second switch
$P(\eta)^{-1}u'=u_1\bar P(\eta_1)$, put $u_2=Cu_1C^{-1}$, and
apply the first switch
$\bar P(\xi)^{-1}u_2=u_0P(\xi_1)$. This gives the claimed puncture map,
with word coordinates $(\xi_1,\eta_1)$ and $n_+$ unused. The torus acts
on $n_+$ by conjugation through $w^{-1}(s)$; the right torus
arguments on $P,\bar P$ are $w^{-1}(s),s$. The automorphism of $U$
used here is triangular by root height, with nonzero diagonal entries
because $C$ is regular. Thus its inverse is polynomial.
\end{proof}

\begin{corollary}[Literal stabilized charts]\label{cor:H-charts}
Iterating the preceding constructions and the rank-one walk induction
gives, on the stabilized representation variety $R\times U_B$, charts
$V\times q^{-1}(c)$, where $V$ is a diagonal $T$-representation,
$q:T^{2g}\times\Gm^k\to T_{\rm der}$ is a homomorphism, and $c$ includes
all fixed Cartan factors. For a closed walk with $a$ up steps,
$\dim V=b(2g+s-1)-L/2+a$ and $L=2a+k$.
\end{corollary}
\begin{proof}
Each operation leaves the subsequent factor unchanged and retains a
leading $U$, so it can be iterated. Moving Cartan factors left uses
$F_w(z)t=w(t)F_w(t^{-1}\cdot z)$. Changing inverse lifts to pinned
reversed words introduces only fixed Cartan factors. The walk induction
is triangular: up coordinates are free, a down coordinate is a regular
function of its predecessors, and at a stay the difference from that
function is inverted. The final Cartan factor is independent of all
free affine coordinates. Invariance of the down locus makes the
recursive difference homogeneous. Logarithm linearizes conjugation
on each unused unipotent factor. The final surface relation uniquely determines the leading
$U$-coordinate. This proves the literal product and its equivariant
zero section.
\end{proof}

\begin{lemma}[Closed filtration on the stabilized space]\label{lem:H-filtration}
These charts admit a finite filtration by closed $T$-invariant unions;
the filtration descends to the rigidified quotient stack and, under
intrinsic genericity, to its coarse space.
\end{lemma}
\begin{proof}
First order the outer product Bruhat strata by a linear extension of
the product Bruhat order. For punctures use $G/T\to G/B$ and the
Schubert strata. Initial lower ideals are closed; intersecting with the
final-relation locus gives closed invariant $X_i$ whose differences
$S_i=X_i\setminus X_{i-1}$ are the outer charts on that locus.

On each $S_i$, the coordinate construction gives regular prefix
maps to $(G/B)^L$. Inverse images of initial lower ideals for its product
Bruhat order form a relatively closed invariant filtration of $S_i$,
with walk strata as differences. If $Y\subset S_i$ is one such relatively
closed subset, $X_{i-1}\cup Y$ is closed in $X_i$, since its complement
is open in the open subset $S_i$. Concatenation gives the required
closed filtration of $X$. Components of a fixed torus fiber are clopen
and are preserved by connected $T$, so finite component refinements are
allowed.

Closed invariant subsets give closed substacks. For coarse spaces,
$X$ is affine and all $H=T/Z(G)$-orbits are closed. Indeed the direct
last-puncture change of variables identifies $X$ with $R\times U$.
An $H$-orbit in $R$ is the intersection of a closed $G$-orbit with the
fixed regular last-monodromy fiber, hence is closed; its stabilizer $F$
is finite by intrinsic genericity. Over this orbit the product is
$H\times^F U$, and its $H$-orbits correspond to finite, hence closed,
$F$-orbits in $U$. Thus the affine GIT quotient is geometric.
Every invariant closed member of the filtration is saturated and has
closed image. No $B$-action on the charts or descent of their partition
along stabilization is required.
\end{proof}

The quotient of these equivariant charts is given by
Lemma~\ref{lem:H-quotient}; its finite group is the kernel of the gauge
map.

\Needspace{5\baselineskip}
\subsubsection{Tracking the surface form}\label{subsec:H-bars}

The bar potential $\Omega$ of Lemma~\ref{lem:H-block-bars} is twice
the conventional fusion form;
the factor $1/2$ is kept uniformly for handles, punctures and word
coordinates. This computes the original Goldman class $E$, without a
later rescaling.

\begin{lemma}[Bar compatibility of the Bruhat switches]\label{lem:H-switch-bars}
Both equivariant switches preserve the local bar potential:
\begin{align*}
 \Omega(\bar P^{-1},d,u)&=\Omega(u',w(d),P'),\\
 \Omega(P^{-1},w(d),u)&=\Omega(u',d,\bar P').
\end{align*}
The identities hold for variable $d\in T$. A Cartan slide through a
reduced-word block also preserves the local potential.
\end{lemma}
\begin{proof}
For $f=u n_w t v\in Bn_wB$, define
$\mathcal B_w(f)=B(u,n_wtv)$. This is independent of the factorization.
Indeed, write another factorization as
$u'=uc$, $k=n_wt$, $d_0=k^{-1}ck\in U$ and
$v'=d_0^{-1}v$, with $c\in U\cap n_wUn_w^{-1}$.
The block $(c,k,d_0^{-1})$ has product $k$ and zero internal potential:
$B(c,k)=0$ and
$B(ck,d_0^{-1})=B(kd_0,d_0^{-1})=0$ by the cocycle and inverse
identities and $B(k,d_0)=0$. Also
$B(u,c)=B(d_0^{-1},v)=0$. Lemma~\ref{lem:H-block-bars} therefore proves
$B(u,kv)=B(uc,kd_0^{-1}v)$.

In the first switch, expand
$\bar P^{-1}=r_{w^{-1}}^{-1}a$ on the left and $P'=ar'_w$ on the right.
Since $w(d)a=ad$, both local bar sums equal the primitive
$\mathcal B_w$ of their common product. All internal bars inserted in
these expansions vanish. The second switch is the same argument with
$w^{-1}$.

For a Cartan slide $F(z)t=w(t)F(t^{-1}\cdot z)$ with
$F=ar(z)$, $r(z)\in U$, both $B(F,t)$ and
$B(w(t),F(t^{-1}\cdot z))$ vanish: their root-valued Maurer--Cartan
forms pair trivially with $\mathfrak t$. This remains true when the
new root coordinates depend on $t$.
\end{proof}

\begin{lemma}[Exact compatibility of the handle chart]\label{lem:H-handle-bars}
Use the notation of Lemma~\ref{lem:H-charts}, and write $b_v=b$ for the chosen lift of $v$. The coordinate
isomorphism of that lemma satisfies
\begin{align*}
 &\Omega(x,y,x^{-1},y^{-1},u)\\
 &\qquad=\Omega_H+
 \Omega(u_0,D_{w,v},P(\xi_2),Q(\theta_2),
                      \bar P(\zeta_2),\bar Q(\eta_1)),
\end{align*}
where
\[
 \Omega_H=
 \Omega\bigl(w(t_1),wv(t_2),wv(t_1^{-1}),wvw^{-1}(t_2^{-1})\bigr).
\]
Moreover, for $n_A=at_1$ and $n_B=bt_2$,
\[
 \Omega_H=\Omega(n_A,n_B,n_A^{-1},n_B^{-1}).
\]
\end{lemma}
\begin{proof}
Expand $x=\bar P(\xi)^{-1}t_1r_1$ and
$y=r_2v(t_2)Q(\eta)$ together with their inverses.
Their internal insertion corrections cancel by the block lemma.
The substitutions $h,k,l$ of the coordinate construction then give
\[
 (\bar P(\xi)^{-1},t_1,k,b,t_2,h^{-1},t_1^{-1},
     \bar P(\xi),Q(\eta)^{-1},v(t_2^{-1}),l)
\]
as a bar word with the same potential. The following zero-cost moves
justify every unipotent simplification in this statement and in the
rest of the construction:
\begin{itemize}
\item Products of $U$-factors have zero internal bars.
\item The splittings $h^{-1}=P(\zeta)^{-1}ah_+$ and
$k=\bar Q(\theta)^{-1}b^{-1}k_+$ have zero internal bars, since
$P^{-1}a$ and $\bar Q^{-1}b^{-1}$ lie in $U$.
\item The replacements $Qr_1^{-1}=bh^{-1}$ and
$\bar P l_1=a^{-1}u_1$ have zero internal bars, since $Q$ and
$\bar P$ have left Maurer--Cartan forms in $\mathfrak u$.
\item The absorptions of $h_+$ and $k_+$ preserve the bar sum:
$aU_w^+a^{-1}\subset U$ and
$b^{-1}U_{v^{-1}}^+b\subset U$ reduce each relevant block to
$U$-factors and one variable Cartan factor, and the internal bars of
such a block vanish.
\end{itemize}
Every Bruhat switch and every Cartan slide is covered by the preceding
lemma. The only remaining costs are the three variable-Cartan mergers.
Put $d=t_1^{-1}w^{-1}(t_2^{-1})$ and $D_0=v(t_2d)$. They retain
\[
 K_1=B(t_1^{-1},w^{-1}(t_2^{-1})),\qquad
 K_2=B(v(t_2),v(d)),\qquad
 K_3=B(w(t_1),w(D_0)),
\]
respectively when forming $d$, $D_0$, and
$D_{w,v}=w(t_1)w(D_0)$. The block lemma determines the signs:
the initial potential is the final potential plus these three terms.
Their sum is $\Omega_H$. To see this, conjugate its four Cartan entries
by $w^{-1}$ and successively merge its last two and then last three
entries. The three retained bars are precisely $K_1,K_2,K_3$.

Finally, expand $at_1$, $bt_2$ and their inverses in the normalizer
commutator. Move all constant lifts to the right. The remaining four
Cartan entries are those in $\Omega_H$ and the remaining constant is
$[a,b]$. Constant entries add no bar. This proves the second identity.
\end{proof}

\begin{lemma}[Exact compatibility of the puncture chart]\label{lem:H-puncture-bars}
For a fixed regular $C\in T$, the $T$-only puncture chart satisfies
\[
 \Omega(g,C,g^{-1},u)
 =\Omega(u_0,P(\xi_1),C,\bar P(\eta_1)).
\]
The left side is twice the fixed-conjugacy-class form fused with the
following unipotent factor.
\end{lemma}
\begin{proof}
If $\theta=g^{-1}dg$, direct differentiation gives
\[
 \Omega(g,C,g^{-1})
 =-\kappa(\operatorname{Ad}_{C^{-1}}\theta\wedge\theta)
 =\kappa(\operatorname{Ad}_C\theta\wedge\theta),
\]
which is twice the conjugacy-class form. Write
$g=\bar P(\xi)^{-1}h$ with $h\in U$ and put
$n=C^{-1}hCh^{-1}$. Expanding $g$ together with its inverse cancels
internal insertion corrections. The blocks $(h,C,h^{-1})$ and $(C,n)$
have the same product and both have zero internal bar sum. Consequently
\[
 \Omega(g,C,g^{-1},u)
 =\Omega(\bar P(\xi)^{-1},C,n,\bar P(\xi),u).
\]
Thus the Lang substitution $h\mapsto n=C^{-1}hCh^{-1}$ is compatible
with the conjugacy-class form, including its differential terms.

The splitting $n=P(\eta)^{-1}an_+$ costs zero since $P^{-1}a\in U$.
The absorption $u'=an_+\bar P(\xi)u$ costs zero because
$an_+a^{-1}\in U$ and $\bar P=a^{-1}r$ with $r\in U$.
The two switches cost zero by the switch lemma; moving the fixed $C$
through the intermediate unipotent also costs zero. Their resulting
bar word is $(u_0,P(\xi_1),C,\bar P(\eta_1))$, proving the assertion.
\end{proof}

\Needspace{5\baselineskip}
\subsubsection{Iteration and the terminal level}\label{subsec:H-iteration}
The fused-double form is
$\frac12\Omega(x,y,x^{-1},y^{-1})$, as follows directly from the bar
cocycle and the usual double formula. The fixed conjugacy-class form
was checked above. On the final $U C_s$ factor the conjugacy form is
zero, since its conjugator lies in $U$. These facts identify the initial
full potential with twice the pullback of the Goldman form.
The block lemma proves compatibility with all subsequent factors, so
the two local identities can be iterated.

Collect the handle Cartan factors to the left by sliding through each
reduced-word block. Keep all four Cartan atoms of every handle, or retain
the bar term each time adjacent Cartan factors are merged. The preceding
Weyl product conjugates the four atoms of the $i$-th handle by $C_{<i}$.
Thus their potential is exactly the potential of the corresponding
normalizer handles, including their inter-handle fusion terms.
The total map becomes $u_0D F_{\rm word}$. Apply the walk induction to
$F_{\rm word}$ alone: its initial Cartan factor is $h_0=1$, not $D$.
On the final closed-walk locus, $F_{\rm word}\in B$ has Cartan part
$\mu$. The remaining two fusion bars are
\[
 B(u_0,D F_{\rm word})=0,\qquad
 B(D,F_{\rm word})=\kappa(d\log D\wedge d\log\mu)=0,
\]
where the second equality uses $D\mu=\text{constant}$.
The raw/rescaled word calculation therefore identifies the restricted
Goldman form itself with the normalizer handle-and-reflection form, so
its parabolic cup-product interpretation and nondegeneracy apply to the
Goldman form, not merely to an abstract complex of the same dimension.

\Needspace{5\baselineskip}
\subsubsection{The word recursion}\label{subsec:H-word}

\begin{lemma}[The word bar recursion]\label{lem:H-word-bars}
Let a partial product be $F=u n_p h v$ in Bruhat form, with the
canonical left unipotent factor $u$, and suppose
\[
 \Omega_{\rm partial}=\Omega_{\log}+B(u,n_phv),
 \qquad B(f,g)=\kappa(f^{-1}df\wedge dg\,g^{-1}).
\]
Prepend $b_\alpha(z)=n_\alpha x_\alpha(z)$ for a simple root $\alpha$.
An up or down step does not change $\Omega_{\log}$. At a stay, with
raw coordinate $a$ and $\beta=p^{-1}\alpha$, the recursion is
\[
 h_{\rm new}=\beta^\vee(a^{-1})h,\qquad
 \Omega_{\log,\rm new}
 =\Omega_{\log}+
   \kappa(d\log\beta^\vee(a^{-1})\wedge d\log h),
\]
up to constant Cartan factors from the chosen lifts.
\end{lemma}
\begin{proof}
Let $U_{\rm rest}$ be generated by the positive root groups other than
$U_\alpha$. The group $U_{\rm rest}$ is normal in $U$, and ordered root multiplication gives
$U=U_\alpha U_{\rm rest}$. Since a simple reflection permutes the
positive roots other than $\alpha$,
$n_\alpha U_{\rm rest}n_\alpha^{-1}=U_{\rm rest}\subset U$.
Root orthogonality gives
\[
 \kappa(\mathfrak u_{\rm rest},\mathfrak b)=0,
 \qquad
 \kappa(\mathfrak u_{\rm rest},\operatorname{Ad}_{n_\alpha}
                                      \mathfrak b)=0;
\]
the only negative root in the second Borel is $-\alpha$.
Pairings with an arbitrary later $n_p$ are not being set to zero;
they are retained in the Bruhat primitive.

The cocycle identity and $B(b_\alpha(z),u)=0$ give the exact expression
\[
 \Omega_{\rm new}
 =\Omega_{\log}+
   \Omega(n_\alpha,x_\alpha(z),u,n_p,h,v).
\]
All changes of the bar word below use the block-potential replacement
of Lemma~\ref{lem:H-block-bars}, so equality of local products and local bar sums suffices even
in the presence of arbitrary later factors.

If $s_\alpha p>p$, the canonical left Bruhat factor $u$ lies in
$U_{\rm rest}$. Put
\[
 u_1=x_\alpha(z)u x_\alpha(-z),\qquad
 u_2=n_\alpha u_1n_\alpha^{-1}\in U.
\]
Then
\[
 n_\alpha x_\alpha(z)u n_p
 =u_2(n_\alpha n_p)(n_p^{-1}x_\alpha(z)n_p).
\]
The intervening changes have zero internal bar cost: they merge or split
positive-unipotent factors or move constant lifts. The last factor is
in the positive root group for $p^{-1}\alpha$, and can be absorbed
past $h$ into $v$ with zero internal cost. The result is the new Bruhat
primitive with $\Omega_{\log}$ unchanged. The same proof works when
$u$ is any $U_{\rm rest}$-valued map, whether or not it is canonical.

For a down or stay step write $p=s_\alpha p'$ and temporarily choose
$n_p=n_\alpha n_{p'}$. Other lift choices differ only by constant
Cartan factors. Write $u=x_\alpha(c)u_0$ with
$u_0\in U_{\rm rest}$, put $a=z+c$, and set
$u_0'=n_\alpha^{-1}u_0n_\alpha\in U_{\rm rest}$.
The exact bar word becomes
\[
 (n_\alpha,x_\alpha(a),n_\alpha,u_0',n_{p'},h,v)
\]
with zero insertion corrections. On the down locus $a=0$, its first
three factors have constant product $n_\alpha^2$. Absorbing this constant
changes only the left unipotent and a constant Cartan multiplier.
Thus $\Omega_{\log}$ is unchanged. Here the forms are pulled back to
the locus $a=0$, so $da=0$.

On the stay locus, use the rank-one identity
\[
 n_\alpha x_\alpha(a)n_\alpha
 =t'(a)x_\alpha(-a)n_\alpha x_\alpha(-a^{-1}),
 \qquad t'(a)=\alpha^\vee(a^{-1}),
\]
for the convention
$n=\left(\begin{smallmatrix}0&-1\\1&0\end{smallmatrix}\right)$.
The other sign convention changes $t'$ by a constant.
Both sides, viewed as blocks, have zero internal bar sum because they
depend on the single variable $a$. The new bar word is therefore
\[
 (t',x_\alpha(-a),n_\alpha,x_\alpha(-a^{-1}),
                       u_0',n_{p'},h,v).
\]
Since $p'^{-1}\alpha$ is positive, the generalized up-step argument
applies to the subword beginning with $n_\alpha$. That argument replaces the
subword, with equal product and equal internal potential, by
$(u_1,n_p,h,v_1)$ for $u_1,v_1\in U$.
Merge the two leading positive-unipotent factors and move $t'$ across
their product; these operations have zero internal bars. Moving $t'$
through the constant $n_p$ yields
\[
 (u_2,n_p,\phi,h,v_1),\qquad
 \phi=p^{-1}(t')=\beta^\vee(a^{-1}).
\]
The only nonzero merger is that of the adjacent variable Cartan factors
$\phi,h$; this merger retains $B(\phi,h)$ with positive sign. The remaining
term is the new Bruhat primitive $B(u_2,n_p\phi h v_1)$, which proves
the asserted recursion. All interactions of $U_{\rm rest}$ with the
later Bruhat factor have remained in that primitive throughout.
\end{proof}

\paragraph{Universal rank-one sign check.}
Let $\kappa=\operatorname{Tr}$ and take independent coordinates
$c,z,t,v$ with $t\ne0$. Put $a=z+c$ and
\[
 F_{\rm old}=x(c)n h(t)x(v),\quad
 n=\begin{pmatrix}0&-1\\1&0\end{pmatrix},\quad
 x(q)=\begin{pmatrix}1&q\\0&1\end{pmatrix},\quad
 h(t)=\begin{pmatrix}t&0\\0&t^{-1}\end{pmatrix}.
\]
On $a\ne0$ exact multiplication gives
\[
 nx(z)F_{\rm old}
 =x(-a^{-1})n h(at)x\bigl(v-(at^2)^{-1}\bigr).
\]
The old primitive plus the fusion bar, the new primitive, and the
Cartan correction are respectively
\[
 -t^2\,da\wedge dv,\qquad
 -t^2\,da\wedge dv-\frac{2}{at}\,da\wedge dt,\qquad
 \frac{2}{at}\,da\wedge dt.
\]
Here $\beta=-\alpha$, so
$\beta^\vee(a^{-1})=\operatorname{diag}(a,a^{-1})$; this verifies the
sign and the coefficient on the full universal chart. On the down
locus $z=-c$, pullback gives $da=0$ and the product is $h(-t)x(v)$.
On the up chart $F_{\rm old}=h(t)x(v)$ the new product is
$nh(t)x(v+z/t^2)$ and its fusion bar is zero.
For arbitrary invariant $\kappa$, its restriction to the rank-one
algebra is a scalar multiple of the trace, while the component of
$d\log h$ in $\ker\alpha$ is orthogonal to that algebra. Thus the
invariant correction above applies to arbitrary $h\in T$ without
choosing an algebraic splitting of $T$.

\paragraph{Endpoint, initial factor, and normalization.}
After the surface transport, all handle Cartan factors are in the
separate factor $D$ of $u_0D F_{\rm word}$. The word induction is applied
to $F_{\rm word}$ alone, so $h_0=1$. At the closed-walk endpoint the
remaining Bruhat primitive is $U$--$B$ and vanishes. Hence
\[
 \Omega_{\rm word}
 =\sum_j\kappa(-\beta_j^\vee d\log a_j
                         \wedge d\log h_{j-1}).
\]
The subsequent substitution $a_j=e_j\beta_j(h_{j-1})$ gives
$v_j=R_jv_{j-1}-q_j$ and
$d\log\mu=W\sum_jR_{<j}q_j$, with no initial $D$ term.
The Goldman form is one half of this bar convention at every
stage: double, puncture, word, normalizer, and parabolic pairing.
There is no later rescaling of $E$. Constant Tits factors change level
constants only and do not alter this normalization or later operator
identities involving $E$.

\subsection*{Acknowledgments \& AI disclosure}
The author has benefited greatly from talks given on various occasions by Davesh Maulik, Anton Mellit, Junliang Shen, and Qizheng Yin on proofs of the $P=W$ conjecture for $\mathrm{GL}_n$.

The project has been developed with substantial AI assistance since 20 August 2026. The main results were obtained after weeks of discussion with AI, in which the author supplied the main ideas and directed the work while the AI served as a tool to check feasibility and help develop technical arguments. 
The source draft records the author's independent checking and substantial revision of the arguments and their references.

\addtocontents{toc}{\protect\setcounter{tocdepth}{1}}


\begin{thebibliography}{dCMS22b}

\bibitem[AHJR14]{AHJR14}
P.~N.~Achar, A.~Henderson, D.~Juteau and S.~Riche,
\emph{Weyl group actions on the Springer sheaf},
Proc. Lond. Math. Soc. (3) \textbf{108} (2014), no.~6, 1501--1528.
\url{https://doi.org/10.1112/plms/pdt055}.

\bibitem[AHLH23]{AHLH23}
J.~Alper, D.~Halpern-Leistner and J.~Heinloth,
\emph{Existence of moduli spaces for algebraic stacks},
Invent. Math. \textbf{234} (2023), no.~3, 949--1038.
\url{https://doi.org/10.1007/s00222-023-01214-4}.

\bibitem[Alp13]{Alp13}
J.~Alper, \emph{Good moduli spaces for Artin stacks},
Ann. Inst. Fourier (Grenoble) \textbf{63} (2013), no.~6, 2349--2402.
\url{https://doi.org/10.5802/aif.2833}.

\bibitem[AMM98]{AMM98}
A.~Alekseev, A.~Malkin and E.~Meinrenken,
\emph{Lie group valued moment maps},
J. Differential Geom. \textbf{48} (1998), no.~3, 445--495.
\url{https://doi.org/10.4310/jdg/1214460860}.

\bibitem[AOV08]{AOV08}
D.~Abramovich, M.~Olsson and A.~Vistoli,
\emph{Tame stacks in positive characteristic},
Ann. Inst. Fourier (Grenoble) \textbf{58} (2008), no.~4, 1057--1091;
Corrigendum, ibid.\ \textbf{64} (2014), no.~3, 945--946.
\url{https://doi.org/10.5802/aif.2378}.
Corrigendum DOI: \url{https://doi.org/10.5802/aif.2869}.

\bibitem[Bal23]{Bal23}
M.~Ballandras,
\emph{Intersection cohomology of character varieties for punctured
Riemann surfaces},
J. \'Ec. polytech. Math. \textbf{10} (2023), 141--198.
\url{https://doi.org/10.5802/jep.215}.

\bibitem[BBD82]{BBD82}
A.~A.~Beilinson, J.~Bernstein and P.~Deligne,
\emph{Faisceaux pervers},
Ast\'erisque \textbf{100} (1982), 5--171.
\url{https://numdam.org/item/AST_1982__100__1_0/}.

\bibitem[BBP17]{BBP17}
V.~Balaji, I.~Biswas and Y.~Pandey,
\emph{Connections on parahoric torsors over curves},
Publ. Res. Inst. Math. Sci. \textbf{53} (2017), no.~4, 551--585.
\url{https://doi.org/10.4171/PRIMS/53-4-3}.

\bibitem[BGPM20]{BGPM20}
O.~Biquard, O.~Garc\'{\i}a-Prada and I.~Mundet i Riera,
\emph{Parabolic Higgs bundles and representations of the fundamental
group of a punctured surface into a real group},
Adv. Math. \textbf{372} (2020), 107305. Locators refer to
arXiv:1510.04207v4.
\url{https://doi.org/10.1016/j.aim.2020.107305}.

\bibitem[Biq97]{Biq97}
O.~Biquard,
\emph{Fibr\'es de Higgs et connexions int\'egrables: le cas logarithmique
(diviseur lisse)},
Ann. Sci. \'Ecole Norm. Sup. (4) \textbf{30} (1997), no.~1, 41--96.
\url{https://doi.org/10.1016/s0012-9593(97)89915-6}.

\bibitem[BM83]{BM83}
W.~Borho and R.~MacPherson,
\emph{Partial resolutions of nilpotent varieties},
Ast\'erisque \textbf{101--102} (1983), 23--74.
\url{https://www.numdam.org/item/AST_1983__101-102__23_0/}.

\bibitem[BMW12]{BMW12}
I.~Biswas, S.~Majumder and M.~L.~Wong,
\emph{Root stacks, principal bundles and connections},
Bull. Sci. Math. \textbf{136} (2012), no.~4, 369--398. Locators refer to
arXiv:1203.6556v1.
\url{https://doi.org/10.1016/j.bulsci.2012.03.006}.

\bibitem[Boa11]{Boa11}
P.~Boalch,
\emph{Riemann--Hilbert for tame complex parahoric connections},
Transform. Groups \textbf{16} (2011), no.~1, 27--50.
\url{https://doi.org/10.1007/s00031-011-9121-1}.

\bibitem[Bor61]{Bor61}
A.~Borel, Sous-groupes commutatifs et torsion des groupes de Lie compacts
connexes, \emph{T\^ohoku Math. J. (2)} \textbf{13} (1961), no.~2, 216--240.
\url{https://doi.org/10.2748/tmj/1178244298}.

\bibitem[BS06]{BS06}
I.~Biswas and G.~Schumacher,
\emph{Geometry of moduli spaces of Higgs bundles},
Comm. Anal. Geom. \textbf{14} (2006), no.~4, 765--793. Locators refer to
arXiv:math/0605589v1.
\url{https://doi.org/10.4310/CAG.2006.v14.n4.a7}.

\bibitem[BS15]{BS15}
V.~Balaji and C.~S.~Seshadri,
\emph{Moduli of parahoric $\mathcal G$-torsors on a compact Riemann
surface}, J. Algebraic Geom. \textbf{24} (2015), no.~1, 1--49.
\url{https://doi.org/10.1090/S1056-3911-2014-00626-3}.

\bibitem[BT71]{BT71}
A.~Borel and J.~Tits,
\emph{\'El\'ements unipotents et sous-groupes paraboliques de groupes
r\'eductifs.~I}, Invent. Math. \textbf{12} (1971), no.~2, 95--104.
\url{https://doi.org/10.1007/BF01404653}.

\bibitem[BV12]{BV12}
N.~Borne and A.~Vistoli,
\emph{Parabolic sheaves on logarithmic schemes},
Adv. Math. \textbf{231} (2012), no.~3--4, 1327--1363.
\url{https://doi.org/10.1016/j.aim.2012.06.015}.

\bibitem[Cad07]{Cad07}
C.~Cadman,
\emph{Using stacks to impose tangency conditions on curves},
Amer. J. Math. \textbf{129} (2007), no.~2, 405--427.
\url{https://doi.org/10.1353/ajm.2007.0007}.

\bibitem[CM93]{CM93}
D.~H.~Collingwood and W.~M.~McGovern,
\emph{Nilpotent Orbits in Semisimple Lie Algebras},
Van Nostrand Reinhold, New York, 1993.

\bibitem[Cor88]{Cor88}
K.~Corlette, \emph{Flat $G$-bundles with canonical metrics},
J. Differential Geom. \textbf{28} (1988), no.~3, 361--382.
\url{https://doi.org/10.4310/jdg/1214442469}.

\bibitem[dA25]{dA25}
L.~de~Amorin, \emph{Motivic mirror symmetry for character stacks},
arXiv:2508.12995v2 (2025).
\url{https://arxiv.org/abs/2508.12995v2}.

\bibitem[dC21]{dC21}
M.~A.~de~Cataldo,
\emph{Projective compactification of Dolbeault moduli spaces},
Int. Math. Res. Not. IMRN \textbf{2021}, no.~5, 3543--3570.
\url{https://doi.org/10.1093/imrn/rnaa069}.

\bibitem[dCHM12]{dCHM12}
M.~A.~de~Cataldo, T.~Hausel and L.~Migliorini,
\emph{Topology of Hitchin systems and Hodge theory of character
varieties: the case $A_1$},
Ann. of Math. (2) \textbf{175} (2012), no.~3, 1329--1407.
\url{https://doi.org/10.4007/annals.2012.175.3.7}.

\bibitem[dCM09]{dCM09}
M.~A.~de~Cataldo and L.~Migliorini,
\emph{The decomposition theorem, perverse sheaves and the topology of
algebraic maps},
Bull. Amer. Math. Soc. (N.S.) \textbf{46} (2009), no.~4, 535--633.
\url{https://doi.org/10.1090/S0273-0979-09-01260-9}.

\bibitem[dCMa20]{dCMa20}
M.~A.~de~Cataldo and D.~Maulik,
\emph{The perverse filtration for the Hitchin fibration is locally
constant},
Pure Appl. Math. Q. \textbf{16} (2020), no.~5, 1441--1464.
\url{https://doi.org/10.4310/PAMQ.2020.v16.n5.a4}.

\bibitem[dCMS22a]{dCMS22}
M.~A.~de~Cataldo, D.~Maulik and J.~Shen,
\emph{Hitchin fibrations, abelian surfaces, and the $P=W$ conjecture},
J. Amer. Math. Soc. \textbf{35} (2022), no.~3, 911--953.
\url{https://doi.org/10.1090/jams/989}.

\bibitem[dCMS22b]{dCMS22SL}
M.~A.~de~Cataldo, D.~Maulik and J.~Shen,
\emph{On the $P=W$ conjecture for $\mathrm{SL}_n$},
Selecta Math. (N.S.) \textbf{28} (2022), no.~5, Paper No.~90.
Locators refer to arXiv:2002.03336.
\url{https://doi.org/10.1007/s00029-022-00803-0}.

\bibitem[Del71]{Del71}
P.~Deligne, \emph{Th\'eorie de Hodge.~II},
Publ. Math. IH\'ES \textbf{40} (1971), 5--57.
\url{https://doi.org/10.1007/BF02684692}.

\bibitem[Del74]{Del74}
P.~Deligne, \emph{Th\'eorie de Hodge.~III},
Publ. Math. IH\'ES \textbf{44} (1974), 5--77.
\url{https://doi.org/10.1007/BF02685881}.

\bibitem[Deo85]{Deo85}
V.~V.~Deodhar,
\emph{On some geometric aspects of Bruhat orderings.~I. A finer
decomposition of Bruhat cells}, Invent. Math. \textbf{79} (1985),
no.~3, 499--511.
\url{https://doi.org/10.1007/BF01388520}.

\bibitem[DH82]{DH82}
J.~J.~Duistermaat and G.~J.~Heckman,
\emph{On the variation in the cohomology of the symplectic form of the
reduced phase space}, Invent. Math. \textbf{69} (1982), no.~2, 259--268.
\url{https://doi.org/10.1007/BF01399506}.

\bibitem[DLP88]{DLP88}
C.~De~Concini, G.~Lusztig and C.~Procesi,
\emph{Homology of the zero-set of a nilpotent vector field on a flag
manifold}, J. Amer. Math. Soc. \textbf{1} (1988), no.~1, 15--34.
\url{https://doi.org/10.1090/S0894-0347-1988-0924700-2}.

\bibitem[Don87]{Don87}
S.~K.~Donaldson,
\emph{Twisted harmonic maps and the self-duality equations},
Proc. London Math. Soc. (3) \textbf{55} (1987), no.~1, 127--131.
\url{https://doi.org/10.1112/plms/s3-55.1.127}.

\bibitem[FHZ26]{FHZ26}
A.~Fernandez Herrero and S.~Zhang,
\emph{Meromorphic Hodge moduli spaces for reductive groups in arbitrary
characteristic}, Michigan Math. J., advance publication (2026), 1--62,
doi:10.1307/mmj/20246514. Locators refer to arXiv:2307.16755v2.

\bibitem[FM22]{FM22}
C.~Felisetti and M.~Mauri,
\emph{$P=W$ conjectures for character varieties with symplectic
resolution}, J. \'Ec. polytech. Math. \textbf{9} (2022), 853--905.
\url{https://doi.org/10.5802/jep.196}.

\bibitem[GM83]{GM83}
M.~Goresky and R.~MacPherson, \emph{Intersection homology.~II},
Invent. Math. \textbf{72} (1983), no.~1, 77--129.
\url{https://doi.org/10.1007/BF01389130}.

\bibitem[Gol84]{Gol84}
W.~M.~Goldman,
\emph{The symplectic nature of fundamental groups of surfaces},
Adv. Math. \textbf{54} (1984), no.~2, 200--225.
\url{https://doi.org/10.1016/0001-8708(84)90040-9}.

\bibitem[Hit87a]{Hit87a}
N.~J.~Hitchin,
\emph{The self-duality equations on a Riemann surface},
Proc. London Math. Soc. (3) \textbf{55} (1987), no.~1, 59--126.
\url{https://doi.org/10.1112/plms/s3-55.1.59}.

\bibitem[Hit87b]{Hit87b}
N.~J.~Hitchin, \emph{Stable bundles and integrable systems},
Duke Math. J. \textbf{54} (1987), no.~1, 91--114.
\url{https://doi.org/10.1215/S0012-7094-87-05408-1}.

\bibitem[HLFH23]{HLFH23}
D.~Halpern-Leistner and A.~Fernandez Herrero,
\emph{The structure of the moduli of gauged maps from a smooth curve},
arXiv:2305.09632v1 (2023), 90 pp.,
\url{https://arxiv.org/abs/2305.09632v1}.

\bibitem[HLRV11]{HLRV11}
T.~Hausel, E.~Letellier and F.~Rodriguez-Villegas,
\emph{Arithmetic harmonic analysis on character and quiver varieties},
Duke Math. J. \textbf{160} (2011), no.~2, 323--400.
\url{https://doi.org/10.1215/00127094-1444258}.

\bibitem[HMMS]{HMMS}
T.~Hausel, A.~Mellit, A.~Minets and O.~Schiffmann,
\emph{$P=W$ via $\mathcal H_2$}, to appear in J. Amer. Math. Soc.
The final preprint is arXiv:2209.05429v3 (3 November 2026).
All numbered locators here refer to arXiv:2209.05429v2
(17 January 2025).
\url{https://arxiv.org/abs/2209.05429v3};
\url{https://arxiv.org/abs/2209.05429v2}.

\bibitem[HRV08]{HRV08}
T.~Hausel and F.~Rodriguez-Villegas,
\emph{Mixed Hodge polynomials of character varieties}, with an appendix
by N.~M.~Katz, Invent. Math. \textbf{174} (2008), no.~3, 555--624.
\url{https://doi.org/10.1007/s00222-008-0142-x}.

\bibitem[HRV15]{HRV15}
T.~Hausel and F.~Rodriguez-Villegas,
\emph{Cohomology of large semiprojective hyperk\"ahler varieties},
Ast\'erisque \textbf{370} (2015), 113--156.
\url{https://doi.org/10.24033/ast.967}.

\bibitem[Kna23]{Kna23}
A.~W.~Knapp, \emph{Lie Groups Beyond an Introduction}, Digital Second
Edition, published by the author, East Setauket, New York, 2023.
\url{https://www.math.stonybrook.edu/~aknapp/download/Beyond2-clickable.pdf}.

\bibitem[Kos63]{Kos63}
B.~Kostant, \emph{Lie group representations on polynomial rings},
Amer. J. Math. \textbf{85} (1963), no.~3, 327--404.
\url{https://doi.org/10.2307/2373130}.

\bibitem[KS90]{KS90}
M.~Kashiwara and P.~Schapira, \emph{Sheaves on Manifolds},
Grundlehren der mathematischen Wissenschaften \textbf{292}, Springer,
Berlin, 1990.
\url{https://doi.org/10.1007/978-3-662-02661-8}.

\bibitem[KSZ24]{KSZ24}
G.~Kydonakis, H.~Sun and L.~Zhao,
\emph{Logahoric Higgs torsors for a complex reductive group},
Math. Ann. \textbf{388} (2024), no.~3, 3183--3228. Locators refer to
arXiv:2107.01977v3.
\url{https://doi.org/10.1007/s00208-023-02605-x}.

\bibitem[Let15]{Let15}
E.~Letellier,
\emph{Character varieties with Zariski closures of $\mathrm{GL}_n$-conjugacy
classes at punctures}, Selecta Math. (N.S.) \textbf{21} (2015), no.~1,
293--344.
\url{https://doi.org/10.1007/s00029-014-0163-9}.

\bibitem[Lop23]{Lop23}
R.~Lopez, \emph{Picard groups of stacky curves}, arXiv:2306.08227v1
(2023).
\url{https://arxiv.org/abs/2306.08227v1}.

\bibitem[Mel25]{Mel25}
A.~Mellit, \emph{Toric stratifications of character varieties},
Publ. Math. IH\'ES \textbf{142} (2025), 153--240.
\url{https://doi.org/10.1007/s10240-025-00158-0}.

\bibitem[Moc02]{Moc02}
T.~Mochizuki, \emph{Asymptotic behaviour of tame nilpotent harmonic
bundles with trivial parabolic structure}, J. Differential Geom.
\textbf{62} (2002), no.~3, 351--559. Locators refer to
arXiv:math/0212232v1. \url{https://doi.org/10.4310/jdg/1090426286}.

\bibitem[Moc05]{Moc05}
T.~Mochizuki,
\emph{Kobayashi--Hitchin correspondence for tame harmonic bundles and an
application}, Ast\'erisque \textbf{309} (2006), viii+117~pp. Locators refer to the
numbering of Kyoto University preprint 2005/15.
\url{https://doi.org/10.24033/ast.716}.

\bibitem[Moc07]{Moc07}
T.~Mochizuki,
\emph{Asymptotic behaviour of tame harmonic bundles and an application to
pure twistor $D$-modules.~I}, Mem. Amer. Math. Soc. \textbf{185} (2007),
no.~869, xii+324~pp. Locators refer to arXiv:math/0312230v2.
\url{https://doi.org/10.1090/memo/0869}.

\bibitem[Moc09]{Moc09}
T.~Mochizuki,
\emph{Kobayashi--Hitchin correspondence for tame harmonic bundles.~II},
Geom. Topol. \textbf{13} (2009), no.~1, 359--455.
\url{https://doi.org/10.2140/gt.2009.13.359}.

\bibitem[MS24]{MS24}
D.~Maulik and J.~Shen, \emph{The $P=W$ conjecture for $\mathrm{GL}_n$},
Ann. of Math. (2) \textbf{200} (2024), no.~2, 529--556.
\url{https://doi.org/10.4007/annals.2024.200.2.3}.

\bibitem[MSY25]{MSY25}
D.~Maulik, J.~Shen and Q.~Yin,
\emph{Perverse filtrations and Fourier transforms},
Acta Math. \textbf{234} (2025), no.~1, 1--69.
See also the updated proof version arXiv:2308.13160v3 (8 October 2025).
\url{https://doi.org/10.4310/ACTA.2025.v234.n1.a1}.

\bibitem[Nit91]{Nit91}
N.~Nitsure, \emph{Moduli space of semistable pairs on a curve},
Proc. London Math. Soc. (3) \textbf{62} (1991), no.~2, 275--300.
\url{https://doi.org/10.1112/plms/s3-62.2.275}.

\bibitem[Ric88]{Ric88}
R.~W.~Richardson,
\emph{Conjugacy classes of $n$-tuples in Lie algebras and algebraic
groups}, Duke Math. J. \textbf{57} (1988), no.~1, 1--35.
\url{https://doi.org/10.1215/S0012-7094-88-05701-8}.



\bibitem[Saf17]{Saf17}
P.~Safronov,
\emph{Symplectic implosion and the Grothendieck--Springer resolution},
Transform. Groups \textbf{22} (2017), no.~3, 767--792.
Locators refer to arXiv:1411.2962v2 (28 April 2016).
\url{https://doi.org/10.1007/s00031-016-9398-1}.

\bibitem[Sai88]{Sai88}
M.~Saito, \emph{Modules de Hodge polarisables},
Publ. Res. Inst. Math. Sci. \textbf{24} (1988), no.~6, 849--995.
\url{https://doi.org/10.2977/PRIMS/1195173930}.

\bibitem[Sai90]{Sai90}
M.~Saito, \emph{Mixed Hodge modules},
Publ. Res. Inst. Math. Sci. \textbf{26} (1990), no.~2, 221--333.
\url{https://doi.org/10.2977/PRIMS/1195171082}.

\bibitem[Sai90e]{SaiExt}
M.~Saito, \emph{Extension of mixed Hodge modules},
Compositio Math. \textbf{74} (1990), no.~2, 209--234.
\url{https://www.numdam.org/item/CM_1990__74_2_209_0/}.

\bibitem[Sai91]{SaiMix}
M.~Saito, \emph{On the formalism of mixed sheaves},
RIMS preprint 784 (1991), arXiv:math/0611597 (2006).
\url{https://arxiv.org/abs/math/0611597}.

\bibitem[San26]{San23}
T.~Santens,
\emph{The Brauer--Manin obstruction for stacky curves},
arXiv:2210.17184v4 (23 September 2026).
\url{https://arxiv.org/abs/2210.17184v4}.



\bibitem[She17]{She17}
V.~Shende,
\emph{The weights of the tautological classes of character varieties},
Int. Math. Res. Not. IMRN \textbf{2017}, no.~22, 6832--6840.
The unnumbered theorem locator refers to arXiv:1411.4975v2, p.~2.
\url{https://doi.org/10.1093/imrn/rnv363}.

\bibitem[Sch73]{Sch73}
W.~Schmid, \emph{Variation of Hodge structure: the singularities of the
period mapping}, Invent. Math. \textbf{22} (1973), 211--319.
\url{https://doi.org/10.1007/BF01389674}.

\bibitem[Sho88]{Sho88}
T.~Shoji,
\emph{Geometry of orbits and Springer correspondence},
in \emph{Orbites unipotentes et repr\'esentations. I.
Groupes finis et alg\`ebres de Hecke},
Ast\'erisque \textbf{168} (1988), 61--140.
\url{https://numdam.org/item/AST_1988__168__61_0/}.

\bibitem[Sim88]{Sim88}
C.~T.~Simpson,
\emph{Constructing variations of Hodge structure using Yang--Mills
theory and applications to uniformization},
J. Amer. Math. Soc. \textbf{1} (1988), no.~4, 867--918.
\url{https://doi.org/10.1090/S0894-0347-1988-0944577-9}.

\bibitem[Sim90]{Sim90}
C.~T.~Simpson, \emph{Harmonic bundles on noncompact curves},
J. Amer. Math. Soc. \textbf{3} (1990), no.~3, 713--770.
\url{https://doi.org/10.1090/S0894-0347-1990-1040197-8}.

\bibitem[Sim92]{Sim92}
C.~T.~Simpson, \emph{Higgs bundles and local systems},
Publ. Math. IH\'ES \textbf{75} (1992), 5--95.
\url{https://doi.org/10.1007/BF02699491}.


\bibitem[Sim97]{Sim97}
C.~T.~Simpson, \emph{The Hodge filtration on nonabelian cohomology}, in:
\emph{Algebraic Geometry---Santa Cruz 1995}, Proc. Sympos. Pure Math.
\textbf{62}, Part~2, Amer. Math. Soc., Providence, RI, 1997, 217--281.
\url{https://doi.org/10.1090/pspum/062.2/1492538}.

\bibitem[Spr76]{Spr76}
T.~A.~Springer,
\emph{Trigonometric sums, Green functions of finite groups and
representations of Weyl groups}, Invent. Math. \textbf{36} (1976), no.~1,
173--207.
\url{https://doi.org/10.1007/BF01390009}.

\bibitem[Spr78]{Spr78}
T.~A.~Springer, \emph{A construction of representations of Weyl groups},
Invent. Math. \textbf{44} (1978), no.~3, 279--293.
\url{https://doi.org/10.1007/BF01403165}.

\bibitem[Spr98]{Spr98}
T.~A.~Springer, \emph{Linear Algebraic Groups}, 2nd ed., Progress in
Mathematics \textbf{9}, Birkh\"auser, Boston, 1998.
\url{https://doi.org/10.1007/978-0-8176-4840-4}.

\bibitem[Ste65]{Ste65}
R.~Steinberg, \emph{Regular elements of semi-simple algebraic groups},
Publ. Math. Inst. Hautes \'Etudes Sci. \textbf{25} (1965), 49--80.
\url{https://doi.org/10.1007/BF02684397}.

\bibitem[Stacks]{Stacks}
The Stacks Project Authors, \emph{The Stacks Project},
Lemma~59.69.3 (cohomology of a smooth affine curve), Tag~03RR, and
Lemma~106.8.7 (infinitesimal lifting criterion), Tag~0DP0.
\url{https://stacks.math.columbia.edu/tag/03RR},
\url{https://stacks.math.columbia.edu/tag/0DP0}.

\bibitem[Ste75]{Ste75}
R.~Steinberg, \emph{Torsion in reductive groups},
Adv. Math. \textbf{15} (1975), no.~1, 63--92.
\url{https://doi.org/10.1016/0001-8708(75)90125-5}.

\bibitem[Sti42]{Sti42}
E.~Stiefel, \"Uber eine Beziehung zwischen geschlossenen Lie'schen Gruppen
und diskontinuierlichen Bewegungsgruppen euklidischer R\"aume und ihre
Anwendung auf die Aufz\"ahlung der einfachen Lie'schen Gruppen,
\emph{Comment. Math. Helv.} \textbf{14} (1941/42), 350--380.
\url{https://doi.org/10.1007/BF02565625}.

\bibitem[Su26]{Su26}
T.~Su, \emph{Cell decomposition and dual boundary complexes of character
varieties}, Adv. Math. \textbf{498} (2026), 111043.
Corrected statements and locators refer to arXiv:2307.16657v5
(11 October 2026).
\url{https://doi.org/10.1016/j.aim.2026.111043};
\url{https://arxiv.org/abs/2307.16657v5}.

\bibitem[Tse33]{Tse33}
Ch.~C.~Tsen, \emph{Divisionsalgebren \"uber Funktionenk\"orpern},
Nachr. Ges. Wiss. G\"ottingen, Math.-Phys. Kl. (1933), 335--339.
\url{https://gdz.sub.uni-goettingen.de/id/PPN252457811_1933}.

\bibitem[Tub25]{Tub25}
S.~Tubach, \emph{Mixed Hodge modules on stacks},
Forum Math. Sigma \textbf{13} (2025), Paper No.~e175,
doi:10.1017/fms.2025.10122.

\bibitem[Zuc79]{Zuc79}
S.~Zucker, \emph{Hodge theory with degenerating coefficients: $L_2$
cohomology in the Poincar\'e metric}, Ann. of Math. (2) \textbf{109}
(1979), no.~3, 415--476.
\url{https://doi.org/10.2307/1971221}.

\end{thebibliography}
\end{document}